\documentclass[reqno,10pt]{amsart}

\usepackage{tikz}
\usepackage{external/takodachi}

\usepackage{external/css}

\title
{
	The Cauchy problem for Manton's Chern-Simons-Schr\"odinger equation
} 
\author{Jason Zhao}
\address{Department of Mathematics, University of California, Berkeley, 94720}
\email{zhao.j@berkeley.edu}
\date{\today}

\begin{document}

\begin{abstract}
	The Chern-Simons-Schr\"odinger equation (in the temporal gauge) arises as the Hamiltonian flow of the abelian Higgs energy on $\R^2$. Manton \cite{Manton1997} introduced the equation as a model for the dynamics of the critical points of the energy, known as vortices. He conjectured that, for a certain range of coupling constants, the vortex motion under the Chern-Simons-Schr\"odinger flow can be effectively captured by a first-order ODE on the moduli space of self-dual vortices constructed by Jaffe-Taubes \cite{JaffeTaubes1980} and Samols \cite{Samols1992}. 
	
	As a first step towards a rigorous proof of Manton's conjecture, we formulate the Cauchy problem in DeTurck gauge within the natural energy space and prove global well-posedness. We also obtain, as a corollary of the well-posedness theory and the stability results in our previous work \cite{Zhao2026}, orbital stability of the self-dual vortices under the Chern-Simons-Schr\"odinger flow near self-dual coupling. The heart of our analysis lies in developing a geometric Littlewood-Paley theory based on the covariant heat equation in caloric gauge, which we use to perform a paradifferential-style decomposition of the Chern-Simons-Schr\"odinger equation. 
\end{abstract}
\maketitle

\tableofcontents

\section{Introduction}

In \cite{Manton1997}, Manton introduced the Chern-Simons-Schr\"odinger equation as a model for vortex motion in thin superconductors. The model describes the dynamics of a Higgs field and an electromagnetic potential, represented by a complex scalar field $\phi : \R \times \R^2 \to \C$ and a real-valued connection $1$-form $A$ on $\R \times \R^2$ respectively. The electromagnetic field associated to the potential is given by the curvature $2$-form $F = dA$; written in coordinates $F_{\mu\nu} = \partial_\mu A_\nu - \partial_\nu A_\mu$, the components $F_{tj}$ constitute the electric field, while $F_{12}$ is the magnetic field. Denoting the covariant derivatives by $\bfD_\mu = \partial_\mu - i A_\mu$, Manton defined the Lagrangian
	\begin{align}\label{eq:lagrange}\tag{L}
		\cL[A, \phi] 
			:= \tfrac12 \int_{\R^{1 + 2}} A \wedge F  - \tfrac12 \int_{\R^{1 + 2}} \Im(\overline \phi \bfD_t \phi) - A_t + |\bfD_1 \phi|^2 +  |\bfD_1 \phi|^2 + |F_{12}|^2 + \tfrac\lambda4 (1 - |\phi|^2)^2 \, dt dx,
	\end{align}
where $\lambda > 0$ is a coupling constant. The first integral is the Chern-Simons action, and in the second integral, the last four terms constitute the energy density for the abelian Higgs model. The Euler-Lagrange equation is given by
	\begin{equation}
		\begin{split}
			i \bfD_t \phi + \bfD_1 \bfD_1 \phi + \bfD_2 \bfD_2 \phi
				&= - \tfrac{\lambda}{2}(1 - |\phi|^2) \phi, \\
			F_{t1}
				&= - \Im (\overline \phi \bfD_2 \phi) - \partial_1 F_{12}, \\
			F_{t2}
				&= + \Im (\overline \phi \bfD_1 \phi) - \partial_2 F_{12},\\
			F_{12}
				&= \tfrac12 (1  - |\phi|^2).
		\end{split}\label{eq:CSS}\tag{L-CSS}
	\end{equation}	
The first equation is a covariant Schr\"odinger equation, the second and third equations constitute Amp\'ere's law, relating the electric field to the charge current and magnetic flux, and the last equation is Gauss's law, relating the magnetic field to the charge density. 

In this article, we consider a generalisation of Manton's model, namely the system obtained by omitting Gauss's law from \eqref{CSS},
	\begin{equation}\label{eq:CSS-amp-schro}\tag{CSS}
        \begin{split}
            \Bigl( i \bfD_t + \bfD^j \bfD_j + \tfrac{\lambda}2 (1 - |\phi|^2) \Bigr) \phi 
                &= 0, \\
            F_{tj} 
                &= - \epsilon_{jk} \Bigl( \Im(\overline \phi \bfD^k \phi) + \partial^\ell F\indices{_\ell^k} \Bigr),
        \end{split}
    \end{equation}
where we have rewritten Amp\'ere's law in the abstract Einstein notation of summing over repeated indices. We refer to this system as the \textit{Chern-Simons-Schr\"odinger equation}. To recover Manton's model \eqref{CSS}, define the \textit{Gauss tension field} as the difference of the charge density and magnetic field, i.e. $\Gauss := \tfrac12(1 - |\phi|^2) - F_{12}$. A formal calculation (see Proposition \ref{prop:CSS-Gauss}) shows that the Gauss tension field is transported under the flow by \eqref{CSS-amp-schro}, 
	\begin{equation}\tag{G}\label{eq:gauss-transport-intro}
		\partial_t \Gauss
			= 0,
	\end{equation}
so Gauss's law is propagated provided that it is satisfied at an initial time. 

The system \eqref{CSS-amp-schro} exhibits a $\mathsf U(1)$-gauge symmetry, that is, the class of configurations $(A, \phi)$ on $\R \times \R^2$ solving the system is invariant under the gauge transformations\footnote{Rather unusually, the Lagrangian \eqref{lagrange} is not fully gauge-invariant due to the $A_t$-term, see \cite[Section 2]{Manton1997} and the references therein for further discussion. }
	\begin{align*}
		A_\mu 
			&\mapsto A_\mu + \partial_\mu \chi, \\
		\phi 
			&\mapsto e^{i \chi} \phi, 
	\end{align*}
for any real scalar field $\chi : \R \times \R^2 \to \R$. The system is also invariant under Galilean boosts, phase rotation, space translation, time reversal\footnote{Unlike the time reversal symmetry for the linear Schr\"odinger equation, the transformation for \eqref{CSS-amp-schro} consists of, in addition to the usual conjugation, conjugation of the gauge potential and reflection of the domain, 
	\begin{align*}
		A_t (t, x^1, x^2)
			&\mapsto A_t (-t, -x^1, x^2),\\
		(A_1 + i A_2) (t, x^1, x^2) 
			&\mapsto (\overline{A_1 + i A_2}) (-t, -x^1, x^2),\\
		\phi(t, x^1, x^2) 
			&\mapsto \overline{\phi(-t, -x^1, x^2)}.
	\end{align*}
}, and time translation. 

The Chern-Simons-Schr\"odinger equation in \textit{temporal gauge}, $A_t = 0$, formally arises as the Hamiltonian flow of the \textit{abelian Higgs energy} (also referred to as the {gauged Ginzburg-Landau energy}), a gauge-invariant functional on configurations which takes the form
	\begin{align*}
		\cE_{\textup{AH}}[A, \phi] 
			&:= \tfrac12 \int_{\R^2}  |\bfD_1 \phi|^2 + |\bfD_2 \phi|^2 + |F_{12}|^2 + \tfrac\lambda4 (1 -|\phi|^2)^2 \, dx.
	\end{align*}	
The critical points of the abelian Higgs energy are examples of topological defects known as vortices, termed as such in analogy with point vortices from classical hydrodynamics. The topological nature of these critical points can be seen from the quantisation of the \textit{topological degree},
	\[
		\deg [A, \phi]  
			:= \tfrac1{2\pi} \int_{\R^2} (1 - |\phi|^2) F_{12} + 2 \Im(\overline{\bfD_1 \phi} \bfD_2 \phi) \, dx,
	\]
which is well-defined, gauge-invariant, and integer-valued if $(A, \phi)$ has finite abelian Higgs energy \cite{CzubakJerrard2014}. Thus, the space of finite-energy configurations can be separated into topological classes. In the case $\lambda = 1$, the minimisers, known as self-dual vortices, within a topological class $N$ form an $N$-complex-dimensional K\"ahler manifold $\cM^N$ \cite{JaffeTaubes1980,Samols1992}. Manton conjectured that, within a low-energy regime and $|\lambda - 1| \ll 1$, one could effectively describe the motion of vortices under \eqref{CSS} by a reduced Lagrangian on the moduli space of self-dual vortices. With an eye towards understanding Manton's conjecture, it is natural to begin by studying the Cauchy problem for \eqref{CSS-amp-schro} in a class of configurations with finite abelian Higgs energy. 

\subsection{The Cauchy problem}

In order to interpret the Chern-Simons-Schr\"odinger equation \eqref{CSS-amp-schro} as a well-posed evolutionary problem, we need to eliminate the gauge freedom. Taking cue from previous work on related models \cite{Demoulini2007,LiuEtAl2014}, we impose the \textit{DeTurck gauge}, 
	\begin{equation}\tag{D}\label{eq:DeTurck}
        A_t 
            = \partial^j A_j.
    \end{equation}
Under the DeTurck gauge, we reformulate \eqref{CSS-amp-schro} (see Lemma \ref{lem:CSS-DT}) as a system consisting of an electromagnetic Schr\"odinger equation for the scalar field and heat equations for the magnetic potential, 
	\begin{equation}\tag{CSS-DT}\label{eq:CSS-DT-intro}
	\begin{split}
		\Bigl( i \partial_t + \bfD^j \bfD_j + \partial^j A_j + \tfrac\lambda2 (1 - |\phi|^2) \Bigr) \phi 
			&= 0,\\
		(\partial_t - \Delta) A_j - \Bigl( \Re(\overline \phi \bfD_j \phi) + \partial_j \Gauss \Bigr) + \epsilon_{jk} \Bigl( \Im(\overline \phi \bfD^k \phi) + \Re(\overline \phi \bfD^k \phi) + \partial^k \Gauss \Bigr)
			&= 0,
	\end{split}
	\end{equation}
for which we study the (forward-in-time)\footnote{The parabolicity of the equations for $A$ fixes the forward direction of time. To study the backwards-in-time Cauchy problem, one can replace the DeTurck gauge \eqref{DeTurck} with
	\begin{equation}\tag{D-}
		A_t = - \partial^j A_j,
	\end{equation}
under which the magnetic potential would solve a backwards heat equation.
} Cauchy problem with initial data
	\[
		(A, \phi)_{|t = 0}
			= (A^{\mathrm{in}}, \phi^{\mathrm{in}}).
	\]
There remains a residual gauge freedom in $t$-independent gauge-transformations, which we eliminate by imposing the \textit{Coulomb gauge} for the initial data, 
    \begin{equation}\tag{{Coul}}\label{eq:coulomb-data}
		\partial^j {A_j}_{|t = 0} = 0.
	\end{equation}

To set in stone the notion of well-posedness of the Cauchy problem for the Chern-Simons-Schr\"odinger equation in DeTurck gauge \eqref{CSS-DT-intro}, we introduce the phase space and a topology thereon with respect to which we will show continuity of the flow map.

\begin{definition}\label{def:topology-1}
    We say that $(A, \phi) \in H^1_\loc (\R^2)$ is an \textit{finite-energy configuration} on $\R^2$ if
		\[
			 \| A \|_{L^4 \, \cap \, \dot H^1}+ \| \bfD \phi \|_{L^2} + \| 1 - |\phi|^2 \|_{L^2} + \| \Gauss \|_{L^\infty \, \cap \, H^1} 
                < \infty.
		\]
	We denote the space of finite-energy configurations by $\frE$, and the sub-class thereof satisfying the Coulomb gauge \eqref{coulomb-data} by $\frE_{\mathrm{df}}$. Furthermore, we endow $\frE$ with the topology induced by the following maps:
		\begin{enumerate}[label=({\roman*})]
			\item the covariant gradient,
				\begin{align*}
					\frE 
						&\longrightarrow L^2 (\R^2)\\
					(A, \phi) 
						&\longmapsto \bfD \phi,
				\end{align*}

			\item the Gauss tension field, 
				\begin{align*}
					\frE 
						&\longrightarrow (H^1 \cap L^\infty) (\R^2)\\
					(A, \phi) 
						&\longmapsto \Gauss,
				\end{align*}

			\item the inclusion map\footnote{By Lemma \ref{lem:gerard}, the scalar field resides in the ambient space $(L^2 + L^\infty) (\R^2)$ whenever $\tfrac12(1 - |\phi|^2) \in L^2 (\R^2)$. },
				\begin{align*}
					\frE 
						&\longrightarrow (L^4 \cap \dot H^1)(\R^2) \times (L^2 + L^\infty) (\R^2)\\
					(A, \phi) 
						&\longmapsto (A, \phi).
				\end{align*}
		\end{enumerate}
	We also introduce the (slightly weaker\footnote{Continuity of the map 
		\begin{align*}
			\frE 
				&\to L^\infty (\R^2),\\
			(A, \phi)
				&\mapsto \Gauss,
		\end{align*}
	where $L^\infty (\R^2)$ is endowed with the norm topology, plays almost no role in any of our analysis. We mainly want (for convenience of statements) a topology which is characterised by its convergent sequences, and the Gauss tension fields of these sequences to be uniformly bounded. An alternative approach would be to work with a closed ball in $L^\infty (\R^2)$ with the weak-$*$ topology. 
	}) metric on the space of finite-energy configurations,
		\begin{align*}
			\dist_\frE \bigl( (A, \phi), (A', \phi') \bigr) 
				&:= \bigl\| \delta \bfD \phi \bigr\|_{L^2} + \bigl\| \delta \Gauss \bigr\|_{H^1} + \bigl\| \delta A \bigr\|_{L^4 \, \cap\, \dot H^1} + \bigl\| \delta \phi \bigr\|_{L^2 + L^\infty},
		\end{align*}
	using the notation\footnote{Refer to Section \ref{subsec:prelim-notation} for further explication.} $\delta \sfU := \sfU - \sfU'$ for differences. Given an interval $I \subseteq \R$, we write $C_t (I \to \frE)$ for the space of continuous maps $I \to \frE$ endowed with the standard topology. 
\end{definition}

The topology induced by the maps (i), (ii), (iii) manifestly render the abelian Higgs energy and topological degree as continuous functionals. We will need the additional regularity assumption $\Gauss \in (L^\infty \cap \dot H^1) (\R^2)$ to handle the appearance of the Gauss tension field in the energy estimates (see, for example, Theorem \ref{thm:smooth-LWP}), as it is merely transported a l\'a \eqref{gauss-transport-intro} along the flow. One should think of the inclusion map (iii) as a ``low frequency" part of the metric, which serves mostly to distinguish configurations.

With the stage set, we are ready to state the main result of this article, 

\begin{theorem}\label{thm:gwp}
	The Cauchy problem for the Chern-Simons-Schr\"odinger equation \eqref{CSS-amp-schro} in DeTurck gauge \eqref{DeTurck} with initial data in Coulomb gauge \eqref{coulomb-data} is globally well-posed in the space of finite-energy configurations in the following sense: 
		\begin{itemize}
			\item \textup{(Existence).} For each finite-energy configuration $(A^{\mathrm{in}}, \phi^{\mathrm{in}}) \in \frE_\df$, there exists $(A, \phi) \in C_t ([0, \infty) \to \frE)$ solving \eqref{CSS-DT-intro} with initial data $(A, \phi)_{|t = 0} = (A^{\mathrm{in}}, \phi^{\mathrm{in}})$. 
			
			\item \textup{(Unconditional uniqueness).} If $(A', \phi') \in C_t ([0, \infty) \to \frE)$ is another solution to \eqref{CSS-DT-intro} with initial data $(A, \phi)_{|t = 0} = (A^{\mathrm{in}}, \phi^{\mathrm{in}})$, then $(A, \phi) = (A', \phi')$.

			\item \textup{(Continuous dependence).} For each $T > 0$, the data-to-solution map,
				\begin{align*}
				\frE_{\mathrm{df}}
					&\longrightarrow C_t ([0, T] \to \frE),\\
				(A^{\mathrm{in}}, \phi^{\mathrm{in}}) 
					&\longmapsto (A, \phi),
				\end{align*}
			is continuous. 
		\end{itemize}
	Furthermore, these solutions satisfy the following conservation laws:
		\begin{itemize}
			\item \textup{(Conservation of abelian Higgs energy).} For all $t \geq 0$, 
				\begin{align*}
					\cE_{\textup{AH}} [A(t), \phi(t)] 
						&= \cE_{\textup{AH}} [A^{\mathrm{in}}, \phi^{\mathrm{in}}].
				\end{align*}

			\item \textup{(Conservation of topological degree).} For all $t \geq 0$, 
				\begin{align*}
					\deg [A(t), \phi(t)] 
						&= \deg [A^{\mathrm{in}}, \phi^{\mathrm{in}}]. 
				\end{align*}
		\end{itemize}
\end{theorem}

\begin{remark}
	One could probably also establish well-posedness in the affine Sobolev space $(\mathring A, \mathring{\vphantom{A}\phi}) + H^1 (\R^2)$ for each finite-energy configuration $(\mathring A, \mathring{\vphantom{A}\phi}) \in \frE^{\mathrm{df}}$. On the one hand, this gives more information about the flow $t \mapsto (A, \phi)(t)$, since the $H^1$-metric is stronger than the $\frE$-metric (see Proposition \ref{prop:magnetic-lip-1}), and suggests that the flow ``preserves'' the spatial asymptotics of the scalar field $\phi$. On the other hand, the $H^1$-metric is ill-suited for comparing configurations with different spatial asymptotics (see Appendix \ref{subapp:example} for examples), so in this sense the uniqueness and continuous dependence statement in Theorem \ref{thm:gwp} can be regarded as stronger.  
\end{remark}

The main conceptual difficulty with studying the Chern-Simons-Schr\"odinger equation (or any geometric flow arising from the abelian Higgs energy, for that matter) is that the phase space of the flow $\frE$ is manifestly non-linear, and consists of fairly complicated objects, see Appendix \ref{subapp:example} for examples. Another difficulty is in estimating the low-frequency part of the electromagnetic potential. For example, one typically writes the divergence-free part of the magnetic potential using the Biot-Savart law -- schematically, $\PP^{\df} A_x = \nabla^{-1} F_{12}$. However, the magnetic field is only $L^2$-integrable, so the most we can read-off is a, say, $\mathrm{BMO}$-bound for $A_x$. It turns out that both difficulties are manifestations of the ``topological'' nature of the problem. On the other hand, these in some sense stem from the ``low-frequency'' behaviour of a configuration, and thus should be rather inconsequential when attempting to settle (local) well-posedness. 

With these considerations in mind, our strategy for proving Theorem \ref{thm:gwp}, in a word, revolves around a smoothing scheme which separates out frequency scales in a geometric fashion. We will use this scheme to construct an auxiliary system of equations and gauge which serve to ``paralinearise'' the Chern-Simons-Schr\"odinger equation, and conduct the bulk of our analysis therein. In particular, the smoothing scheme will allow us to cleanly decompose the configuration into a ``low-frequency'' part, which is smooth and captures the ``non-linear'' or ``topological'' aspects of the configuration, and a ``high-frequency'' part, which is rough yet ``linearised'', and thus amenable to functional-analytic tools. Consequently, we can handle the difficulties of the problem, namely the novel ``low-frequency'' considerations discussed earlier and standard ``high-frequency'' ones which are more reminiscent of the analysis of similar models, in a fairly modular fashion. We outline the argument in more detail in Section \ref{subsec:overview}.

\begin{remark}
	We can also pose the Cauchy problem for the Chern-Simons-Schr\"odinger equation in other gauges. In fact, we conduct the bulk of our analysis in an auxiliary gauge with better regularity properties, albeit being more difficult to state. One can then read-off the regularity in another gauge by estimating the gauge transformation. For comparison, the transformations to the DeTurck, temporal, and Coulomb gauges are given respectively by
		\begin{align*}
			\chi^{\text{D}}
				&:= - (\partial_t - \Delta)^{-1} (A_t - \partial^j A_j),\\
			\chi^{\text{T}}
				&:= - \partial_t^{-1} A_t, \\
			\chi^{\text{C}}
				&:= - \Delta^{-1} \partial^j A_j. 
		\end{align*}
	The temporal gauge, which may seem natural due to revealing the Hamiltonian structure of the equations, is obtained by solving an ODE, so we can expect it to be a full derivative rougher compared to the Coulomb and DeTurck gauges, where one can use elliptic and parabolic smoothing. At the other end of the spectrum, the Coulomb gauge, while favourable at high-frequencies, is difficult to estimate for low-frequencies, since the electromagnetic field is not well-localised. As a concrete example, the electric potential under Coulomb gauge takes the form
		\[
			\Delta A_t
				= \partial^j \Bigl( \epsilon_{jk} \Im(\overline \phi \bfD^k \phi) + \partial_j F_{12} \Bigr).
		\]
	Heuristically, $|\phi| \to 1$ at spatial infinity, so the most one can read-off here is $A_t \in \dot H^1 (\R^2) \subseteq \mathrm{BMO} (\R^2)$. The lack of control over low-frequencies unfortunately causes several technical issues.
\end{remark}

As a relatively simple corollary of Theorem \ref{thm:gwp}, we can show orbital stability of self-dual vortices under the Chern-Simons-Schr\"odinger flow near self-dual coupling $|\lambda - 1| \ll 1$. Recall the following quantitative stability result adapted from the our previous work \cite[Corollary 1.2]{Zhao2026},

\begin{theorem}\label{thm:stability}
    Let $N \in \N_0$ be a non-negative integer, and suppose that $(A, \phi) \in \frE$ is a finite-energy configuration with topological degree $N$ such that
        \begin{equation}\label{eq:small-self-dual}
            \int_{\R^2} |(\bfD_1 + i \bfD_2) \phi|^2 + |\Gauss|^2 \, dx 
                \ll_N 1. 
        \end{equation}
    Then there exists a self-dual vortex $(\mathring A, \mathring{\vphantom{A}\phi}) \in \cM^N$ such that 
        \[
            \bigl\| (A, \phi) - (\mathring A, \mathring{\vphantom{A}\phi}) \bigr\|_{H^1 \times L^2}^2 
                \lesssim_{N} \int_{\R^2} |(\bfD_1 + i \bfD_2) \phi|^2 + |\Gauss|^2 \, dx. 
        \]
\end{theorem}

By conservation of the abelian Higgs energy and topological degree under the Chern-Simons-Schr\"odinger flow, the \textit{Bogomoln'yi identity} \cite{Bogomolny1976}, \cite[Proposition 4.1]{Zhao2026},
	\begin{equation} \label{eq:bogomolnyi-identity}\tag{Id}
		\cE_{\textup{AH}} [A, \phi] = \tfrac{\lambda - 1}{8} \int_{\R^2} (1 - |\phi|^2)^2 \, dx  + \tfrac12\int_{\R^2} |(\bfD_1 + i \bfD_2) \phi|^2 + |\Gauss|^2 \, dx  + \pi \deg[A, \phi],
	\end{equation}
and Theorem \ref{thm:stability}, we obtain

\begin{corollary}[Orbital stability of self-dual vortices]\label{cor:orbital}
    Let $N \in \N_0$ be a non-negative integer, $|1 - \lambda| \ll_N 1$, and suppose that $(A^{\mathrm{in}}, \phi^{\mathrm{in}}) \in \frE_\df$ is a finite-energy configuration in Coulomb gauge \eqref{coulomb-data} with topological degree $N$ satisfying \eqref{small-self-dual}. If $(A, \phi) \in C_t ([0, \infty) \to \frE)$ solves the Chern-Simons-Schr\"odinger equation in DeTurck gauge \eqref{CSS-DT-intro} with initial data $(A, \phi)_{|t = 0} = (A^{\mathrm{in}}, \phi^{\mathrm{in}})$, then
        \[
            \sup_{t \geq 0} \int_{\R^2} |(\bfD_1 + i \bfD_2) \phi|^2 + |\Gauss|^2 \, dx 
                \ll_N 1,
        \]
    and consequently, for each $t \geq 0$ there exists a self-dual vortex $(\mathring A(t), \mathring{\vphantom{A}\phi}(t)) \in \cM^N$ such that 
        \[
            \bigl\| (A(t), \phi(t)) - (\mathring A(t), \mathring{\vphantom{A}\phi}(t)) \bigr\|_{H^1 \times L^2}^2 
                \lesssim_{N}  \int_{\R^2} |(\bfD_1 + i \bfD_2) \phi|^2 + |\Gauss|^2 \, dx.
        \]
\end{corollary}

\begin{remark}
    Inspecting the statement of \cite[Corollary 1.2]{Zhao2026} more closely, one can upgrade the distance between the solution and self-dual vortex from $H^1 \times L^2$ to $H^1 \times H^1$, at the cost of the implicit constant depending on $\| A \|_{L^4}$. In DeTurck gauge \eqref{DeTurck}, this unfortunately means the constant can grow in $t$, though one can remove the growth by working in a different gauge. We plan on returning to this point in future work. 
\end{remark}

\begin{remark}
    The analogous results for negative topological degrees follow by reflection.  
\end{remark}

\subsection{Background and comparison to related models}\label{sec:intro-back}

The Chern-Simons-Schr\"odinger equation belongs to a broader collection of equations modeling the dynamics of vortices. We give a brief account of the subject, particularly Manton's original motivations for introducing \eqref{CSS} and comparing the Cauchy problems of related models. The interested reader should consult the monograph of Manton-Sutcliffe \cite[Chapter 7]{MantonSutcliffe2004} for a more in-depth survey.

\subsubsection*{Abelian Higgs model on $\R^2$ and Manton's conjectures}

In his seminal paper \cite{Bogomolny1976}, Bogomoln'yi derived the following identity for the abelian Higgs energy,
	\begin{equation} \tag{Id}
		\cE_{\textup{AH}} [A, \phi] = \cV_\lambda [\phi] + \cE_{\pm} [A, \phi] \pm \pi \deg[A, \phi],
	\end{equation}
where 
	\begin{align*}
		 \cV_\lambda [\phi]
		 	&:= \tfrac{\lambda - 1}{8} \int_{\R^2} (1 - |\phi|^2)^2 \, dx,\\
		\cE_\pm [A, \phi] 
			&:= \tfrac12 \int_{\R^2} |(\bfD_1 \pm i \bfD_2) \phi|^2 + \bigl|\tfrac12(1 - |\phi|^2) \mp F_{12}\bigr|^2 \, dx.
	\end{align*}
The identity implies that, in the case of \textit{self-dual} coupling, $\lambda = 1$, the minimisers within each topological class solve the first-order system of equations
	\begin{equation}\tag{B}\label{eq:bogomolnyi}
	\begin{split}
		(\bfD_1 \pm i \bfD_2) \phi
			&= 0,\\
		F_{12} 
			&= \pm \tfrac12(1 - |\phi|^2).
	\end{split}
	\end{equation}
Jaffe-Taubes \cite[Chapter III]{JaffeTaubes1980} proved that every critical point of the abelian Higgs energy in the self-dual case is a solution to \eqref{bogomolnyi}, and completely constructed the moduli space of solutions $\mathcal M^N$ as a smooth manifold diffeomorphic to $\C^{|N|}$ for each topological degree $N \in \Z$. Later Samols \cite{Samols1992} identified a natural K\"ahler structure $(\cM^N, \bfg, \mathbb J)$ on the moduli space. Much less is known about vortices away from self-dual coupling, $\lambda \neq 1$ -- we refer the reader to \cite[Chapter 7]{MantonSutcliffe2004} and \cite[Chapter 14.3]{SandierEtAl2007} and the references therein for a survey of the known results. 

While the self-dual vortices are not critical points of the abelian Higgs energy when $\lambda \neq 1$, they should nevertheless serve as reasonable approximate critical points near self-dual coupling, i.e. $|1 - \lambda| \ll 1$. Therefore, in this regime, one may use the moduli space of self-dual vortices $\cM^N$, with all its rich structure, as a proxy for the non-self-dual vortices. More precisely, we state the following conjecture, generally attributed to Manton, 

\begin{conjecture}
	In the low-energy regime and near self-dual coupling, i.e. $|1 - \lambda| \ll 1$, the dynamics of vortices with topological degree $N$ under a geometric flow arising from the abelian Higgs energy $\cE_{\mathrm{AH}} [A, \phi]$ can be modeled for long time by a corresponding flow arising from an effective energy on the moduli space of self-dual vortices $\cM^N$. 
\end{conjecture}

\begin{table}[!htbp]
\centering
\renewcommand{\arraystretch}{2}
\resizebox{\textwidth}{!}{%
    \begin{tabular}{|l|c|c|c|c|}
    \hline
    \textbf{Flow} & \textbf{Energy} & \textbf{Equation for $(A, \phi)$} & \textbf{Effective energy} & \textbf{Equation for $q$}
        \\ \hline
    Gradient & $\cE_{\mathrm{AH}} [A, \phi]$ & $\partial_t \binom{A}{\phi} = -(\nabla \cE_{\mathrm{AH}}) \binom{A}{\phi}$ & $\cV_\lambda [q]$ & $\bfdot q = -(\nabla \cV_\lambda) q$ 
        \\ \hline
    Schr\"odinger & $\cE_{\mathrm{AH}} [A, \phi]$ &$\partial_t \binom{A}{\phi} = -i (\nabla \cE_{\mathrm{AH}})\binom{A}{\phi}$ & $\cV_\lambda [q]$ &$\bfdot q = -\mathbb J (\nabla \cV_\lambda) q$ 
        \\ \hline
    Hyperbolic & $ \int_\R \cE_{\mathrm{AH}} [A, \phi] - \tfrac12 \| \partial_t (A, \phi) \|_{L^2}^2 \, dt$ &$\partial_t^2 \binom{A}{\phi} = -(\nabla \cE_{\mathrm{AH}}) \binom{A}{\phi}$ & $\int_\R \cV_\lambda [q] - |\bfdot q|_{\bfg}^2 \, dt$ & $\nabla_{\bfdot q} \bfdot q = - (\nabla \cV_\lambda) q$ 
        \\ \hline
    \end{tabular}
}

\vspace{1em}

    \caption{A summary of Manton's conjecture for the abelian Higgs model on $\R^2$, namely for the abelian Higgs gradient flow, the Chern-Simons-Schr\"odinger equation, and the hyperbolic abelian Higgs equation, the Euler-Lagrange equation for the abelian Higgs energy on Minkowski space $\R^{1 + 2}$. In view of the Bogomoln'yi identity \eqref{bogomolnyi-identity} and equations \eqref{bogomolnyi}, the energy reduces to an effective energy on $\cM^N$ dictating the dynamics of vortices $q(t) \in \cM^N$. When $\lambda = 1$, the dynamics are (effectively) trivial for the gradient and Schr\"odinger flows, and are governed by the geodesic equation on $\cM^N$ for the hyperbolic flow. }\label{table:conjecture}
\end{table} 

Stuart initiated the rigorous study of Manton's conjecture in the contexts of the hyperbolic abelian Higgs equations \cite{Stuart1994} and abelian Higgs gradient flow \cite{Stuart1996}. As for the Chern-Simons-Schr\"odinger equation, not much was known in terms of rigorous results prior to our work. Manton's original article \cite{Manton1997} gives a formal derivation of the vortex dynamics from a effective reduction of the Lagrangian \eqref{lagrange}. Formulated as a concrete mathematical conjecture, his argument suggests a low-energy solution to the Chern-Simons-Schr\"odinger equation \eqref{CSS} near self-dual coupling, i.e. $|1 - \lambda| \ll 1$, admits the decomposition\footnote{We are slightly abusive of notation here, since the vortex $q$ represents a gauge-equivalence class, and the addition operation does not interact well with gauge-transformation. Some care needs to be taken when interpretting the decomposition of a configuration $(A, \phi)$ into a vortex plus error. }
	\begin{align*}
		(A, \phi)(t) 
			&= q(t) + \text{error}, 
	\end{align*}
where $q : [0, T] \to \cM^N$, to leading-order, obeys the Hamiltonian flow of the reduced energy $\cV_\lambda [q]$,
	\begin{align*}
		\bfdot q
			&= - \mathbb J (\nabla \cV_\lambda) q + \text{error}.
	\end{align*}
There is some numerical evidence for the conjecture by Krusch-Sutcliffe \cite{KruschSutcliffe2006}. As far as rigorous results go, Corollary \ref{cor:orbital} shows the existence of a decomposition of the solution $(A, \phi)(t)$ into a self-dual vortex $q(t)$ plus error, albeit without any detailed information on the dynamics of $q(t)$. 

Restricting to the self-dual case, a stronger conjecture is expected in the mathematical folklore, namely asymptotic stability of the self-dual vortices under geometric flows arising from the self-dual abelian Higgs energy. For the abelian Higgs gradient flow, Demoulini-Stuart \cite{DemouliniStuart1997} in fact established an even stronger result, namely that $\cM^N$ forms a global attractor for the flow, though only for a class of smooth and localised data. The author \cite{Zhao2026} recently showed asymptotic stability in the natural energy class. In the same year, L\"uhrmann et. al. \cite{luhrmann2026a, luhrmann2026c, luhrmann2026b} resolved the asymptotic stability of a $1$-vortex, i.e. $q \in \cM^1$, under the hyperbolic abelian Higgs equation for a class of smooth, localised and equivariant perturbations. The analogous problem for the Chern-Simons-Schr\"odinger equation remains open, and seems difficult, as it amounts to proving global dispersive bounds for a variable-coefficient Schr\"odinger equation with quadratic non-linearity on $\R^2$. 

Of course, to avoid putting the cart before the horse, one should first address the Cauchy problem before studying dynamics of vortices. In the hyperbolic setting, one can pass $H^k$-well-posedness in Lorenz gauge \cite{Moncrief1980} to $H^k_\loc$-well-posedness via finite-speed of propagation. As for the gradient flow, Demoulini-Stuart \cite[Section 4]{DemouliniStuart1997} proved well-posedness in $C^{k,\alpha}$-spaces using parabolic smoothing. An alternate approach is to decompose the solution into a fixed smooth background and a perturbation, and then pose the Cauchy problem for the perturbation. This idea is applied in Burzlaff-Moncrief \cite{BurzlaffMoncrief1985} for the hyperbolic flow, the work of the author \cite[Appendix A]{Zhao2026} for the gradient flow, and Appendix \ref{app:picard} for \eqref{CSS-amp-schro}. 

\subsubsection*{Manton's Chern-Simons-Schr\"odinger equation on Riemann surfaces}

In addition to the Euclidean space $\R^2$, the abelian Higgs energy and Chern-Simons-Schr\"odinger equation can be defined on other Riemann surfaces \cite[Chapter 7.14]{MantonSutcliffe2004}, such as the torus $\TT^2$, the sphere $\SS^2$, and the hyperbolic plane $\HH^2$, and one can make analogous conjectures. On closed Riemann surfaces, Demoulini \cite{Demoulini2007} proved global existence for initial data $(A^{\mathrm{in}}, \phi^{\mathrm{in}}) \in (H^1 \times H^2) (\Sigma)$ -- the local result follows from the energy method, while global regularity follows from conservation of energy and a covariant Brezis-Gallouet inequality. This was followed by work of Demoulini-Stuart \cite{DemouliniStuart2009} on the adiabatic limit as $\lambda \to 1$, verifying a soft version of Manton's conjecture. On the hyperbolic plane $\HH^2$, Landoulsi-Shahshahani \cite{LandoulsiShahshahani2025} recently proved asymptotic stability of the $N$-vortex for equivariant $H^1$-perturbations in Coulomb gauge. Under their symmetry and gauge assumptions, the equation may be recast as a semi-linear variable-coefficient Schr\"odinger equation, without derivatives in the non-linearity, for the perturbation. They prove Strichartz estimates for the linear flow, which play a crucial role in the asymptotic stability analysis, and allow them to close a contraction mapping argument for local well-posedness \cite[Appendix A]{LandoulsiShahshahani2025}. 

\subsubsection*{Gross-Pitaevskii equation}

As alluded in the name, the gauged Ginzburg-Landau energy has a non-gauged analogue in the {Ginzburg-Landau energy}, 
	\[
		\cE_{\mathrm{GP}} [\phi] 
			:= \tfrac12 \int_{\R^2} |\nabla \phi|^2 + \tfrac{1}{2} (1 - |\phi|^2)^2 \, dx.
	\]
The critical points are also known as vortices; of particular interest is the $1$-equivariant vortex. The corresponding Schr\"odinger flow is the {Gross-Pitaevskii equation}, 
	\begin{equation}\label{eq:GP} \tag{GP}
		i \partial_t \phi + \Delta \phi = - (1 - |\phi|^2) \phi.
	\end{equation}	
The $1$-equivariant vortex is conjectured to be asymptotically-stable under \eqref{GP}. In terms of rigorous results, we highlight spectral stability \cite{WeinsteinXin1996}, orbital stability \cite{GravejatEtAl2022}, and linear stability \cite{CollotEtAl2025, LuhrmannEtAl2025}. Let us also mention the recent work \cite{PinoEtAl2025} justifying the Kirkhoff-Onsager dynamics of interacting $1$-vortices under \eqref{GP}. 

As for well-posedness, Gallo \cite{Gallo2004} proved local well-posedness for \eqref{GP} in Zhidkov spaces by energy arguments. Using Strichartz estimates, global well-posedness was established by Gerard \cite{Gerard2006,Gerard2008} for the energy space, and by 
Bethuel-Smets \cite{BethuelSmets2007} and Gallo \cite{Gallo2008} for $H^1$-perturbations of certain smooth backgrounds. In some sense, Theorem \ref{thm:gwp} is an analogue of Gerard's result for \eqref{CSS-amp-schro}, however the energy space for the abelian Higgs admits topologically non-trivial configurations, while the Ginzburg-Landau energy space does not. For example, consider $n$-equivariant scalar fields $\phi (r, \theta)= \rho(r) e^{in \theta}$, and compute 
	\begin{align*}
		| \nabla \phi|^2 
			&= |\partial_r \rho|^2 + \frac{n^2}{r^2} |\rho|^2, \\ 
		| \bfD \phi|^2 
			&= |\bfD_r \rho|^2 + \frac1{r^2} |\rho|^2 |A_\theta - n|^2.
	\end{align*}
If $1 - |\phi|^2 \in L^2 (\R^2)$, then, formally,
	\[
		\lim_{r \to \infty} |\rho(r)| = 1. 
	\]	
It follows that necessarily $n = 0$ if $\cE_{\text{GP}} [\phi] < \infty$, as otherwise the gradient would have an $r^{-1}$-tail arising from the angular derivative of the phase, leading to $\nabla \phi \not\in L^2 (\R^2)$. In the magnetic setting $\cE_{\text{AH}} [A, \phi] < \infty$, topologically non-trivial scalar fields, $n \neq 0$, are consistent with $\bfD \phi \in L^2 (\R^2)$, as the angular component of the potential can cancel out the tail when $A_\theta \to n$ in an appropriate sense.

\subsubsection*{Jackiw-Pi's Chern-Simons-Schr\"odinger equation} 

In \cite{JackiwPi1990, JackiwPi1990a, JackiwPi1991}, Jackiw-Pi introduced the model
	\begin{equation}\label{eq:JP-CSS}\tag{JP-CSS}
	\begin{split}
		i \bfD_t \phi + \bfD_1 \bfD_1 \phi + \bfD_2 \bfD_2 \phi
			&= - \lambda |\phi|^2 \phi, \\
		F_{0j}
			&= - \epsilon_{jk}\Im (\overline \phi \bfD^k \phi), \\
		F_{12}
				&= - \tfrac12 |\phi|^2,
	\end{split}
	\end{equation}
where $\lambda \in \R$ is a coupling constant. It is also referred in the mathematical literature as the Chern-Simons-Schr\"odinger equation; to avoid confusion with \eqref{CSS-amp-schro}, we will refer to \eqref{JP-CSS} as the Jackiw-Pi model. The following energy and charge are conserved under the flow, 
	\begin{align*}
		\cE_{\mathrm{JP}} [A, \phi]
			&:= \tfrac12 \int_{\R^2}  |\bfD_1 \phi|^2 + |\bfD_2 \phi|^2 - \tfrac{\lambda}{2} |\phi|^4 \, dx,\\
		\cM_{\mathrm{JP}} [A, \phi]
			&:= \tfrac{1}{2} \int_{\R^2} |\phi|^2  \, dx.
	\end{align*}
The Jackiw-Pi model obeys a scaling symmetry which preserves the charge. Hence, one may view \eqref{JP-CSS} as a gauge-covariant analogue of the mass-critical non-linear Schr\"odinger equation on $\R^2$. It is natural then to study the Cauchy problem for \eqref{JP-CSS} on the scale of Sobolev spaces $H^s (\R^2)$, where, in contrast to \eqref{CSS-amp-schro}, one has access to standard PDE and harmonic analysis tools off-the-shelf. We highlight two low-regularity local-wellposedness results; for large-data in $H^1 (\R^2)$ by Lim \cite{Lim2018}, and for small-data in $H^{0+}(\R^2)$ by Liu-Smith-Tataru \cite{LiuEtAl2014}. To give a brief summary of the proofs of each respective result, Lim argues by energy estimates and Strichartz estimates with derivative loss for the linear paradifferential equation, while Liu-Smith-Tataru use local smoothing estimates and null structure. We should also mention the result of Huh \cite{Huh2013}, who proved unconditional uniqueness of $H^1$-solutions to \eqref{JP-CSS} in Coulomb gauge. 

\subsection{Overview of the argument}\label{subsec:overview}

\subsubsection*{Geometric Littlewood-Paley theory}

In Section \ref{sec:littlewoodpaley}, we set-up the main functional-analytic tools by developing a Littlewood-Paley-type theory adapted to the space of finite-energy configurations. We have three goals in mind: (\texttt{i}) a covariant notion of frequency-localisation, (\texttt{ii}) a Besov-type characterisation of the topology, (\texttt{iii}) approximation of rough finite-energy configurations by smooth configurations. The analysis is based on the covariant heat equation, 
	 \begin{equation}\tag{P}\label{eq:intro-heat}
	\begin{split}
        \bigl( \bfD_s - \bfD^\ell \bfD_\ell \bigr) \phi
				&= 0 
                \hspace{6em}
                \raisebox{-.7\normalbaselineskip}[0pt][0pt]
                {%
                    on $\R^2 \times [0, \underline s]$.
                } \\
		F_{sj}
				&= \partial^\ell F_{\ell j}
	\end{split}
	 \end{equation}
To fix the gauge (modulo $s$-independent transformations), we take inspiration from work on the hyperbolic Yang-Mills equation, e.g. \cite{Oh2014a,Oh2015a,Gavrus2022}, which in turn draws from Tao's work on wave maps \cite{Tao2004}, and impose the \textit{caloric gauge},
	 \begin{equation}\tag{Cal} \label{eq:intro-cal}
		A_s = 0.
	 \end{equation}
     
By well-posedness theory, there is a one-to-one correspondence between the space of finite-energy configurations $\frE$ and their flows under the equation \eqref{intro-heat} in \eqref{intro-cal}, which we refer as $\frE$-caloric extensions and denote the space thereof by $\pzcE^1$. In a word, we want to justify the heuristic that the covariant Hessian of the caloric extension may be thought of as a Littlewood-Paley projection,
	 \[
	 	(s \bfD^{(2)}) \phi(s)
	 		\approx \text{projection of $\phi(s = 0)$ to frequency $s^{-1/2}$}.
	 \]
This heuristic will be encoded by a menagerie of parabolic smoothing estimates, which follow from abstract parabolic theory developed in Section \ref{subsec:abstract-parabolic}. Our covariant analysis is similar to the geometric Littlewood-Paley theory of Klainerman-Rodnianski \cite{KlainermanRodnianski2006}, though our approach introduces a notion of frequency-envelopes to capture the distribution of energy across frequency-scales. Our characterisation of the energy space takes inspiration from Tataru's \cite[Section 3]{Tataru2005} and Tao's \cite{Tao2009} works on wave maps at critical regularity (c.f. also Oh-Tataru's work on hyperbolic Yang-Mills \cite{OhTataru2019, OhTataru2020,OhTataru2021,OhTataru2022}), though since our setting is ``sub-critical" in some sense, we make a distinguished choice of ``low-frequency" scale, similar to the local caloric gauge in \cite{Oh2014a, Oh2015a, Gavrus2022}.

\subsubsection*{Continuation of smooth solutions}

In Section \ref{sec:high-LWP}, we outline the high-regularity well-posedness theory for the Chern-Simons-Schr\"odinger equation, which consists of (\texttt{i}) local existence of a solution arising from smooth initial data, (\texttt{ii}) higher-order energy estimates, (\texttt{iii}) a quantitative continuation criterion, and (\texttt{iv}) conservation of the abelian Higgs energy and topological degree. To this end, we construct covariant higher-order energy functionals which obey the following energy estimate\footnote{We write $\lessapprox$ and $\approx$ to indicate these are just caricatures of the actual estimates and definitions. },
	\[
		\frac{d}{dt} \cE^N [A, \phi]
			\lessapprox \| \phi \|_{L^\infty_x}^2 \cdot \cE^N [A, \phi].
	\]
We claim that smooth solutions exist as long as the $N$-th order covariant energy function remains bounded for $N \gg 1$. Consequently, the energy estimate furnishes the continuation criterion,  
	\begin{equation}\label{eq:intro-cont}\tag{$\star$}
		\int_0^T \| \phi(t) \|_{L^\infty_x}^2 \, dt < \infty.
	\end{equation}
If we can show that $\| \phi \|_{L^2_t L^\infty_x (I)} \leq 2$ on time-scales $I \subseteq [0, T]$ which depend only on the abelian Higgs energy and the Gauss tension field, then the conservation laws yield global regularity of smooth solutions.

At this point, an astute reader should object that the energy estimate in itself does not suffice for showing existence of a solution. Nonetheless, we can design a Picard iteration scheme which mimics the covariant analysis and establishes existence of a class of smooth solutions to the Chern-Simons-Schr\"odinger equation in DeTurck gauge \eqref{CSS-DT-intro} which may be continued as long as \eqref{intro-cont} holds. The convergence of the iteration follows from a technical adaptation of the energy method in the affine Sobolev space 
	\[
        (\mathring A, \mathring{\vphantom{A}\phi}) + H^N (\R^2) 
            := \Bigl\{ (A, \phi) \in H^N_{\loc} (\R^2) : (A, \phi) - (\mathring A, \mathring \phi) \in H^N (\R^2) \Bigr\},
    \]
where $(\mathring A, \mathring{\vphantom{A}\phi})$ is a fixed smooth background configuration. As the details are standard yet cumbersome, we defer them to Appendix \ref{app:picard}. 

\subsubsection*{The paradifferential formulation of the Chern-Simons-Schr\"odinger equation}

Following in the tradition of smoothing out a geometric dispersive equation by an associated heat flow (c.f. \cite{Tao2004,BejenaruEtAl2011,Oh2014a, Oh2015a,Gavrus2022,OhTataru2019, OhTataru2020, OhTataru2021, OhTataru2022} and the references therein), we extend a solution on $[0, T] \times \R^2$ to the Chern-Simons-Schr\"odinger equation to a configuration on $[0, T] \times \R^2 \times [0, \underline s]$ solving the system
	\begin{equation*}\tag{dP-CSS}\label{eq:dP-CSS-intro}
        \begin{split}
        \bigl( \bfD_s - \bfD^\ell \bfD_\ell \bigr) \phi
				&= 0 
                \hspace{12em}
                \raisebox{-.7\normalbaselineskip}[0pt][0pt]
                {%
                    on $[0, T] \times \R^2 \times [0, \underline s]$,
                } \\
		F_{s\mu}
				&= \partial^\ell F_{\ell \mu} \\
        \Bigl( i \bfD_t + \bfD^j \bfD_j + \tfrac\lambda2 (1 - |\phi|^2) \Bigr) \phi 
				&= 0
                \hspace{12em}
                \raisebox{-.7\normalbaselineskip}[0pt][0pt]
                {%
                    on $[0, T] \times \R^2 \times \{0\}.$
                }\\
		F_{tj}
				&= - \epsilon_{jk} \Big(\mathfrak \Im(\overline \phi \bfD^k \phi) + \partial^\ell F\indices{_\ell^k}\Big)
        \end{split}
     \end{equation*}
The first two equations in \eqref{dP-CSS-intro} form an augmentation of the covariant heat equation to configurations on $[0, T] \times \R^2 \times [0, \underline s]$ which induces a paradifferential-style decomposition of \eqref{CSS-amp-schro} in the form of a covariant Schr\"odinger equation for the covariant Hessian, 
	\begin{equation} \tag{$\star\star$}\label{eq:Ds-eq}
		\bigl( i \bfD_t + \bfD^j \bfD_j \bigr) (s \bfD^{(2)}) \phi (s)
			= \text{perturbative}, \qquad \text{on $[0, T] \times \R^2$}.
	\end{equation}
We dedicate Section \ref{sec:paradiff} to deriving this paradifferential-type equation for the covariant Hessian, along with companion equations, namely covariant parabolic equations for the electromagnetic field and tension fields. The equation \eqref{Ds-eq} has two advantages over the original Chern-Simons-Schr\"odinger equation \eqref{CSS-amp-schro} -- first, it separates out the main non-perturbative parts of the evolution, and second, it is an equation for a linearised object and thus can be analysed in conventional functional spaces, in contrast to the original dynamical variables $(A, \phi)$ which evolve in a non-linear phase space. 

\subsubsection*{The caloric-temporal-Coulomb gauge}

Rather than directly analysing the Chern-Simons-Schr\"odinger equation in DeTurck gauge, in Section \ref{sec:gauge}, we introduce an auxiliary gauge with better regularity properties, namely the caloric-temporal gauge (termed in deference to \cite{Oh2014a, Oh2015a}),
	\begin{equation}\label{eq:CTC-intro}\tag{C-T}
    \begin{split}
        A_s 
            &= 0
            \qquad \text{on $[0, T] \times \R^2 \times [0, \underline s]$,}
            \notag 
            \\
        {A_t}
            &= 0
            \qquad \text{on $[0, T] \times \R^2 \times \{\underline s\}$.}
    \end{split}
    \end{equation}
Morally, our goal is to show that the electromagnetic potential is ``smoother'' in this gauge compared to the covariant Hessian, and so one can think of the left-hand side of \eqref{Ds-eq} as having a ``low-high paraproduct'' structure. 

Thanks to parabolic smoothing alotted by the covariant heat equation, the geometric variables at $s = \underline s$ are smooth, and so it is favourable to propagate estimates in $t$ here. Accordingly, we recover the electromagnetic potential $A_\mu (t, s)$ by integrating the curvature $F_{s\mu} (s')$ on $s' \in [s, \underline s]$, followed by integrating the electric field $F_{tj} (t', \underline s)$ on $t' \in [0, t]$, and finally solving an elliptic-type equation for $A_j (0, \underline s)$. This last step is rather novel. First, we note that the curl-free part of the magnetic potential is transported by \eqref{intro-heat} in \eqref{intro-cal}, so in particular, the Coulomb gauge condition \eqref{coulomb-data} is propagated in $s$ at $t = 0$, 
	 \begin{equation}\tag{$\underline{\mathrm{Coul}}$}\label{eq:coulomb-sbar-intro}
		\partial^j {A_j} = 0 \qquad \text{on $\{0\} \times \R^2 \times \{\underline s\}$}.
	\end{equation}
From here, the traditional approach to writing $A_j (0, \underline s)$ is via the Biot-Savart law; schematically $A = \nabla^{-1} F_{12}$. However, the magnetic field is only $L^2$-integrable, so we can only read-off, say, a bounded mean oscillation bound for the low-frequency part of $A_j$. To remedy this deficiency, we derive a modified Biot-Savart law,
	\[
		A = \nabla^{-1} \omega(A, \phi) + \phi \cdot \bfD \phi,
	\]
where $\omega(A, \phi)$ is the vorticity, the density defining the topological degree. The vorticity is $L^1$-integrable, and so this representation has better behavior at low-frequencies. The derivation is found in Appendix \ref{app:conservation}. 

\subsubsection*{The linear paradifferential equation}

We will view the equation for the covariant Hessian \eqref{Ds-eq} from the lens of the abstract linear paradifferential equation, 
	\begin{equation}\tag{Paralin}\label{eq:paralinear-intro}
		(i \bfD_t + \bfD^j \bfD_j) u(s) 
			= \sfN(s). 
	\end{equation}
In Section \ref{sec:strichartz}, we show that the equation satisfies suitable energy estimates and Strichartz-type estimates. To prove the latter, we would like to regard \eqref{paralinear-intro} as a perturbation of the constant-coefficient Schr\"odinger equation. However, expanding out the left-hand side, the main enemy is the magnetic potential term $2i A^j \partial_j$, which incurs a loss of derivative if naively moved to the right. Nevertheless, the energy estimates we prove will imply that, roughly speaking, \eqref{paralinear-intro} propagates localisation to frequency $s^{-1/2}$, so, when estimated on the right in a time-interval $J$, 
	\[
		\bigl\| A(s) \nabla u(s) \bigr\|_{L^1_t L^2_x (J)}
			\lessapprox |J| \cdot s^{-\frac12} \cdot \| A (s) \|_{L^\infty_{t, x} (J)} \bigl\| u(s) \bigr\|_{L^\infty_t L^2_x (J)}.
	\]
Under caloric-temporal gauge \eqref{CTC-intro}, the contribution of the magnetic field is $O(1)$, so taking the time-scale $|J| \sim s^{1/2}$, this term may be regarded as perturbative, and we can show that $u$ satisfies a classical Strichartz-type estimate on $J$. For unit time-scale, we simply partition $[0, T]$ into $O(s^{-1/2})$-many sub-intervals of length $s^{1/2}$ and sum up these individual contributions, yielding a Strichartz estimate with ``derivative-loss",
	\[
		\| s^{\frac1{2p}} u (s) \|_{L^p_t L^q_x ([0, T])}
			\lessapprox \| u^{\mathrm{in}} (s) \|_{L^2_x} + \| \sfN (s) \|_{L^1_t L^2_x ([0, T])}. 
	\]

\subsubsection*{Covariant energy and Strichartz estimates}

With the prior machinery at our disposal, in Section \ref{sec:energy}, we establish global regularity of smooth solutions to the Chern-Simons-Schr\"odinger equation. In fact, we prove a more refined statement, namely that the energy profile of the initial data propagates under the flow, and that the solution obeys covariant Strichartz-type bounds. Phrased in the language of frequency envelopes, we show that there exists $T_* \ll 1$ depending only on the abelian Higgs energy and the Gauss tension field such that, if $T \leq T_*$ and $\ttc(s)$ is a frequency envelope for the initial data, the following bounds hold for \eqref{dP-CSS-intro}, 
\begin{itemize}
	\item \textup{(Control parameter bound).}
		\[
			\int_0^T \| \phi (t, s) \|_{L^\infty_x}^2 \, dt 
				\lesssim 1. 
		\]
	\item \textup{(Covariant energy estimate).}
		\[
			\sup_{t \in [0, T]} \bigl\| s^{-\frac12} (s \bfD^{(2)}) \phi(s) \bigr\|_{L^2_x}
				\lesssim \ttc(s).
		\]
	
	\item \textup{(Strichartz estimate).}
		\[
			\bigl\| s^{\frac1{2p} - \frac12} (s \bfD^{(2)}) \phi(s) \bigr\|_{L^p_t L^q_x} 
				\lesssim \ttc(s).
		\]
\end{itemize}
Since $T_*$ depends only on the conserved quantities, the analysis from Section \ref{sec:high-LWP} implies global existence of smooth solutions. To prove the bounds, we break gauge-covariance by passing to caloric-temporal gauge \eqref{CTC-intro}, and view the equation for the covariant Hessian \eqref{Ds-eq} as a perturbation of the linear paradifferential equation \eqref{paralinear-intro}. The task at hand reduces to showing that the right-hand side is indeed perturbative, i.e. 
	\[
		\bigl\| \text{R.H.S.\eqref{Ds-eq}} (s)\bigr\|_{L^2_x}
			\lessapprox \| \phi(s) \|_{L^\infty_x}^2 \cdot \ttc(s). 
	\]
This follows from the parabolic smoothing bounds developed in Sections \ref{sec:littlewoodpaley} and \ref{sec:tension}.

\subsubsection*{Difference estimates and uniqueness}

We approach the difference estimates for the Chern-Simons-Schr\"odinger equation from the persepective of quasi-linear hyperbolic equations, aiming to propagate a weaker $L^2$-type distance functional defined in Section \ref{subsec:weak},
	\[
		\dist_{\pzcE^0} \bigl((A, \phi) , (A', \phi') \bigr)
			\approx \bigl\| \delta (s^\frac12 \bfD) \phi \bigr\|_{L^2_{\frac{ds}{s}} L^2_x}.
	\]
Interpolating a bit with energy estimates in the high-frequency regime, we can show that, roughly speaking, the non-linear interactions are log-Lipschitz with respect to the $L^2$-type distance functional for finite-energy solutions in caloric-temporal gauge. Putting these together in Section \ref{sec:difference}, we (morally, at least) arrive at the differential inequality,
	\[
		\frac{d}{dt} |\dist_{\pzcE^0}|^2 \lessapprox \langle t \rangle \cdot \log\Bigl( \frac{1}{|\dist_{\pzcE^0}|^2} \Bigr) \cdot |\dist_{\pzcE^0}|^2,
	\]
where the implicit constant depends only on the Gauss tension field and the abelian Higgs energy. A non-linear Gronwall inequality implies a weak H\"older-type continuity of the data-to-solution map in caloric-temporal gauge. Notably, the weak distance functional propagates at the level of $\frE$-regularity, so this argument also yields unconditional uniqueness in $C_t ([0, T] \to \frE)$. 

\subsubsection*{Continuity of the flow map}

In Section \ref{sec:cts}, we conclude the proof of Theorem \ref{thm:gwp} by constructing solutions to the Chern-Simons-Schr\"odinger equation in DeTurck gauge arising from generic initial data in the space of finite-energy configurations, and showing continuity of the ensuing data-to-solution map with respect to the energy topology. Formally, there is a one-to-one correspondence between solutions in DeTurck gauge and caloric-temporal gauge via an appropriate gauge transformation, so we argue by showing the analogue of Theorem \ref{thm:gwp} for \eqref{CTC}, and then establishing regularity of the gauge-transformation, i.e. the maps
	\begin{align*}
		\begin{matrix}
			\frE_\df \\
			\text{initial data in \eqref{coulomb-data}}
		\end{matrix}
		\qquad\overset{\text{data-to-sol.}}{\longrightarrow}\qquad
		\begin{matrix}
			C_t ([0, T] \to \frE)\\
			\textup{solution in \eqref{CTC}} 
		\end{matrix}
		\qquad\overset{\text{gauge trans.}}{\longrightarrow}\qquad
		\begin{matrix}
			C_t ([0, T] \to \frE)\\
			\textup{solution in \eqref{DeTurck}}.
		\end{matrix}
	\end{align*}
are well-defined and continuous. To accomplish the former, we implement the well-posedness roadmap of Ifrim-Tataru \cite{IfrimTataru2022}, which works by interpolating between the ``compactness'' alotted by frequency envelope bounds and and ``uniqueness'' from the weak difference bounds. Moreover, this argument produces rough solutions as limits of smooth solutions, so, by continuity of the abelian Higgs energy and topological degree with respect to the $\frE$-topology, these rough solutions retain the same conservation laws. 

\subsection{Outline of the article}

The following two sections are dedicated to developing functional analytic tools: 
	\begin{itemize}
		\item In Section \ref{sec:preliminaries}, we collect the basic notation, covariant identities, magnetic Sobolev and interpolation inequalities, and develop a cook-book of abstract parabolic estimates.
		
		\item In Section \ref{sec:littlewoodpaley}, we develop the geometric Littlewood-Paley theory, which consists of a covariant notion of frequency-localisation, a Besov-type characterisation of the topology, and a smooth approximation scheme. 
	\end{itemize}
The bulk of the analysis of the Chern-Simons-Schr\"odinger equation will be conducted on an auxiliary system and choice of gauge, which is detailed in the following sections: 
	\begin{itemize}
		\item In Section \ref{sec:paradiff}, we introduce the paradifferential formulation of the Chern-Simons-Schr\"odinger equation and derive covariant equations of motion induced by the system.
		
		\item In Section \ref{sec:tension}, we estimate the tension fields and electromagnetic field using the parabolic equations of motion derived in the previous section. 
		
		\item In Section \ref{sec:gauge}, we introduce the caloric-temporal gauge and derive a system of transport and elliptic equations for reconstructing and then estimating the electromagnetic potential. 
		
		\item In Section \ref{sec:strichartz}, we introduce the linear paradifferential equation, for which we prove energy estimates and Strichartz estimates (with derivative-loss). 
	\end{itemize}
With the preceding analysis at hand, the proof of global well-posedness, Theorem \ref{thm:gwp}, follows a more-or-less standard argument:
	\begin{itemize}
		\item In Appendix \ref{app:picard}, we prove the existence of smooth solutions via a Picard iteration argument in high-regularity affine Sobolev spaces. 

		\item In Appendix \ref{app:conservation}, we prove conservation of energy and topological degree by deriving and integrating density-flux identities.

		\item In Section \ref{sec:high-LWP}, we prove higher-order energy estimates and deduce a continuation criterion for smooth solutions. 
		
		\item In Section \ref{sec:energy}, we close the control parameter bound and propagate frequency envelope and Strichartz bounds for smooth solutions on unit time-scale depending on the abelian Higgs energy and Gauss tension field. By the continuation criterion and conservation laws, this implies global regularity.
		
		\item In Section \ref{sec:difference}, we prove difference estimates with respect to the weak distance functional, and deduce unconditional uniqueness.  
		
		\item In Section \ref{sec:cts}, we construct rough solutions as limits of smooth solutions and proving continuity of the ensuing data-to-solution map in caloric-temporal gauge. We pass these results back to DeTurck gauge by estimating the gauge-transformation, completing the proof of Theorem \ref{thm:gwp}. 
	\end{itemize}

\subsection*{Acknowledgements}

The author thanks Sung-Jin Oh for his guidance and encouraging discussions over the course of conducting this work. Part of the work was done under the generous hospitality of the Korean Institute for Advanced Study and the Simons Center for Geometry and Physics. The author was partially supported by the NSF CAREER Grant (\texttt{NSF-DMS-1945615}) and NSF DMS Analysis Grant (\texttt{NSF-DMS-2452760}).

\subsection*{Disclosure of AI usage}

\begin{figure}[ht]
    \begin{center}
    \begin{tikzpicture}[scale=3, thick]
  
    \definecolor{lightyellow}{RGB}{255,250,205}

        \fill[lightyellow!80!white] (0,0) rectangle (2,2);
        \draw[->] (-0.05,0) -- (2.2,0) node[right] {$s$};
        \draw[red, ->] (0,-0.05) -- (0,2.2) node[above] {$t$} node[left] at (0, 1.1) {$\CSS = 0$} node[left] at (0, 0.9) {$\Amp = 0$}; 
        \node[red, below] at (0,-0.05) {$s = 0$};

        \draw[blue, very thick] (2,0) -- (2,2);
        \node[blue, below] at (2,-0.05) {$s = \underline s$};

        \node[yellow!50!orange, above] at (1,1) {(dP-$\phi$)-(dP-$A$)};

    \end{tikzpicture}
    \caption{The system \eqref{dP-CSS} arises from imposing a configuration $(\overline A, \overline\phi)$ on $[0, T] \times \R^2$ solving the Chern-Simons-Schr\"odinger system \eqref{CSS-amp-schro} as heat-temporal initial data and then evolving in $s$ under the dynamic covariant heat equation \eqref{d-covariant-heat}-\eqref{d-covariant-heat2}.}\label{fig:dP-CSS}
    \end{center}
\end{figure}
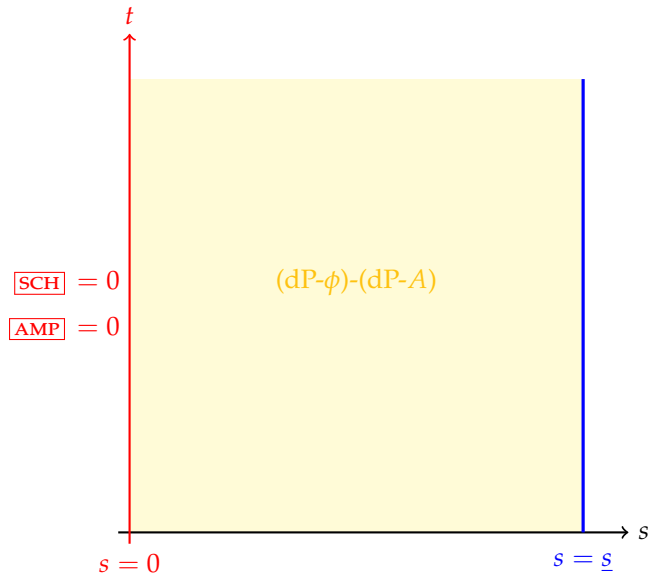

\begin{figure}[!htbp]
    \begin{center}
    \begin{tikzpicture}[scale=3, thick]
  
    \definecolor{lightyellow}{RGB}{255,250,205}

        \fill[lightyellow!80!white] (0,0) rectangle (2,2);
        \draw[->] (-0.05,0) -- (2.2,0) node[right] {$s$};
        \draw[->] (0,-0.05) -- (0,2.2) node[above] {$t$};
        \node[red, below] at (0,-0.05) {$s = 0$};

        \draw[blue, very thick] (2,0) -- (2,2);
        \node[blue, right] at (2,1) {$\underline A_t = 0$};
        \node[blue, below] at (2,-0.05) {$s = \underline s$};

        \fill[green!60!black] (2,0) circle (1pt);
        \node[green!60!black, above left] at (1.95,0.05) {$\partial^\ell A_\ell = 0$};

        \draw[yellow!80!orange, very thick] (0.5,0.5) -- (2,0.5);
        \fill[red] (0.5,0.5) circle (1pt);
        \node[red, above] at (0.45,0.5) {$A_\mu (t,s)$};
        \fill[blue] (2,0.5) circle (1pt);
        \node[blue, right] at (2.05,0.5) {$\underline A_j (t)$};
        \node[yellow!50!orange, above] at (1,1) {$A_s = 0$};
    
    \end{tikzpicture}
    \caption{The caloric-temporal gauge \eqref{CTC} conists of the caloric gauge $A_s = 0$ for all $t, s$, and the temporal gauge $A_t = 0$ at $s = \underline s$. We recover the electromagnetic potential $A_\mu (t, s)$ by integrating the curvature component $F_{s\mu}$, which obeys the (dynamic) Yang-Mills heat flow \eqref{d-covariant-heat2}, on $[s, \underline s]$, and integrating the electric field $\underline F_{tj}$ on $[0,t]$. Since the divergence of $A_x$ is transport in $s$ under caloric gauge, we will view \eqref{coulomb-data} as imposed at $s = \underline s$ and $t = 0$. }\label{fig:CTC}
    \end{center}
\end{figure}
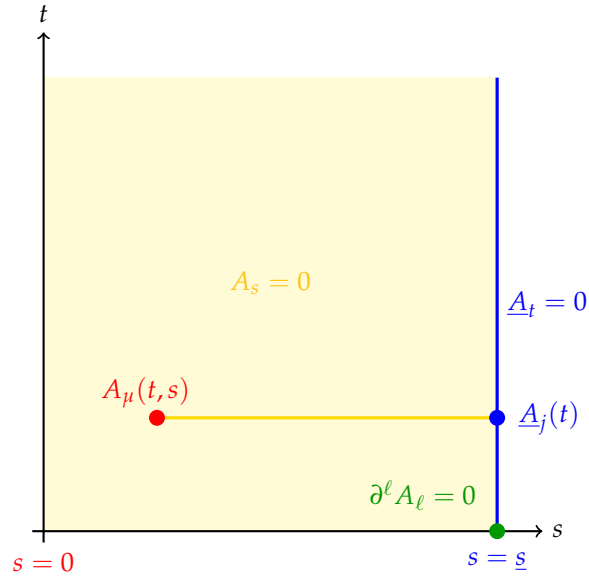

Figures \ref{fig:dP-CSS} and \ref{fig:CTC} were converted from the author's hand-drawings by ChatGPT. The connection of Lemma \ref{lem:gronwall} to Osgood's uniqueness criterion was pointed out by Gemini 3.1 Pro. 

\section{Preliminaries}\label{sec:preliminaries}

\subsection{Notation}\label{subsec:prelim-notation}

\subsubsection*{Indices}

We employ the usual Einstein summation convention of summing repeated upper and lower indices. The Levi-Civita symbol $\epsilon_{jk}$ is the anti-symmetric tensor
    \[
        \epsilon_{jk} 
            := 
            \begin{cases}
                +1 & \text{if $(j, k) = (1, 2)$}, \\
                -1 & \text{if $(j, k) = (2, 1)$},\\
                0
                    & \text{otherwise.}
            \end{cases}
    \]
Latin indices, e.g. $j, k, \ell$, will be reserved for spatial coordinates $(x^1, x^2)$, Greek indices, e.g. $\alpha, \beta, \mu, \nu$, for space-time coordinates $(t, x^1, x^2)$, and bold indices, e.g. $\bfa, \bfb, \bfc$, for space-Schr\"odinger-heat-time coordinates $(t, x^1, x^2, s)$. The $x$-subscript, e.g. $A_x$ and $F_{t x}$, refers to an arbitrary spatial index, i.e. $x \in \{1, 2\}$. 

\subsubsection*{Derivatives}

We will use $\bfD^{(n)}$ and $\nabla^{(n)}$ for $n$-fold iterated spatial derivatives. Similarly, we write $\bfD^{(\leq n)}$ and $\nabla^{(\leq n)}$ for iterated spatial derivatives up to order $n$. More precisely, in equations, this will be a placeholder for any $n$-fold or up-to-$n$-fold combination of derivatives, e.g.
    \begin{align*}
        \bfD^{(n)} 
            &= \bfD_{j_1} \cdots \bfD_{j_n} 
            \qquad 
            \text{for some indices $j_1, \dots, j_n \in \{1, 2\}$}, \\
        \bfD^{(\leq n)}
            &= \bfD_{j_1} \cdots \bfD_{j_k}, \qquad 
            \text{for some indices $j_1, \dots, j_k \in \{1, 2\}$ and $k \leq n$},
    \end{align*}
while the norms of these expressions will represent the sum total of all such combinations, e.g. 
    \begin{align*}
        \bigl\| \bfD^{(n)} \Phi \bigr\|_{L^p_x} 
            &:= \sum_{j_1, \dots, j_n \in \{1, 2\} } \bigl\| \bfD_{j_1} \cdots \bfD_{j_n} \Phi \bigr\|_{L^p_x}, \\
        \bigl\| \bfD^{(\leq n)} \Phi \bigr\|_{L^p_x} 
            &:= \sum_{\substack{j_1, \dots, j_k \in \{1, 2\} \\ k \leq n} } \bigl\| \bfD_{j_1} \cdots \bfD_{j_k} \Phi \bigr\|_{L^p_x}.
    \end{align*}
The sans-serif font $\sfD$ will sometimes be used as a placeholder for an arbitrary derivative, either magnetic $\sfD = \bfD$ or non-magnetic $\sfD = \nabla$. We will refer to $\bfD^{(2)}$ as the covariant Hessian, and $\bfD^j \bfD_j$ as the covariant Laplacian. 

\subsubsection*{Multi-linear forms}

We use dots as catch-all notation for multi-linear expressions which are gauge-covariant or gauge-invariant, depending on the context. For example, 
    \[
        \Phi \cdot \Psi \cdot \Gamma 
            = a \overline \Phi \Psi \Gamma + b\Phi \overline\Psi \Gamma + c\Phi \Psi \overline\Gamma , \qquad \text{for some $a, b, c \in \C$}
    \]
The reader should simply keep in mind that these expressions behave well with respect to the product rule and its covariant analogues. 

\subsubsection*{Differences}

The difference operator $\delta$ maps a field to a difference of its inputs. For example, if $\Phi$ and $\Phi'$ are scalar fields, and $A$ and $A'$ are connection $1$-forms, then 
    \begin{align*}
        \delta \Phi 
            &= \Phi - \Phi', \\
        \delta \bfD \Phi 
            &= (\nabla - i A) \Phi - (\nabla - i A') \Phi', \\
        \delta F_{\mu\nu} 
            &= F_{\mu\nu} - F_{\mu\nu}'
    \end{align*}
The difference operator commutes with the usual derivatives $[\delta, \partial_\bfa] = 0$ (in fact, any translation-invariant linear operator), and satisfies the product rule, 
    \[
        \delta (\Phi \cdot \Psi) 
            = \delta \Phi \cdot \Psi + \Phi \cdot \delta \Psi, 
    \]
where the ``undifferenced'' factors are understood to represent either input, e.g. the $\Psi$ in $\delta \Phi \cdot \Psi$ is actually a placeholder for both $\Psi$ and $\Psi'$.

\subsubsection*{Fundamental theorem of calculus}

We write $(\partial_t)^{-1}$ for integration with respect to $dt$ from $t = 0$, and $(\partial_s)^{-1}$ for integration with respect to $ds$ from $s = \underline s$, and $(s \partial_s)^{-1}$ for integration with respect to $ds/s$ from $s = \underline s$, 
    \begin{align*}
        (\partial_t)^{-1} \Phi (t) 
            &:= \int_0^t \Phi(t') \, dt', \\
        (\partial_s)^{-1} \Phi (s) 
            &:= \int_s^{\underline s} \Phi(s') \, ds', \\
        (s \partial_s)^{-1} \Phi(s) 
            &:= \int_s^{\underline s} \Phi(s') \, \frac{ds'}{s'}. 
    \end{align*}

\subsubsection*{Linear heat equation}

For $\Phi_0 \in L^1_\loc (\R^2)$, define the \textit{linear heat propagator} by
    \begin{align*}
        (e^{s \Delta} \Phi_0) (x) 
            &:= \Bigl(\tfrac{1}{4 \pi s}  e^{-\frac{|\bullet|^2}{4s}} * \Phi_0\Bigr)(x),
    \end{align*}
which is the unique distributional solution to the homogeneous heat equation $(\partial_s - \Delta) \Phi = 0$ with initial data $\Phi_{|s = 0} = \Phi_0$. Using the linear heat propagator and Duhamel's principle, we can solve the heat equation with forcing $\sfN \in L^1_{s, x\loc} (\R^2 \times I)$, defining the \textit{(heat) Duhamel integral} by  
    \begin{align*}
        (\partial_s - \Delta)^{-1} \sfN (s) 
            &:= \int_0^s e^{(s - s') \Delta} \sfN(s') \, ds'.
    \end{align*}
This is the unique distributional solution to the inhomogeneous heat equation $(\partial_s - \Delta) \Phi = \sfN$ with zero initial data $\Phi_{|s = 0} = 0$. 

\subsubsection*{Helmholtz decomposition}

Given a vector field $A \in L^p (\R^2 \to \R^2)$ for $1 < p < \infty$, we can decompose
    \[
        A 
             = \PP^{\df} A + \PP^{\cf} A,
    \]
where $\PP^{\df}$ is the projection to divergence-free vector fields, and $\PP^{\cf}$ is the projection to curl-free vector fields, defined respectively by 
    \begin{align*}
        \PP^{\df}_j A
            &:= A_j - \partial_j \Delta^{-1} \partial^\ell A_\ell,\\
        \PP^{\cf}_j A  
            &:= A_j + \epsilon_{jk} \partial^k \Delta^{-1} F_{12}.
    \end{align*}
Rearranging the Helmholtz decomposition yields
    \begin{align*}
        \PP^{\df}_j A
            &= - \epsilon_{jk} \partial^k \Delta^{-1} F_{12},\\
        \PP^{\cf}_j A
            &= \partial_j \Delta^{-1} \partial^\ell A_\ell.
    \end{align*}
The former is commonly known as the Biot-Savart law. 

\subsubsection*{Asymptotic notation}

The notations $X \lesssim Y$ denote the inequality $X \leq C Y$ for an implicit constant $C > 0$. Moreover, $X \sim Y$ is an abbreviation for $X \lesssim Y \lesssim X$. We will write $X \ll Y$ for $CX \leq Y$ and conversely $X \gg Y$ for $X \geq C Y$ for a sufficiently large constant $C > 0$. We indicate the dependence of implicit constants on parameters by subscripts, e.g. $X \lesssim_s Y$ denotes $X \lesssim C(s) Y$.

\subsection{Covariant identities and inequalities}

Throughout, $A = A_\bfa \, dx^\bfa$ and $A' = A_\bfa' \, dx^\bfa$ will be connection $1$-formss on a (subset of) Euclidean space, e.g. $\R^2$, $[0, T] \times \R^2$, and $[0, T] \times \R^2 \times [0, \underline s]$, and $\Phi, \Psi$ will be complex scalar fields with the appropriate domain. 

\begin{lemma}[Gauge covariance]
    Let $\chi$ be a sufficiently regular real-scalar field, and set 
        \begin{align*}
            \widetilde \Phi 
                &:= e^{i \chi} \Phi,\\
            \widetilde A 
                &:= A + d \chi.
        \end{align*}
    Then we have
    \begin{enumerate}
        \item gauge-invariance of the magnetic field, 
            \begin{equation}
                F_{\bfa \bfb} = \widetilde F_{\bfa \bfb},
            \end{equation}

        \item gauge-covariance of the covariant derivative 
            \begin{equation}
                \bfD_\bfa \Phi 
                    = e^{-i \chi} \widetilde \bfD_\bfa \widetilde \Phi. 
            \end{equation}
    \end{enumerate}
\end{lemma}

\begin{proof}
    \leavevmode
    \begin{enumerate}
        \item Expanding the right-hand side and using the famous identity $\partial_\bfa \partial_\bfb = \partial_\bfb \partial_\bfa$ yields
            \begin{align*}
                \widetilde F_{\bfa \bfb} 
                    &= \partial_\bfa \widetilde A_\bfb - \partial_\bfb \widetilde A_\bfa \\
                    &= \partial_\bfa (A_\bfb + \partial_\bfb \chi) - \partial_\bfb (A_\bfa + \partial_\bfa \chi) \\
                    &= (\partial_\bfa A_\bfb - \partial_\bfb A_\bfa) + (\partial_\bfa \partial_\bfb \chi - \partial_\bfb \partial_\bfa \chi) = F_{\bfa \bfb}. 
            \end{align*}

        \item Expanding the right-hand side, it follows from the product rule that
            \begin{align*}
                \widetilde \bfD_\bfa \widetilde \Phi 
                    &= (\partial_\bfa - i \widetilde A_\bfa) (e^{i \chi} \Phi) \\
                    &= e^{i \chi} \bigl((\partial_\bfa - i \widetilde A_\bfa) \Phi + i \partial_\bfa \chi \, \Phi\bigr) \\
                    &= e^{i \chi} \bigl( (\partial_\bfa - i A_\bfa) \Phi + i(\partial_\bfa \chi - \partial_\bfa \chi) \Phi\bigr)= e^{i \chi} \bfD_\bfa \Phi.
            \end{align*}
    \end{enumerate}
\end{proof}

\begin{lemma}
    \leavevmode
    \begin{enumerate}
        \item \textup{(Covariant derivative commutator identity).} 
            \begin{equation}\label{eq:commute}
                [\bfD_\bfa, \bfD_\bfb] 
                    = - i F_{\bfa \bfb}. 
            \end{equation}

        \item \textup{(Bianchi identity).}
            \begin{equation}\label{eq:bianchi}
                \partial_\bfa F_{\bfb \bfc} 
                    = \partial_\bfb F_{\bfa \bfc} - \partial_{\bfc} F_{\bfa \bfb}. 
            \end{equation}

        \item \textup{(Covariant Laplacian commutator identity).}
            \begin{equation}\label{eq:laplace-commute}
                [\bfD_\bfa, \bfD^j \bfD_j]
                    = -2i F\indices{_\bfa^j} \bfD_j - i \partial^j F_{\bfa j} .
            \end{equation}

        \item \textup{(Covariant product rule).}
            \begin{align}
            \partial_\bfa (\overline \Phi \Psi) 
                &= \overline{\bfD_\bfa \Phi} \Psi + \overline \Phi \bfD_\bfa \Psi. \label{eq:product}
            \end{align}
    \end{enumerate}
    
\end{lemma}

\begin{proof}
    \leavevmode
    \begin{enumerate}
        \item Applying the product rule and commuting mixed partial derivatives, 
            \begin{align*}
                [\bfD_\bfa, \bfD_\bfb] \Phi
                    &= (\partial_\bfa - i A_\bfa) (\partial_\bfb - i A_\bfb)\Phi - (\partial_\bfb - i A_\bfb)(\partial_\bfa - i A_\bfa)\Phi \\
                    &= - i \left(\partial_\bfa (A_\bfb \Phi)  + A_\bfa \partial_\bfb \Phi - \partial_\bfb(A_\bfa \Phi) - A_\bfb \partial_\bfa \Phi \right) \\
                    &= - i (\partial_\bfa A_\bfb - \partial_\bfb A_\bfa) \Phi = - i F_{\bfa \bfb} \Phi.
            \end{align*}
        \item Commuting mixed partial derivatives, 
            \begin{align*}
            \partial_\bfa F_{\bfb \bfc} 
                &= \partial_\bfa (\partial_\bfb A_\bfc - \partial_\bfc A_\bfb) + (\partial_\bfc \partial_\bfb A_\bfa - \partial_\bfb \partial_\bfc A_\bfa)\\
                &= \partial_\bfb (\partial_\bfa A_\bfc - \partial_\bfc A_\bfa) - \partial_\bfc (\partial_\bfa A_\bfb - \partial_\bfb A_\bfa) = \partial_\bfb F_{\bfa \bfc} - \partial_{\bfc} F_{\bfa \bfb}.
        \end{align*}

        \item Commuting covariant derivatives a l\'a \eqref{commute} gives 
        \begin{align*}
            \bfD_\bfa \bfD^j \bfD_j \Phi
                &= \bfD^j \bfD_\bfa \bfD_j \Phi - i F\indices{_\bfa^j} \bfD_j \Phi \\
                &= \bfD^j \bfD_j \bfD_\bfa \Phi - i \bfD^j \big( F\indices{_\bfa_j} \Phi\big)- i F\indices{_\bfa^j} \bfD_j \Phi \\
                &= \bfD^j \bfD_j \bfD_\bfa \Phi - 2i F\indices{_\bfa^j} \bfD_j \Phi - i \partial^j F_{\bfa j} \Phi,
        \end{align*}
    which, upon rearranging terms, proves the identity. 

        \item By the usual product rule and remarking $\overline{(-i A_\bfa \Phi)} \Psi + \overline \Phi (-i A_\bfa \Psi) = 0$, we compute
        \begin{align*}
            \partial_\bfa (\overline \Phi \Psi) 
                &= \overline{\partial_\bfa \Phi} \Psi + \overline \Phi \partial_\bfa \Psi \\
                &=  \overline{\partial_\bfa \Phi} \Psi + \overline \Phi \partial_\bfa \Psi + \left( \overline{(-i A_\bfa \Phi)} \Psi + \overline \Phi (-i A_\bfa \Psi) \right) = \overline{\bfD_\bfa \Phi} \Psi + \overline \Phi \bfD_\bfa \Psi,
        \end{align*}
    which gives the result. 
    \end{enumerate}
\end{proof}

\begin{lemma}
    The following inequalities holds in distribution:
    \begin{enumerate}
        \item \textup{(Diamagnetic inequality).}
            \begin{equation}
                |\partial_\bfa |\Phi|| \leq |\bfD_\bfa \Phi|.\label{eq:diamagnetic} 
            \end{equation}
        
        \item \textup{(Parabolic Bochner inequality).}
            \begin{equation}\label{eq:bochner-heat}
            (\partial_s - \Delta) |\Phi| \leq |(\bfD_s - \bfD^j \bfD_j) \Phi|.
            \end{equation}
    \end{enumerate} 
        
\end{lemma}

\begin{proof}
    \leavevmode
    \begin{enumerate}
        \item Writing $|\Phi|^2 = \overline \Phi \Phi$, we compute using the product rule \eqref{product} and Cauchy-Schwarz 
        \begin{align*}
            |\partial_\bfa |\Phi||
                &= \frac{1}{\sqrt{\overline \Phi \Phi}} |\Re(\overline \Phi \bfD_\bfa \Phi)|\\
                &\leq |\bfD_\bfa \Phi|, 
        \end{align*}
    as desired.

        \item Writing $|\Phi|^2 = \overline\Phi \Phi$, we compute using the product rule \eqref{product},
        \begin{align*}
            \partial_s |\Phi|
                &= \frac{1}{\sqrt{\overline \Phi \Phi}} \Re(\overline \Phi \bfD_s \Phi), \\
            \Delta |\Phi|
                &= \frac{1}{\sqrt{\overline \Phi \Phi}} \left(\Re(\overline \Phi \bfD^j \bfD_j \Phi) + \overline{\bfD^j \Phi} \bfD_j \Phi - \frac{\Re(\overline \Phi \bfD^j \Phi) \Re(\overline \Phi \bfD_j \Phi)}{\overline \Phi \Phi} \right).
        \end{align*}
        The main observation is that last two terms in the identity for $\Delta |\Phi|$ combined are non-negative by Cauchy-Schwarz. Upon another application of Cauchy-Schwarz, we conclude that 
            \begin{align*}
                (\partial_s - \Delta) |\Phi| 
                    &= \frac{1}{\sqrt{\overline \Phi \Phi}} \left( \Re\left(\overline \Phi (\bfD_s - \bfD^j \bfD_j) \Phi\right) \right) - \frac{1}{\sqrt{\overline\Phi \Phi}} \left( \overline{\bfD^j \Phi} \bfD_j \Phi - \frac{\Re(\overline \Phi \bfD^j \Phi) \Re(\overline \Phi \bfD_j \Phi)}{\overline \Phi \Phi}  \right)\\
                    &\leq \frac{1}{\sqrt{\overline \Phi \Phi}} \left( \Re\left(\overline \Phi (\bfD_s - \bfD^j \bfD_j) \Phi\right) \right) \\
                    &\leq |(\bfD_s - \bfD^j \bfD_j) \Phi|,
            \end{align*}
        as desired. 
    \end{enumerate}
   
\end{proof}

\begin{remark}
    Strictly speaking, the proofs above are merely formal arguments, ignoring technical issues when $\Phi(x) = 0$. One can easily make the arguments rigorous by a standard regularisation, replacing $|\Phi|$ with $\sqrt{|\Phi|^2 + \epsilon}$, testing against non-negative test functions, and taking $\epsilon \to 0$.
\end{remark}

\begin{corollary}[Parabolic maximum principle]
    We have the following pointwise inequality
        \begin{equation}\label{eq:maximum-principle}
            |\Phi(x,s)| 
                \leq e^{s \Delta} |\Phi_{|s = 0} (x)| + (\partial_s - \Delta)^{-1} \big|\big(\bfD_s - \bfD^j \bfD_j \big) \Phi (x)\big| .
        \end{equation}
\end{corollary}

\begin{proof}
    This follows from the parabolic Bochner inequality \eqref{bochner-heat}, Duhamel's formula, and the fact that the heat kernel is non-negative. 
\end{proof}

\begin{lemma}[Difference identities]
    \leavevmode
    \begin{enumerate}
        \item covariant derivative commutator,
            \begin{equation}\label{eq:commute-diff}
                [\delta, \bfD_\bfa]
                    = \delta A_\bfa,
            \end{equation}
        
        \item covariant Laplacian commutator,
            \begin{equation}\label{eq:commute-diff-laplace}
                [\delta, \bfD^\ell \bfD_\ell] 
                    = \delta A \cdot \bfD + \delta \partial^\ell A_\ell + \delta A \cdot A,
            \end{equation}

        \item composition of difference with magnetic derivative,
            \begin{equation}\label{eq:diff-magn}
                \delta \bfD = \nabla \delta + \delta A + A \delta.
            \end{equation}
    \end{enumerate}
\end{lemma}

\begin{proof}
    Fix an arbitrary complex scalar field $\Phi$. 
    \begin{enumerate}
        \item We compute 
            \begin{align*}
                \delta \bfD_\bfa \Phi 
                    &= (\partial_\bfa - i A_\bfa) \Phi - (\partial_\bfa - i A_\bfa') \Phi' \\
                    &= (\partial_\bfa - i A_\bfa) (\Phi - \Phi') - i (A_\bfa - A_\bfa') \Phi' = \bfD_\bfa \delta \Phi - i \delta A_\bfa \Phi' .
            \end{align*}
            Rearranging gives the result. 

        \item We compute 
            \begin{align*}
                \delta \bfD^\ell \bfD_\ell \Phi
                    &= (\partial^\ell - i A^\ell) (\partial_\ell - i A_\ell) \Phi -(\partial^\ell - i {A'}^\ell) (\partial_\ell - i A_\ell') \Phi'\\
                    &= (\partial^\ell - i A^\ell) (\partial_\ell - i A_\ell) (\Phi - \Phi') + \bigl(  (\partial^\ell - i A^\ell) (\partial_\ell - i A_\ell) - (\partial^\ell - i {A'}^\ell) (\partial_\ell - i A_\ell') \bigr) \Phi'\\
                    &= \bfD^\ell \bfD_\ell \, \delta \Phi - 2i (A_\ell - A_\ell') (\partial^\ell - i{A'}^\ell) \Phi'  - i \partial^\ell (A_\ell - A_\ell') \Phi' - (A^\ell A_\ell - {A'}^\ell A_\ell') \Phi' \\
                    &= \bfD^\ell \bfD_\ell \, \delta \Phi - 2i \delta A_\ell {\bfD'}^\ell \Phi' - 2 \delta A_\ell  {A'}^\ell \Phi' - i \delta \partial^\ell A_\ell \Phi'  - (A^\ell + {A'}^\ell) \delta A_\ell \Phi' .
            \end{align*}
        Rearranging gives the result.
        
        \item We compute 
            \begin{align*}
                \delta \bfD \Phi 
                    &= (\nabla - i A) \Phi - (\nabla - i A') \Phi'\\
                    &= \nabla (\Phi - \Phi') - i A (\Phi - \Phi') - i (A - A') \Phi' ,
            \end{align*}
        which is of the desired form. 
    \end{enumerate}
\end{proof}

\begin{remark}
    Symmetrisation the notation and calculations would yield cleaner schematic formulas. For example, the difference of covariant Laplacians can be written as 
        \[
            \bfD^\ell \bfD_\ell \Phi - {\bfD'}^\ell {\bfD'}_\ell \Phi'
                = -\tfrac12 i (A - A')^\ell (\bfD + \bfD')_\ell (\Phi + \Phi') - \tfrac12 i \partial^\ell (A - A')_\ell (\Phi + \Phi').
        \]
    The cost is that we need to symmetrise the definition of $\bfD$ in the schematic notation. Taking the ``wrong'' magnetic derivative can be quite dangerous (see Appendix \ref{app:example}), so we avoid the symmetrised notation.
\end{remark}

\subsection{Magnetic Sobolev estimates}

We record some analogues of Sobolev inequalities for magnetic derivatives $\bfD = \nabla - i A$. Define the corresponding \textit{magnetic Sobolev spaces} and \textit{magnetic Zhidkov spaces} by 
    \begin{align*}
        H^N_A (\R^2) 
            &:= \Bigl\{ \Phi \in H^N_{x,\loc} (\R^2 \to \C) : \sum_{n = 0}^N \bigl\| \bfD^{(n)} \Phi \bigr\|_{L^2_x} < \infty  \Bigr\},\\
        X^N_A (\R^2) 
            &:= \Bigl\{ \Phi \in H^N_{x,\loc} (\R^2 \to \C) : \| \Phi \|_{L^\infty_x} + \sum_{n = 1}^N \bigl\| \bfD^{(n)} \Phi \bigr\|_{L^2_x} < \infty  \Bigr\}.
    \end{align*}
For these spaces, we have the following interpolation inequalities, 

\begin{lemma}
\leavevmode
    \begin{enumerate}
        \item \textup{(Magnetic Bethuel-Smets inequality\footnote{The reader should compare with \cite[Lemma A.1]{BethuelSmets2007}.}).} Let $\Phi \in X^2_A (\R^2)$, then 
            \begin{equation}
                 \|\bfD \Phi \|_{L^4_x} 
                    \lesssim \|\Phi\|_{L^\infty_x}^{\frac12} \ \|\bfD^{(2)} \Phi\|_{L^2_x}^{\frac12}.\label{eq:GN4} 
            \end{equation}

        \item \textup{(Magnetic Agmon inequality).} Let $\Phi \in H^2_A (\R^2)$, then 
            \begin{equation}
            \|\Phi\|_{L^\infty_x} 
                \lesssim \| \Phi \|_{L^2_x}^{\frac12} \ \bigl\| \bfD^{(2)} \Phi\bigr\|_{L^2_x}^{\frac12}. \label{eq:GNinfty}
            \end{equation}
    \end{enumerate}
\end{lemma}

\begin{proof}
\leavevmode
    \begin{enumerate}
        \item Integrating-by-parts and applying Cauchy-Schwartz yields 
            \begin{align*}
                \|\bfD \Phi\|_{L^4_x}^4
                    &\approx \int_{\R^2} \bfD \Phi \cdot \bfD \Phi \cdot \bfD \Phi \cdot \bfD \Phi \, dx \\
                    &\approx \int_{\R^2} \Phi \cdot \bfD \Phi \cdot \bfD \Phi \cdot \bfD^{(2)} \Phi \, dx \lesssim \|\Phi\|_{L^\infty_x} \ \|\bfD \Phi\|_{L^4_x}^2 \ \|\bfD^{(2)}\Phi \|_{L^2_x}.
            \end{align*}
        Rearranging gives the result.

        \item Consider the following instance of the usual Gagliardo-Nirenberg inequality,
            \[
                \| w \|_{L^\infty_x} \lesssim \|w\|_{L^1_x}^{\frac13} \|\nabla^{(2)} w\|_{L^2_x}^{\frac23}.
            \]
        Applying this to the gauge-invariant quantity $|\Phi|^2$ yields
        \begin{align*}
            \|\Phi\|_{L^\infty_x}^2 
                &\lesssim \| |\Phi|^2 \|_{L^1_x}^{\frac13} \ \| \nabla^{(2)} |\Phi|^2 \|_{L^2_x}^{\frac23} \\
                &\lesssim \|\Phi\|_{L^2_x}^{\frac23} \left( \| \bfD \Phi\|_{L^4_x}^2 + \|\Phi\|_{L^\infty_x} \|\bfD^{(2)} \Phi\|_{L^2_x} \right)^{\frac23} \\
                &\lesssim \|\Phi\|_{L^\infty_x}^{\frac23} \ \|\Phi\|_{L^2_x}^{\frac23} \ \| \bfD^{(2)} \Phi\|_{L^2_x}^{\frac23},
        \end{align*}
    using the covariant product rule \eqref{product} in the second line, and the magnetic interpolation inequality \eqref{GN4} in the third line. Rearranging and taking appropriate roots gives the result. 
    \end{enumerate}
\end{proof}

\begin{lemma}[Magnetic Gagliardo-Nirenberg interpolation inequalities]
\leavevmode
    \begin{enumerate}
        \item Let $n \in \N_0$ be a non-negative integer, $j = 0, 1, \dots, n$, and suppose $\Phi \in H^{n}_A (\R^2)$, then 
            \begin{equation}\label{eq:L2-interpolation}
                \bigl\| \bfD^{(j)} \Phi \bigr\|_{L^2_x} 
                    \leq \bigl\| \Phi \bigr\|_{L^2_x}^{1 - \frac{j}{n}} \bigl\| \bfD^{(n)} \Phi \bigr\|_{L^2_x}^{\frac{j}{n}},
            \end{equation}
        \item Let $n \in \N$ be a positive integer, $j = 1, \dots, n$, and suppose $\bfD\Phi \in H^{n - 1}_A (\R^2)$, then 
        \begin{equation}\label{eq:GN-interpolation}
            \bigl\| \bfD^{(j)} \Phi \bigr\|_{L^{2n/j}_x} 
                \lesssim \| \bfD \Phi \|_{L^2_x}^{1 - \frac{j}{n}}  \bigl\| \bfD^{(n)} \Phi \bigr\|_{L^2_x}^{\frac{j}{n}} .
        \end{equation}      
    \end{enumerate}
\end{lemma}

\begin{proof}
   \leavevmode
    \begin{enumerate}
        \item The first non-trivial case is $n = 2$ and $j = 1$, which follows immediately from integration-by-parts and Cauchy-Schwartz. 
            \begin{align*}
                \int_{\R^2} \bfD\Phi \cdot \bfD \Phi \, dx 
                    &= - \int_{\R^2} \Phi \cdot \bfD^{(2)} \Phi \, dx \\
                    &\leq \bigl\| \Phi \bigr\|_{L^2_x} \, \bigl\| \bfD^{(2)} \Phi \bigr\|_{L^2_x}.
            \end{align*}
        With this strategy in mind, we proceed inductively on $n$, assuming the result holds for $n$, and aim to prove the result for $n + 1$. 
        
        We divide into two cases; first, when $j = n$, we have
            \begin{align*}
                \bigl\| \bfD^{(n)} \Phi \bigr\|_{L^2_x} 
                    &\leq \bigl\| \bfD^{(n - 1)} \Phi\bigr\|_{L^2_x}^{\frac12} \,\bigl\| \bfD^{(n + 1)} \Phi \bigr\|_{L^2_x}^{\frac12}\\
                    &\leq \bigl\| \Phi \bigr\|_{L^2_x}^{\frac12 (1 - \frac{n - 1}{n})} \, \bigl\| \bfD^{(n)}  \Phi\bigr\|_{L^2_x}^{\frac12 \frac{n - 1}{n}} \,\bigl\| \bfD^{(n + 1)} \Phi \bigr\|_{L^2_x}^{\frac12},
            \end{align*}
        using the base case in the first line, and the induction hypotheses in the second line. Rearranging yields the desired result. 

        Second, when $j = 1, \dots, n - 1$, we have
            \begin{align*}
                \bigl\| \bfD^{(j)} \Phi \bigr\|_{L^2_x}
                    &\leq \bigl\| \Phi \bigr\|_{L^2_x}^{1 - \frac{j}{n}} \, \bigl\| \bfD^{(n)} \Phi \bigr\|_{L^2_x}^{\frac{j}{n}} \\
                    &\leq \bigl\| \Phi \bigr\|_{L^2_x}^{1 - \frac{j}{n}} \, \bigl\| \Phi \bigr\|_{L^2_x}^{(1 - \frac{n}{n + 1})\frac{j}{n}} \, \bigl\| \bfD^{(n + 1)} \Phi \bigr\|_{L^2_x}^{\frac{n}{n + 1}\frac{j}{n}},
            \end{align*}
        directly applying the inductive hypothesis and the $j = n$ case. Simplifying the exponents, we see that this is precisely the desired result. 

        \item Consider the following instance of the usual Gagliardo-Nirenberg inequality,
            \[
                \| w \|_{L^{2n/j}_x} 
                \lesssim \| w \|_{L^2_x}^{\frac{j}{n}} \, \| \nabla w \|_{L^2_x}^{1 - \frac{j}{n}}.
            \]
        Applying this to the gauge-invariant quantity $|\bfD^{(j)} \Phi|$, using the diamagnetic inequality \eqref{diamagnetic}, and then interpolating the $L^2$-magnetic Sobolev norms via \eqref{L2-interpolation}, we obtain
        \begin{align*}
            \bigl\| \bfD^{(j)} \Phi \bigr\|_{L^{2n/j}_x} 
                &\lesssim  \bigl\| \bfD^{(j)} \Phi \bigr\|_{L^2_x}^{\frac{j}{n}} \, \bigl\| \nabla |\bfD^{(j)} \Phi|\bigr\|_{L^2_x}^{1 - \frac{j}{n}} \\
                &\lesssim \bigl\| \bfD^{(j)} \Phi \bigr\|_{L^2_x}^{\frac{j}{n}} \, \bigl\| \bfD^{(j+1)} \Phi\bigr\|_{L^2_x}^{1 - \frac{j}{n}}\\
                &\lesssim \| \bfD \Phi \|_{L^2_x}^{\frac{j}{n} (1 - \frac{j - 1}{n - 1}) + (1 - \frac{j}{n})(1 - \frac{j}{n - 1})} \, \bigl\| \bfD^{(n)} \Phi \bigr\|_{L^2_x}^{\frac{j}{n} \frac{j - 1}{n - 1} + (1 - \frac{j}{n}) \frac{j}{n - 1}}. 
        \end{align*}
        Simplifying the exponents, we see that this is precisely the desired result.
    \end{enumerate}
\end{proof}

We will also need an $L^\infty_x$-bound for scalar fields obeying the topological boundary condition $|\Phi| \to 1$ as $|x| \to \infty$. Given $\bfD \Phi \in L^2 (\R^2)$, it follows from the diamagnetic inequality \eqref{diamagnetic} that $|\Phi| \in \dot H^1 (\R^2)$; this is unfortunately not enough to deduce boundedness as the Sobolev inequality 
    \[
        {\dot H}^1 (\R^2)\not\hookrightarrow L^\infty (\R^2),
    \]
fails by a logarithm. We will use the following inequality as a suitable replacement, 

\begin{lemma}[Sobolev-type inequality]
    Let $\Phi \in X^2_A (\R^2)$, and suppose that $1 - |\Phi|^2 \in L^2 (\R^2)$, then
        \begin{equation}\label{eq:inefficientLinfty}
            \|\Phi\|_{L^\infty_x} 
                \lesssim 1 + \big\|1 - |\Phi|^2\big\|_{L^2_x}^{\frac13} \ \big\|\bfD^{(2)} \Phi\big\|_{L^2_x}^{\frac13}.
        \end{equation}
\end{lemma}

\begin{proof}
    Consider the usual Agmon inequality,
        \[
            \| w\|_{L^\infty_x} 
                \lesssim \bigl\|w\bigr\|_{L^2_x}^{\frac12} \bigl\|\nabla^{(2)} w\bigr\|_{L^2_x}^{\frac12}.
        \]
    Applying this to the gauge-invariant quantity $1 - |\Phi|^2$ yields
        \begin{align*}
            \big\| 1 - |\Phi|^2 \big\|_{L^\infty_x} 
                &\lesssim \big\|1 - |\Phi|^2 \big\|_{L^2_x}^{\frac12} \ \big\| \nabla^{(2)} |\Phi|^2 \big\|_{L^2_x}^{\frac12} \\
                &\lesssim \big\|1 - |\Phi|^2 \big\|_{L^2_x}^{\frac12} \left( \|\bfD \Phi\|_{L^4_x}^2 + \|\Phi\|_{L^\infty_x}  \ \big\|\bfD^{(2)} \Phi\big\|_{L^2_x} \right)^{\frac12}\\
                &\lesssim  \|\Phi\|_{L^\infty_x}^{\frac12} \ \big\|1 - |\Phi|^2 \big\|_{L^2_x}^{\frac12} \ \big\| \bfD^{(2)} \Phi\big\|_{L^2_x}^{\frac12},
        \end{align*}
   using the covariant product rule \eqref{product} in the second line, and the magnetic interpolation inequality \eqref{GN4} in the third line. Inserting this into a trivial triangle inequality bound gives
        \[
            \|\Phi\|_{L^\infty_x}^2 
                \lesssim 1 + \|\Phi\|_{L^\infty_x}^{\frac12} \ \big\|1 - |\Phi|^2 \big\|_{L^2_x}^{\frac12} \ \big\| \bfD^{(2)} \Phi\big\|_{L^2_x}^{\frac12}.
        \]
   Interpolating, we can absorb the $L^\infty_x$-term to the left-hand side, completing the proof. 
\end{proof}

For scalar fields with $L^2$-integrable charge density, it will be convenient to record a decomposition into a bounded part and a decaying part. Following the idea of G\'erard \cite[Lemma 1]{Gerard2006} in the non-magnetic case, we have 

\begin{lemma}[Scalar field decomposition]\label{lem:gerard}
    Let $\Phi \in L^1_\loc (\R^2 \to \C)$ be a scalar field such that the charge density satisfies $1 - |\Phi|^2 \in L^2 (\R^2)$. Then there exist $f \in L^\infty (\R^2)$ and $g \in L^2 (\R^2)$ such that 
        \[
            \Phi 
                = f + g,
        \]
    with
        \begin{align}
            \| f \|_{L^\infty_x} 
                &\leq 3,\label{eq:Linfty-part}\\
            \| g \|_{L^2_x} 
                &\leq \big\| 1 - |\Phi|^2 \big\|_{L^2_x} \label{eq:L2-part}.
        \end{align}
    Furthermore, if $\bfD \Phi \in L^2 (\R^2)$ for some magnetic potential $A$, then $f \in X^1_A (\R^2)$ and $g \in H^1_A (\R^2)$ with
        \begin{align}
            \| f \|_{X^1_A}
                &\lesssim 1 + \| \bfD \Phi \|_{L^2_x}, \label{eq:Lp-part1}\\
            \| g \|_{H^1_A} 
                &\lesssim \| 1 - |\Phi|^2 \|_{L^2_x} + \| \bfD \Phi \|_{L^2_x}.\label{eq:Lp-part2}
        \end{align}
\end{lemma}

\begin{proof}
    We cut-off $\Phi$ according to its amplitude, fixing a smooth cut-off $\chi \in C^\infty_c (\C \to [0, 1])$ such that 
        \[
            \chi(\Phi) 
                := 
                \begin{cases}
                    1,  
                        &\text{if $|\Phi| \leq 2$,}\\
                    0, 
                        &\text{if $|\Phi| > 3$},
                \end{cases}
        \]
    and setting 
        \[
            f:= \chi(\Phi) \Phi, \qquad g := (1 - \chi(\Phi)) \Phi.
        \]
    Then \eqref{Linfty-part} is clear from construction, and \eqref{L2-part} follows from the pointwise bound $|g| \leq |1 - |\Phi|^2|$, which holds since $|\Phi| \geq 2$ on the support of $g$. For higher regularity, we compute
        \begin{align*}
            \bfD f 
                &= \chi(\Phi) \bfD \Phi + \chi'(\Phi) \nabla |\Phi| \Phi, \\
            \bfD g 
                &= (1 - \chi(\Phi)) \bfD \Phi - \chi' (\Phi) \nabla |\Phi| \Phi,
        \end{align*}
    then \eqref{Lp-part1}-\eqref{Lp-part2} follows from the diamagnetic inequality \eqref{diamagnetic} and $|\Phi| \leq 3$ on the support of $\chi'(\Phi)$. 
\end{proof}

\subsection{Abstract parabolic theory}\label{subsec:abstract-parabolic}

We will write $s$ for the \textit{heat-temporal} variable (in contrast to the \textit{Schr\"odinger temporal} variable $t$). We typically integrate in $s$ in a dimension-less fashion using the scale-invariant measure $\frac{ds}{s}$, defining for $1 \leq p < \infty$ the \textit{heat-temporal Lebesgue norms},
    \[
        \| \mathsf U \|_{L^p_{\frac{ds}{s}} (I)} := \left( \int_I | \mathsf U (s) |^p \, \frac{ds}{s} \right)^{\frac1p}. 
    \]
To set expectations for the ensuing estimates, recall that parabolic scaling dictates that spatial derivatives scale relative to the heat-time as $\sfD \approx s^{-1/2}$, and the Lebesgue measure on $\R^2$ scales as $dx \approx s^1$. All-in-all, a Sobolev-type norm for a function $\mathsf U(x)$ scales like 
    \[
        \|  \sfD^{(n)} \mathsf U \|_{L^p_x} 
            \approx s^{-\frac{n}{2}} s^{\frac1p}  \cdot \mathsf U. 
    \]
Accordingly, we will often work with parabolic-normalised derivatives $(s^\frac12 \sfD)$ and $(s \sfD^{(2)})$, which approximately act like bounded operators. In line with the dimensional analysis heuristics, we have the following interpolation lemmas, which are essentially parabolically-normalised versions of the magnetic interpolation inequalities from the previous section,

\begin{lemma}[Parabolic interpolation lemma]\label{lem:interpolate}
    Let $\cK > 0$ and $2 \leq p, q \leq \infty$ and $n_0 \in \N_0$. If $\sfU : \R^2 \times [0, \underline s] \to \C$ is a sufficiently regular complex-scalar field satisfying
        \begin{equation}\label{eq:L2-base}
            \bigl\| (s^\frac12 \sfD)^{(n)} \sfU \bigr\|_{L^q_{\frac{ds}{s}} L^2_x}
                \lesssim_n \cK,
        \end{equation}
    for each non-negative integer $n \in \N_0$, then
        \begin{equation}\label{eq:interpolate}
            \bigl\| s^{\frac12 - \frac1p} (s^\frac12 \sfD)^{(n)} \sfU \bigr\|_{L^q_{\frac{ds}{s}} L^p_x}
                \lesssim_n \cK.
        \end{equation}
\end{lemma}

\begin{proof}
    By the interpolation inequality \eqref{GNinfty} when $\sfD = \bfD$, or its non-magnetic counter-part when $\sfD = \nabla$, 
        \begin{align*}
            \bigl\| s^{\frac12} (s^\frac12 \sfD)^{(n)} \sfU \bigr\|_{L^q_{\frac{ds}{s}} L^\infty_x}
                &\lesssim \bigl\|  (s^\frac12 \sfD)^{(n)} \sfU \bigr\|_{L^q_{\frac{ds}{s}} L^2_x}^\frac12 \bigl\| (s^\frac12 \sfD)^{(n + 2)} \sfU \bigr\|_{L^q_{\frac{ds}{s}} L^2_x}^\frac12 \\
                &\lesssim_n \cK,
        \end{align*}
    using \eqref{L2-base} in the second line. This yields \eqref{interpolate} for $p = \infty$. For the non-endpoint cases, we interpolate
        \begin{align*}
            \bigl\| s^{\frac12 - \frac1p} (s^\frac12 \sfD)^{(n)} \sfU \bigr\|_{L^q_{\frac{ds}{s}} L^p_x}
                &\lesssim \bigl\| (s^\frac12 \sfD)^{(n)} \sfU \bigr\|_{L^q_{\frac{ds}{s}} L^2_x}^{\frac2p} \bigl\| s^{\frac12} (s^\frac12 \sfD)^{(n)} \sfU \bigr\|_{L^q_{\frac{ds}{s}} L^\infty_x}^{1 - \frac2p} \\
                &\lesssim_n \cK,
        \end{align*}
    using \eqref{L2-base} and the $p = \infty$ case in the second line. 
\end{proof}

\begin{lemma}[Abstract bilinear interpolation]\label{lem:bilinear}
    For sufficiently regular scalar fields $\sfU, \sfV : \R^2 \times [0, \underline s] \to \C$, we have  
        \begin{equation}
            \bigl\| s^\frac12 \cdot \sfU \cdot \sfV \bigr\|_{L^2_x}
                \lesssim \bigl\| (s^\frac12 \sfD) \sfU \bigr\|_{L^2_x} \bigl\| \sfV \bigr\|_{L^2_x} + \bigl\| \sfU \bigr\|_{L^2_x} \bigl\| (s^\frac12 \sfD) \sfV \bigr\|_{L^2_x},
        \end{equation}
    where $\sfD \in \{ \nabla, \bfD\}$. 
\end{lemma}

\begin{proof}
    We compute 
        \begin{align*}
            \bigl\| s^\frac12 \cdot \sfU \cdot \sfV \bigr\|_{L^2_x}
                &\lesssim \bigl\| (s^\frac12 \sfD) (\sfU \cdot \sfV) \bigr\|_{L^1_x} \\
                &\lesssim \bigl\| (s^\frac12 \sfD) \sfU \bigr\|_{L^2_x} \bigl\| \sfV \bigr\|_{L^2_x} + \bigl\| \sfU \bigr\|_{L^2_x} \bigl\| (s^\frac12 \sfD) \sfV \bigr\|_{L^2_x},
        \end{align*}
    using the Gagliardo-Nirenberg-Sobolev inequality (along with the diamagnetic inequality \eqref{diamagnetic} in the case $\sfD = \bfD$) in the first line, and the product rule in the second line and Cauchy-Schwarz in the second line. 
\end{proof}

Our main order of business herein is to write a cook-book of smoothing estimates for abstract parabolic equations. Let $A$ be a connection $1$-form on $\R^2 \times [0, \infty)$ solving the \textit{(abelian) Yang-Mills heat flow}
    \begin{equation}\tag{P-$A$}\label{eq:YMH}
        F_{sj}
            = \partial^\ell F_{\ell j},
    \end{equation}
in \textit{caloric gauge}, 
    \begin{equation}\tag{Cal}\label{eq:caloric-prelim}
        A_s = 0.
    \end{equation}
Fix $\sfD_\bfa \in \{ \bfD_\bfa, \partial_\bfa \}$, and suppose $\sfU : \R^2 \times [0, \underline s] \to \C$ solves the \textit{abstract parabolic equation},
    \begin{equation}\label{eq:abstract-heat-eq}\tag{Par}
        \bigl(\sfD_s - \sfD^\ell \sfD_\ell \bigr) \sfU
            = \sfN.
    \end{equation}
We begin with some basic calculations for the Yang-Mills heat flow, 

\begin{lemma}
    Let $A$ be a connection $1$-form on $\R^2 \times [0, \underline s]$ solving the Yang-Mills heat flow \eqref{YMH}. 
    \begin{enumerate}
        \item \textup{(Linear heat equation).} The magnetic field satisfies 
        \begin{equation}\label{eq:F12-heat}
            (\partial_s - \Delta) F_{12} 
                = 0.
        \end{equation}
        In particular, the divergence-free part of the magnetic potential satisfies 
        \begin{equation}
            (\partial_s - \Delta) \PP^\df A_x
                = 0.\label{eq:A-heat}
        \end{equation}

        \item \textup{(Linear transport equation).} $F_{sj}$ is divergence-free,
        \begin{equation}\label{eq:Fsj-transport}
            \partial^j F_{sj} 
                = 0.
        \end{equation}
        In particular, if $A$ is in caloric gauge \eqref{caloric-prelim}, then the curl-free part of the magnetic potential satisfies,  
            \begin{equation}
                 \partial_s \PP^\cf A_x 
                    = 0. \label{eq:A-transport}
            \end{equation}

        \item \textup{(Commutator identity).} The covariant derivative satisfies the identity
        \begin{equation}\label{eq:parabolic-commute}
                \bigl[ \bfD_j, \bfD_s - \bfD^\ell \bfD_\ell \bigr]
                    = 2i F\indices{_j^\ell} \bfD_\ell.
            \end{equation}
    \end{enumerate}

\end{lemma}

\begin{proof}
    \leavevmode
    \begin{enumerate}
        \item By the Bianchi identity \eqref{bianchi} and the Yang-Mills heat flow equation \eqref{YMH}, 
        \begin{align*}
            \partial_s F_{12} 
                &= \partial_1 F_{s2} - \partial_2 F_{s 1} \\
                &= \partial^\ell ( \partial_1 F_{\ell 2}- \partial_2 F_{\ell 1}) \\
                &= \Delta F_{12}.
        \end{align*}

        \item By the Yang-Mills heat flow \eqref{YMH} and anti-symmetry, 
        \begin{align*}
            \partial^j F_{sj} 
                &= \partial^j \partial^\ell F_{\ell j} \\
                &= \epsilon_{\ell j} \partial^j \partial^\ell F_{12} \\
                &= 0.
        \end{align*}
        Imposing the caloric gauge, this identity reduces to the transport equation for the divergence,
            \[
                \partial_s \partial^j A_j = 0.
            \]

        \item Using the commutator identities \eqref{commute}, \eqref{laplace-commute}, and the Yang-Mills heat flow \eqref{YMH},
            \begin{align*}
                \bigl[ \bfD_j, \bfD_s - \bfD^\ell \bfD_\ell \bigr]
                    &= i \bigl( F_{sj} + 2 F\indices{_j^\ell} \bfD_\ell + \partial^\ell F_{j \ell}  \bigr) \\
                    &= 2i F\indices{_j^\ell} \bfD_\ell.
            \end{align*}
    \end{enumerate}
\end{proof}

The basic ingredients for the smoothing estimates consist of energy estimates, and properties of the linear heat propagator, which we can transfer to solutions to the abstract parabolic equation using Duhamel's principle and the maximum principle. 

\begin{lemma}[Properties of linear heat propagator]
    \leavevmode
    \begin{itemize}
        \item \textup{(Parabolic smoothing).} For $1 \leq q \leq p \leq \infty$ and non-negative integers $n \in \N_0$, 
            \begin{equation} \label{eq:standard-heat}
                \big\| s^{\frac1q - \frac1p} (s^\frac12\nabla)^{(n)} e^{s \Delta} u_0 \big\|_{L^p_x} 
                    \lesssim_{n, p, q} \|u_0\|_{L^q_x}.
            \end{equation} 

        \item \textup{(Parabolic ``Strichartz").} For $2 \leq p \leq \infty$ and non-negative integers $n \in \N_0$ such that $(p, n) \neq (2, 0)$, 
            \begin{equation}\label{eq:parabolic-strichartz}
            \bigl\| s^{\frac12 - \frac1p} (s^\frac12 \nabla)^{(n)} e^{s \Delta} u_0 \bigr\|_{L^2_{\frac{ds}s} L^p_x}
                \lesssim_{n, p} \| u_0 \|_{L^2_x}. 
            \end{equation}

    \end{itemize}
\end{lemma}

\begin{proof}
    The former is a standard consequence of Young's convolution inequality and the explicit form of the heat kernel. For the latter, when $p = 2$, parabolic smoothing reduces the estimate to showing the $n = 1$ case, which is precisely the parabolic energy estimate, c.f. \eqref{abs-heat-energy}. For $2 < p \leq \infty$, it suffices to prove the result for $n = 0$ by parabolic smoothing. Then one can employ a $T$-$T^*$-argument, see \cite[Lemma 2.5]{Tao2009}.
\end{proof}

As a consequence of the standard heat equation estimates \eqref{standard-heat}-\eqref{parabolic-strichartz} applied to the equation for the magnetic field \eqref{F12-heat}, we have

\begin{proposition}[Magnetic field smoothing]
    Let $A$ be a connection $1$-form on $\R^2 \times [0, \infty)$ solving the Yang-Mills heat flow \eqref{YMH}. Then for $2 \leq p \leq \infty$ and non-negative integers $n \in \N_0$,
        \begin{equation}\label{eq:F12-sobolev}
        \begin{split}
            \bigl\| s^{\frac12 - \frac1p} (s^\frac12 \nabla)^{(n)}  F_{12} \bigr\|_{L^\infty_s L^p_x} 
                &\lesssim_{p, n} \bigl\|  \overline F_{12} \bigr\|_{L^2_x},\\
            \bigl\| s^{\frac12 - \frac1p} (s^\frac12 \nabla)^{(n)}  F_{12} \bigr\|_{L^2_{\frac{ds}{s}} L^p_x} 
                &\lesssim_{p, n} \bigl\|  \overline F_{12} \bigr\|_{L^2_x}, \qquad (p, n) \neq (2, 0). 
        \end{split}
        \end{equation}
\end{proposition}

To prove smoothing estimates for the abstract parabolic equation \eqref{abstract-heat-eq}, we proceed by the energy method. Starting from the equation \eqref{abstract-heat-eq}, we can inductively use the commutator identity \eqref{parabolic-commute} and the product rule to derive 
    \begin{equation}\label{eq:abstract-heat-eq-n}\tag{$\mathrm{Par}_n$}
        \bigl( \sfD_s - \sfD^\ell \sfD_\ell \bigr) \sfD^{(n)} \sfU 
            = \sum_{\substack{a + b = n \\ b \neq 0}} \nabla^{(a)} F_{12} \cdot \sfD^{(b)} \sfU + \sfD^{(n)} \sfN.
    \end{equation}

\begin{lemma}
    Let $A$ be a connection $1$-form on $\R^2 \times [0, \underline s]$ solving the Yang-Mills heat flow \eqref{YMH}, and suppose that  
         \begin{equation}\label{eq:heat-data-0}
            {\bigl\| \overline F_{12} \bigr\|_{L^2_x}}
                \leq \underline s^{-\frac12}.
        \end{equation}
    If $\sfU : \R^2 \times [0, \underline s] \to \C$ is a sufficiently regular solution to the abstract linear parabolic equation \eqref{abstract-heat-eq}, then we have the following estimates: 
    \begin{itemize}   
        \item \textup{(Energy estimate).} 
        \begin{equation}\label{eq:abs-heat-energy}
            \bigl\| \sfU \bigr\|_{L^\infty_s L^2_x} + \bigl\| (s^\frac12 \sfD) \sfU \bigr\|_{L^2_{\frac{ds}{s}} L^2_x}  
                \lesssim \bigl\| \overline \sfU \bigr\|_{L^2_x} + \bigl\| s \, \sfN \bigr\|_{L^1_{\frac{ds}{s}} L^2_x}.
        \end{equation}
        
        \item \textup{(Weighted energy estimate).} For each $\gamma > 0$, 
        \begin{equation}
             \big\| s^\gamma \sfU \big\|_{L^\infty_s L^2_x} + \big\| s^\gamma (s^{\frac12} \sfD) \sfU \big\|_{L^2_{\frac{ds}{s}} L^2_x} 
                \lesssim \gamma \big\|s^\gamma \sfU\big\|_{L^2_{\frac{ds}{s}} L^2_x} + \big\| s^{\gamma + 1} \sfN \big\|_{L^1_{\frac{ds}{s}} L^2_x}. \label{eq:abs-heat-weight}
        \end{equation}
        
        \item \textup{(Smoothing estimate).} For each non-negative integer $n \in \N_0$,  
            \begin{equation}\label{eq:abs-smooth}
            \begin{split}
            \bigl\| (s^\frac12 \sfD)^{(n)}  \sfU \bigr\|_{L^\infty_s L^2_x} &+ \bigl\|  (s^\frac12 \sfD)^{(n + 1)}  \sfU \bigr\|_{L^2_{\frac{ds}{s}} L^2_x}
                \lesssim_n \bigl\| \overline\sfU\bigr\|_{L^2_x} +  \bigl\| s (s^\frac12 \sfD)^{(\leq n)} \sfN \bigr\|_{L^1_{\frac{ds}{s}} L^2_x}.
            \end{split}
            \end{equation}

        \item \textup{(Dyadic smoothing estimate).} For each non-negative integer $n \in \N_0$ and $\gamma \in \R$, 
            \begin{equation}\label{eq:dyadic-smooth}
            \bigl\| {s'}^\gamma ({s'}^\frac12 \sfD)^{(n)} \sfU (s') \bigr\|_{L^2_x}
                \lesssim_n \bigl\| {s'}^\gamma \sfU (\tfrac12 s') \bigr\|_{L^2_x} + \bigl\| s^{1 + \gamma} (s^\frac12 \sfD)^{(\leq n)} \sfN \bigr\|_{L^2_{\frac{ds}{s}} L^2_x ([\frac12 s', s'])}.
            \end{equation}
            
    \end{itemize}
\end{lemma}

\begin{proof}[Proof of energy estimate \eqref{abs-heat-energy}]
    Multiplying the abstract parabolic equation \eqref{abstract-heat-eq} by $\sfU$, integrating-by-parts on $\R^2$ yields the energy identity 
        \[
            \tfrac12 \partial_s\bigl\| \sfU \bigr\|_{L^2_x}^2 +  \bigl\| \sfD \sfU \bigr\|_{L^2_x}^2  
                = \int_{\R^2} \sfN \cdot \sfU \, dx.
        \]
    Cauchy-Schwarz allows us to bound the second term on the right by $\| \sfN \|_{L^2_x} \| \sfU \|_{L^2_x}$. Integrating in $s$ on the interval $[0, s']$ and taking supremum over $s' \in [0, \underline s]$ furnishes 
        \begin{align*}
            \tfrac12 \bigl\| \sfU \bigr\|_{L^\infty_s L^2_x}^2 + \bigl\| (s^\frac12 \sfD) \sfU \bigr\|_{L^2_{\frac{ds}{s}} L^2_x}^2 
                &\leq \tfrac12 \bigl\| \overline \sfU \bigr\|_{L^2_x}^2 + \bigl\| s \, \sfN \bigr\|_{L^1_{\frac{ds}{s}} L^2_x}  \bigl\| \sfU \bigr\|_{L^\infty_s L^2_x} .
        \end{align*}
    Interpolating, the third term on the right-hand side may be controlled by $\| s \sfN \|_{L^1_{ds/s} L^2_x}^2$, plus a term which can be absorbed into the left-hand side. This completes the proof. 
\end{proof}

\begin{proof}[Proof of weighted energy estimate \eqref{abs-heat-weight}]
    Commuting $s^\gamma$ into the abstract parabolic equation \eqref{abstract-heat-eq}, 
        \[
            \big(\sfD_s - \sfD^j \sfD_j\big) (s^\gamma \sfU) 
                = \gamma s^{\gamma - 1} \sfU + s^\gamma \sfN.
        \]
    Multiplying by $s^\gamma \sfU$, integrating on $\R^2$, and integrating-by-parts furnishes the energy identity
        \[
            \tfrac12 \partial_s \big\|s^\gamma \sfU\big\|_{L^2_x}^2 + \frac1s \big\| s^\gamma (s^{\frac12} \sfD) \sfU \big\|_{L^2_x}^2 
                = \frac{\gamma}{s} \big\| s^\gamma \sfU \big\|_{L^2_x}^2 + \int_{\R^2} s^\gamma \sfN \cdot s^\gamma \sfU \, dx. 
        \]
    Cauchy-Schwartz allows us to bound the second term on the right by $\|s^\gamma \sfN\|_{L^2_x} \|s^\gamma \sfU\|_{L^2_x}$. Integrating in $s$ on the interval $[0, s']$ and taking supremum over $s' \in [0, \underline s]$ yields the inequality
        \begin{align*}
            \tfrac12 \big\|s^\gamma \sfU\big\|_{L^\infty_s L^2_x}^2 + \big\| s^\gamma (s^{\frac12} \sfD) \sfU \big\|_{L^2_{\frac{ds}{s}} L^2_x}^2 
                &\leq \gamma \big\|s^\gamma \sfU\big\|_{L^2_{\frac{ds}{s}} L^2_x}^2 + \big\| s^{\gamma + 1} \sfN \big\|_{L^1_{\frac{ds}{s}} L^2_x} \big\| s^\gamma \sfU \big\|_{L^\infty_s L^2_x}.
        \end{align*}
    Interpolating, the third term on the right-hand side may be controlled by $\|s^{\gamma + 1} \sfN\|^2_{L^1_{ds/s} L^2_x}$ plus a term which can be absorbed by the first term on the left-hand side. This completes the proof. 
\end{proof}

\begin{proof}[Proof of smoothing estimate \eqref{abs-smooth}]
    We will consider the case $\sfD_\bfa = \bfD_\bfa$; the case $\sfD_\bfa = \partial_\bfa$ follows a similar but easier argument due to the lack of commutator terms. For brevity, denote $\mathrm E_{n}$ for the energy-type quantity, and $\mathrm B_{n}$ for the bulk-type quantity, 
        \begin{align*}
            \mathrm E_{n} 
                &:= \bigl\| (s^\frac12 \bfD)^{(n - 1)} \sfU \bigr\|_{L^\infty_s L^2_x}, \\
                \mathrm B_{n} 
                &:= \bigl\| (s^\frac12 \bfD)^{(n)} \sfU \bigr\|_{L^2_{\frac{ds}{s}} L^2_x}.
        \end{align*}
    Observe the $n = 0$ case is precisely the energy estimate \eqref{abs-heat-energy}. Proceeding inductively, we assume 
        \[
            \mathrm{E}_b + \mathrm{B}_b 
                \lesssim \bigl\| \overline \sfU \bigr\|_{L^2_x} + \bigl\| s (s^\frac12 \bfD)^{(\leq b)} \sfN \bigr\|_{L^1_{\frac{ds}{s}} L^2_x}
        \]
    for $b = 0, \dots, n$. Applying the weighted energy estimate \eqref{abs-heat-weight} with $\gamma = \frac{n}{2}$ to the equation for the $n$-th order derivatives \eqref{abstract-heat-eq-n}
        \begin{align*}
            \mathrm E_{n + 1} + \mathrm B_{n + 1}
                &\lesssim \mathrm B_{n} + \bigl\| s (s^\frac12 \bfD)^{(n)} \sfN \bigr\|_{L^1_{\frac{ds}{s}} L^2_x} + \sum_{\substack{a + b = n \\ b \neq 0}} \bigl\| s (s^\frac12 \nabla)^{(a)} F_{12} \cdot (s^\frac12 \bfD)^{(b)} \sfU \bigr\|_{L^1_{\frac{ds}{s}} L^2_x}\\
                &\lesssim \bigl\| \overline \sfU \bigr\|_{L^2_x} + \bigl\| s (s^\frac12 \bfD)^{(\leq n)} \sfN \bigr\|_{L^1_{\frac{ds}{s}} L^2_x} + \sum_{\substack{a + b = n \\ b \neq 0}} \underline s^\frac12 \bigl\| s^\frac12 (s^\frac12 \nabla)^{(a)} F_{12} \bigr\|_{L^\infty_{s, x}} \cdot \mathrm{B}_{b} \\
                &\lesssim \bigl\| \overline \sfU \bigr\|_{L^2_x} + \bigl\| s (s^\frac12 \bfD)^{(\leq n)} \sfN \bigr\|_{L^1_{\frac{ds}{s}} L^2_x} ,
        \end{align*}
    using H\"older to estimate the commutator terms, and the inductive hypotheses for the bulk terms, parabolic smoothing of the magnetic field \eqref{F12-sobolev} and \eqref{heat-data-0}. 
\end{proof}

\begin{proof}[Proof of dyadic smoothing estimate \eqref{dyadic-smooth}]
    Since $s' \sim s$ in the interval $s \in [\frac12 s', s']$, it suffices to show the result for $\gamma = 0$. Denoting heat-temporal translation by $\tfrac12 s'$ with tildes, e.g.
            \begin{align*}
                \widetilde \sfU (s) 
                    &:= \sfU (\tfrac12 s' + s),\\
                \widetilde \sfN (s) 
                    &:= \sfN (\tfrac12 s' + s), 
            \end{align*}
        we apply the abstract parabolic smoothing \eqref{abs-smooth} to $\widetilde \sfU$ on the interval $[0, \tfrac12 s']$. On the left hand-side of the inequality, we have 
            \begin{align*}
                \bigl\| ({s'}^\frac12 \sfD)^{(n)} \sfU (s') \bigr\|_{L^2_x} 
                    \lesssim \bigl\| (s^\frac12 \widetilde\sfD)^{(n)} \widetilde\sfU \bigr\|_{L^\infty_s L^2_x ([0, \frac12 s'])} .
            \end{align*}
    On the right-hand side, it is convenient to first convert from the scale-invariant measure $\frac{ds}{s}$ to the translation-invariant measure $ds$. Applying Cauchy-Schwarz and observing that $s \sim s'$ within $s \in [\tfrac12 s', s']$, 
        \begin{align*}
            \bigl\| s (s^\frac12 \widetilde\sfD)^{(n)} \widetilde \sfN \bigr\|_{L^1_{\frac{ds}{s}} L^2_x ([0, \frac12 s'])}
                &\lesssim \bigl\| s^{\frac{n}{2}} \bigr\|_{L^2_{s} ([0, \frac12 s'])} \bigl\| \sfD^{(n)} \sfN \bigr\|_{L^2_{s, x} ([\frac12 s', s'])} \\
                &\lesssim |s'|^{\frac{n + 1}{2}} \bigl\| \sfD^{(n)} \sfD \sfN \bigr\|_{L^2_{s, x} ([\frac12 s', s'])}  \\
                &\lesssim \bigl\| s (s^\frac12\sfD)^{(n)} \sfD \sfN \bigr\|_{L^2_{\frac{ds}{s}} L^2_x ([\frac12 s', s'])} .
        \end{align*}
    Collecting these observations gives the desired inequality.
\end{proof}

We also record an analogue of the smoothing estimate \eqref{abs-smooth} for differences. Let $A$ and $A'$ be connection $1$-forms on $\R^2 \times [0, \underline s]$ solving the Yang-Mills heat flow \eqref{YMH} under caloric gauge \eqref{caloric-prelim}, and suppose that $\sfU, \sfU' : \R^2 \times [0, \underline s] \to \C$ solve the abstract linear parabolic equation \eqref{abstract-heat-eq} with the respective background connection $1$-form. Commuting the difference operator $\delta$ through the higher-order equation \eqref{abstract-heat-eq-n} gives
    \begin{equation}\label{eq:abstract-heat-diff}\tag{$\delta\text{Par}_n$}
        \begin{split}
            \bigl(\partial_s - \bfD^\ell \bfD_\ell\bigr) \delta \bfD^{(n)} \sfU 
                &= \bigl[\delta, \bfD^\ell \bfD_\ell\bigr] \bfD^{(n)} \sfU +  \sum_{\substack{a + b = n \\ b \neq 0}} \nabla^{(a)} \delta F_{12} \cdot \bfD^{(b)} \sfU  + \nabla^{(a)} F_{12} \cdot \delta \bfD^{(b)} \sfU + \delta \bfD^{(n)} \sfN  \\
                &=: \mathrm{I}_n + \mathrm{II}_n + \mathrm{III}_n + \delta \bfD^{(n)} \sfN .
    \end{split}
    \end{equation}

\begin{lemma}[Abstract parabolic difference estimate]\label{lem:diff-smooth}
    Let $A$ and $A'$ be connection $1$-forms on $\R^2 \times [0, \underline s]$ solving the Yang-Mills heat flow \eqref{YMH} in caloric gauge \eqref{caloric-prelim} with heat-temporal initial data satisfying \eqref{heat-data-0}, and suppose that $\sfU, \sfU' : \R^2 \times [0, \underline s] \to \C$ solve respective abstract linear parabolic equations \eqref{abstract-heat-eq} and satisfy the bounds
        \begin{equation}\label{eq:U-diff-K}
            \bigl\| s^\frac18 (s^\frac12 \bfD)^{(n)} \sfU \bigr\|_{L^\infty_s L^2_x} + \bigl\| s^\frac18 (s^\frac12 \bfD)^{(n)} \sfU' \bigr\|_{L^\infty_s L^2_x}  
                \lesssim_n \cK,
        \end{equation}
    for some $\cK > 0$. Then 
        \begin{equation}
        \begin{split}
            \bigl\| \delta (s^\frac12 \bfD)^{(n)} \sfU \bigr\|_{L^\infty_s L^2_x} &+ \bigl\| \delta (s^\frac12 \bfD)^{(n + 1)} \sfU \bigr\|_{L^2_{\frac{ds}{s}} L^2_x}\\
                    &\lesssim \bigl\| \delta \overline \sfU \bigr\|_{L^2_x}  + \bigl\| s \, \delta (s^\frac12 \bfD)^{(\leq n)} \sfN \bigr\|_{L^1_{\frac{ds}{s}} L^2_x} \\
                    &\qquad +  \underline s^\frac18 \cdot \cK \cdot \Bigl( \bigl\| \delta \overline A \bigr\|_{L^4_x} + \underline s^\frac14 \cdot \bigl\| \overline A \bigr\|_{L^4_x} \cdot \bigl\| \delta \overline A \bigr\|_{L^4_x} + \underline s^\frac14 \cdot \bigl\| \nabla \delta \overline A \bigr\|_{L^2_x}\Bigr).
        \end{split}
        \end{equation}
\end{lemma}

\begin{proof}   
    Define the $\mathrm{E}_n$ the energy-type quantity, $\mathrm{B}_n$ and $\widetilde{\mathrm{B}}_n$ for the bulk-type quantities,
        \begin{align*}
            \delta \mathrm{E}_n 
                &:= \bigl\| \delta (s^\frac12 \bfD)^{(n - 1)} \sfU  \bigr\|_{L^\infty_s L^2_x},\\
            \delta \mathrm{B}_n  
                &:= \bigl\| (s^\frac12 \bfD) \delta (s^\frac12 \bfD)^{(n - 1)} \sfU \bigr\|_{L^2_{\frac{ds}{s}} L^2_x},\\
            \delta \widetilde{\mathrm{B}}_n 
                &:= \bigl\| \delta (s^\frac12 \bfD)^{(n)} \sfU \bigr\|_{L^2_{\frac{ds}{s}} L^2_x} ,
        \end{align*}
    Commuting $\delta$ with $\bfD$ a l\'a \eqref{commute-diff}, applying H\"older, interpolation from Lemma \ref{lem:interpolate} with \eqref{U-diff-K},
        \begin{align*}
            \delta \widetilde{\mathrm{B}_n} 
                &\lesssim \delta\mathrm{B}_n + \bigl\| [s^\frac12 \bfD, \delta] (s^\frac12 \bfD)^{(n - 1)} \sfU \bigr\|_{L^2_{\frac{ds}{s}} L^2_x} \\
                &\lesssim \delta\mathrm{B}_n + \bigl\| \delta \overline A \bigr\|_{L^4_x} \bigl\| s^\frac12 (s^\frac12 \bfD)^{(n - 1)} \sfU \bigr\|_{L^2_{\frac{ds}{s}} L^4_x} \\
                &\lesssim \delta \mathrm{B}_n + \underline s^\frac18 \cdot \cK \cdot  \bigl\| \delta \overline A \bigr\|_{L^4_x}.
        \end{align*}
    This will be useful on the left. Next, we estimate the terms on the right of \eqref{abstract-heat-diff}; the first term consists of commutators with the covariant Laplacian, so expanding via the identity \eqref{commute-diff-laplace}, applying H\"older and the assumption \eqref{U-diff-K} appropriately, 
        \begin{align*}
            \bigl\| s \cdot s^{\frac{n}{2}}\cdot \mathrm{I}_n \bigr\|_{L^1_{\frac{ds}{s}} L^2_x}
                &\lesssim \bigl\| \delta \overline A \bigr\|_{L^4_x} \bigl\| s^\frac12 (s^\frac12 \bfD)^{(n + 1)} \sfU \bigr\|_{L^\infty_s L^4_x} +  \Bigl( \bigl\| \delta \partial^\ell \overline A_\ell \bigr\|_{L^2_x} + \bigl\| \delta \overline A \bigr\|_{L^4_x} \bigl\| \overline A \bigr\|_{L^4_x} \Bigr) \bigl\| s(s^\frac12 \bfD)^{(n)} \sfU \bigr\|_{L^\infty_{s, x}} \\
                &\lesssim \underline s^\frac18 \cdot \cK \cdot \bigl\| \delta \overline A \bigr\|_{L^4_x} + \underline s^{\frac38} \cdot \cK \cdot \Bigl( \bigl\| \delta \partial^\ell \overline A_\ell \bigr\|_{L^2_x} + \bigl\| \overline A \bigr\|_{L^4_x} \bigl\| \delta \overline A \bigr\|_{L^4_x} \Bigr).
        \end{align*}
    The second term consists of differences falling on the magnetic field. Applying parabolic smoothing \eqref{F12-sobolev} and the assumption \eqref{U-diff-K}, 
        \begin{align*}
            \bigl\| s \cdot s^{\frac{n}{2}} \cdot\mathrm{II}_n \bigr\|_{L^1_{\frac{ds}{s}} L^2_x} 
                &\lesssim  \bigl\| s^\frac12 (s^\frac12 \nabla)^{(\leq n)} \delta F_{12} \bigr\|_{L^2_{\frac{ds}{s}} L^\infty_x} \cdot \bigl\| s^\frac12 (s^\frac12 \bfD)^{(\leq n)} \sfU \bigr\|_{L^2_{\frac{ds}{s}} L^2_x} \\
                &\lesssim \underline s^{\frac38} \cdot \cK \cdot \bigl\| \delta \overline F_{12} \bigr\|_{L^2_x} 
        \end{align*}
    The third term consists of differences falling on derivatives of $\sfU$. By Cauchy-Schwarz and the assumption \eqref{heat-data-0} on $\overline F_{12}$, 
        \begin{align*}
            \bigl\|  s \cdot s^{\frac{n}{2}} \cdot \mathrm{III}_n \bigr\|_{L^1_{\frac{ds}{s}} L^2_x} 
                &\lesssim \sum_{\substack{a + b = n \\ b \neq 0}} \underline s^\frac12 \bigl\| s^\frac12 (s^\frac12 \nabla)^{(a)} F_{12} \bigr\|_{L^\infty_{s, x}} \cdot \delta \widetilde{\mathrm{B}}_b \\
                &\lesssim \sum_{b = 1}^n \delta \widetilde{\mathrm{B}}_b.
        \end{align*}
    
    Armed with the previous calculations, we now apply the energy estimate \eqref{abs-heat-energy} to the difference equation \eqref{abstract-heat-diff} with $n = 0$,
        \begin{align*}
            \delta\mathrm{E}_1 + \delta\mathrm{B}_1
                &\lesssim \bigl\| \delta \overline \sfU \bigr\|_{L^2_x} + \bigl\| s \cdot (\mathrm{I}_0 + \mathrm{II}_0) \bigr\|_{L^1_{\frac{ds}{s}} L^2_x} + \bigl\| s \, \delta \sfN \bigr\|_{L^1_{\frac{ds}{s}} L^2_x} \\
                &\lesssim \bigl\| \delta \overline \sfU \bigr\|_{L^2_x}  +  \underline s^\frac18 \cdot \cK \cdot \bigl( 1 + \underline s^\frac12 \bigl\| \overline A \bigr\|_{L^4_x} \bigr) \cdot \bigl\| \delta \overline A \bigr\|_{L^4_x} + \underline s^{\frac38} \cdot \cK \cdot \bigl\| \nabla \delta \overline A \bigr\|_{L^2_x} + \bigl\| s \, \delta \sfN \bigr\|_{L^1_{\frac{ds}{s}} L^2_x} .
        \end{align*}
    Recall $\delta \widetilde{\mathrm B}_1$ can be controlled by $\delta \mathrm B_1$, so the result for $n = 0$ holds. Proceeding inductively, we assume the result holds up to $n$. Applying the abstract weighted energy estimte \eqref{abs-heat-weight} with $\gamma = \tfrac{n}{2}$ to \eqref{abstract-heat-diff}, 
        \begin{align*}
             \delta\mathrm{E}_{n + 1} + \delta\mathrm{B}_{n + 1} 
                &\lesssim \delta \widetilde{\mathrm{B}}_n + \bigl\| s \cdot s^\frac{n}{2} \cdot (\mathrm{I}_n + \mathrm{II}_n + \mathrm{III}_n ) \bigr\|_{L^1_{\frac{ds}{s}} L^2_x} + \bigl\| s \cdot \delta (s^\frac12 \bfD)^{(n)} \sfN \bigr\|_{L^1_{\frac{ds}{s}} L^2_x} \\
                &\lesssim \bigl\| \delta \overline \sfU \bigr\|_{L^2_x}  +  \underline s^\frac18 \cdot \cK \cdot \bigl( 1 + \underline s^\frac12 \bigl\| \overline A \bigr\|_{L^4_x} \bigr) \cdot \bigl\| \delta \overline A \bigr\|_{L^4_x} + \underline s^{\frac38} \cdot \cK \cdot \bigl\| \nabla \delta \overline A \bigr\|_{L^2_x} + \bigl\| s \, \delta (s^\frac12 \bfD)^{(\leq n)} \sfN \bigr\|_{L^1_{\frac{ds}{s}} L^2_x} ,
        \end{align*}
    using the inductive hypothesis to handle the appearance of the bulk term on the right. As $\delta \widetilde{\mathrm B}_{n + 1}$ can be controlled by $\delta \mathrm{B}_n$, we can complete the induction. 
\end{proof}

\begin{remark}
    The choice $s^{1/8}$ in the assumption \eqref{U-diff-K} is a bit arbitrary. Similarly, we could have replaced $\nabla \delta A$ in $L^2$ with the weaker quantity
        \[
            \|s^\frac14 \delta \nabla A \|_{L^\infty_s L^2_x}
        \]
    on the right-hand side. The point is that when commuting $\delta$ through the equation, there is quite a bit of room in terms of powers of $s$. 
\end{remark}

\section{Geometric Littlewood-Paley theory}\label{sec:littlewoodpaley}

The analytic heart of this article rests in constructing a smoothing scheme which (\texttt i) separates frequency scales, (\texttt{ii}) is gauge-covariant, and (\texttt{iii}) preserves the class of finite-energy configurations. One's first-take might be to use the standard Littlewood-Paley theory, as is ubiquitous in the study of PDEs and by design satisfies (\texttt{i}). However, it manifestly fails to satisfy (\texttt{ii}), and it is not obvious whether (\texttt{iii}) holds. To remedy these deficiencies, we develop a geometric analogue of the Littlewood-Paley theory based on the \textit{covariant heat equation} on $\R^2 \times [0, \underline s]$, 
    \begin{align}
        \big(\bfD_s - \bfD^\ell \bfD_\ell\big) \phi 
            &= 0, \label{eq:covariant-heat} \tag{P-$\phi$}\\
        F_{sj} 
            &= \partial^\ell F_{\ell j} \label{eq:covariant-heat2} \tag{P-$A$},
    \end{align}
The system of equations consists of the (abelian) Yang-Mills heat flow \eqref{covariant-heat2}, coupled with an electromagnetic heat equation for the scalar field \eqref{covariant-heat}. It is covariant under the gauge transformation 
    \begin{align*}
        A_\bfa 
            &\mapsto A_\bfa + \partial_\bfa \chi, \\
        \phi 
            &\mapsto e^{i \chi} \phi, 
    \end{align*}
for any real scalar field $\chi : \R^2 \times [0, \underline s] \to \R$. To fix the gauge (modulo $s$-independent transformations), we take inspiration from work on the hyperbolic Yang-Mills equation, e.g. \cite{Oh2014a,Oh2015a,Gavrus2022}, which in turn draws from Tao's work on wave maps \cite{Tao2004}, and impose the \textit{caloric gauge},
    \begin{equation} \tag{Cal}\label{eq:caloric-gauge}
        A_s 
            = 0.
    \end{equation}

We refer to $s \in [0, \underline s]$ as the \textit{heat-temporal} variable (in contrast to $t \in [0, T]$ the \textit{Schr\"odinger-temporal} variable). Given a fixed positive heat-time $\underline s > 0$, denote $s$-dependent functions and operators restricted to $s = 0$ by overlining, and similarly the restriction to $s = \underline s$ by underlining, e.g.
    \begin{alignat*}{2}
        \overline{\mathsf U} (x) 
            &:= \mathsf U(x, 0), \qquad
        \overline \bfD_j 
            &:= \partial_j - i A_j (x, 0),\\
        \underline{\mathsf U} (x) 
            &:= \mathsf U(x, \underline s),  \qquad
        \underline \bfD_j 
            &:= \partial_j - i A_j (x, \underline s).
    \end{alignat*}
We will often refer to $s$ close to $s = 0$ as ``high-frequency'', and $s$ close to $s = \underline s$ as ``low-frequency''. 

\begin{table}[!htbp]
\centering
\renewcommand{\arraystretch}{1.4}
    \resizebox{\textwidth}{!}{%
    \begin{tabular}{|l|c|c|c|}
    \hline
    & \textbf{Fourier multiplier} & \textbf{Linear heat flow} & \textbf{Covariant heat flow}
        \\ \hline
    Frequency scale & $N$ & $s^{-1/2}$ & $s^{-1/2}$
        \\ \hline
    Low frequency & $P_{\leq N} \overline\phi$ & $e^{s\Delta} \overline\phi$ & $\phi(s)$
        \\ \hline
    Medium frequency & $P_N \overline\phi$ & $-s \partial_s e^{s\Delta} \overline\phi$ & $-s \bfD_s \phi (s)$
        \\ \hline
    Scale-invariant measure & discrete measure on $2^\N$ &  $\frac{ds}{s}$ & $\frac{ds}{s}$
        \\ \hline
    Resolution of the identity & $\overline\phi = P_{\leq \underline N} \phi + \sum_{N > \underline N} P_N \phi$ & $\overline\phi = e^{\underline s \Delta} \phi +\int_0^{\underline s} (-s\partial_s e^{s \Delta} \phi) \, \frac{ds}{s}$ & $\overline\phi = \underline \phi + \int_0^{\underline s} (-s \bfD_s \phi) \, \frac{ds}{s} $
        \\ \hline
    \end{tabular}
    }

\vspace{1em}

    \caption{The reader may find the table above as a useful dictionary for translating between the Littlewood-Paley theory based on Fourier multipliers (see e.g. \cite[Appendix A]{Tao2006}), its formulation through the linear heat flow, and our approach via the covariant heat flow \eqref{covariant-heat}-\eqref{covariant-heat2} in caloric gauge \eqref{caloric-gauge}. }\label{table:dictionary-basic}
\end{table} 

The covariant heat equation plays three roles in our analysis. First, in analogy with the Littlewood-Paley theory, we will develop a covariant notion frequency-localisation for finite-energy configurations. Second, it is difficult to directly analyse the space of finite-energy configurations as it consists of rough non-linear objects. Instead, leveraging a Littlewood-Paley-type decomposition, we split the difficulty by identifying a finite-energy configuration with a smooth non-linear object and a superposition of rough linear objects. Third, we conduct the bulk of the analysis of the Chern-Simons-Schr\"odinger equation in a class of smooth solutions arising from smooth initial data. With suitable stability bounds at hand, we can construct rough solutions by smoothly approximating rough initial data using the covariant heat equation. Let us elaborate on the first two roles; the last one is fairly standard. 

To see how the covariant heat equation gives rise to a Littlewood-Paley-type theory, it is instructive to use the linear heat equation as a caricature. In particular, one should think of the covariant Hessian as 
    \[
        (s \bfD^{(2)}) \phi(s)
            \approx (s\Delta) e^{s \Delta} \overline \phi.
    \]
Viewing the right-hand side on the Fourier-side, we see that frequencies above $s^{-1/2}$ are damped exponentially by the symbol of the linear heat propagator, while frequencies below $s^{-1/2}$ are damped polynomially by the symbol of the Laplacian. Thus, one can regard it as Littlewood-Paley projection to frequency $s^{-1/2}$. This viewpoint goes back to the work of Stein \cite{Stein1970a}; see also \cite{Tao2004, KlainermanRodnianski2006,Oh2014a,Oh2015a} and the references therein. 

As evidenced by the examples in Appendix \ref{app:example}, a finite-energy configuration is a non-linear object, so, in some sense, one should view the space of finite-energy configurations as an infinite-dimensional manifold rather than a linear function space. It will be more convenient to identify finite-energy configurations with their flows under the covariant heat equation in caloric gauge. More precisely, we integrate the equation \eqref{covariant-heat} in caloric gauge \eqref{caloric-gauge} in $s$, which by the fundamental theorem of calculus yields the decomposition
    \[
        \overline\phi 
            = \underline \phi + \int_0^{\underline s} (s \bfD^{(2)}) \phi(s) \frac{ds}{s}.
    \]
This decomposes the scalar field into a ``low-frequency" part, which is smooth and captures the ``non-linear'' aspects of the configuration, and a superposition of ``high-frequency'' pieces, which are linear, and thus amenable to standard functional analytic tools, and capture the ``roughness''. Abstractly, the Littlewood-Paley decomposition ``linearises'' the non-linear function space -- it is instructive to compare with the setting of maps into manifolds, c.f. Tataru's lecture notes \cite{KochEtAl2014} and the references therein. 

To identify finite-energy configurations with solutions of the covariant heat equation, fix\footnote{The constant may also depend on the number of derivatives we are considering, i.e. $\epsilon \ll_n 1$ when proving an $n$-th order parabolic smoothing bound. In practice, we will only need finitely many derivatives, e.g. $n \leq 100000$, so we suppress this dependence. } $\epsilon \ll 1$ sufficiently small, and for\footnote{We will not need to use the fact that $\underline s^{-1}$ controls the Gauss tension field within this section. However, it will be convenient to keep it in the definition of $\frE(\underline s)$ to abbreviate some implicit constants when analysing the Chern-Simons-Schr\"odinger equation. } $\underline s > 0$, define the open neighborhood of the space of finite-energy configuration, 
    \begin{align*}
        \frE(\underline s) 
            := \Bigl\{ (A, \phi) \in \frE : 1 + \cE_{\text{AH}} [A, \phi] + \| \nabla \Gauss \|_{L^2_x}^2 < \epsilon \cdot \underline s^{-1} \Bigr\}.
    \end{align*}
These neighborhoods manifestly cover the space of finite-energy configurations, i.e. $\frE = \cup_{\underline s > 0} \frE(\underline s)$. By the well-posedness theory for the covariant heat equation in caloric gauge (Lemma \ref{lem:lwp-heat}), there is a one-to-one correspondence between the neighborhood $\frE(\underline s)$ and the \textit{space of $\frE$-caloric extensions} to $\R^2 \times [0, \underline s]$, 
    \begin{align*}
        \pzcE^1 ([0, \underline s])
            &:= \Bigl\{ (A, \phi) : \text{solution on $\R^2 \times [0, \underline s]$ to \eqref{covariant-heat}-\eqref{covariant-heat2} in \eqref{caloric-gauge} with $(A, \phi)_{|s = 0} \in \frE(\underline s)$} \Bigr\}.
    \end{align*}
Our primary goal herein is to further develop the correspondence in line with the previous discussion. The contents of each subsection are as follows: 
    \begin{enumerate}[label=\bf\arabic*.]
        \item[\ref{subsec:LP-sobolev}.] Capture a gauge-invariant notion of frequency localisation through covariant parabolic smoothing.
        
        \item[\ref{subsec:LP-topology}.] Give a Besov-type characterisation of the topology on the space of finite-energy configurations.
        
        \item[\ref{subsec:weak}.] Introduce a weaker distance functional on the space of finite-energy configurations, which we will show propagates under the Chern-Simons-Schr\"odinger flow in Section \ref{sec:difference}. 

        \item[\ref{subsec:env}.] Introduce the notion of a frequency envelope as a robust tool for capturing the distribution of energy across frequency-scales. 

        \item[\ref{subsec:LP-regularise}.] Approximate generic finite-energy configurations by a family of smooth configurations by showing that the covariant heat equation serves as a non-linear ``smooth approximation to the identity''. 
        
    \end{enumerate}

\subsection{Covariant energy and smoothing estimates}\label{subsec:LP-sobolev}

We put the abstract parabolic theory to work in proving two types of estimates for the covariant heat equation, namely energy estimates and parabolic smoothing. The energy estimates will show that \eqref{covariant-heat}-\eqref{covariant-heat2} propagates regularity; that is, if the heat-temporal initial data is a finite-energy configuration, then the flow remains finite-energy uniformly in $s$, and similarly if the heat-temporal initial data is a smooth configuration. The reader should view the parabolic smoothing as quantifying the heuristic that the heat-flow ``projects'' the heat-temporal initial data to ``low-frequencies'',
    \[
        \bfD\phi(s) = \text{projection of $\overline{\bfD\phi}$ to frequency $\lesssim s^{-1/2}$} \longleftrightarrow \bigl\| s^{\frac12 - \frac1p} (s^\frac12 \bfD)^{(n)} \bfD\phi(s) \bigr\|_{L^p_x} 
            \lesssim \bigl\| \overline{\bfD\phi} \bigr\|_{L^2_x}.
    \]
The estimates herein are completely gauge-invariant, and so goes for the notion of ``frequency-localisation''. It is instructive to compare with the geometric Littlewood-Paley theory of Klainerman-Rodnianski \cite{KlainermanRodnianski2006} and  the covariant smoothing estimates of Gavrus \cite[Section 6]{Gavrus2022} for the Yang-Mills heat flow.

To facilitate the analysis, we derive the parabolic equations of motion for the relevant geometric variables, 

\begin{lemma}[Covariant equations of motion]
    Let $(A, \phi)$ be a configuration obeying the covariant heat equation \eqref{covariant-heat}-\eqref{covariant-heat2}. Then the following parabolic equations hold for
    \begin{enumerate}
        \item the (spatially) differentiated scalar field:
            \begin{equation}\label{eq:Dphi-heat}
            \big(\bfD_s - \bfD^\ell \bfD_\ell \big)  \bfD \phi
                = F_{12} \cdot\bfD \phi,
            \end{equation}

        \item the covariant Hessian:
            \begin{equation}\label{eq:Dphi-heat-2}
                \bigl( \bfD_s - \bfD^\ell \bfD_\ell \bigr) \bfD^{(2)} \phi
                    = F_{12} \cdot \bfD^{(2)} \phi + \nabla F_{12} \cdot \bfD \phi,
            \end{equation}

        \item the charge density: 
            \begin{align}
                (\partial_s - \Delta) \tfrac12(1 - |\phi|^2) 
                &= \bfD \phi \cdot \bfD \phi ,\label{eq:charge-heat}
            \end{align}

        \item the Gauss tension field:
            \begin{equation}\label{eq:Gauss-heat-prelim}
                (\partial_s - \Delta) \Gauss 
                    = \bfD \phi \cdot \bfD \phi.
            \end{equation}
    \end{enumerate}
    
\end{lemma}

\begin{proof}
\leavevmode
    \begin{enumerate}
        \item This is an immediate consequence of the commutator identity \eqref{parabolic-commute}.
        
        \item This is an immediate consequence of the commutator identity \eqref{parabolic-commute} and the product rule. 

        \item By the covariant heat equation \eqref{covariant-heat} and product rule \eqref{product},
            \begin{align*}
                \partial_s \tfrac12(1 - |\phi|^2) 
                    &= - \Re(\overline \phi \cdot \bfD_s \phi) \\
                    &= - \Re(\overline \phi \bfD^\ell \bfD_\ell \phi) \\
                    &= \Delta \tfrac12(1 - |\phi|^2) + |\bfD \phi|^2.
            \end{align*}

        \item Subtracting the equations for the magnetic field \eqref{F12-heat} and the charge density \eqref{charge-heat} gives the result. 
    \end{enumerate}
    
\end{proof}

\begin{proposition}
    Let $(A, \phi)$ be a configuration on $\R^2 \times [0, \underline s]$ solving the covariant heat equation \eqref{covariant-heat}-\eqref{covariant-heat2} such that $(\overline A, \overline \phi) \in \frE(\underline s)$. Then for each non-negative integer $n \in \N_0$ and $2 \leq p \leq \infty$, we have: 
    \begin{itemize}
        \item \textup{(Parabolic smoothing for $\bfD \phi$).} 
        \begin{equation}\label{eq:Dphi-sobolev}
            \bigl\| s^{\frac12 - \frac1p} (s^\frac12 \bfD)^{(n)} \bfD \phi \bigr\|_{L^\infty_s L^p_x} + \bigl\| s^{\frac12 - \frac1p} (s^\frac12 \bfD)^{(n + 1)} \bfD \phi \bigr\|_{L^2_{\frac{ds}{s}} L^p_x} 
                \lesssim_{n} \bigl\| \overline{\bfD \phi} \bigr\|_{L^2_x}.
        \end{equation}

        \item \textup{(Parabolic Strichartz estimate for $\bfD \phi$).}
            \begin{equation}\label{eq:Dphi-strichartz-heat}
            \bigl\| s^{\frac12 - \frac1p} \bfD \phi  \bigr\|_{L^2_{\frac{ds}{s}} L^p_x}
                \lesssim_p \bigl\| \overline{\bfD \phi} \bigr\|_{L^2_x}, \qquad p \neq 2.
            \end{equation}

         \item \textup{(Gauss tension field parabolic smoothing).} 
            \begin{equation} \label{eq:Gauss-sobolev}
                \bigl\| s^{\frac12 - \frac1p} (s^\frac12 \nabla)^{(n)} \nabla \Gauss \bigr\|_{L^\infty_s L^p_x} +  \bigl\| s^{\frac12 - \frac1p} (s^\frac12 \nabla)^{(n + 1)} \nabla \Gauss \bigr\|_{L^2_{\frac{ds}{s}} L^p_x}
                    \lesssim_n \bigl\| \nabla \overline \Gauss \bigr\|_{L^2_x} + \bigl\| \overline{\bfD \phi} \bigr\|_{L^2_x}^2.
            \end{equation}

        \item \textup{(Uniform Gauss tension field bound).} 
            \begin{equation}\label{eq:Gauss-Linfty}
                \| \Gauss \|_{L^\infty_{s, x}}
                    \lesssim \bigl\| \overline \Gauss \bigr\|_{L^\infty_x} + \bigl\| \overline{\bfD \phi} \bigr\|_{L^2_x}^2.
            \end{equation}
    \end{itemize}

\end{proposition}

\begin{proof}[Proof of parabolic smoothing \eqref{Dphi-sobolev}]
    By the interpolation Lemma \ref{lem:interpolate}, it suffices to prove the case $p = 2$. We compute the derivatives of the right-hand side of the equation \eqref{Dphi-heat} for $\bfD \phi$, 
        \[
            (s^\frac12 \bfD)^{(n)} (\text{R.H.S.\eqref{Dphi-heat}}) 
                = \sum_{a + b = n} (s^\frac12 \nabla)^{(a)} F_{12} \cdot (s^\frac12 \bfD)^{(b)} \bfD \phi.
        \]
       Then, by H\"older, parabolic smoothing for the magnetic field \eqref{F12-sobolev}, and the assumption on the heat-temporal initial data $(\overline A, \overline{\vphantom{A}\phi}) \in \frE(\underline s)$,
         \begin{align*}
            \bigl\| s \cdot (s^\frac12 \bfD)^{(\leq n)} (\text{R.H.S.\eqref{Dphi-heat}}) \bigr\|_{L^1_{\frac{ds}{s}} L^2_x} 
                &\lesssim \sum_{a + b \leq n} \underline s^\frac12 \bigl\| s^\frac12 (s^\frac12 \nabla)^{(a)} F_{12} \bigr\|_{L^\infty_{s, x}} \bigl\| (s^\frac12 \bfD)^{(b)}  \bfD \phi \bigr\|_{L^\infty_s L^2_x}\\
                &\lesssim \underline s^{\frac12} \bigl\| \overline F_{12} \bigr\|_{L^2_x} \bigl\| (s^\frac12 \bfD)^{(\leq n)} \bfD \phi \bigr\|_{L^\infty_s L^2_x} \\
                &\lesssim \epsilon \cdot \bigl\| (s^\frac12 \bfD)^{(\leq n)} \bfD \phi \bigr\|_{L^\infty_s L^2_x} ,
            \end{align*}
    Inserting the calculation above into the abstract smoothing estimate \eqref{abs-smooth} applied to the equation \eqref{Dphi-heat}, 
        \begin{align*}
            \bigl\| (s^\frac12 \bfD)^{(\leq n)} \bfD \phi \bigr\|_{L^\infty_s L^2_x} + \bigl\| (s^\frac12 \bfD) (s^\frac12 \bfD)^{(\leq n)} \bfD \phi \bigr\|_{L^2_{\frac{ds}{s}} L^2_x}
                &\lesssim \bigl\| \overline{\bfD \phi}\bigr\|_{L^2_x} +  \bigl\| s \cdot (s^\frac12 \bfD)^{(\leq n)} (\text{R.H.S.\eqref{Dphi-heat}}) \bigr\|_{L^1_{\frac{ds}{s}} L^2_x} \\
                &\lesssim  \bigl\| \overline{\bfD \phi}\bigr\|_{L^2_x}  + \epsilon \cdot \bigl\| (s^\frac12 \bfD)^{(\leq n)} \bfD \phi \bigr\|_{L^\infty_s L^2_x}.
        \end{align*}
    If $\epsilon_n \ll 1$, then the second term on the right may be absorbed into the left-hand side. 
\end{proof}

\begin{proof}[Proof of parabolic Strichartz bound \eqref{Dphi-strichartz-heat}]
    Applying the maximum principle \eqref{maximum-principle} to the equation \eqref{Dphi-heat} for $\bfD \phi$ gives
            \begin{align*}
                \bigl\| \bfD \phi (s)\bigr\|_{L^p_x} 
                    &\leq \bigl\| e^{s \Delta} | \overline{\bfD \phi} | \bigr\|_{L^p_x} + \int_0^s  \Big\| e^{(s - s') \Delta} |F_{12} \cdot \bfD \phi| (s') \Big\|_{L^p_x} \, ds'  \\
                    &\lesssim \bigl\| \overline{\bfD \phi} \bigr\|_{L^2_x} + \Bigl\|\int_0^s  \Big\| e^{(s - s') \Delta} |F_{12} \cdot \bfD \phi| (s') \Big\|_{L^p_x} \, ds' \Bigr\|_{L^2_{\frac{ds}{s}}},
            \end{align*}
        Multiplying by $s^{1/2 - 1/p}$ and taking $L^2_{ds/s}$, we estimate the homogeneous flow by parabolic Strichartz \eqref{parabolic-strichartz},
            \[
                \bigl\| s^{\frac12 - \frac1p} e^{s \Delta} |\overline{\bfD \phi}| \bigr\|_{L^2_{\frac{ds}{s}} L^p_x}
                    \lesssim \bigl\| \overline{\bfD \phi} \bigr\|_{L^2_x}
            \]
        For the inhomogeneous flow, we apply the $(L^2 \to L^p)$-smoothing estimate \eqref{standard-heat} and the parabolic smoothing \eqref{F12-sobolev}, \eqref{Dphi-sobolev} to obtain 
            \begin{align*}
                \Bigl\| s^{\frac12 - \frac1p} \int_0^s  \Big\| e^{(s - s') \Delta} |F_{12} \cdot \bfD \phi| (s') \Big\|_{L^p_x} \, ds' \Bigr\|_{L^2_{\frac{ds}{s}}}
                    &\lesssim \Bigl\|   \int_0^s \Bigl( \frac{s}{s - s'}\Bigr)^{\frac12 - \frac1p} \, \bigl\| F_{12} (s') \bigr\|_{L^\infty_x} \bigl\| \bfD \phi(s') \bigr\|_{L^2_x} \, ds' \Bigr\|_{L^2_{\frac{ds}{s}}} \\
                    &\lesssim  \Bigl\| \int_0^s \Bigl( \frac{s}{s - s'}\Bigr)^{\frac12 - \frac1p} |s'|^{\frac12} \, \bigl\| \overline F_{12} \bigr\|_{L^2_x} \bigl\| \overline{\bfD \phi} \bigr\|_{L^2_x} \,  \frac{ds'}{s'} \Bigr\|_{L^2_{\frac{ds}{s}}}\\
                    &\lesssim \underline s^\frac12 \bigl\| \overline F_{12} \bigr\|_{L^2_x} \bigl\| \overline{\bfD \phi} \bigr\|_{L^2_x} \\
                    &\lesssim\bigl\| \overline{\bfD \phi} \bigr\|_{L^2_x}
            \end{align*}
    recalling the size of the heat-temporal initial data $(\overline A, \overline{\vphantom{A}\phi}) \in \frE(\underline s)$ in the last line. 
\end{proof}

\begin{proof}[Proof of Gauss tension field parabolic smoothing \eqref{Gauss-sobolev}]
    By the interpolation Lemma \ref{lem:interpolate}, it suffices to prove the case $p = 2$. We compute the derivatives of the right-hand side of the equation \eqref{Gauss-heat-prelim} for $\Gauss$ using the product rule \eqref{product}, 
        \[
            (s^\frac12 \nabla)^{(n)} \nabla (\text{R.H.S.\eqref{Gauss-heat-prelim}})
                = \sum_{a + b = n} (s^\frac12 \bfD)^{(a)} \bfD \phi \cdot (s^\frac12 \bfD)^{(b + 1)} \bfD \phi. 
        \]  
    Then, by H\"older, parabolic smoothing \eqref{Dphi-sobolev} and parabolic Strichartz \eqref{Dphi-strichartz-heat}, 
        \begin{align*}
            \bigl\|s \cdot  (s^\frac12 \nabla)^{(\leq n)} \nabla (\text{R.H.S.\eqref{Gauss-heat-prelim}}) \bigr\|_{L^1_{\frac{ds}{s}} L^2_x}
                &\lesssim \sum_{a + b \leq n} \bigl\| s^\frac12 (s^\frac12 \bfD)^{(a)} \bfD \phi \bigr\|_{L^2_{\frac{ds}{s}} L^\infty_x} \bigl\| (s^\frac12 \bfD)^{(b + 1)} \bfD \phi\bigr\|_{L^2_{\frac{ds}{s}} L^2_x} \\
                &\lesssim \bigl\| \overline{\bfD \phi} \bigr\|_{L^2_x}^2.
        \end{align*}
    Inserting this calculation into the abstract smoothing estimate \eqref{abs-smooth} applied to \eqref{Gauss-heat-prelim} gives \eqref{Gauss-sobolev}.
\end{proof}

\begin{proof}[Proof of Gauss tension field uniform bound \eqref{Gauss-Linfty}]
    Estimating Duhamel's formula to the equation for the Gauss tension field \eqref{Gauss-heat-prelim} in $L^\infty_x$, 
        \begin{align*}
            \| \Gauss (s) \|_{L^\infty_x}
                &\lesssim \|  \overline\Gauss \|_{L^\infty_x} + \int_0^s \| {s'}^\frac12 \bfD \phi (s') \|_{L^\infty_x}^2 \, \frac{ds'}{s'} \\
                &\lesssim \| \overline \Gauss \|_{L^\infty_x} + \| \overline{\bfD \phi} \|_{L^2_x}^2,
        \end{align*}
    using the parabolic Strichartz estimate \eqref{Dphi-strichartz-heat} in the second line. 
\end{proof}

\begin{remark}
    To illustrate some analogies with the standard Littlewood-Paley theory, one should compare the smoothing estimates \eqref{Dphi-sobolev} to the Sobolev-Bernstein estimates, c.f. \cite[Appendix A]{Tao2006}, and the parabolic Strichartz bound \eqref{Dphi-strichartz-heat} with the Sobolev embedding 
        \[
            \dot H^1 (\R^2) \hookrightarrow \dot B^{-2/p, p}_2 (\R^2) ,
        \]
    under the interpretation that $(s^\frac12 \bfD) \phi (s)$ is localised to medium frequencies $s^{-1/2}$.
\end{remark}

\begin{corollary}
    Let $(A, \phi)$ be a configuration on $\R^2 \times [0, \underline s]$ solving the covariant heat equation \eqref{covariant-heat}-\eqref{covariant-heat2} such that $(\overline A, \overline \phi) \in \frE(\underline s)$. Then we have the following bounds:
        \begin{itemize}
            \item \textup{(Charge density bound).} 
                \begin{equation}\label{eq:charge-sobolev}
                    \big\| 1 - |\phi|^2 \big\|_{L^\infty_s L^2_x} 
                        \lesssim \big\| 1 - |\overline\phi|^2 \big\|_{L^2_x} + \big\| \overline{\bfD \phi} \big\|_{L^2_x}.
                \end{equation}
                
            \item \textup{(Brezis-Gallou\"et-type bound).} For any $s_1 < s_2$, 
                \begin{equation}\label{eq:brezis-0}
                    \bigl\| \phi(s_1) \bigr\|_{L^\infty_x} 
                        \lesssim \bigl\| \phi(s_2) \bigr\|_{L^\infty_x} + \log^\frac12 \Bigl( \frac{s_2}{s_1} \Bigr) \cdot \bigl\| \overline{\bfD \phi} \bigr\|_{L^2_x}.
                \end{equation}
            In particular, 
                \begin{align}\label{eq:brezis}
                    \bigl\|\phi (s) \bigr\|_{L^\infty_{x}}
                        &\lesssim \underline s^{-\frac12} \Bigl( 1 + \log^\frac12 \Bigl( \frac{\underline s}{s} \Bigr) \Bigr)   .
                \end{align}
            and, for any $\gamma > 0$, 
                \begin{align}\label{eq:inefficient-smooth}
                    \bigl\| s^\gamma \phi(s) \bigr\|_{L^\infty_{s, x}} 
                        \lesssim \underline s^{\gamma - \frac12}. 
                \end{align}
        \end{itemize}
    
\end{corollary}

\begin{proof}[Proof of charge density bound \eqref{charge-sobolev}]
    We estimate the right-hand side of the equation \eqref{charge-heat}, 
        \begin{align*}
        \bigl\| s \cdot \text{R.H.S.\eqref{charge-heat}} \bigr\|_{L^1_{\frac{ds}{s}} L^2_x}
            &\lesssim \bigl\| s^\frac12 \bigr\|_{L^1_{\frac{ds}{s}}} \bigl\| s^\frac14 \bfD \phi \bigr\|_{L^\infty_s L^4_x}^2 \\
            &\lesssim \underline s^\frac12 \bigl\| \overline{\bfD \phi} \bigr\|_{L^2_x}^2 \\
            &\lesssim \bigl\| \overline{\bfD \phi} \bigr\|_{L^2_x}, 
        \end{align*}
    using smoothing \eqref{Dphi-sobolev} in the second line and the size of the heat-temporal initial data $(\overline A, \overline{\vphantom{A}\phi}) \in \frE(\underline s)$ in the third line. The estimate \eqref{charge-sobolev} follows from the energy estimate \eqref{abs-heat-energy} applied to the equation \eqref{charge-heat}. 
\end{proof}

\begin{proof}[Proof of Brezis-Gallou\"et-type bound \eqref{brezis-0}-\eqref{inefficient-smooth}]
    Rewriting $|\phi|$ using the fundamental theorem of calculus in $s$ and using the diamagnetic inequality \eqref{diamagnetic}, 
        \begin{align*}
            |\phi(s_1)|  
                \leq |\phi(s_2)| + \int_{s_1}^{s_2} |s \bfD_s \phi(s)| \frac{ds}{s}.
        \end{align*}
    For the ``high-frequency'' part, we use parabolic smoothing \eqref{Dphi-sobolev} and Cauchy-Schwarz, 
        \begin{align*}
             \int_{s_1}^{s_2} \|s \bfD^{(2)} \phi(s)\|_{L^\infty_s} \frac{ds}{s}
                &\lesssim \log^\frac12 \Bigl( \frac{s_2}{s_1} \Bigr) \bigl\| s^\frac12 ({s}^\frac12 \bfD) \bfD \bigr\|_{L^2_{\frac{ds}{s} } L^\infty_x} \\
                &\lesssim  \underline s^{-\frac12}  \log^\frac12 \Bigl( \frac{s_2}{s_1} \Bigr) .
        \end{align*}
    This proves \eqref{brezis-0}. 
    
    For \eqref{brezis}, take $s_1 = s$ and $s_2 = \underline s$; we estimate the ``low-frequency'' part using the Sobolev-type estimate \eqref{inefficientLinfty} and parabolic smoothing \eqref{Dphi-sobolev} and the charge density bound \eqref{charge-sobolev}, 
        \begin{align*}
            \bigl\|\underline \phi\bigr\|_{L^\infty_x} 
                &\lesssim  1 + \bigl\| 1 - |\underline \phi \bigr\|_{L^2_x}^\frac13 \bigl\| \underline \bfD^{(2)} \underline\phi \bigr\|_{L^2_x}^{\frac13} \\
                &\lesssim \underline s^{-\frac12}, 
        \end{align*}
    recalling the size of $\underline s^{-1}$ in the second line. 

    Multiplying \eqref{brezis} by $s^\gamma$ and recalling basic calculus, one easily sees \eqref{inefficient-smooth}. 
\end{proof}

\begin{proposition}[Higher-order energy estimates]
    Let $(A, \phi)$ be a sufficiently regular configuration on $\R^2 \times [0, \underline s]$ solving the covariant heat equation \eqref{covariant-heat}-\eqref{covariant-heat2} such that $(\overline A, \overline \phi) \in \frE(\underline s)$. Then, for each positiver integer $n \in \N$, we have the following bounds, 
    \begin{itemize}
        \item \textup{(Higher-order covariant derivative energy estimate).}
            \begin{equation}\label{eq:Dphi-sobolev-high}
            \bigl\| \bfD^{(n)} \bfD \phi \bigr\|_{L^\infty_s L^2_x} + \bigl\| (s^\frac12 \bfD) \bfD^{(n)} \bfD \phi \bigr\|_{L^2_{\frac{ds}{s}} L^2_x} 
                \lesssim  \bigl\| \overline\bfD^{(n)} \overline{\bfD \phi} \bigr\|_{L^2_x} + \bigl\| \nabla^{(n)} \overline F_{12} \bigr\|_{L^2_x} .
        \end{equation}

        \item \textup{(Higher-order Gauss tension field energy estimate).}
            \begin{equation}\label{eq:gauss-sobolev-high}
            \bigl\| \nabla^{(n)} \nabla \Gauss \bigr\|_{L^\infty_s L^2_x}
                \lesssim \bigl\| \nabla^{(n)} \nabla \overline \Gauss \bigr\|_{L^2_x} +  \bigl\| \overline{\bfD \phi} \bigr\|_{L^2_x} \Bigl(\bigl\| \overline\bfD^{(n)} \overline{\bfD \phi} \bigr\|_{L^2_x} + \bigl\| \nabla^{(n)} \overline F_{12} \bigr\|_{L^2_x}\Bigr) .
        \end{equation}
    \end{itemize}
    
\end{proposition}

\begin{proof}
    Writing the equations \eqref{F12-heat}, \eqref{Dphi-heat} for $\sfU = \{ F_{12}, \bfD \phi\}$ schematically, it follows from the product rule and the commutator identity \eqref{commute} that 
        \[
            (\sfD_s - \sfD^\ell \sfD_\ell) \sfD^{(n)} \sfU 
                = \sum_{a + b = n} \sfD^{(a)} \sfU \cdot \sfD^{(b)}\sfU. 
        \]
    Applying the abstract energy estimate \eqref{abs-heat-energy} to the equation above, 
        \begin{align*}
            \bigl\| \sfD^{(n)} \sfU \bigr\|_{L^\infty_s L^2_x} + \bigl\|  (s^\frac12 \sfD)\sfD^{(n)} \sfU \bigr\|_{L^2_{\frac{ds}{s}} L^2_x} 
                &\lesssim \bigl\| \overline \sfD^{(n)} \overline \sfU \bigr\|_{L^2_x} + \sum_{a + b = n} \bigl\| s \, \sfD^{(a)} \sfU \cdot \sfD^{(b)} \sfU  \bigr\|_{L^1_{\frac{ds}{s}} L^2_x}\\ 
                &\lesssim \bigl\| \overline \sfD^{(n)} \overline \sfU \bigr\|_{L^2_x} + \sum_{a + b = n} \underline s^{\frac12} \bigl\| s^{\frac{a}{2n}}\sfD^{(a)} \sfU \bigr\|_{L^\infty_s L^{2n/a}_x} \bigl\| s^{\frac{b}{2n}} \sfD^{(b)} \sfU  \bigr\|_{L^\infty_s L^{2n/b}_x} \\
                &\lesssim \bigl\| \overline \sfD^{(n)} \overline \sfU \bigr\|_{L^2_x} + \underline s^{\frac12} \Bigl( \bigl\| (s^\frac12 \sfD) \sfU \bigr\|_{L^\infty_s L^2_x} + \bigl\| s^\frac12 \sfU \bigr\|_{L^\infty_{s, x}} \Bigr) \bigl\| \sfD^{(n)} \sfU \bigr\|_{L^\infty_s L^2_x} \\
                &\lesssim \bigl\| \overline \sfD^{(n)} \overline \sfU \bigr\|_{L^2_x} + \underline s^{\frac12} \bigl\| \overline \sfU \bigr\|_{L^2_x} \bigl\| \sfD^{(n)} \sfU \bigr\|_{L^\infty_s L^2_x},
        \end{align*}
    applying H\"older in the second line, interpolating appropriately in the third line, and the parabolic smoothing estimates \eqref{F12-sobolev}, \eqref{Dphi-sobolev} in the last line. Since $(\overline A, \overline \phi) \in \frE(\underline s)$, the second term in the last line can be absorbed into the left-hand side, furnishing \eqref{Dphi-sobolev-high}.  

    Applying the product rule to the equation for the Gauss tension field \eqref{Gauss-heat-prelim}, 
        \[
            (\partial_s - \Delta) \nabla^{(n)} \nabla \Gauss 
                = \sum_{a + b = n + 1} \bfD^{(a)} \bfD \phi \cdot \bfD^{(b)} \bfD \phi. 
        \]
    Then 
        \begin{align*}
            \bigl\| \nabla^{(n)} \nabla \Gauss  \bigr\|_{L^\infty_s L^2_x} &+  \bigl\| (s^\frac12 \nabla)\nabla^{(n)} \nabla \Gauss  \bigr\|_{L^2_{\frac{ds}{s}} L^2_x} 
                \\
                &\lesssim \bigl\| \nabla^{(n)} \nabla \overline \Gauss \bigr\|_{L^2_x}  + \sum_{a + b = n + 1} \bigl\| s \, \bfD^{(a)} \bfD \phi \cdot \bfD^{(b)} \bfD \phi \bigr\|_{L^1_{\frac{ds}{s}} L^2_x} \\
                &\lesssim \bigl\| \nabla^{(n)} \nabla \overline \Gauss \bigr\|_{L^2_x} + \Bigl( \bigl\| s^\frac12 \bfD \phi \bigr\|_{L^2_{\frac{ds}{s}} L^\infty_x} + \bigl\| (s^\frac12 \bfD) \bfD \phi \bigr\|_{L^2_{\frac{ds}{s}} L^2_x} \Bigr) \bigl\| (s^\frac12 \bfD) \bfD^{(n)} \bfD \phi \bigr\|_{L^2_{\frac{ds}{s}} L^2_x} \\
                &\lesssim \bigl\| \nabla^{(n)} \nabla \overline \Gauss \bigr\|_{L^2_x} + \bigl\| \overline{\bfD \phi} \bigr\|_{L^2_x} \Bigl( \bigl\| \overline \bfD^{(n)} \overline{\bfD \phi} \bigr\|_{L^2_x} + \bigl\| \nabla^{(n)} \overline F_{12} \bigr\|_{L^2_x} \Bigr),
        \end{align*}
    using a similar circle of ideas as in the previous proof, along with the estimate \eqref{Dphi-sobolev-high} in the last line. This furnishes \eqref{gauss-sobolev-high}. 
\end{proof}

\subsection{Littlewood-Paley characterisation of the energy topology}\label{subsec:LP-topology}

While the topology on the space of finite-energy configurations introduced in Definition \ref{def:topology-1} conveniently renders the abelian Higgs energy a continuous functional, its defining metric will be cumbersome to directly propagate under the Chern-Simons-Schr\"odinger flow. Instead, we give a more refined description of the metric by identifying the configuration on $\R^2$ with its caloric extension to $\R^2 \times [0, \underline s]$. 

By the analysis from Sections \ref{subsec:LP-sobolev}-\ref{subsec:env}, the covariant heat equation induces a separation of frequency-scales -- at the level of the scalar field, which is a rough non-linear object, we can integrate its equation \eqref{covariant-heat} in $s$ under caloric gauge \eqref{caloric-gauge}, then the fundamental theorem of calculus yields the decomposition 
    \[
        \overline \phi 
            = \underline \phi + \int_0^{\underline s} (s \bfD^{(2)}) \phi (s) \frac{ds}{s}.
    \]
One should view $\underline \phi$ as a smooth non-linear part of $\overline \phi$, and $(s \bfD^{(2)}) \phi(s)$ as high-frequency linearised objects which can be estimated in $L^2$-based spaces. We can similarly integrate the equation \eqref{Dphi-heat} for $\bfD \phi$,  
    \[
        \overline{\bfD \phi} = \underline{\bfD \phi} + \int_0^{\underline s}  (s \bfD^{(2)}) \bfD \phi (s) \frac{ds}{s} + \text{lower-order terms}.
    \]
By the analysis from Sections \ref{subsec:LP-sobolev}-\ref{subsec:env}, one should think of $\underline{\bfD\phi}$ as the localisation of $\overline{\bfD \phi}$ to frequencies below $\underline s^{-1/2}$, and $(s \bfD^{(2)}) \bfD \phi(s)$ as the localisation of $\overline{\bfD \phi}$ to frequencies at $s^{-1/2}$. Thus, at least heuristically, 
    \[
        \bigl\| \delta \overline{\bfD \phi} \bigr\|_{L^2_x}^2
            \approx \bigl\| \delta \underline{\bfD \phi} \bigr\|_{L^2_x}^2 + \int_0^{\underline s} \bigl\|s^{-\frac12} \delta (s \bfD^{(2)}) \phi (s)\bigr\|_{L^2_x}^2 \frac{ds}{s}.
    \]
This motivates the following definition, 

\begin{definition}\label{def:dist-2}
    Given $(\overline A, \overline \phi) \in \frE(\underline s)$ a finite-energy configuration on $\R^2$, we say that a configuration $(A ,\phi)$ on $\R^2 \times [0, \underline s]$ is its \textit{$\frE$-caloric extension} if it is the unique solution to the covariant heat equation \eqref{covariant-heat}-\eqref{covariant-heat2} in caloric gauge \eqref{caloric-gauge} with heat-temporal initial data $(A, \phi)_{|s = 0} = (\overline A, \overline{\vphantom{A}\phi})$. We denote $\pzcE^1 ([0, \underline s])$ for the space of $\frE$-caloric extensions, endowed with the topology induced by the following maps:  
        \begin{enumerate}[label=({\roman*})]
            \item (derivatives of) the covariant Hessian, 
                \begin{align*}
                    \pzcE^1 ([0, \underline s]) 
                        &\longrightarrow L^2_{\frac{ds}{s}} L^2_x (\R^2 \times [0, \underline s]),\\
                    (A, \phi) 
                        &\longmapsto  s^{-\frac12} (s^\frac12 \bfD)^{(\leq 2026)} (s \bfD^{(2)}) \phi, 
                \end{align*}

            \item the Gauss tension field at $s = 0$, 
                \begin{align*}
                    \pzcE^1 ([0, \underline s]) 
                        &\longrightarrow (L^\infty \cap H^1) (\R^2), \\
                    (A, \phi) 
                        &\longmapsto \overline \Gauss,
                \end{align*}

            \item (derivatives of) the covariant gradient at $s = \underline s$, 
                \begin{align*}
                    \pzcE^1 ([0, \underline s]) 
                        &\longmapsto L^2 (\R^2),\\
                    (A, \phi) 
                        &\longmapsto\underline{\bfD}^{(\leq 100)} \underline{\bfD\phi},
                \end{align*}

            \item the magnetic potential at $s = 0$ and scalar field at $s = \underline s$,
                \begin{align*}
                    \pzcE^1 ([0, \underline s]) 
                        &\longrightarrow (\dot H^1 \cap L^4) (\R^2) \times L^\infty (\R^2) , \\
                    (A, \phi)
                        &\longmapsto(\overline A, \underline \phi).
                \end{align*}
        \end{enumerate}
    Furthermore, define the distance functional
        \begin{align*}
            \dist_{\pzcE^1([0, \underline s])} \bigl( (A, \phi), (A', \phi') \bigr)
                &:=  \bigl\| s^{-\frac12}\delta (s^\frac12 \bfD)^{(\leq 2026)} (s \bfD^{(2)})\phi \bigr\|_{L^2_{\frac{ds}{s}} L^2_x} + \bigl\| \delta \overline \Gauss \bigr\|_{H^1_x} \\
                &\qquad + \bigl\| \delta \underline \bfD^{(\leq 100)} \underline{\bfD \phi} \bigr\|_{L^2_x} +  \bigl\| \delta \overline A \bigr\|_{(\dot H^1 \cap L^4)_x}  + \bigl\| \delta \underline \phi \bigr\|_{L^\infty_x}.
        \end{align*}

\end{definition}

\begin{remark}
    The idea of taking smooth extensions of rough functions and constructing function spaces which abstractly capture frequency localisation originates, to our knowledge, with Tataru's work on wave maps \cite[Section 3]{Tataru2005}. In this case, we have a distinguished smooth extension via the covariant heat equation, so our description of the energy space is also reminiscent of Tao's construction in \cite{Tao2009} of the energy space for maps into hyperbolic space $\R^2 \to \HH^m$ via the harmonic map heat flow. 
\end{remark}

\begin{lemma}[Well-posedness of the covariant heat equation]\label{lem:lwp-heat}
    Let $(\overline A, \overline{\vphantom{A}\phi}) \in \frE$ be a finite-energy configuration. Then there exists a unique configuration $(A, \phi)$ on $\R^2 \times [0, \infty)$ solving the covariant heat equation \eqref{covariant-heat}-\eqref{covariant-heat2} in caloric gauge \eqref{caloric-gauge} with heat-temporal initial data $(A, \phi)_{|s = 0} = (\overline A, \overline{\vphantom{A}\phi})$. 
\end{lemma}

\begin{proof}
    It is not difficult to see that $A := A^\df + A^\cf$, where $A^\df$ solves the heat equation \eqref{A-heat} and $A^\cf$ solves the transport equation \eqref{A-transport}, is a solution to the Yang-Mills heat flow \eqref{YMH} in caloric gauge \eqref{caloric-gauge}, and this solution is unique in $(L^4 \cap \dot H^1) (\R^2)$. Fixing then the magnetic potential $A$, we regard \eqref{covariant-heat} as a linear heat equation for $\phi$. By the maximum principle and duality, one can show that the equation is well-posed in $(L^2 + L^\infty) (\R^2)$. We leave the technical verifications to the reader. 
\end{proof}

By the well-posedness of the covariant heat equation, there is a one-to-one correspondence between the space of $\frE$-caloric extensions $\pzcE^1 ([0, \underline s])$ to $\R^2 \times [0, \underline s]$ and the open neighborhood of the space of finite-energy configurations $\frE(\underline s)$. Our goal herein is to show that the correspondence is also locally bi-Lipschitz, and so one can regard the metric on $\pzcE^1 ([0, \underline s])$ as a suitable substitute for that on $\frE(\underline s)$, 

\begin{proposition}[Equivalence of topologies]\label{thm:topology}
    The map
        \begin{align*}
            \pzcE^1 ([0, \underline s])
                &\longrightarrow \frE(\underline s),\\
            (A, \phi)
                &\longmapsto (\overline A, \overline{\vphantom{A}\phi}),
        \end{align*}
    is bijective, and satisfies the (local) bi-Lipschitz bound 
        \[
            \dist_{\pzcE^1([0, \underline s])} \bigl( (A, \phi), (A', \phi') \bigr)
                \sim_{\underline s, \| \overline A \|_{L^4}} \dist_{\frE} \bigl( (\overline A, \overline{\vphantom{A} \phi}), (\overline{A}', \overline{\vphantom{A}\phi}') \bigr).
        \]
\end{proposition}

We begin by studying the caloric extension map $\frE(\underline s) \to \pzcE^1 ([0, \underline s])$. Comparing Definitions \ref{def:topology-1} and \ref{def:dist-2}, we only need to estimate the contributions 
    \[
        \bigl\| s^{-\frac12}\delta (s^\frac12 \bfD)^{(\leq 2026)} (s \bfD^{(2)})\phi \bigr\|_{L^2_{\frac{ds}{s}} L^2_x} , \qquad \bigl\| \delta \underline \phi \bigr\|_{L^\infty_x},
    \]
from the metric on $\pzcE^1 ([0, \underline s])$. This map is the data-to-solution map for the covariant heat equation in caloric gauge, so we proceed by parabolic smoothing and energy arguments applied to the difference equations, 

\begin{proposition}[Lipschitz continuity of caloric extension map]
   The caloric extension map 
        \begin{align*}
            \frE(\underline s)
                &\longrightarrow \pzcE^1 ([0, \underline s]) , \\
            (\overline A, \overline{\vphantom{A}\phi})
                &\longmapsto (A, \phi),
        \end{align*}
    is locally Lipschitz continuous. More precisely, let $(A, \phi), (A' \phi') \in \pzcE^1([0, \underline s])$ be $\frE$-caloric extensions to $\R^2 \times [0, \underline s]$, then the following difference bounds hold:
    \begin{itemize}
        \item \textup{(Scalar field difference at $s = \underline s$).}
        \begin{equation}\label{eq:phi-extend}
            \bigl\| \delta \underline \phi \bigr\|_{L^\infty_x} 
                \lesssim_{\underline s, \| (\overline A, \overline A') \|_{L^4_x}} \bigl\| \delta \overline \phi \bigr\|_{(L^2 + L^\infty)_x} + \bigl\| \delta \overline A \bigr\|_{(\dot H^1 \cap L^4)_x}.
        \end{equation}

        \item \textup{(Parabolic smoothing differences).} For each non-negative integer $n \in \N_0$, 
            \begin{equation}\label{eq:heat-lip}
            \begin{split}
            \bigl\| \delta (s^\frac12 \bfD)^{(n)} \bfD \phi\bigr\|_{L^\infty_s L^2_x} + \bigl\| \delta (s^\frac12 \bfD)^{(n + 1)} \bfD \phi \bigr\|_{L^2_{\frac{ds}{s}}L^2_x} 
                &\lesssim_{n, \underline s, \| (\overline A, \overline A') \|_{L^4_x}} \bigl\| \delta \overline{\bfD \phi} \bigr\|_{L^2_x} +  \bigl\| \delta \overline A \bigr\|_{(\dot H^1 \cap L^4)_x}.
            \end{split}
        \end{equation}  

    \end{itemize}
\end{proposition}

\begin{proof}[Proof of scalar field at $s = \underline s$ difference bound \eqref{phi-extend}]
    By the maximum principle \eqref{maximum-principle},
        \begin{align*}
            \bigl\| \delta \underline \phi \bigr\|_{L^\infty_x} 
                &\lesssim \bigl\| e^{\underline s \Delta} |\delta \overline \phi| \bigr\|_{L^\infty_x} + \int_0^{\underline s} \bigl\| e^{(\underline s - s) \Delta} (\partial_s - \bfD^\ell \bfD_\ell) \delta \phi (s) \bigr\|_{L^\infty_x} \, ds \\
                &\lesssim \underline s^{-\frac12} \bigl\| \delta \overline \phi \bigr\|_{L^2 + L^\infty} + \int_0^{\underline s} |\underline s - s|^{-\frac12}  \bigl\| \bigl[ \delta, \bfD^\ell \bfD_\ell \bigr] \phi (s) \bigr\|_{L^2_x} \, ds  \\
                &\lesssim \underline s^{-\frac12} \bigl\| \delta \overline \phi \bigr\|_{L^2 + L^\infty} + \int_0^{\underline s} \Bigl( \frac{s}{\underline s - s} \Bigr)^\frac12 \cdot s^\frac12 \cdot\Bigl( \| A \|_{L^4} \| \bfD \phi \|_{L^4} +\bigl( \| \delta \partial^\ell A_\ell \|_{L^2} + \| \delta A \|_{L^4} \| A \|_{L^4}\bigr) \|\phi \|_{L^\infty} \Bigr)\, \frac{ds}{s}\\
                &\lesssim_{\underline s, \|\overline A \|_{L^4}}  \bigl\| \delta \overline \phi \bigr\|_{L^2 + L^\infty}  + \bigl\| \overline A \bigr\|_{\dot H^1 \cap L^4} ,
        \end{align*}
    using $(L^2 \to L^\infty)$-parabolic smoothing \eqref{standard-heat} in the second line, the commutator identity \eqref{commute-diff-laplace} and H\"older in the third line, and parabolic smoothing \eqref{Dphi-sobolev} and \eqref{inefficient-smooth} in the last line. 
\end{proof}

\begin{proof}[Proof of smoothing difference bound \eqref{heat-lip}]
    Denote 
        \begin{align*}
            \delta \mathrm{E}_n 
                &:= \bigl\| \delta (s^\frac12 \bfD)^{(n - 1)} \bfD \phi \bigr\|_{L^\infty_s L^2_x},\\
            \delta \widetilde{\mathrm{B}}_n
                &:= \bigl\| \delta (s^\frac12 \bfD)^{(n)} \bfD \phi \bigr\|_{L^2_{\frac{ds}{s}} L^2_x}.
        \end{align*}
    We compute the differences of the equations for $\bfD \phi$ and derivatives thereof -- by the product rule, 
        \begin{align*}
            \delta (s^\frac12 \bfD)^{(n)} \bigl( \text{R.H.S.\eqref{Dphi-heat}}\bigr)
                &= \sum_{a + b = n} (s^\frac12 \nabla)^{(a)} \delta F_{12} \cdot (s^\frac12 \bfD)^{(b)} \bfD \phi + \sum_{a + b = n} (s^\frac12 \nabla)^{(a)}  F_{12} \cdot \delta (s^\frac12 \bfD)^{(b)} \bfD \phi\\
                &=: \mathrm{I}_n + \mathrm{II}_n. 
        \end{align*}
    Using parabolic smoothing \eqref{F12-sobolev}, \eqref{Dphi-sobolev}, 
        \begin{align*}
            \bigl\| s \cdot \mathrm{I}_{\leq n} \bigr\|_{L^1_{\frac{ds}{s}} L^2_x}
                &\lesssim \underline s^\frac12 \bigl\| s^\frac12 (s^\frac12 \nabla)^{(\leq n)} \delta F_{12}  \bigr\|_{L^\infty_{s, x}} \bigl\| (s^\frac12 \bfD)^{(\leq n)} \bfD \phi \bigr\|_{L^\infty_s L^2_x} \\
                &\lesssim \underline s^\frac12  \bigl\| \delta \overline F_{12} \bigr\|_{L^2_x} \bigl\| \overline{\bfD \phi} \bigr\|_{L^2_x} \\
                &\lesssim \bigl\| \delta \overline F_{12} \bigr\|_{L^2_x}, 
        \end{align*}
    and similarly,
        \begin{align*}
            \bigl\| s \cdot \mathrm{II}_n \bigr\|_{L^1_{\frac{ds}{s}} L^2_x} 
                &\lesssim \underline s^\frac12 \bigl\| s^\frac12 (s^\frac12 \nabla)^{(\leq n)} F_{12}  \bigr\|_{L^\infty_{s, x}} \bigl\| \delta (s^\frac12 \bfD)^{(\leq n)} \bfD \phi \bigr\|_{L^\infty_s L^2_x} \\
                &\lesssim \underline s^\frac12  \bigl\| \delta \overline F_{12} \bigr\|_{L^2_x} \bigl\| \delta (s^\frac12 \bfD)^{(\leq n)} \bfD \phi \bigr\|_{L^\infty_s L^2_x} \\
                &\lesssim \epsilon \cdot \mathrm{E}_{\leq n + 1},
        \end{align*}
    using the definition of $\frE(\underline s)$ in the last line. Furthermore, parabolic smoothing \eqref{Dphi-heat} implies
        \begin{align*}
            \bigl\| s^{\frac18} (s^\frac12 \bfD)^{(n)} \bfD \phi \bigr\|_{L^\infty_s L^2_x}
                \lesssim \underline s^{-\frac38}. 
        \end{align*}
    Inserting these calculations into the abstract difference estimates from Lemma \ref{lem:diff-smooth} applied to \eqref{Dphi-heat}, 
        \begin{align*}
            \mathrm{E}_{\leq n + 1} + \widetilde{\mathrm{B}}_{\leq n + 1} 
                &\lesssim \bigl\| \delta \overline{\bfD \phi} \bigr\|_{L^2_x} + \bigl\| s \, \delta (s^\frac12 \bfD)^{(\leq n - 1)} \bigl( \text{R.H.S.\eqref{Dphi-heat}} \bigr) \bigr\|_{L^1_{\frac{ds}{s}} L^2_x} \\
                    &\qquad +  \underline s^{-\frac14} \cdot \bigl\| \delta \overline A \bigr\|_{L^4_x} + \bigl\| \overline A \bigr\|_{L^4_x} \cdot \bigl\| \delta \overline A \bigr\|_{L^4_x} + \bigl\| \nabla \delta \overline A \bigr\|_{L^2_x} \\
                &\lesssim \bigl\| \delta \overline{\bfD \phi} \bigr\|_{L^2_x}  + \epsilon \cdot \mathrm E_{\leq n + 1} \\
                    &\qquad +  \underline s^{-\frac14} \cdot \bigl\| \delta \overline A \bigr\|_{L^4_x} + \bigl\| \overline A \bigr\|_{L^4_x} \cdot \bigl\| \delta \overline A \bigr\|_{L^4_x} + \bigl\| \nabla \delta \overline A \bigr\|_{L^2_x} .
        \end{align*}
    Since $\epsilon \ll 1$, the energy-type term on the right may be absorbed into the left. This completes the proof. 
\end{proof}

We now look at the caloric trace map $\pzcE^1 ([0, \underline s]) \to \frE(\underline s)$. Again, comparing Definitions \ref{def:topology-1} and \ref{def:dist-2}, we only need to estimate the contributions 
    \[
        \bigl\| \delta \overline{\bfD \phi} \bigr\|_{L^2_x}, \qquad \bigl\| \delta \overline \phi \bigr\|_{(L^2 + L^\infty)_x},
    \]
from the metric on $\frE(\underline s)$. We proceed by almost-orthogonality, leveraging the ``separation of scales'' implicit in the definition of the metric on the space of caloric extensions,

\begin{proposition}[Lipschitz continuity of caloric trace map]
    The caloric trace map 
        \begin{align*}
            \pzcE^1 ([0, \underline s]) 
                &\longrightarrow \frE(\underline s), \\
            (A, \phi) 
                &\longmapsto (\overline A, \overline{\vphantom{A}\phi}),
        \end{align*}
    is locally Lipschitz continuous. More precisely, let $(A, \phi), (A' \phi') \in \pzcE^1([0, \underline s])$ be $\frE$-caloric extensions to $\R^2 \times [0, \underline s]$, then the following difference bounds hold: 
    \begin{itemize}
        \item \textup{(Scalar field difference).}
            \begin{equation}\label{eq:phi-trace}
            \bigl\| \delta \overline \phi \bigr\|_{L^2 + L^\infty} 
                \lesssim \bigl\| \delta \underline \phi \bigr\|_{L^\infty_x} + \bigl\| s^{-\frac12} \delta (s \bfD^{(2)}) \phi\bigr\|_{L^2_{\frac{ds}{s}} L^2_x}.
            \end{equation}

        \item \textup{(Covariant gradient difference).}
            \begin{equation}\label{eq:Dphi-diff-plancharel}
                \bigl\| \delta \overline{\bfD \phi} \bigr\|_{L^2_x} 
                    \lesssim_{\underline s} \bigl\| \delta \overline A \bigr\|_{(\dot H^1\cap L^4)_x} + \bigl\| \delta \underline{\bfD \phi} \bigr\|_{L^2_x} +  \bigl\| s^{-\frac12} \delta (s^\frac12 \bfD)^{(\leq 2)} (s \bfD^{(2)}) \phi \bigr\|_{L^2_{\frac{ds}{s}} L^2_x}.
            \end{equation}

    \end{itemize}

\end{proposition}

\begin{proof}[Proof of scalar field trace difference \eqref{phi-trace}]
    Integrating the covariant heat equation \eqref{covariant-heat} under \eqref{caloric-gauge}, 
        \[
            \delta \overline \phi 
                = \delta \underline \phi + (s \partial_s)^{-1} \delta (s \bfD^{(2)}) \phi
        \]
    Placing the high-frequency part in $L^2$ and the low-frequency part in $L^\infty$, we obtain 
        \begin{align*}
            \bigl\| \delta \overline \phi \bigr\|_{L^2 + L^\infty} 
                &\lesssim \bigl\| \delta \underline \phi \bigr\|_{L^\infty_x} + \int_0^{\underline s} \bigl\| \delta (s\bfD^{(2)}) \phi \bigr\|_{L^2_x} \, \frac{ds}{s} \\
                &\lesssim \bigl\| \delta \underline \phi \bigr\|_{L^\infty_x} + \bigl\| s^{-\frac12} \delta (s \bfD^{(2)}) \phi\bigr\|_{L^2_{\frac{ds}{s}} L^2_x},
        \end{align*}
    using Cauchy-Schwarz in the second line. 
\end{proof}

Before proceeding to the proof of the covariant gradient difference \eqref{Dphi-diff-plancharel}, it will be convenient to record the following two lemmas.  

\begin{lemma}
    Let $(A, \phi), (A', \phi') \in \pzcE^1 ([0, \underline s])$ be $\frE$-caloric extensions, then for any $s' \in [0, \underline s]$,
        \begin{equation}\label{eq:Dphi-s-diff-1}
            \bigl\|\delta \bfD \phi (s') \bigr\|_{L^2_x} 
                \lesssim \bigl\| \delta \underline{\bfD \phi} \bigr\|_{L^2_x} + \Bigl\| \int_{s'}^{\underline s} s^{-\frac12} \delta (s^\frac12 \bfD) (s \bfD^{(2)}) \phi \frac{ds}{s} \Bigr\|_{L^2_x} + \bigl\| \delta \overline F_{12} \bigr\|_{L^2_x}.
        \end{equation}
    In particular, 
        \begin{equation}\label{eq:Dphi-s-diff-2}
            \bigl\| \delta \bfD \phi (s') \bigr\|_{L^2_x} 
                \lesssim \bigl\| \delta \underline{\bfD \phi} \bigr\|_{L^2_x} + \log^\frac12 \Bigl( \frac{\underline s}{s'} \Bigr) \bigl\| s^{-\frac12} \delta (s^\frac12 \bfD) (s \bfD^{(2)}) \phi \bigr\|_{L^2_{\frac{ds}{s}}L^2_x} + \bigl\| \delta \overline F_{12} \bigr\|_{L^2_x}.
        \end{equation}
\end{lemma}

\begin{proof}
    The second estimate \eqref{Dphi-s-diff-2} is an immediate consequence of the first \eqref{Dphi-s-diff-1} and Cauchy-Schwarz in $s$. To prove the first estimate, take the differences of the equations for the covariant gradients \eqref{Dphi-heat} in caloric gauge \eqref{caloric-gauge} and integrate in $s$. It follows from the fundamental theorem of calculus that 
        \begin{align*}
            \delta \bfD \phi (s') 
                &= \delta \underline{\bfD \phi} + \int_{s'}^{\underline s} s^{-\frac12} \delta (s^\frac12 \bfD) (s \bfD^{(2)}) \phi \frac{ds}{s} + \int_{s'}^{\underline s} s^\frac12 \cdot s^\frac12 \delta F_{12} \cdot \bfD \phi \, \frac{ds}{s}\\
                    &\qquad + \int_{s'}^{\underline s} s^\frac12 \cdot s^\frac12 F_{12} \cdot \delta \bfD \phi \, \frac{ds}{s}. 
        \end{align*}
    Estimating in $L^2_x$, we obtain 
        \begin{align*}
        \bigl\| \delta \bfD \phi \bigr\|_{L^\infty_{s} L^2_x ([s', \underline s]) } 
            &\lesssim \bigl\| \delta \underline{\bfD \phi} \bigr\|_{L^2_x} +\Bigl\| \int_{s'}^{\underline s} s^{-\frac12} \delta (s^\frac12 \bfD) (s \bfD^{(2)}) \phi \frac{ds}{s} \Bigr\|_{L^2_x} + \int_{s'}^{\underline s} s^\frac12 \cdot \bigl\|  s^\frac12 \delta F_{12} \bigr\|_{L^\infty_x} \bigl\| \bfD \phi \bigr\|_{L^2_x} \frac{ds}{s}  \\
                &\qquad + \int_{s'}^{\underline s} s^\frac12 \cdot \bigl\| s^\frac12 F_{12} \bigr\|_{L^\infty_x} \bigl\| \delta \bfD \phi \bigr\|_{L^2_x} \frac{ds}{s} \\
            &\lesssim \bigl\| \delta \underline{\bfD \phi} \bigr\|_{L^2_x} + \Bigl\| \int_{s'}^{\underline s} s^{-\frac12} \delta (s^\frac12 \bfD) (s \bfD^{(2)}) \phi \frac{ds}{s} \Bigr\|_{L^2_x} + \underline s^\frac12 \bigl\| \delta \overline F_{12} \bigr\|_{L^2_x} \bigl\| \overline{\bfD\phi} \bigr\|_{L^2_x} \\
                &\qquad + \underline s^\frac12 \bigl\| \overline F_{12} \bigr\|_{L^2_x} \bigl\| \delta \bfD \phi \bigr\|_{L^2_s L^2_x ([s', \underline s])}  \\
            &\lesssim \bigl\| \delta \underline{\bfD \phi} \bigr\|_{L^2_x} + \Bigl\| \int_{s'}^{\underline s} s^{-\frac12} \delta (s^\frac12 \bfD) (s \bfD^{(2)}) \phi \frac{ds}{s} \Bigr\|_{L^2_x} + \bigl\| \delta \overline F_{12} \bigr\|_{L^2_x} \\
                &\qquad + \epsilon \cdot \bigl\| \delta \bfD \phi \bigr\|_{L^\infty_s L^2_x ([s', s])},
    \end{align*}
    using Cauchy-Schwarz and parabolic smoothing \eqref{F12-sobolev}, \eqref{Dphi-sobolev} in the second line, and the assumption on the heat-temporal initial data $(\overline A, \overline{\vphantom{A}\phi}) \in \frE(\underline s)$ in the last line. As $\epsilon \ll 1$, we absorb the last term into the left-hand side, completing the proof.
\end{proof}

One should think of boundedness of the $s$-weighted derivatives $(s^\frac12 \bfD)^{(\leq n)} \sfU(s)$ as conferring frequency-localisation to $\sfU(s)$. Accordingly, we have the following lemma which captures ``almost-orthogonality'' of these frequency-localised objects. 

\begin{lemma}[Almost orthogonality lemma]\label{lem:orthogonal}
    Let $A$ be a connection $1$-form on $\R^2 \times [0, \underline s]$ solving the Yang-Mills heat flow \eqref{covariant-heat2} in caloric gauge \eqref{caloric-gauge} such that 
        \[
            \| \overline F_{12} \|_{L^2_x} \leq \underline s^{-1/2}.
        \]
    Then for sufficiently regular scalar fields $\sfU, \sfV : \R^2 \times [0, \underline s] \to \C$, we have  
        \[
            \int_0^{\underline s} \int_{s"}^{\underline s} \bigl\langle \sfU(s'), ({s"}^\frac12 \bfD) \sfV(s") \bigr\rangle_{L^2_x} \frac{ds'}{s'} \frac{ds"}{s"} 
                \lesssim \bigl\| (s^\frac12 \bfD)^{(\leq 1)} \sfU \bigr\|_{L^2_{\frac{ds}{s}} L^2_x} \bigl\| \sfV \bigr\|_{L^2_{\frac{ds}{s}} L^2_x}.
        \]
\end{lemma}

\begin{proof}
    Our strategy is to ``integrate-by-parts'' the derivative from the ``high-frequency'' factor at $s"$ to the ``low-frequency'' factor at $s'$, 
        \[
            \langle \sfU(s'), \bfD \sfV(s") \rangle_{L^2_x} 
                \approx - \langle \bfD \sfU(s'), \sfV(s")  \rangle_{L^2_x}.
        \]
    An astute reader will note that the formula above is not quite accurate, as on the left, the magnetic derivative is evaluated at $s"$, while on the right it is evaluated at $s'$. Nevertheless, the difference is merely a lower-order factor; indeed, integrating the Yang-Mills heat flow \eqref{covariant-heat2} in caloric gauge \eqref{caloric-gauge}, 
        \begin{align*}
            \bfD \sfV(s") 
                &= (\nabla - i A(s")) \sfV(s")\\
                &= (\nabla - i A(s')) \sfV(s") - i (A(s") - A(s')) \sfV(s")\\
                &= (\nabla - i A(s"))\sfV(s") + \Bigl( \int_{s"}^{s'} s^\frac12 (s^\frac12 \nabla) F_{12} \,  \frac{ds}{s} \Bigr)\cdot \sfV(s").
        \end{align*}
    Thus,
        \begin{align*}
            \bigl\langle \sfU(s') , ({s"}^\frac12 \bfD) \sfV(s") \bigr\rangle_{L^2_x} 
                &=  \Bigl( \frac{s"}{s'} \Bigr)^\frac12 \cdot \bigl\langle ({s'}^\frac12 \bfD) \sfV(s'), \sfV(s") \bigr\rangle_{L^2_x}\\
                    &\qquad +  \bigl\langle \sfU(s'), \Bigl( \int_{s"}^{s'} s^\frac12 (s^\frac12 \nabla) F_{12} \,  \frac{ds}{s} \Bigr)  \cdot \sfV(s") \bigr\rangle_{L^2_x} \\
                &\lesssim \Bigl(\frac{s"}{s'} \Bigr)^\frac12 \bigl\| ({s'}^\frac12\bfD) \sfU(s') \bigr\|_{L^2_x} \bigl\| \sfV (s") \bigr\|_{L^2_x} \\
                    &\qquad + |s"|^{\frac12} \log \Bigl( \frac{s'}{s"} \Bigr) \bigl\| \overline F_{12} \bigr\|_{L^2_x} \bigl\| \sfU(s') \bigr\|_{L^2_x} \bigl\| \sfV(s") \bigr\|_{L^2_x} =: \mathrm{I} + \mathrm{II},
        \end{align*}
    using the previous calculation and integrating-by-parts in the first line, and applying Cauchy-Schwarz and parabolic smoothing for the magnetic field \eqref{F12-sobolev} in the second line. Integrating in $s'$ and $s"$, 
        \begin{align*}
            \int_0^{\underline s} \int_{s"}^{\underline s} \mathrm{I} \, \frac{ds'}{s'} \frac{ds"}{s"} 
                &\lesssim \Bigl\| \int_{s"}^{\underline s} \Bigl(\frac{s"}{s'} \Bigr)^\frac12 \bigl\| ({s'}^\frac12\bfD) \sfU(s') \bigr\|_{L^2_x} \frac{ds'}{s'} \Bigr\|_{L^2_{\frac{ds"}{s"}}} \cdot \bigl\| \sfV(s") \bigr\|_{L^2_{\frac{ds"}{s"}} L^2_x} \\
                &\lesssim \bigl\| ({s'}^\frac12 \bfD) \sfU(s') \bigr\|_{L^2_{\frac{ds'}{s'}} L^2_x}  \bigl\| \sfV(s") \bigr\|_{L^2_{\frac{ds"}{s"}} L^2_x} ,
        \end{align*}
    by Cauchy-Schwarz in $s"$ and Schur's test, 
        \begin{align*}
             \int_0^{\underline s} \int_{s"}^{\underline s} \mathrm{II} \, \frac{ds'}{s'} \frac{ds"}{s"}
                &\lesssim \underline s^\frac12 \cdot \bigl\| \overline F_{12} \bigr\|_{L^2_x} \cdot \bigl\| \sfU(s') \bigr\|_{L^2_{\frac{ds'}{s'}} L^2_x} \bigl\| \sfV(s") \bigr\|_{L^2_{\frac{ds"}{s"} }L^2_x} \\
                &\lesssim \bigl\| \sfU(s') \bigr\|_{L^2_{\frac{ds'}{s'}} L^2_x} \bigl\| \sfV(s") \bigr\|_{L^2_{\frac{ds"}{s"}} L^2_x} ,
        \end{align*}
    by Cauchy-Schwarz in $s'$ and $s"$, basic calculus, and the assumption on $\overline F_{12}$. 
\end{proof}

\begin{proof}[Proof of covariant gradient difference \eqref{Dphi-diff-plancharel}]
    By \eqref{Dphi-s-diff-1},
        \[
             \bigl\|\delta \overline{\bfD \phi} \bigr\|_{L^2_x} 
                \lesssim \bigl\| \delta \underline{\bfD \phi} \bigr\|_{L^2_x} + \Bigl\| \int_{0}^{\underline s} \bfD^{(3)} \phi \frac{ds}{s} \Bigr\|_{L^2_x} + \bigl\| \delta \overline F_{12} \bigr\|_{L^2_x}.
        \]
    It remains to estimate the second term on the right. 
        \begin{align*}
             \Bigl\| \int_{0}^{\underline s}  \delta \bfD^{(3)} \phi \frac{ds}{s} \Bigr\|_{L^2_x}^2
                &= 2 \int_0^{\underline s} \int_{s"}^{\underline s} 
                    \Bigl\langle  \delta \bfD^{(3)} \phi (s'),  \delta \bfD^{(3)} \phi (s") \Bigr\rangle_{L^2_x} \, ds'\, ds"\\
                &= 2 \int_0^{\underline s} \int_{s"}^{\underline s}  \Bigl\langle  \delta \bfD^{(3)} \phi (s'),  \bfD \delta \bfD^{(2)} \phi (s") \Bigr\rangle_{L^2_x} \, ds' \, ds" \\
                    &\qquad + 2 \int_0^{\underline s} \int_{s"}^{\underline s}  \Bigl\langle  \delta \bfD^{(3)} \phi (s'),  \delta A \cdot \bfD^{(2)} \phi (s") \Bigr\rangle_{L^2_x} \, ds' \, ds" =: \mathrm{I} + \mathrm{II}. 
        \end{align*}

    We estimate the first term on the right using almost orthogonality, namely Lemma \ref{lem:orthogonal},
        \begin{align*}
            \mathrm{I} 
                &\lesssim \bigl\| s^{-\frac12} (s^\frac12 \bfD)^{(\leq 1)} \delta (s^\frac12 \bfD)^{(3)} \phi  \bigr\|_{L^2_{\frac{ds}{s}} L^2_x} \bigl\| s^{-\frac12} \delta (s \bfD^{(2)}) \phi \bigr\|_{L^2_{\frac{ds}{s}} L^2_x}\\
                &\lesssim  \bigl\| [\delta, \bfD] (s^\frac12 \bfD)^{(3)} \phi  \bigr\|_{L^2_{\frac{ds}{s}} L^2_x}^2 + \bigl\| s^{-\frac12} \delta (s^\frac12 \bfD)^{(\leq 2)} (s \bfD^{(2)}) \phi \bigr\|_{L^2_{\frac{ds}{s}} L^2_x}^2 \\
                &\lesssim_{\underline s} \bigl\| \delta \overline A \bigr\|_{L^4_x}^2 + \bigl\| s^{-\frac12} \delta (s^\frac12 \bfD)^{(\leq 2)} (s \bfD^{(2)}) \phi \bigr\|_{L^2_{\frac{ds}{s}} L^2_x}^2,
        \end{align*}
    using Cauchy-Schwarz, the commutator identity \eqref{commute-diff}, H\"older, and parabolic smoothing \eqref{Dphi-sobolev}. 

    The second term on the right can be handled directly by Cauchy-Schwarz in $s'$ and $s"$, 
        \begin{align*}
            \mathrm{II} 
                &\lesssim  \bigl\| s^{-\frac12} \delta (s^\frac12 \bfD)^{(3)} \phi \bigr\|_{L^2_{\frac{ds}{s}} L^2_x} \cdot \bigl\| \delta \overline A \bigr\|_{L^4_x} \cdot \bigl\| s^\frac14 \log(\underline s/s) \bigr\|_{L^2_{\frac{ds}{s}}} \cdot \bigl\| s^{-\frac14} (s^\frac12\bfD)^{(2)} \phi \bigr\|_{L^\infty_s L^4_x} \\
                &\lesssim \bigl\| s^{-\frac12} \delta (s^\frac12 \bfD) (s \bfD)^{(2)} \phi \bigr\|_{L^2_{\frac{ds}{s}} L^2_x} \cdot \bigl\| \delta \overline A \bigr\|_{L^4_x} \cdot \underline s^\frac14 \cdot \bigl\| \overline{\bfD \phi} \bigr\|_{L^2_x} \\
                &\lesssim_{\underline s} \bigl\| s^{-\frac12} \delta (s^\frac12 \bfD) (s \bfD)^{(2)} \phi \bigr\|_{L^2_{\frac{ds}{s}} L^2_x} \cdot \bigl\| \delta \overline A \bigr\|_{L^4_x} . 
        \end{align*}
    This completes the proof. 
\end{proof}

\subsection{Weak distance functional}\label{subsec:weak}

We will approach difference estimates for the Chern-Simons-Schr\"odinger flow in a quasi-linear fashion, that is, by propagating a weaker metric than the energy topology. To motivate the construction, one should view the top-order terms in the definition of the $\pzcE^1$-metric,
    \[
         \bigl\| s^{-\frac12} \delta (s^\frac12 \bfD)^{(\leq 2026)} (s \bfD^{(2)})\phi \bigr\|_{L^2_{\frac{ds}{s}} L^2_x} ,\qquad \bigl\| \delta \nabla \overline \Gauss \bigr\|_{L^2_x},
    \]
as encoding an $\dot H^1$-type distance. Analogously, our weak distance functional primarily consists of the terms, 
     \[
         \bigl\| \delta (s^\frac12 \bfD)^{(\leq 2027)} (s^\frac12 \bfD)\phi \bigr\|_{L^2_{\frac{ds}{s}} L^2_x} ,\qquad \bigl\| \delta \overline \Gauss \bigr\|_{L^2_x},
    \]
which encode an $L^2$-type distance. It will also be convenient to measure the differences of the magnetic potentials at $s = \underline s$ rather than $s = 0$, to create a better separation of low- and high-frequency scales. Putting these considerations together, we are led to introduce

\begin{definition}\label{def:weak-dist}
    Let $(A, \phi), (A', \phi') \in \pzcE^1 ([0, \underline s])$ be $\frE$-caloric extensions, then define the \textit{$\pzcE^0$-metric} by 
       \begin{align*}
            \dist_{\pzcE^0 ([0, \underline s])} \bigl( (A, \phi), (A', \phi') \bigr)
                &:=  \Big( \bigl\| \delta (s^\frac12 \bfD)^{(\leq 2027)} (s^\frac12 \bfD)\phi \bigr\|_{L^2_{\frac{ds}{s}} L^2_x}^2 + \bigl\| \delta \overline \Gauss \bigr\|_{L^2_x}^2 \\
                &\qquad\qquad + \bigl\| \delta \underline \bfD^{(\leq 100)} \underline{\bfD \phi} \bigr\|_{L^2_x}^2 +  \bigl\| \delta \underline A \bigr\|_{L^4_x}^2 + \bigl\| \nabla \delta \underline A \bigr\|_{H^{100}_x}^2 + \bigl\| \delta \underline \phi \bigr\|_{L^\infty_x}^2 \Big)^{\frac12}.
        \end{align*}
\end{definition}

We put $\delta A$ in a high Sobolev norm as it will be convenient later on to estimate it in $L^\infty$. The downside is that, for a generic $\frE$-caloric extension, this part of the metric need not be well-defined, as the curl-free part is merely transported by the Yang-Mills heat flow and thus remains much rougher than desired here. On the other hand, if we restrict to divergence-free configurations, which is precisely the class of initial data considered, this becomes a non-issue, as the divergence-free part is smoothed out by the linear heat flow.   

\begin{lemma}[Weakness of $\pzcE^0$-topology]\label{lem:weak}
    Let $(A, \phi), (A', \phi') \in \pzcE^1 ([0, \underline s])$ be $\frE$-caloric extensions to $\R^2 \times [0, \underline s]$, and suppose that $(\overline A, \overline{\vphantom{A} \phi}), (\overline A', \overline{\vphantom{A} \phi}')$ satisfy the Coulomb gauge \eqref{coulomb-data}. Then 
        \[
            \dist_{\pzcE^0 ([0, \underline s])} \bigl( (A, \phi), (A', \phi') \bigr)
                \lesssim_{\underline s} \dist_{\pzcE^1 ([0, \underline s])} \bigl( (A, \phi), (A', \phi') \bigr).
        \]
\end{lemma}

\begin{proof}
    Comparing the two metrics, the non-obvious terms to estimate are
        \[
             \bigl\| \delta (s^\frac12 \bfD)\phi\bigr\|_{L^2_{\frac{ds}{s}} L^2_x} , \qquad \bigl\|\nabla \delta \underline A \bigr\|_{H^{100}_x}, \qquad \bigl\|\delta\underline A\bigr\|_{L^4_x}.
        \]
    The first is controlled by the $\pzcE^1$-metric in view of \eqref{Dphi-s-diff-2} and $s^\frac12 \log^{1/2} (\underline s/s) \in L^2_{ds/s}$. For the second and third terms, note that the divergence-free condition is transported along the caloric gauge via \eqref{A-transport}, so $A = \PP^\df A$ solves the linear heat equation \eqref{A-heat}. Thus these can be estimate by $\overline A$ via parabolic smoothing. 
\end{proof}

The remainder of this subsection is dedicated to some Lipschitz and $\log^\frac12$-Lipschitz-type estimates with respect to the $\pzcE^0$-metric, which will be useful for difference estimates for the Chern-Simons-Schr\"odinger flow. The heuristic one should keep in mind is the Rellich-Kondrachov lemma for Sobolev spaces, that is, assuming additional regularity, one gains compactness in weaker topologies. The ensuing estimates are a manifestation of the idea in this non-linear setting. 

We begin with a Lipschitz-type bound for the Gauss tension field,

\begin{lemma}[Gauss tension field difference bound]
    Let $(A, \phi), (A', \phi') \in \pzcE^1 ([0, \underline s])$ be $\frE$-caloric extensions to $\R^2 \times [0, \underline s]$, then 
        \begin{equation}\label{eq:gauss-diff}
            \bigl\|  \delta \Gauss \bigr\|_{L^\infty_{s} L^2_x} + \bigl\| (s^\frac12 \nabla) \delta \Gauss \bigr\|_{L^2_{\frac{ds}{s}} L^2_x}
                \lesssim \bigl\| \delta \overline \Gauss \bigr\|_{L^2_x} + \bigl\| \overline{\bfD \phi} \bigr\|_{L^2_x} \bigl\| \delta (s^\frac12 \bfD) \phi \bigr\|_{L^2_{\frac{ds}{s}} L^2_x}.
        \end{equation}
\end{lemma}

\begin{proof}
    The right-hand side of the heat equation for the difference of $\Gauss$ and $\Gauss'$ is 
        \begin{align*}
            \delta(\text{R.H.S.\eqref{Gauss-heat-prelim}})
                 =  \bfD \phi \cdot \delta \bfD \phi.
        \end{align*}
    Estimating by Cauchy-Schwarz, 
        \begin{align*}
            \bigl\| s \delta(\text{R.H.S.\eqref{Gauss-heat-prelim}}) \bigr\|_{L^1_{\frac{ds}{s}} L^2_x}
                &\lesssim  \bigl\| s^\frac12 \bfD \phi \bigr\|_{L^2_{\frac{ds}{s}} L^\infty_x} \bigl\| \delta (s^\frac12 \bfD)\phi \bigr\|_{L^2_{\frac{ds}{s}} L^2_x} \\
                &\lesssim \bigl\| \overline{\bfD \phi} \bigr\|_{L^2_x} \bigl\| (s^\frac12 \bfD) \phi \bigr\|_{L^2_{\frac{ds}{s}} L^2_x},
        \end{align*}
    by parabolic smoothing \eqref{Dphi-strichartz-heat}. Inserting this into the abstract parabolic energy estimate \eqref{abs-heat-energy} applied to the differences of \eqref{Gauss-heat-prelim} gives the result. 
\end{proof}

\begin{remark}
    One can think of \eqref{gauss-diff} as analogous to the Lipschitz continuity of the high $\times$ high paraproduct, 
        \[
            \phi \mapsto \sum_N P_N (P_{> N} \phi \cdot P_{> N} \phi)
        \]
    with respect to the $L^2$-metric on bounded subsets of $H^1 (\R^2)$. 
\end{remark}

Recall that the $\pzcE^1$-metric had a term arising from the magnetic potential at $s = 0$, 
    \[
         \bigl\| \delta \overline A \bigr\|_{(\dot H^1 \cap L^4)_x}.
    \]
We claim that this term is genuinely ``lower-order'', in the sense that it is $\log^{\frac12}$-Lipschitz with respect to the weaker $\pzcE^0$-metric. To prove the claim, we argue by interpolating between ``high-frequency" decay alotted by the parabolic smoothing estimates from Section \ref{subsec:LP-sobolev}, and summability of low-frequencies. This will be useful throughout the article when proving difference estimates, 

\begin{lemma}[Low-high interpolation lemma]\label{lem:low-high-inter}
    Let $\mathtt f(s) : [0, \underline s] \to [0, \infty]$ be a non-negative function, and suppose that $\gamma > 0$ and $1 \leq p \leq q \leq \infty$, then for all $0 < s_* < \underline s$,
        \[
            \bigl\| \mathtt f(s) \bigr\|_{L^p_{\frac{ds}{s}}} 
                \lesssim_\gamma |s_*|^\gamma \bigl\| s^{-\gamma} \mathtt f(s) \bigr\|_{L^\infty_s ([0, s_*])} + \log^{\frac1p - \frac1q} \Bigl( \frac{\underline s}{s_*} \Bigr) \bigl\| \mathtt f(s) \bigr\|_{L^q_{\frac{ds}{s}} ([s_*, \underline s])}.
        \]
\end{lemma}

\begin{proof}
    By the triangle inequality and H\"older, 
        \begin{align*}
            \bigl\| \mathtt f(s) \bigr\|_{L^p_{\frac{ds}{s}}}
                &\leq \bigl\| \mathtt f(s) \bigr\|_{L^p_{\frac{ds}{s}} ([0, s_*])} + \bigl\| \mathtt f(s) \bigr\|_{L^p_{\frac{ds}{s}} ([s_*, \underline s])} \\
                &\leq \bigl\| s^\gamma \bigr\|_{L^p_{\frac{ds}{s}} ([0, s_*])}  \bigl\| s^{-\gamma} \mathtt f(s) \bigr\|_{L^\infty_s ([0, s_*])} + \bigl\| 1 \bigr\|_{L^r_{\frac{ds}{s}} ([s_*, \underline s])} \bigl\| \mathtt f(s) \bigr\|_{L^q_{\frac{ds}{s}}([s_*, \underline s])} \\
                &\lesssim_\gamma |s_*|^\gamma \bigl\| s^{-\gamma} \mathtt f(s) \bigr\|_{L^\infty_s ([0, s_*])} + \log^{\frac1r} \Bigl( \frac{\underline s}{s_*} \Bigr) \bigl\| \mathtt f(s) \bigr\|_{L^q_{\frac{ds}{s}}([s_*, \underline s])},
        \end{align*}
    where $\tfrac1r = \tfrac1p - \tfrac1q$.
\end{proof}

\begin{proposition}\label{prop:Ax-diff}
    Let $(A, \phi), (A', \phi') \in \pzcE^1 ([0, \underline s])$ be $\frE$-caloric extensions to $\R^2 \times [0, \underline s]$, and fix $0 < s_* < \tfrac{1}{1000} \underline s$. Then we have the following bounds: 
    \begin{itemize}
        \item \textup{(Scalar field difference bound).}
            \begin{equation}\label{eq:phi-diff}
                \bigl\| \delta \phi (s) \bigr\|_{(L^2 + L^\infty)_x} 
                    \leq \bigl\| \delta \underline \phi \bigr\|_{L^\infty_x} + \log^\frac12 \Bigl( \frac{\, \underline s\, }{s} \Bigr) \bigl\| \delta (s^\frac12 \bfD)^{(2)} \phi \bigr\|_{L^2_{\frac{ds}{s}} L^2_x}.
            \end{equation}

        \item \textup{(Magnetic field difference bound, I).} For each non-negative integer $n \in \N_0$, 
            \begin{equation}\label{eq:F12-diff}
            \begin{split}
                \bigl\| (s^\frac12 \nabla)^{(n)} \delta F_{12}\bigr\|_{L^\infty_s L^2_x} &+ \bigl\| (s^\frac12 \nabla)^{(n + 1)} \delta F_{12}  \bigr\|_{L^2_{\frac{ds}{s}} L^2_x} 
                    \\
                    &\lesssim_{n, \underline s} \log^\frac12 \Bigl( \frac{\underline s}{s_*} \Bigr) \cdot \Bigl( |s_*|^\frac12 + \dist_{\pzcE^0 ([0, \underline s])} \bigl( (A, \phi), (A', \phi') \bigr) \Bigr) .
            \end{split}
            \end{equation}

        \item \textup{(Magnetic field difference bound, II).}
            \begin{equation}\label{eq:F12-diff-lot}
                \bigl\| s^\frac18 \delta F_{12} \bigr\|_{L^\infty_s L^2_x}
                    \lesssim_{\underline s} \dist_{\pzcE^0 ([0, \underline s])} \bigl( (A, \phi), (A', \phi') \bigr).
            \end{equation}

        \item \textup{(Magnetic potential $L^4$-difference bound).}
            \begin{equation}\label{eq:Ax-diff}
                \bigl\| \delta A_x \bigr\|_{L^\infty_s L^4_x}
                    \lesssim_{\underline s} \dist_{\pzcE^0 ([0, \underline s])} \bigl( (A, \phi), (A', \phi') \bigr) .
            \end{equation}

        \item \textup{(Magnetic potential $L^\infty$-difference bound).}
            \begin{equation}\label{eq:Ax-diff-linfty}
                \bigl\| \delta A_x \bigr\|_{L^\infty_{s, x}}
                    \lesssim \log\Bigl( \frac{\underline s}{s_*} \Bigr) \cdot\Bigl( |s_*|^\frac12 + \dist_{\pzcE^0} \Bigr) .
            \end{equation}
    \end{itemize}
\end{proposition}

\begin{proof}[Proof of scalar field difference bound \eqref{phi-diff}]
    Taking the difference of \eqref{covariant-heat} in caloric gauge \eqref{caloric-gauge}, 
        \[
            \delta \phi(s) 
                = \delta \underline \phi + (s \partial_s)^{-1} \delta (s^\frac12 \bfD)^{(2)} \phi (s).
        \]
    Then by Cauchy-Schwarz, 
        \begin{align*}
            \| \delta \phi(s) \|_{(L^2 + L^\infty)_x}
                &\leq \bigl\| \delta \underline \phi \bigr\|_{L^\infty_x} + \int_s^{\underline s} \bigl\| \delta ({s'}^\frac12 \bfD)^{(2)} \phi (s') \bigr\|_{L^2_x} \frac{ds'}{s'} \\
                &\leq  \bigl\| \delta \underline \phi \bigr\|_{L^\infty_x}  + \log^\frac12 \Bigl( \frac{\underline s}{s} \Bigr) \bigl\| \delta (s^\frac12 \bfD)^{(2)} \phi \bigr\|_{L^2_{\frac{ds}{s}} L^2_x},
        \end{align*}
    as desired. 
\end{proof}

\begin{proof}[Proof of magnetic field difference bound \eqref{F12-diff}]
    By parabolic smoothing \eqref{standard-heat} for \eqref{F12-heat}, it suffices to consider $n = 0$. The energy identity for the linear heat equation \eqref{F12-heat} implies 
        \begin{align*}
            \tfrac12 \bigl\| \delta F_{12} \bigr\|_{L^\infty_s L^2_x}^2 
                \leq \tfrac12 \bigl\| \delta \underline F_{12} \bigr\|_{L^2_x}^2 + \bigl\| (s^\frac12 \nabla) \delta F_{12} \bigr\|_{L^2_{\frac{ds}{s}} L^2_x}^2,
        \end{align*}
    so it remains to estimate the $(s^\frac12 \nabla) \delta F_{12}$ in $L^2_{ds/s} L^2_x$. Recalling the definition of the Gauss tension field, 
        \begin{align*}
            (s^\frac12 \nabla) \delta F_{12}
                =   \delta (\phi \cdot (s^\frac12\bfD) \phi) + (s^\frac12 \nabla) \delta \Gauss .
        \end{align*}
    We estimate the contribution of the Gauss tension field via \eqref{gauss-diff}. For the remaining terms, we interpolate between high-frequencies and low-frequencies via Lemma \ref{lem:low-high-inter}. The high-frequencies decay; for $s \leq s_*$, we have
        \begin{align*}
            \bigl\| \delta (\phi \cdot (s^\frac12\bfD) \phi) \bigr\|_{L^2_x} 
                &\lesssim \bigl\| s^\frac12 \phi \bigr\|_{L^\infty} \bigl\| \bfD \phi \bigr\|_{L^2_x} \\
                &\lesssim \Bigl( \frac{\,s \,}{\underline s}\Bigr)^\frac12 \cdot \log^\frac12 \Bigl( \frac{\underline s}{\,s\,}\Bigr) \cdot \bigl\| \overline{\bfD \phi} \bigr\|_{L^2_x} ,
        \end{align*}
    using the Brezis-Gallou\"et-type bound \eqref{brezis-0} and parabolic smoothing \eqref{Dphi-sobolev}. For the low-frequencies $s \geq s_*$, we estimate
        \begin{align*}
             \bigl\| \delta (\phi \cdot (s^\frac12\bfD) \phi) \bigr\|_{L^2_{\frac{ds}{s}} L^2_x ([s_*, \underline s])} 
                &\leq \bigl\| \delta \phi \cdot (s^\frac12\bfD) \phi \bigr\|_{L^2_{\frac{ds}{s}} L^2_x ([s_*, \underline s])}  + \bigl\| \phi \cdot  \delta  (s^\frac12\bfD) \phi \bigr\|_{L^2_{\frac{ds}{s}} L^2_x ([s_*, \underline s])} \\
                &\lesssim  \bigl\| \delta \phi \bigr\|_{L^\infty_s (L^2 + L^\infty)_x ([s_*, \underline s])} \bigl\| (s^\frac12\bfD) \phi \bigr\|_{L^2_{\frac{ds}{s}} (L^\infty \cap L^2)_x ([s_*, \underline s])}\\  
                    &\qquad + \bigl\| \phi \bigr\|_{L^\infty_{s, x} ([s_*, \underline s])} \bigl\|  \delta  (s^\frac12\bfD) \phi \bigr\|_{L^2_{\frac{ds}{s}}L^2_x ([s_*, \underline s])} \\
                &\lesssim \bigl\| \overline{\bfD \phi} \bigr\|_{L^2_x}  \bigl\| \delta \underline \phi \bigr\|_{L^\infty_x} + \log^\frac12 \Bigl( \frac{\underline s}{s_*} \Bigr) \cdot \bigl\| \overline{\bfD \phi} \bigr\|_{L^2_x}  \cdot \bigl\|  \delta (s^\frac12 \bfD)^{(\leq 1)} (s^\frac12\bfD) \phi \bigr\|_{L^2_{\frac{ds}{s}}L^2_x ([s_*, \underline s])},
        \end{align*}
    estimating $\phi$ with the Brezis-Gallou\"et-type bound \eqref{brezis-0}, its difference via the analogous bound \eqref{phi-diff}, and using parabolic smoothing \eqref{Dphi-sobolev}, \eqref{Dphi-strichartz-heat} for $\bfD \phi$. Inserting this into Lemma \ref{lem:low-high-inter} completes the proof. 
\end{proof}

\begin{proof}[Proof of magnetic field difference bound, II \eqref{F12-diff-lot}]
    Following the proof of the previous estimate \eqref{F12-diff}, we have 
        \begin{align*}
            \| F_{12} (s) \|_{L^2_x} 
                \lesssim_{\underline s} \Bigl( 1 + \log^\frac12 \Bigl( \frac{\underline s}{s} \Bigr) \Bigr) \cdot \dist_{\pzcE^0 ([0, \underline s])} \bigl( (A, \phi), (A', \phi') \bigr).
        \end{align*}
    Multiplying by $s^\frac18$ damps the growth of the logarithmic factor near $s = 0$.
\end{proof}

\begin{proof}[Proof of magnetic potential $L^4$-difference bound \eqref{Ax-diff}]
    Integrating the transport equation \eqref{Ax-transport-df},
        \[
            \PP^\df \delta (A_x(s) - \underline A_x)
                = - (s \partial_s)^{-1} (s \, \delta \, \nabla F_{12}) (s). 
        \]
    Estimating the right-hand side, 
        \begin{align*}
            \bigl\| (s \partial_s)^{-1} (s \, \delta \, \nabla F_{12}) \bigr\|_{L^\infty_s L^4_x}
                &\lesssim \int_0^{\underline s} \bigl\| s\, \nabla \delta F_{12} \bigr\|_{L^4_x} \frac{ds}{s} \\
                &\lesssim \bigl\| s^\frac18 \delta F_{12} \bigr\|_{L^\infty_s L^2_x}\\
                &\lesssim \dist_{\pzcE^0 ([0, \underline s])} \bigl( (A, \phi), (A' , \phi') \bigr). 
        \end{align*}
    by parabolic smoothing \eqref{standard-heat} and \eqref{F12-diff-lot}.
\end{proof}

\begin{proof}[Proof of magnetic potential $L^\infty$-difference bound \eqref{Ax-diff-linfty}]
    Writing in caloric gauge, we arrive at the transport equation
        \begin{align*}
            \delta A_x (s)
                = \delta \underline A_x + (s \partial_s)^{-1} (s \nabla \delta F_{12}). 
        \end{align*}
    The first term on the right may be estimated using Sobolev embedding, 
        \begin{align*}
            \bigl\| \delta \underline A_x \bigr\|_{L^\infty_x} 
                &\lesssim \bigl\| \delta \underline A_x \bigr\|_{L^4_x} +  \bigl\| \nabla^{(2)}\delta \underline A_x \bigr\|_{L^2_x}.
        \end{align*}
    The second term, we use parabolic smoothing \eqref{standard-heat},
        \begin{align*}
            \bigl\|  (s \partial_s)^{-1} (s \nabla \delta F_{12}) \bigr\|_{L^\infty_x} 
                &\lesssim \bigl\| s^\frac12 (s^\frac12 \nabla) \delta F_{12} \bigr\|_{L^1_{\frac{ds}{s}} L^\infty_x}   \\
                &\lesssim \bigl\|  (s^\frac12 \nabla) \delta F_{12} \bigr\|_{L^1_{\frac{ds}{s}} L^2_x}.
        \end{align*}
    To estimate, we interpolate between low- and high-frequencies via Lemma \ref{lem:low-high-inter}. For the high-frequencies $s \leq s_*$, we use parabolic smoothing \eqref{Dphi-sobolev}, \eqref{Gauss-sobolev} and \eqref{brezis}, 
        \begin{align*}
            \bigl\| (s^\frac12 \nabla) \delta F_{12} \bigr\|_{L^2_x}
                &\lesssim \bigl\| \phi \bigr\|_{L^\infty_x} \bigl\| (s^\frac12 \bfD) \phi \bigr\|_{L^2_x} + \bigl\| (s^\frac12 \nabla) \Gauss \bigr\|_{L^2_x} \\
                &\lesssim_{\underline s} s^\frac12 \log^\frac12 \Bigl( \frac{\underline s}{s} \Bigr). 
        \end{align*}
    For the low-frequencies $s \geq s_*$, recall the magnetic field difference bound \eqref{F12-diff}, 
        \begin{align*}
           \bigl\| (s^\frac12 \nabla) \delta F_{12} \bigr\|_{L^2_{\frac{ds}{s}} L^2_x}
                &\lesssim_{\underline s} \log^\frac12 \Bigl( \frac{\underline s}{s_*} \Bigr) \cdot \Bigl( |s_*|^\frac12 +  \dist_{\pzcE^0} \Bigr).
        \end{align*}
    This completes the proof. 
\end{proof}

\begin{remark}
    One should think of the magnetic field difference bound \eqref{F12-diff} as analogous to the $\log^{\frac12}$-Lipschitz continuity of the map 
        \[
            \phi \mapsto \phi \cdot \phi,
        \]
    with respect to the $L^2$-metric on bounded subsets of $H^1 (\R^2)$. The main obstruction to a Lipschitz bound is the low-high part of the product, particularly the failure of Sobolev embedding $H^1 \not\hookrightarrow L^\infty$ when attempting to estimate the low-frequency factor. On the other hand, $B^{1, 2}_1 \hookrightarrow L^\infty$, so the Sobolev embedding only misses by half a logarithm after accounting for the difference between $\ell^1$-summability and $\ell^2$-summability of the Littlewood-Paley pieces. 
\end{remark}

As a straightforward corollary, we can show that if a sequence of $\frE$-caloric extensions converges with respect to the $\pzcE^0$-metric, then the ``low-frequency'' part of the $\pzcE^1$-metric also converges. 

\begin{corollary}\label{cor:converge-low}
     Let $\{ (\smoothing{n}{A}, \smoothing{n}{\vphantom{A}{\phi}}) \}_n \subseteq \pzcE^1 ([0, \underline s])$ be a sequence of $\frE$-caloric extensions such that 
        \[
            \lim_{n \to \infty} \dist_{\pzcE^0 ([0, \underline s])} \bigl( (A, \phi),  (\smoothing{n}{A}, \smoothing{n}{\vphantom{A}{\phi}}) \bigr)
                = 0,
        \]
    for some $\frE$-caloric extension $(A, \phi) \in \pzcE^1 ([0, \underline s])$. Then 
        \[
            \lim_{n \to \infty} \bigl\| \overline A - \smoothing{n}{\overline A} \bigr\|_{L^4 \, \cap \, \dot H^1}
                = 0. 
        \]
\end{corollary}

\begin{proof}
    By \eqref{F12-diff} and \eqref{Ax-diff}, 
        \[
             \bigl\| \overline A - \smoothing{n}{\overline A} \bigr\|_{L^4 \, \cap \, \dot H^1}
                \lesssim_{\underline s} \log^\frac12 \Bigl( \frac{\underline s}{s_*} \Bigr) \cdot \Bigl( |s_*|^\frac12 + \dist_{\pzcE^0 ([0, \underline s])} \bigl( (A, \phi),  (\smoothing{n}{A}, \smoothing{n}{\vphantom{A}{\phi}}) \bigr) \Bigr).
        \]
    Taking $n \to \infty$ and then $s_* \to 0$ gives the result. 
\end{proof}

\subsection{Frequency envelopes}\label{subsec:env}

Following Tao \cite{Tao2001a, Tao2009}, we introduce the notion of a frequency envelope as a robust tool for capturing the distribution of energy across frequency-scales. Before doing so, we need to identify a suitable analogue of a ``Littlewood-Paley projection'' for a configuration. Previous work on other geometric equations, e.g. \cite{Tao2004, BejenaruEtAl2011,Oh2014a,Oh2015a,Gavrus2022}, suggests that $(s F_{sx}, (s \bfD_s) \phi)$ is a reasonable candidate. However, the Gauss tension field is more regular than its defining variables in the definition of finite-energy configuration (Definition \ref{def:topology-1}), so, at least in $\frE$, it is somewhat more natural to view $(\Gauss, \phi)$ as the fundamental variables rather than the configuration $(A, \phi)$ itself. Accordingly, we interpret
    \[
        \bigl( (s \nabla^{(2)}) e^{s \Delta} \overline\Gauss,  (s \bfD^{(2)}) \phi (s)\bigr)
            \approx \text{projection of $(A, \phi)$ to frequency $s^{-1/2}$}.
    \]
and define the \textit{energy profile} of $(A, \phi)$
    \[
        \tte^1 [A, \phi] (s) 
            := \bigl\| s^{-\frac12} (s \nabla^{(2)}) e^{s \Delta} \overline\Gauss \bigr\|_{L^2_x}^2 + \bigl\| s^{-\frac12}  (s \bfD^{(2)}) \phi (s) \bigr\|_{L^2_x}^2.
    \]
Roughly speaking, $\tte^1 [A, \phi] (s)$ describes the ``(high)-frequency profile'' of the configuration $(A, \phi)$ at energy regularity. Indeed, by parabolic smoothing \eqref{parabolic-strichartz}, \eqref{Dphi-sobolev}, we know that
    \[
        \sup_{s \in [0, \underline s]} \tte^1 [A, \phi](s) + \int_0^{\underline s} \tte^1 [A, \phi](s) \frac{ds}{s} 
            \lesssim \bigl\| \nabla \overline \Gauss \bigr\|_{L^2_x}^2 + \bigl\| \overline{\bfD \phi} \bigr\|_{L^2_x}^2.
    \]

In practice, it will be inconvenient to directly work with the energy profile due to the non-linear nature of the setting, which causes non-trivial interactions between various frequency-scales. Similarly, when we turn to the Chern-Simons-Schr\"odinger equation, the flow can transfer energy between scales. Nevertheless, in both settings, we expect the interactions between differing frequency-scales to be mild, so a ``smoothed-out''-version of the energy profile should prove more useful. To this end, we introduce 

\begin{definition}
    We say that $\ttc(s)$ is an \textit{envelope} for a non-negative function $\mathtt f (s) :[0, \underline s] \to [0, \infty)$ if
    \begin{itemize}
        \item $\ttc(s)$ is \textit{slowly-varying}, i.e. for all $s, s' \in [0, \underline s]$, 
        \begin{equation}\label{eq:env-slow}
            \frac{\ttc(s)}{\ttc(s')}
                \leq \max\Bigl\{ \Bigl( \frac{s'}{s} \Bigr)^{\frac1{1000}}, \Bigl( \frac{s}{s'} \Bigr)^{\frac{1}{1000}} \Bigr\},
        \end{equation}

        \item $\mathtt f(s)$ \textit{lies underneath} $\ttc(s)$, i.e. for all $s \in [0, \underline s]$, 
        \begin{equation}\label{eq:under}
            \mathtt f(s) 
                \leq \ttc(s). 
        \end{equation}
    \end{itemize}
    Given a $\frE$-caloric extension $(A, \phi) \in \pzcE^1 ([0, \underline s])$, we say that $\ttc(s) \in L^2_{ds/s} ([0, \underline s])$ is a \textit{$\frE$-frequency envelope} for $(A, \phi)$ if it is an envelope for $\sqrt{\tte^1 [A, \phi](s)}$ and it satisfies the additional properties:
    \begin{itemize}
        \item low-frequency control, 
            \begin{equation}\label{eq:env-low}
            \ttc(\underline s) 
                \sim \sqrt{\cE_{\text{AH}} [\overline A, \overline{\vphantom{A}\phi}]}.
            \end{equation}    
        
        \item finite-energy, 
            \begin{equation}\label{eq:env-sharp}
                \int_0^{\underline s} |\ttc(s)|^2 \, \frac{ds}{s} 
                \lesssim \cE_{\text{AH}} [\overline A, \overline{\vphantom{A}\phi}].
            \end{equation}
    \end{itemize} 
\end{definition}

Observe that the energies $\tte^1 [A, \phi]$ and thus the notion of a $\frE$-frequency envelope are gauge-invariant\footnote{In fact, we could drop the assumption that $(A, \phi)$ satisfies the caloric gauge and simply require that it solve the covariant heat equation with heat-temporal data $(A, \phi)_{|s = 0} \in \frE(\underline s)$.}. The property \eqref{env-low} serves to control the ``low-frequency'' part of $(A, \phi)$, and is stated in a manner which allows us to easily exploit the conservation of the abelian Higgs energy energy. It is a standard exercise to show that $\frE$-frequency envelopes are uniformly bounded, with
        \[
          \ttc(s) 
            \lesssim \sqrt{\cE_{\text{AH}} [\overline A, \overline{\vphantom{A}\phi}]}.
        \]

A frequency envelope always exists, via a maximal operator-type construction, 

\begin{proposition}[Construction of a $\frE$-frequency envelope]\label{prop:env-stable}
    Let $(A, \phi) \in \pzcE^1 ([0, \underline s])$ be a $\frE$-caloric extension, then $\ttc(s)$ defined by 
        \[
        \ttc(s) 
            := \sup_{s" \in [0, \underline s]} \min \Bigl\{ \Bigl( \frac{s"}{s} \Bigr)^{\frac{1}{1000}}, \Bigl( \frac{s}{s"} \Bigr)^{\frac{1}{1000}} \Bigr\} \cdot \sqrt{\tte^1 [A, \phi] (s")} +\Bigl( \frac{s}{\,\underline s \, } \Bigr)^{\frac{1}{1000}} \cdot \sqrt{\cE_{\textup{AH}} [\overline A, \overline{\vphantom{A}\phi}]} ,
        \]
    is a $\frE$-frequency envelope. Furthermore, if $(A', \phi') \in \pzcE^1 ([0, \underline s])$ is also a $\frE$-caloric extension and $\ttc'(s)$ its $\frE$-frequency envelope, constructed as above, then we have the stability bound
        \begin{equation}\label{eq:env-lip}
        \begin{split}
            \| \delta \ttc(s) \|_{L^2_{\frac{ds}{s}}}
                &\lesssim \bigl\|\delta \overline F_{12} \bigr\|_{L^2_x} + \bigl\| \delta \tfrac12(1 - |\overline \phi|^2) \bigr\|_{L^2_x} + \bigl\| \delta \overline{\bfD \phi} \bigr\|_{L^2_x} + \bigl\| \nabla \overline \Gauss \bigr\|_{L^2_x} \\
                    &\qquad +  \bigl\| \delta \underline{\bfD \phi} \bigr\|_{L^2_x} + \bigl\| s^{-\frac12} (s^\frac12 \bfD)^{(\leq 2)} (s \bfD^{(2)}) \phi \bigr\|_{L^2_{\frac{ds}{s}} L^2_x}.
        \end{split}
        \end{equation}
\end{proposition}

As an immediate corollary, the stability bound \eqref{env-lip} along with Proposition \ref{thm:topology} imply that we can choose a $\frE$-frequency envelope for any $\frE$-caloric extension, 
    \begin{align*}
        \pzcE^1 ([0, \underline s])
            &\longrightarrow L^2_{\frac{ds}{s}} ([0, \underline s]),\\
        (A, \phi) 
            &\longmapsto \ttc(s),
    \end{align*}
in a locally Lipschitz-continuous fashion. In particular, given a convergent sequence of $\frE$-caloric extensions, we can choose corresponding $\frE$-frequency envelopes such that their energies are ``tight''. This will be useful in Section \ref{sec:cts} when closing the continuous dependence argument. 

It is easy to check that $\ttc(s)$ is slowly-varying \eqref{env-slow}, lies above the energy profile \eqref{under}, and controls the low-frequencies \eqref{env-low}. To complete the proof of Proposition \ref{prop:env-stable}, it remains then to estimate the size of the envelope \eqref{env-sharp}, and show the stability bound \eqref{env-lip}. The following lemma, adapted from \cite[Lemma 7.3]{Tao2009a}, will help facilitate the proof, 

\begin{lemma}[Abstract $L^2_{ds/s}$-envelope bound]\label{lem:abs-env-dyadic}
    Let $\mathtt f:  [0, \underline s] \to [0, \infty)$ be a sufficiently regular function, and set 
        \[
        \ttc(s) 
            := \sup_{s" \in (0, \underline s]} \min\Bigl\{ \Bigl( \frac{s"}{s} \Bigr)^{\frac{1}{1000}}, \Bigl( \frac{s}{s"} \Bigr)^{\frac{1}{1000}} \Bigr\} \cdot \sqrt{\mathtt f(s")}.
        \]
    Then
        \[
            \int_0^{\underline s} |\ttc(s)|^2 \frac{ds}{s} 
                \lesssim  \int_0^{\underline s} |\mathtt f(s)| \frac{ds}{s} + \int_0^{\underline s}  |(s \partial_s) \mathtt f(s)| \frac{ds}{s}.
        \]
\end{lemma}

\begin{proof}
    It is obvious from the construction that $\ttc(s)$ is an envelope for $\sqrt{\mathtt f(s)}$. Dividing $[0, \underline s]$ into dyadic blocks $J_N := [\frac{1}N \underline s, \frac{1}{2N} \underline s]$ for $N \in 2^{\N_0}$ and leveraging the slowly-varying property, 
        \begin{align*}
            \int_0^{\underline s} |\ttc(s)|^2 \frac{ds}{s}
                &\sim \sum_{N \in 2^{\N_0}} \int_{J_N} |\ttc(s)|^2 \frac{ds}{s} \\
                &\sim \sum_{N \in 2^{\N_0}} \bigl|\ttc\bigl(\tfrac{1}{N} \underline s\bigr)\bigr|^2.
        \end{align*}
    By construction of $\ttc(s)$, we can bound 
        \begin{align*}
            \ttc \bigl( \tfrac{1}{N} \underline s \bigr) 
                &\lesssim \sum_{M \in 2^{\N_0}} \sup_{s" \in J_M} \min \Bigl\{ \Bigl( \frac{s"}{\frac1N \underline s} \Bigr)^{\frac{1}{1000}}, \Bigl( \frac{\frac1N \underline s}{s"} \Bigr)^{\frac{1}{1000}} \Bigr\} \cdot \sqrt{\mathtt f(s")}\\
                &\lesssim \sum_{M \in 2^{\N_0}} \min \Bigl\{ \Bigl( \frac{N}{M} \Bigr)^{\frac{1}{1000}}, \Bigl( \frac{M}{N} \Bigr)^{\frac{1}{1000}} \Bigr\} \cdot \sup_{s" \in J_M} \sqrt{\mathtt f(s")}.
        \end{align*}
    The integral kernel $(N, M) \mapsto \min\{ (\frac{N}{M})^{\frac{1}{1000}}, (\frac{M}{N})^{\frac1{1000}}\}$ decays exponentially off of the diagonal $N \neq M$ and has unit size on the diagonal $N = M$, so by Schur's test, 
        \begin{align*}
            \sum_{N \in 2^{\N_0}} \bigl| \ttc \bigl( \tfrac{1}{N} \underline s \bigr)  \bigr|^2 
                \lesssim \sum_{N \in 2^{\N_0}} \sup_{s \in J_N} |\mathtt f(s)|.  
        \end{align*}
    On each of these dyadic blocks, we can invoke the Poincar\'e-Sobolev inequality, 
        \begin{align*}
            \sup_{s \in J_N} |\mathtt f(s)|
                \lesssim \int_{J_N} |\mathtt f(s)| \frac{ds}{s} + \int_{J_N} |(s \partial_s) \mathtt f(s)| \frac{ds}{s},
        \end{align*}
    so, summing in $N$, we can conclude the result. 
\end{proof}

\begin{proof}[Proof of \eqref{env-sharp}]
    Applying Lemma \ref{lem:abs-env-dyadic} to the definition of $\ttc(s)$,  
        \begin{align*}
            \int_0^{\underline s} |\ttc(s)|^2 \frac{ds}{s}
                &\lesssim \int_0^{\underline s} \tte^1 [A,\phi] (s) \frac{ds}{s} +\int_0^{\underline s} |(s \partial_s) \tte^1 [A, \phi] (s)| \frac{ds}{s} +  \int_0^{\underline s} \Bigl( \frac{s}{\underline s} \Bigr)^{\frac{1}{500}}  \cE_{\mathrm{AH}} [\overline A, \overline{\vphantom{A}\phi}]\frac{ds}{s} .
            \end{align*}
    The first term on the right is bounded by the abelian Higgs energy thanks to parabolic smoothing \eqref{parabolic-strichartz}, \eqref{Dphi-sobolev}, and the third term on the right is easily integrated. It remains to estimate the second term on the right; by the product rule and Cauchy-Schwarz,  
        \begin{align*}
            (s \partial_s) \tte^1 [A, \phi] (s)
                &= (s \partial_s)  \int_{\R^2} \bigl|s^{-\frac12} \sfU (s)\bigr|^2 \, dx\\
                &\lesssim  \int_{\R^2} \bigl|s^{-\frac12}  \sfU (s)\bigr|^2 \, dx + \Bigl| \Bigl\langle s^{-\frac12}  \sfU (s), s^{-\frac12} (s \sfD_s)  \sfU (s)\Bigr\rangle_{L^2_x} \Bigr| \\
                &\lesssim  \int_{\R^2} \bigl|s^{-\frac12}  \sfU (s)\bigr|^2 \, dx + \bigl\| s^{-\frac12} \sfU \bigr\|_{L^2_x} \bigl\| s^{-\frac12} (s \sfD_s) \sfU (s) \bigr\|_{L^2_x}
        \end{align*}
    writing $\sfD \in \{ \nabla, \bfD \}$ and $\sfU \in \{ (s \nabla^{(2)}) e^{s \Delta} \overline \Gauss, (s \bfD^{(2)})\phi \}$. Integrating the against $ds/s$,
        \begin{align*}
            \int_0^{\underline s} |(s \partial_s) \tte^1 [A, \phi]|  \frac{ds}{s}
                &\lesssim \cE_{\mathrm{AH}} [\overline A, \overline{\vphantom{A}\phi}] + \cE_{\mathrm{AH}} [\overline A, \overline{\vphantom{A}\phi}]^\frac12 \bigl\| s^{-\frac12} (s \sfD_s) \sfU \bigr\|_{L^2_{\frac{ds}{s}} L^2_x}.
        \end{align*}
    We estimate the case $\sfU = (s \bfD^{(2)}) \phi$, as the case $\sfU = (s \nabla^{(2)}) e^{s \Delta} \overline \Gauss$ is strictly easier. Multiplying the equation for the covariant Hessian \eqref{Dphi-heat-2} by $s^{-\frac12 + 1 + 1}$ and commuting powers appropriately, we obtain
        \begin{align*}
            s^{-\frac12} (s \bfD_s)  (s \bfD^{(2)}) \phi 
                &= s^{-\frac12}  (s^\frac12 \bfD)^{(2)} (s\bfD^{(2)}) \phi + s^{-\frac12} (s\bfD^{(2)}) \phi  + s\cdot F_{12} \cdot s^{-\frac12}  (s \bfD^{(2)}) \phi  \\
                    &\qquad + s\cdot s^{\frac12} \nabla F_{12} \cdot \bfD \phi . 
            \end{align*}
    Estimating in $L^2_{ds/s} L^2_x$, 
        \begin{align*}
            \bigl\| s^{-\frac12} (s \bfD_s)  (s \bfD^{(2)}) \phi \bigr\|_{L^2_{\frac{ds}{s}}L^2_x} 
                &\lesssim \Bigl( 1 + \underline s^{\frac12} \bigl\| s^\frac12 F_{12} \bigr\|_{L^\infty_{s,x}} \Bigr) \bigl\| s^{-\frac12} (s^\frac12 \bfD)^{(\leq 2)} (s \bfD^{(2)}) \phi \bigr\|_{L^2_{\frac{ds}{s}} L^2_x}  \\
                    &\qquad + \underline s^{\frac12} \bigl\|  (s^\frac12 \nabla) F_{12} \bigr\|_{L^\infty_{s} L^2_x}\bigl\| s^\frac12 \bfD \phi \bigr\|_{L^2_{\frac{ds}{s}} L^\infty_x} \\
                &\lesssim \Bigl( 1 + \underline s^{\frac12} \bigl\| \overline F_{12} \bigr\|_{L^2_x}\Bigr) \bigl\| \overline{\bfD \phi} \bigr\|_{L^2_x} \\
                &\lesssim \cE_{\mathrm{AH}} [\overline A, \overline{\vphantom{A}\phi}]^\frac12,
        \end{align*}
    by H\"older's inequality in the first line, parabolic smoothing \eqref{F12-sobolev}, \eqref{Dphi-sobolev}, \eqref{Dphi-strichartz-heat} in the second line, and the definition of $(\overline A, \overline{\vphantom{A}\phi}) \in \frE(\underline s)$ in the last line. Collecting the previous inequalities, we conclude \eqref{env-sharp}. 
\end{proof}

The proof of the stability bound \eqref{env-lip} follows largely the same circle of ideas. The main difference is that when converting $s$-derivatives to $x$-derivatives using the covariant heat equation, we will see differences of the covariant gradient, i.e. $\delta \bfD \phi$. 

\begin{proof}[Proof of \eqref{env-lip}]
    By the definition of the envelopes $\ttc(s)$ and $\ttc(s')$ and the triangle inequality, 
        \begin{align*}
            | \delta \ttc(s)| 
                &\leq  \sup_{s" \in [0, \underline s]} \min \Bigl\{ \Bigl( \frac{s"}{s} \Bigr)^{\frac{1}{1000}}, \Bigl( \frac{s}{s"} \Bigr)^{\frac{1}{1000}} \Bigr\} \sqrt{ \delta \tte^1 [A, \phi] (s")} \\
                    &\qquad + \Bigl( \frac{\,s\,}{\underline s} \Bigr)^{\frac{1}{1000}} \Bigl| \delta \sqrt{\cE_{\mathrm{AH}} [\overline A, \overline{\vphantom{A}\phi}]} \Bigr| \\
                &\lesssim \sup_{s" \in [0, \underline s]} \min \Bigl\{ \Bigl( \frac{s"}{s} \Bigr)^{\frac{1}{1000}}, \Bigl( \frac{s}{s"} \Bigr)^{\frac{1}{1000}} \Bigr\} \sqrt{ \delta \tte^1 [A, \phi] (s")}\\
                    &\qquad + \Bigl( \frac{\,s\,}{\underline s} \Bigr)^{\frac{1}{1000}} \Bigl(\bigl\|\delta \overline F_{12} \bigr\|_{L^2_x} + \bigl\| \delta \tfrac12(1 - |\overline \phi|^2) \bigr\|_{L^2_x} + \bigl\| \delta \overline{\bfD \phi} \bigr\|_{L^2_x}\Bigr) .
        \end{align*}
    Employing Lemma \ref{lem:abs-env-dyadic} as in the proof of \eqref{env-sharp} and parabolic smoothing, 
        \begin{align*}
            \int_0^{\underline s} |\delta \ttc(s)|^2 \frac{ds}{s}
                &\lesssim   \bigl\| s^{-\frac12} \delta (s \nabla^{(2)}) e^{s \Delta} \overline \Gauss \bigr\|_{L^2_{\frac{ds}{s}} L^2_x} + \bigl\| s^{-\frac12} \delta (s \bfD^{(2)}) \phi \bigr\|_{L^2_{\frac{ds}{s}} L^2_x} + \int_0^{\underline s} |(s \partial_s) \delta \tte^1 [A, \phi] (s)| \frac{ds}{s} \\
                    &\qquad + \bigl\|\delta \overline F_{12} \bigr\|_{L^2_x} + \bigl\| \delta \tfrac12(1 - |\overline \phi|^2) \bigr\|_{L^2_x} + \bigl\| \delta \overline{\bfD \phi} \bigr\|_{L^2_x} \\
                &\lesssim \bigl\|\nabla \overline \Gauss \bigr\|_{L^2_x} + \bigl\| s^{-\frac12} \delta (s \bfD^{(2)}) \phi \bigr\|_{L^2_{\frac{ds}{s}} L^2_x} + \int_0^{\underline s} |(s \partial_s) \delta \tte^1 [A, \phi] (s)| \frac{ds}{s} \\
                    &\qquad + \bigl\|\delta \overline F_{12} \bigr\|_{L^2_x} + \bigl\| \delta \tfrac12(1 - |\overline \phi|^2) \bigr\|_{L^2_x} + \bigl\| \delta \overline{\bfD \phi} \bigr\|_{L^2_x}
        \end{align*}
    It remains to estimate the contribution of the $s$-derivatives. Proceeding analogously to the previous proof, we consider the differences of the equations for the covariant Hessians \eqref{Dphi-heat-2}, 
        \begin{align*}
            s^{-\frac12} \delta (s \bfD_s)  (s \bfD^{(2)}) \phi 
                &= s^{-\frac12}  \delta (s^\frac12 \bfD)^{(2)} (s\bfD^{(2)}) \phi + s^{-\frac12} \delta (s\bfD^{(2)}) \phi  + s\cdot F_{12} \cdot s^{-\frac12} \delta (s \bfD^{(2)}) \phi  \\
                    &\qquad + s \cdot \delta F_{12} \cdot s^{-\frac12} (s \bfD^{(2)}) \phi + s\cdot s^{\frac12} \nabla \delta F_{12} \cdot \bfD \phi \\
                    &\qquad + s \cdot s^{\frac12} \nabla F_{12} \cdot \delta \bfD \phi \\
                &=: \mathrm{I} + \mathrm{II} + \mathrm{III} 
            \end{align*}
    Estimating in $L^2_{ds/s} L^2_x$, the terms in the first line on the right yield
        \begin{align*}
            \bigl\| \mathrm{I} \bigr\|_{L^2_{\frac{ds}{s}}L^2_x} 
                &\lesssim \Bigl( 1 + \underline s^{\frac12} \bigl\| s^\frac12 F_{12} \bigr\|_{L^\infty_{s,x}} \Bigr) \bigl\| s^{-\frac12} \delta(s^\frac12 \bfD)^{(\leq 2)} (s \bfD^{(2)}) \phi \bigr\|_{L^2_{\frac{ds}{s}} L^2_x} \\
                &\lesssim  \Bigl( 1 + \underline s^{\frac12} \bigl\| \overline F_{12} \bigr\|_{L^2_x} \Bigr)\bigl\| s^{-\frac12} \delta(s^\frac12 \bfD)^{(\leq 2)} (s \bfD^{(2)}) \phi \bigr\|_{L^2_{\frac{ds}{s}} L^2_x}  \\
                &\lesssim \bigl\| s^{-\frac12} \delta(s^\frac12 \bfD)^{(\leq 2)} (s \bfD^{(2)}) \phi \bigr\|_{L^2_{\frac{ds}{s}} L^2_x} ,
        \end{align*}
    using parabolic smoothing \eqref{F12-sobolev} and the size of $\underline s$. For the terms in the second line on the right,  
        \begin{align*}
            \bigl\| \mathrm{II} \bigr\|_{L^2_{\frac{ds}{s}} L^2_x}
                &\lesssim \underline s^{\frac12} \bigl\|  (s^\frac12 \nabla)^{(\leq 1)} \delta F_{12} \bigr\|_{L^\infty_{s} L^2_x}\bigl\| (s^\frac12 \bfD)^{(\leq 1)} (s^\frac12 \bfD) \phi \bigr\|_{L^2_{\frac{ds}{s}} L^\infty_x} \\
                &\lesssim \underline s^{\frac12} \bigl\| \delta \overline F_{12} \bigr\|_{L^2_x} \bigl\| \overline{\bfD \phi} \bigr\|_{L^2_x} \\
                &\lesssim \bigl\| \overline F_{12} \bigr\|_{L^2_x}, 
        \end{align*}
    using parabolic smoothing \eqref{F12-sobolev}, \eqref{Dphi-sobolev} and the size of $\underline s$. For the term in the third line on the right, 
        \begin{align*}
            \bigl\| \mathrm{III} \bigr\|_{L^2_{\frac{ds}{s}} L^2_x} 
                &\lesssim \bigl\| s^\frac12 (s^\frac12 \nabla) F_{12} \bigr\|_{L^\infty_{s, x}} \bigl\| s^\frac12 \delta \bfD \phi \bigr\|_{L^2_{\frac{ds}{s}} L^2_x} \\
                &\lesssim \bigl\| \overline F_{12} \bigr\|_{L^2_x} \Bigl( \underline s^\frac12 \cdot \bigl(\bigl\| \delta \underline{\bfD \phi} \bigr\|_{L^2_x} + \bigl\| \delta \overline F_{12} \bigr\|_{L^2_x}\bigr) + \bigl\| s^\frac12 \log\bigl( \tfrac{\underline s}{s} \bigr) \bigr\|_{L^2_{\frac{ds}{s}} L^2_x} \bigl\| s^{-\frac12} \delta (s \bfD) (s \bfD^{(2)}) \phi \bigr\|_{L^2_{\frac{ds}{s}} L^2_x} \Bigr) \\
                &\lesssim \bigl\| \delta \underline{\bfD \phi} \bigr\|_{L^2_x} + \bigl\| \delta \overline F_{12} \bigr\|_{L^2_x} + \bigl\| s^{-\frac12} \delta (s \bfD) (s \bfD^{(2)}) \phi \bigr\|_{L^2_{\frac{ds}{s}} L^2_x},
        \end{align*}
    using parabolic smoothing \eqref{F12-sobolev} and the difference bound on the covariant gradient \eqref{Dphi-s-diff-2}. 
\end{proof}

This completes the proof of Proposition \ref{prop:env-stable}. We now turn to proving dyadic versions of the smoothing estimates from the previous subsection. Here the slowly-varying property is very convenient, as it implies $\ttc(s) \sim \ttc(s')$ whenever $s \sim s'$. Thus, we can abuse notation quite a bit in using the same notation for various variables evaluated at the same scale of heat-times, e.g. $F_{12}$ for both $F_{12}(s)$ and $F_{12} (s/2)$.

\begin{proposition}\label{prop:F12-env}
    Let $(A, \phi) \in \pzcE^1 ([0, \underline s])$ be a $\frE$-caloric extension, and suppose that $\ttc(s)$ is a $\frE$-frequency envelope for $(A, \phi)$. Then, for each non-negative integer $n \in \N_0$, the following hold:
        \begin{itemize}
            \item \textup{(Gauss tension field envelope bound).} 
                \begin{equation}\label{eq:gauss-env}
                    \bigl\| s^{-\frac{1}{2}} (\partial_s - \Delta)^{-1} (\bfD \phi \cdot \bfD \phi) \bigr\|_{L^2_x}
                        \lesssim \bigl\| \overline{\bfD \phi} \bigr\|_{L^2_x} \cdot \ttc (s).
                \end{equation}

            \item \textup{(Magnetic field envelope bound).} 
                \begin{equation}\label{eq:F12-env}
                    \bigl\| s^{-\frac12} (s^\frac12 \nabla)^{(\leq n)} (s \nabla^{(2)}) F_{12} \bigr\|_{L^2_x} 
                        \lesssim_n \Bigl( 1 + \bigl\| \phi \bigr\|_{L^\infty_x} + \bigl\| \overline{\bfD \phi} \bigr\|_{L^2_x} \Bigr) \cdot \ttc (s) .
                \end{equation}
        \end{itemize}
    
\end{proposition}

\begin{proof}[Proof of Gauss tension field bound \eqref{gauss-env}]
    By Duhamel's formula, 
        \begin{align*}
            \bigl\| s^{-\frac12} (\partial_s - \Delta)^{-1} (\bfD \phi \cdot \bfD \phi) \bigr\|_{L^2_x}
                &\lesssim \int_0^s  \Bigl( \frac{s'}{s} \Bigr)^\frac12 \bigl\| {s'}^\frac12 \bfD \phi \cdot \bfD \phi \bigr\|_{L^2_x} \, \frac{ds'}{s'} \\
                &\lesssim \int_0^s \Bigl( \frac{s'}{s} \Bigr)^\frac12 \bigl\| \bfD \phi \bigr\|_{L^2_x} \bigl\| {s'}^{-\frac12} (s' \bfD^{(2)}) \phi \bigr\|_{L^2_x} \frac{ds'}{s'} \\
                &\lesssim \int_0^s \Bigl( \frac{s'}{s} \Bigr)^{\frac12 - \frac{1}{1000}} \bigl\| \overline{\bfD \phi} \bigr\|_{L^2_x} \cdot \ttc(s) \frac{ds'}{s'} \lesssim \bigl\| \overline{\bfD \phi} \bigr\|_{L^2_x} \cdot \ttc(s),
        \end{align*}
    using the abstract bilinear estimate from Lemma \ref{lem:bilinear} in the second line, parabolic smoothing \eqref{Dphi-sobolev} and the slowly-varying property in the third line.
\end{proof}

\begin{proof}[Proof of magnetic field bound \eqref{F12-env}]
    By parabolic smoothing \eqref{standard-heat} for the linear heat equation \eqref{F12-heat}, and the slowly-varying property $\ttc(s) \sim \ttc(s/2)$, it suffices to show the result for $n = 0$. Using the definition of the Gauss tension field, Duhamel's formula for \eqref{Gauss-heat-prelim}, and the product rule \eqref{product}.
        \begin{align*}
            \nabla^{(2)} F_{12} (s) 
                &= \nabla^{(2)} e^{\frac12 s \Delta} \Bigl( \tfrac12(1 - |\phi|^2) - \Gauss \Bigr)(\tfrac12 s) \\
                &= \nabla^{(2)} e^{\frac12 s \Delta} \Bigl( \tfrac12 (1 - |\phi|^2) - e^{s \Delta} \overline \Gauss + (\partial_s - \Delta)^{-1} (\bfD \phi \cdot \bfD \phi) \Bigr)(\tfrac12 s) \\
                &= e^{\frac12 s \Delta}\bigl(\bfD \phi \cdot \bfD \phi + \phi \cdot \bfD^{(2)} \phi\bigr) + \nabla^{(2)} e^{s \Delta} \overline \Gauss + e^{\frac12 s \Delta} \nabla^{(2)} (\partial_s - \Delta)^{-1} (\bfD \phi \cdot \bfD \phi)(\tfrac12 s).
        \end{align*}
    We estimate each term on the right individually; for the first term, 
        \begin{align*}
            \bigl\| s^{\frac12} \cdot e^{\frac12 s \Delta}( \bfD \phi \cdot \bfD \phi) \bigr\|_{L^2_x} 
                &\lesssim \bigl\| \bfD \phi \bigr\|_{L^2_x} \bigl\| s^{-\frac12} (s \bfD)^{(2)} \phi \bigr\|_{L^2_x} \\
                &\lesssim \bigl\| \overline{\bfD \phi} \bigr\|_{L^2_x} \cdot \ttc(s), 
        \end{align*}
    using boundedness of the linear heat propagator on $L^2$ and the bilinear estimate from Lemma \ref{lem:bilinear} in the first line and parabolic smoothing \eqref{Dphi-sobolev} in the second line. The second and third terms are straightforward applications of the definition of the envelope, 
        \begin{align*}
            \bigl\| s^{\frac12} \cdot e^{\frac12 s \Delta} (\phi \cdot \bfD^{(2)}\phi)\bigr\|_{L^2_x}
                &\lesssim \| \phi \|_{L^\infty} \cdot \ttc(s), \\
            \bigl\| s^{-\frac12} (s \nabla^{(2)}) e^{s \Delta} \overline \Gauss \bigr\|_{L^2_x}
                &\lesssim \ttc(s),
        \end{align*}
    an the last term, we use parabolic smoothing \eqref{standard-heat} and the Gauss tension field bound \eqref{gauss-env}, 
        \begin{align*}
            \bigl\| s^{-\frac12} (s^\frac12 \nabla)^{(2)} e^{\frac12 s \Delta} (\partial_s - \Delta)^{-1} (\bfD \phi \cdot \bfD \phi) \bigr\|_{L^2_x}
                &\lesssim \bigl\| s^{-\frac12} (\partial_s - \Delta)^{-1} (\bfD \phi \cdot \bfD \phi) \bigr\|_{L^2_x} \\
                &\lesssim \bigl\| \overline{\bfD \phi} \bigr\|_{L^2_x} \cdot \ttc(s).
        \end{align*}
    Collecting the previous calculations, we conclude the result.  
\end{proof}

The following bootstrap-type argument will be useful in the sequel, 

\begin{lemma}[Abstract bootstrap argument]\label{lem:abs-env-boot}
    Let $\mathtt f(s) : [0, \underline s] \to [0, \infty)$ be a non-negative function, and suppose that $\mathtt g(s): [0, \underline s] \to [0, \infty)$ is a slowly-varying function such that 
        \begin{equation}\label{eq:abs-BS-env}
            \mathtt f(s) 
                \lesssim \mathtt g(s) + \epsilon\cdot \ttc(s),
        \end{equation}
    for $\epsilon \ll 1$ sufficiently small and any envelope $\ttc(s)$ for $\mathtt f(s)$. Then 
        \begin{equation}
            \mathtt f (s) 
                \lesssim \mathtt g(s).
        \end{equation}
\end{lemma}

\begin{proof}
    Set
        \[
        \ttc(s) 
            := \sup_{s" \in (0, \underline s]} \min\Bigl\{ \Bigl( \frac{s"}{s} \Bigr)^{\frac{1}{1000}}, \Bigl( \frac{s}{s"} \Bigr)^{\frac{1}{1000}} \Bigr\} \cdot \mathtt f(s").
        \]
    It is easy to verify that this is an envelope for $\mathtt f(s)$. Then inserting \eqref{abs-BS-env} into the definition and using the slowly-varying property \eqref{env-slow},
        \begin{align*}
            \ttc(s)  
                &\lesssim \sup_{s" \in (0, \underline s]} \min\Bigl\{ \Bigl( \frac{s"}{s} \Bigr)^{\frac{1}{1000}}, \Bigl( \frac{s}{s"} \Bigr)^{\frac{1}{1000}} \Bigr\} \cdot \bigl( \mathtt g(s") + \epsilon \cdot {\ttc} (s") \bigr)   \\
                &\lesssim \sup_{s" \in (0, \underline s]} \min\Bigl\{ \Bigl( \frac{s"}{s} \Bigr)^{\frac{1}{1000}}, \Bigl( \frac{s}{s"} \Bigr)^{\frac{1}{1000}} \Bigr\} \cdot\max\Bigl\{ \Bigl( \frac{s"}{s} \Bigr)^{\frac{1}{1000}}, \Bigl( \frac{s}{s"} \Bigr)^{\frac{1}{1000}} \Bigr\}\bigl( \mathtt g(s) + \epsilon \cdot {\ttc} (s) \bigr) \\
                &\lesssim\mathtt g(s) + \epsilon \cdot {\ttc} (s).
        \end{align*}
    We can absorb the second term in the last line into the left-hand side provided that $\epsilon \ll 1$, as desired. 
\end{proof}

\begin{proposition}[Dyadic parabolic smoothing]\label{prop:no-loss}
    Let $(A, \phi)$ be a $\frE$-caloric extension, and suppose that $\ttc(s)$ is a $\frE$-frequency envelope for $(A, \phi)$. Then for each non-negative integer $n \in \N_0$, we have 
        \begin{equation}\label{eq:no-loss}
            {\bigl\| s^{-\frac12} (s^\frac12 \nabla)^{(\leq n)} (s \nabla^{(2)}) e^{s \Delta} \overline\Gauss \bigr\|_{L^2_x} + \bigl\| s^{-\frac12} (s^\frac12 \bfD)^{(\leq n)} (s \bfD^{(2)}) \phi (s) \bigr\|_{L^2_x}}
                \lesssim_n \ttc(s).
        \end{equation}
\end{proposition}

\begin{proof}
    Set 
        \[
            \tte^1_n [A, \phi] (s)
                := {\bigl\| s^{-\frac12} (s^\frac12 \nabla)^{(\leq n)} (s \nabla^{(2)}) e^{s \Delta} \overline\Gauss \bigr\|_{L^2_x} + \bigl\| s^{-\frac12} (s^\frac12 \bfD)^{(\leq n)} (s \bfD^{(2)}) \phi (s) \bigr\|_{L^2_x}},
        \]
    and let $\ttc_n (s)$ be an envelope for $\sqrt{\tte^1_n [A, \phi](s)}$. By the abstract bootstrap argument from Lemma \ref{lem:abs-env-boot}, it will suffice to show an estimate of the form 
        \[
            \sqrt{\tte^1_n [A, \phi] (s)}
                \lesssim \ttc(s) + \epsilon \cdot \ttc_n (s), \qquad \text{where }\epsilon \ll 1.
        \]

    For the contribution of the Gauss tension field, we can invoke parabolic smoothing \eqref{standard-heat} on $[\frac12s, s]$ and the slowly-varying property for $\ttc(s)$, 
        \begin{align*}
            \bigl\| s^{-\frac12} (s^\frac12 \nabla)^{(\leq n)} (s \nabla^{(2)}) e^{s \Delta} \overline\Gauss \bigr\|_{L^2_x} 
                \lesssim \bigl\| s^{-\frac\sigma2} (s \nabla^{(2)}) e^{s \Delta} \overline\Gauss \bigr\|_{L^2_x} \lesssim \ttc(\tfrac12 s) \sim \ttc(s),
        \end{align*}
    so it remains to check the contribution of the covariant Hessian. Applying the product rule to the right-hand side of the equation for the covariant Hessian,
        \begin{align*}
            (s^\frac12 \bfD)^{(n)} \bigl(\textup{R.H.S.\eqref{Dphi-heat-2}}\bigr)
                &= (s^\frac12 \nabla)^{(n)} \nabla F_{12}\cdot \bfD \phi + \sum_{a + b = n} (s^\frac12\nabla)^{(a)} F_{12} \cdot (s^\frac12\bfD)^{(b)} \bfD^{(2)} \phi\\
                &=: \mathrm{I}_n + \mathrm{II}_n.
        \end{align*}
    For the first term on the right, 
        \begin{align*}
            \bigl\| s^{\frac12 + 1} \cdot \mathrm{I}_n \bigr\|_{L^2_x}
                &\lesssim s^\frac12 \bigl\| s^{-\frac12} (s^\frac12 \nabla)^{(n)} (s \nabla^{(2)}) F_{12} \bigr\|_{L^2_x} \bigl\| \bfD \phi \bigr\|_{L^2_x} + s^\frac12 \bigl\| (s^\frac12 \nabla)^{(n + 1)} F_{12}\bigr\|_{L^2_x} \bigl\| s^{-\frac12} (s \bfD^{(2)}) \phi \bigr\|_{L^2_x} \\
                &\lesssim s^\frac12 \cdot \Bigl( 1 + \|\phi \|_{L^\infty_x} + \bigl\| \overline{\bfD \phi} \bigr\|_{L^2_x} + \bigl\| \overline F_{12} \bigr\|_{L^2_x} \Bigr) \cdot \ttc(s) \\
                &\lesssim \ttc(s), 
        \end{align*}
    using the abstract bilinear estimate from Lemma \ref{lem:bilinear} in the first line, parabolic smoothing \eqref{F12-sobolev}, \eqref{Dphi-sobolev} along with the envelope bound \eqref{F12-env} in the second line, and the $L^\infty$-bound \eqref{inefficient-smooth} along with the definition of $(\overline A, \overline{\vphantom{A}\phi}) \in \frE(\underline s)$ in the last line. For the second term on the right, 
        \begin{align*}
            \bigl\| s^{\frac12 + 1} \cdot \mathrm{II}_n \bigr\|_{L^2_x} 
                &\lesssim \sum_{a + b = n} s^\frac12 \bigl\| s^{\frac12} (s^\frac12 \nabla)^{(a)} F_{12} \bigr\|_{L^\infty_x} \bigl\| s^{-\frac12} (s^\frac12 \bfD)^{(b)} (s \bfD^{(2)}) \phi \bigr\|_{L^2_x} \\
                &\lesssim \underline s^{\frac12} \cdot \bigl\| \overline F_{12} \bigr\|_{L^2_x} \cdot \ttc_n (s) \\
                &\lesssim \epsilon \cdot \ttc_n (s), 
        \end{align*}
    using H\"older's inequality in the first line, parabolic smoothing \eqref{F12-sobolev} in the second line, and the definition of $(\overline A, \overline{\vphantom{A}\phi}) \in \frE(\underline s)$ in the last line. Inserting these into the dyadic smoothing bound \eqref{dyadic-smooth} with $\gamma = \tfrac12$ applied to the equation for the covariant Hessian \eqref{Dphi-heat-2}, we obtain 
        \begin{align*}
            \bigl\| s^{-\frac12} (s^\frac12 \bfD)^{(\leq n)} (s \bfD^{(2)}) \phi(s) \bigr\|_{L^2_x} 
                &\lesssim \bigl\|  s^{-\frac12} (s \bfD^{(2)}) \phi \bigr\|_{L^2_x} + \bigl\| {s'}^{\frac12 + 1} ({s'}^\frac12 \bfD)^{(\leq n)} (\text{R.H.S.\eqref{Dphi-heat-2}})\bigr\|_{L^2_{\frac{ds'}{s'}}L^2_x ([\frac12 s, s])}\\
                    &\lesssim \ttc(s) + \epsilon \cdot \sup_{s' \in [\frac12 s , s]}\ttc_n (s') \\
                    &\lesssim \ttc(s) + \epsilon \cdot \ttc_n (s)
            \end{align*}
        using the slowly-varying property in the second line. This completes the proof. 
\end{proof}

\subsection{Regularisation of initial data} \label{subsec:LP-regularise}

The last role played by the covariant heat equation is regularising rough sets of initial data. Before stating the main result, we should first clarify what we mean by a regular configuration, 

\begin{definition}
    Given a finite-energy configuration $(A, \phi) \in \frE$, we say that it is a \textit{smooth configuration} if for each non-negative integer $n \in \N_0$, 
        \[
            \bigl\| \bfD^{(n)} \bfD \phi \bigr\|_{L^2_x} + \bigl\| \nabla^{(n)} F_{12} \bigr\|_{L^2_x}
                < \infty.
        \] 
    We denote the space of smooth configurations by $\frE^\infty$, and the sub-class thereof satisfying Coulomb gauge \eqref{coulomb-data} by $\frE^\infty_\df$. 
\end{definition}

It is clear that if $(A, \phi) \in \frE^\infty_\df$, then $\nabla A_x \in H^\infty_x (\R^2)$. Furthermore, we can use Sobolev embedding to convert the control of magnetic derivatives into boundedness of usual derivatives, showing $\phi \in C^\infty_x (\R^2)$. These observations will play a mostly qualitative role, so we defer their proofs to Appendix \ref{sec:smooth-configurations}. 

Given a generic finite-energy configuration $(A, \phi) \in \frE$, we smooth it out by flowing along the covariant heat equation. Similar to how the linear heat propagator serves as a smooth approximation to the identity, the covariant heat equation satisfies a similar property with respect to the $\frE$-metric, 

\begin{proposition}[Continuity of the covariant heat flow]
    Let $(A, \phi) \in \pzcE^1 ([0, \underline s])$ be a $\frE$-caloric extension. Then \label{prop:cts-s}
        \[
            \dist_{\frE} \bigl( (A(s_1), \phi(s_1)), (A(s_2), \phi(s_2))\bigr)
                \overset{|s_1 - s_2| \to 0}{\longrightarrow} 0.
        \]
\end{proposition}

\begin{proof}
    The linear heat equation and transport equation generate strongly continuous semi-groups on $L^p (\R^2)$ for every $1 \leq p < \infty$, so the contributions of the magnetic potential clearly vanish as $|s_1 - s_2| \to 0$. Similarly, by Duhamel's formula, we can write for $s_1 < s_2$,
        \begin{align*}
            \Gauss(s_2) - \Gauss(s_1) 
                = \Bigl( e^{(s_2 - s_1) \Delta} \Gauss (s_1) - \Gauss (s_1) \Bigr) + \int_0^{s_2 - s_1} e^{(s_2 - s_1 - s) \Delta} (\bfD \phi \cdot \bfD \phi) \, ds
        \end{align*}
    Thus,
        \begin{align*}
            \lim_{|s_2 - s_1| \to 0}\bigl\| \nabla( \Gauss(s_2) - \Gauss(s_1) ) \bigr\|_{L^2_x} 
                &\leq \lim_{|s_2 - s_1| \to 0} \int_0^{s_2 - s_1} \bigl\| s^\frac12 \bfD \phi \bigr\|_{L^\infty_x} \bigl\| (s^\frac12 \bfD) \bfD \phi \bigr\|_{L^2_x} \, \frac{ds}{s}\\
                &\leq \lim_{|s_2 - s_1| \to 0}\bigl\| s^\frac12 \bfD \phi \bigr\|_{L^2_{\frac{ds}{s}} L^\infty_x ([0, s_2 - s_1])} \bigl\| (s^\frac12 \bfD) \bfD \phi \bigr\|_{L^2_{\frac{ds}{s}} L^2_x ([0, s_2 - s_1])} = 0,
        \end{align*}
    by monotone convergence theorem and \eqref{Dphi-sobolev}, \eqref{Dphi-strichartz-heat}. The estimate for $\Gauss$ in $L^2$ is similar but easier. 

   Rewriting the difference of the scalar fields using the fundamental theorem of calculus and the caloric gauge \eqref{caloric-gauge}, 
        \begin{align*}
            \phi(s_2) - \phi(s_1) 
                = \int_{s_1}^{s_2} \bfD^{(2)} \phi(s') \, ds' .
        \end{align*}
    Taking the $L^2$-norm, 
        \begin{align*}
            \bigl\| \phi(s_2) - \phi(s_1) \bigr\|_{L^2_x}
                &\lesssim \Bigl( \int_{s_1}^{s_2} |s'|^{-\frac12} \, ds' \Bigr) \bigl\| ({s}^\frac12 \bfD) \bfD \phi \bigr\|_{L^\infty_s L^2_x} \overset{|s_1 - s_2| \to 0}{\longrightarrow} 0,
        \end{align*}
    by the parabolic smoothing bounds \eqref{Dphi-sobolev}. 

    For the covariant derivatives of the scalar field, we again rewrite the difference using the fundamental theorem of calculus in $s$ and the caloric gauge \eqref{caloric-gauge}, 
        \begin{align*}
            \bfD \phi (s_2) - \bfD \phi(s_1) 
                &= \int_{s_1}^{s_2} \bfD_s \bfD \phi (s') \, ds' \\
                &= \int_{s_1}^{s_2} ({s'}^\frac12 \bfD)^{(2)} \bfD \phi (s') \, \frac{ds'}{s'} + \int_{s_1}^{s_2} F_{12} (s')\cdot \bfD \phi (s') \, ds',
        \end{align*}
    using the equation \eqref{Dphi-heat} in the second line. It will suffice to show that these two terms on the right vanish in $L^2_x$-norm as $|s_1 - s_2| \to 0$. For the lower-order terms, 
        \begin{align*}
            \Bigl\| \int_{s_1}^{s_2} F_{12} (s') \cdot \bfD \phi(s') \, ds' \Bigr\|_{L^2_x} 
                &\lesssim \Bigl( \int_{s_1}^{s_2} |s'|^{-\frac12} ds' \Bigr) \bigl\| s^\frac12 F_{12} \bigr\|_{L^\infty_{s, x}}  \bigl\| \bfD \phi \bigr\|_{L^\infty_s L^2_x} \overset{|s_1 - s_2| \to 0}{\longrightarrow} 0,
        \end{align*}
    by H\"older's inequality and parabolic smoothing \eqref{F12-sobolev}, \eqref{Dphi-sobolev}. For the top-order term, we argue by almost-orthogonality, expanding out the $L^2$-norm as an inner product, 
        \begin{align*}
            \Bigl\| \int_{s_1}^{s_2} ({s'}^\frac12 \bfD)^{(2)} \bfD\phi (s') \, \frac{ds'}{s'}  \Bigr\|_{L^2_x}^2
                &= 2\int_{s_1}^{s_2} \int_{s"}^{s_2} \Bigl\langle ({s'}^\frac12 \bfD)^{(2)} \bfD \phi (s') , ({s"}^\frac12 \bfD)^{(2)} \bfD \phi (s")\Bigr\rangle_{L^2_x} \frac{ds'}{s'} \frac{ds"}{s"} 
        \end{align*}
    using symmetry to write this as an integral on $s" < s'$. Then by Lemma \ref{lem:orthogonal},
        \begin{align*}
            \Bigl\| \int_{s_1}^{s_2} ({s'}^\frac12 \bfD)^{(2)} \bfD\phi (s') \, \frac{ds'}{s'}  \Bigr\|_{L^2_x}^2
                &\lesssim \bigl\| ({s}^\frac12 \bfD)^{(3)} \bfD \phi \bigr\|_{L^2_{\frac{ds}{s}} L^2_x ([s_1, s_2])} \bigl\| ({s}^\frac12 \bfD) \bfD \phi \bigr\|_{L^2_{\frac{ds}{s}} L^2_x ([s_1, s_2])} \overset{|s_1 - s_2| \to 0}{\longrightarrow} 0,
        \end{align*}
    concluding convergence via parabolic smoothing \eqref{Dphi-sobolev} and the monotone convergence theorem. 
\end{proof}

\begin{corollary}\label{cor:smooth-data}
    Let $(A, \phi) \in \frE(\underline s)$ be a finite-energy configuration. Then there exists a sequence $\{ (\smoothing{n}{A}, \smoothing{n}{{\vphantom{A}\phi}}) \}_n \subseteq \frE^\infty$ of smooth configurations satisfying:
    \begin{itemize}
        \item \textup{(Uniform energy bounds).}
            \[
                (\smoothing{n}{ A}, \smoothing{n}{{\vphantom{A}\phi}})
                    \in \frE(\underline s)
            \]
        
        \item \textup{(Uniform Gauss tension field bounds).}
            \[
               \lim_{n \to \infty}\| \smoothing{n}{\Gauss} \|_{L^\infty_x}
                    \leq \| \Gauss \|_{L^\infty_x}.
            \]

        \item \textup{(Convergence in $\frE$).}
            \[
                \lim_{n \to \infty} \dist_\frE \bigl( (\smoothing{n}{A}, \smoothing{n}{{\vphantom{A}\phi}}), ({A}, {{\vphantom{A}\phi}}) \bigr) = 0.  
            \]
    \end{itemize}
    Furthermore, if $(A, \phi)$ satisfies the Coulomb gauge \eqref{coulomb-data}, then $\{ (\smoothing{n}{A}, \smoothing{n}{{\vphantom{A}\phi}}) \}_n$ also satisfy the Coulmb gauge. 
\end{corollary}

\begin{proof}
    Let $(A, \phi) \in \pzcE^1 ([0, \underline s])$ be the $\frE$-caloric extension of $(A, \phi)$, and set 
        \[
             (\smoothing{n}{A}, \smoothing{n}{{\vphantom{A}\phi}})
                := (A, \phi)_{|s = \frac1n}, \qquad n \gg 1.
        \]
    The covariant estimates from Section \ref{subsec:LP-sobolev} imply that these are smooth configurations and the Gauss tension field is uniformly bounded. Furthermore, by Proposition \ref{prop:cts-s}, we have convergence of the $\frE$-metric, and, since $\frE(\underline s)$ is characterised by an open condition, this sequence can be chosen to reside in $\frE(\underline s)$. Since the curl-free part of the magnetic potential is transported \eqref{A-transport}, the Coulomb gauge condition is also propagated. This completes the proof. 
\end{proof}

\section{Continuation of smooth solutions}\label{sec:high-LWP}

To set some expectations, we outline a local ``well-posedness"\footnote{A proper well-posedness theory should include a notion of continuity of the data-to-solution map. For our purposes herein, we emphasise only the local existence of a smooth solution, though continuous dependence on data with respect to $H^N$-perturbations most likely follows by standard arguments. On the other hand, the examples from Appendix \ref{subapp:example} suggest that $H^1$-perturbations are much more restrictive compared to those allowed by the energy topology as defined in Definition \ref{def:topology-1}, so we elect not to pursue this point further. } theory for smooth solutions to the Chern-Simons-Schr\"odinger equation, which consists of (\texttt{i}) local existence of a solution arising from smooth initial data, (\texttt{ii}) higher-order energy estimates, (\texttt{iii}) a quantitative continuation criterion, and (\texttt{iv}) conservation of the abelian Higgs energy and topological degree. Herein we will emphasis the quantitative aspects, (\texttt{ii}) and (\texttt{iii}), and defer the proofs of the qualitative aspects, (\texttt{i}) and (\texttt{iv}), to the appendices. 

As the phase space for the Chern-Simons-Schr\"odinger equation is rather non-standard, we should clarify the notion of a smooth solution. Taking cue from the abelian Higgs energy, for each positive integer $N \in \N$, define the \textit{$N$-th order covariant energy} of a configuration $(A, \phi)$ on $\R^2$ by 
    \begin{align*}
        \cE^N [A, \phi] 
            &:= \bigl\| 1 - |\phi|^2 \bigr\|_{L^2_x}^2 + \bigl\| \Gauss \bigr\|_{H^N}^2  + \sum_{n = 1}^N \bigl\| \bfD^{(n)} \phi \bigr\|_{L^2_x}^2.
    \end{align*}
With these higher-order energies at hand, we introduce our notion of ``regular'' solutions. One should not read into the definition too much -- the moral of the story is that there exists a class of classical solutions to the Chern-Simons-Schr\"odinger equation wherein we can justify the technical calculations within this article with impunity. 

\begin{definition}\label{def:smooth-1}
    We say that a configuration $(A, \phi)$ on $[0, T] \times \R^2$ is a \textit{smooth solution} to the Chern-Simons-Schr\"odinger equation if it solves \eqref{CSS-amp-schro} and satisfies the following regularity properties:
    \begin{itemize} 
        \item The configuration $(A_x(t), \phi(t)) \in \frE^\infty$ is smooth for each $t \in [0, T]$, and obeys the uniform bounds
            \[
                \sup_{t \in [0, T]} \cE^N [A_x(t), \phi(t)] 
                    < \infty, \qquad \text{for each $N \in \N$.}
            \]

        \item It forms a continuous flow into the space of finite-energy configurations, $(A_x, \phi) \in C_t ([0, T] \to \frE)$. 
        
        \item The electromagnetic potential has regularity $A_t, \nabla A_x \in C^\infty_t H^\infty_x ([0, T] \times \R^2)$. 
        
        \item The scalar field has regularity $\phi \in C^\infty_{t, x} ([0, T] \times \R^2)$. 
    \end{itemize}
\end{definition}

The main protagonist in our analysis will be played by the \textit{control parameter}, a gauge-invariant quantity which serves to bound the growth of energy along the Chern-Simons-Schr\"odinger flow, 
    \[
        \cC[A, \phi]  
            :=  1 + \|\Gauss \|_{L^\infty_x} + \| \bfD \phi \|_{L^2_x}^2  + \| \phi \|_{L^\infty_x}^2.
    \]
The stage is now set for us to state the local existence of smooth solutions, 

\begin{theorem}[Local existence of smooth solutions]\label{thm:smooth-LWP}
    Let $(A^{\mathrm{in}}, \phi^{\mathrm{in}}) \in \frE^\infty$ be a smooth configuration on $\R^2$ in Coulomb gauge \eqref{coulomb-data}, then there exists $T \equiv T(\cE^{9000} [A^{\mathrm{in}}, \phi^{\mathrm{in}}]) > 0$ and a smooth solution $(A, \phi)$ on $[0, T] \times \R^2$ to the Chern-Simons-Schr\"odinger equation \eqref{CSS-amp-schro} in DeTurck gauge \eqref{DeTurck} with initial data $(A, \phi)_{|t = 0} = (A^{\mathrm{in}}, \phi^{\mathrm{in}})$. Furthermore, the solution obeys the following properties: 
    \begin{itemize}        
        \item \textup{(Covariant energy estimates).} For each positive integer $N \in \N$, 
        \begin{equation}\label{eq:EN-bound-1}
            \frac{d}{dt} \cE^N [A, \phi] 
                \lesssim_N \cC[A, \phi] \cdot \cE^N [A, \phi].
        \end{equation}  

        \item \textup{(Continuation criterion).} The solution can be extended to a smooth solution on $[0, T^+]\times \R^2$ solving the Chern-Simons-Schr\"odinger equation for some $T^+ > T$ whenever 
        \begin{equation}\label{eq:high-blow-up}
            \int_0^T \cC[A(t), \phi(t)] \, dt  
                < \infty.
        \end{equation}
        
        \item \textup{(Conservation of energy and degree).} The abelian Higgs energy and topological degree are conserved by the flow, i.e. for all $t \in [0, T]$, 
        \begin{equation}\label{eq:high-conserve}
            \cE_{\textup{AH}}[A(t), \phi(t)] 
                = \cE_{\textup{AH}} [A^{\mathrm{in}}, \phi^{\mathrm{in}}],
        \end{equation}
        and 
        \begin{equation}
            \deg [A(t), \phi(t)] 
                = \deg [A^{\mathrm{in}}, \phi^{\mathrm{in}}].
        \end{equation}
    \end{itemize}    
\end{theorem}

We defer the proofs of local existence and the conservation laws to Appendices \ref{app:picard} and \ref{app:conservation} respectively. To say just a word on both, the latter is a covariant analogue of the usual density-flux calculations, while for the former, we in fact prove a more general result, namely local existence within the affine Sobolev space
    \[
        (\mathring A, \mathring \phi) + H^N (\R^2) 
            := \Bigl\{ (A, \phi) \in H^N_{\loc} (\R^2) : (A, \phi) - (\mathring A, \mathring \phi) \in H^N (\R^2) \Bigr\},
    \]  
where $(\mathring A, \mathring \phi) \in \frE^\infty$ is a smooth configuration on $\R^2$ in Coulomb gauge \eqref{coulomb-data}, provided the initial data satisfies $\Gauss^{\mathrm{in}} \in H^N (\R^2)$ and sufficiently high regularity $N \gg 1$. The proof proceeds by a standard energy argument and Picard iteration -- the main non-trivial takeaway is that the time-interval of existence may be taken to depend only on the $N$-th order covariant energy. The continuation criterion \eqref{high-blow-up} then follows as a corollary of the covariant energy estimate \eqref{EN-bound-1} and Gronwall's inequality. 

The continuation criterion \eqref{high-blow-up} sets the following benchmark for our desired well-posedness theory, 

\begin{proposition}[Control parameter bound]\label{prop:control}
    Let $(A, \phi)$ be a smooth solution on $[0, T] \times \R^2$ to the Chern-Simons-Schr\"odinger equation \eqref{CSS-amp-schro}. Then there exists $T_* \equiv T_* (\| \Gauss^{\mathrm{in}} \|_{L^\infty \cap H^1}, \cE_{\textup{AH}} [A^{\mathrm{in}}, \phi^{\mathrm{in}}]) > 0$ such that if $J \subseteq [0, T]$ is a sub-interval of length $|J| \leq T_*$, then
        \begin{equation}\label{eq:control-bound}
            \int_J \cC[A(t), \phi(t)] \, dt 
                \leq 2. 
        \end{equation}
\end{proposition}

We defer the proof to Section \ref{sec:energy}, though to give some suggestions of the strategy, note that the statement of Proposition \ref{prop:control} is gauge-invariant and time-translation invariant. Thus, we are not necessarily beholden to the choices \eqref{DeTurck}+\eqref{coulomb-data}, and indeed, we will conduct the bulk of our analysis in an auxiliary gauge with better regularity properties, see Section \ref{sec:gauge}. 

Using Proposition \ref{prop:control}, we can close the energy estimate \eqref{EN-bound-1} within the space of finite-energy configurations $\frE$. While this by itself does not constitute a complete argument for well-posedness, it nonetheless serves as compelling evidence for thereof. Moreover, we deduce global existence of smooth solutions by a standard argument combining the continuation criterion \eqref{high-blow-up}, the control parameter bound \eqref{control-bound}, and conservation of the Gauss tension field \eqref{gauss-transport-intro} and abelian Higgs energy \eqref{high-conserve}, 
    
\begin{theorem}[Global regularity of smooth solutions]\label{thm:global-1}
    Let $(A^{\mathrm{in}}, \phi^{\mathrm{in}}) \in \frE^\infty$ be a smooth configuration on $\R^2$ in Coulomb gauge \eqref{coulomb-data}. Then there exists a smooth solution $(A, \phi)$ on $[0, \infty) \times \R^2$ to the Chern-Simons-Schr\"odinger equation \eqref{CSS-amp-schro} in DeTurck gauge \eqref{DeTurck} with the initial data $(A, \phi)_{|t = 0} = (A^{\mathrm{in}}, \phi^{\mathrm{in}})$. 
    
    Furthermore, for each positive integer $N \in \N$, the solution obeys the covariant energy estimate 
        \begin{equation}\label{eq:EN-bound-2}
            \cE^N [A (t), \phi (t)] 
                \leq \cE^N [A^{\mathrm{in}}, \phi^{\mathrm{in}}] \exp (\kappa t),
        \end{equation} 
    for some constant $\kappa \equiv \kappa(N, \| \Gauss^{\mathrm{in}} \|_{L^\infty \cap H^1}, \cE_{\textup{AH}} [A^{\mathrm{in}}, \phi^{\mathrm{in}}])$.
\end{theorem}

The remainder of this section will be dedicated to the proof of the $N$-th order covariant energy estimate \eqref{EN-bound-1}. Since the Gauss tension field obeys the transport equation \eqref{gauss-transport-intro}, it will suffice to prove the differential inequality for the charge density and the covariant derivatives of $\phi$. For the latter, we argue by the energy method for the covariant Schr\"odinger equations, 

\begin{lemma}[Abstract energy estimate]\label{lem:abs-schro-energy}
    Let $A$ be a connection $1$-form on $[0, T] \times \R^2$, and $\sfV : [0, T] \times \R^2 \to \R$ a real-valued potential. Then for sufficiently regular solutions $u : [0, T] \times \R^2 \to \C$ to the electromagnetic Schr\"odinger equation,
        \[
            \big(i \bfD_t + \bfD^j \bfD_j + \sfV \big) u
                = \sfN,
        \]
    the following energy estimate holds, 
        \begin{equation}\label{eq:abstract-energy}
            \frac{d}{dt} \| u \|_{L^2_x}^2
                \leq \| u \|_{L^2_x} \| \sfN \|_{L^2_x}. 
        \end{equation}
\end{lemma}

\begin{proof}
    By the product rule \eqref{product} and integration-by-parts, 
        \begin{align*}
            \frac{d}{dt} \| u \|_{L^2_x}^2
                &= \langle u, \bfD_t u \rangle_{L^2_x}\\
                &= \langle u , i (\bfD^j \bfD_j + \sfV) u \rangle_{L^2_x} - \langle u, i \sfN \rangle_{L^2_x} \\
                &= \langle u, - i \sfN \rangle_{L^2_x}.
        \end{align*}
    We conclude \eqref{abstract-energy} by Cauchy-Schwarz. 
\end{proof}

\begin{lemma}[Higher-order Schr\"odinger equations]\label{lem:higher-schrodinger}
    Let $(A, \phi)$ be a solution to the Chern-Simons-Schr\"odinger equation \eqref{CSS-amp-schro}. Then for each positive integer $n \in \N$, we have the following schematic equation, 
        \begin{equation}\label{eq:higher-schrodinger}
        \begin{split}
            \Bigl(i \bfD_t + \bfD^j \bfD_j + \tfrac\lambda2 (1 - |\phi|^2) \Bigr) \bfD^{(n)} \phi 
                &= \bfD^{(n)} \phi + \sum_{a + b + c = n} \bfD^{(a)} \phi \cdot \bfD^{(b)} \phi \cdot \bfD^{(c)} \phi\\
                &\qquad + \sum_{a + b = n} \nabla^{(a)} \Gauss \cdot \bfD^{(b)} \phi .
        \end{split}
        \end{equation}
\end{lemma}

\begin{proof}
    We have the commutator identity
        \begin{align*}
            \Bigl[\bfD, i \bfD_t + \bfD^j \bfD_j + \tfrac\lambda2(1 - |\phi|^2) \Bigr] 
                &= F_{tx} + F_{12} \cdot \bfD  + \nabla F_{12}  + \phi \cdot \bfD \phi \\
                &= F_{12} \cdot \bfD  + \nabla F_{12}  + \phi \cdot \bfD \phi\\
                &= \bfD + \phi \cdot \phi \cdot \bfD + \Gauss \cdot \bfD + \nabla \Gauss + \phi \cdot \bfD \phi
        \end{align*}
    using the basic commutator identities \eqref{commute}, \eqref{laplace-commute} in the first line, and writing Amp\'ere law's schematically as $F_{tx} = \phi \cdot \bfD \phi + \nabla F_{12}$ in the second line, and re-expressing the magnetic field using the Gauss tension field and charge density in the third line. The desired equation \eqref{higher-schrodinger} then follows from an inductive application of the product rule and the commutator identity above. 
\end{proof}

\begin{proposition}[Higher-order energy estimates]\label{prop:aprioriestimate}
    Let $(A, \phi)$ be a smooth solution on $[0, T] \times \R^2$ to the Chern-Simons-Schr\"odinger equation \eqref{CSS-amp-schro}. Then for each positive integer $n \in \N$, we have 
        \begin{equation}\label{eq:energyestimate}
        \begin{split}
            \frac{d}{dt} \bigl\| \bfD^{(n)} \phi \bigr\|_{L^2_x}^2 
                &\lesssim_n \cC[A, \phi] \cdot \Bigl( \bigl\| \nabla^{(n)} \Gauss \bigr\|_{L^2_x}^2 + \bigl\| \bfD^{(n)} \phi \bigr\|_{L^2_x}^2  \Bigr).
        \end{split}
        \end{equation}
\end{proposition}

\begin{proof}
    Applying the abstract energy estimate \eqref{abstract-energy} to the covariant Schr\"odinger equation \eqref{higher-schrodinger}, it remains to estimate the right-hand side of the equation in $L^2 (\R^2)$. The linear term can be estimated immediately. For the non-linear terms, heuristically, the worst terms on the right-hand side of the equation \eqref{higher-schrodinger} are the interactions with ``unbalanced" derivatives, namely 
    \[
        \phi \cdot \phi \cdot \bfD^{(n)} \phi ,\qquad \phi \cdot \nabla^{(n)} \Gauss, \qquad \Gauss \cdot \bfD^{(n)} \phi , \qquad \phi \cdot \bigl( \bfD^{(a)} \phi \cdot \bfD^{(b)} \phi\bigr) ,
    \]
    where $a, b \neq 0$. We estimate each term by $L^\infty \times L^2 \to L^2$, and the magnetic interpolation inequality \eqref{GN-interpolation} for the last term, 
        \begin{align*}
            \bigl\| \text{unbalanced} \bigr\|_{L^2_x}
                &\lesssim \Bigl( \bigl\| \Gauss \bigr\|_{L^\infty_x} +  \bigl\| \phi \bigr\|_{L^\infty_x}^2  \Bigr) \bigl\| \bfD^{(n)} \phi \bigr\|_{L^2_x} + \bigl\| \phi \bigr\|_{L^\infty_x} \bigl\| \nabla^{(n)} \Gauss \bigr\|_{L^2_x} \\
                    &\qquad + \sum_{\substack{a + b = n \\ a, b \neq 0}} \bigl\| \phi \bigr\|_{L^\infty_x} \bigl\| \bfD^{(a)} \phi \bigr\|_{L^{2n/a}_x} \bigl\| \bfD^{(b)} \phi \bigr\|_{L^{2n/b}_x}  \\
                &\lesssim \cC[A, \phi] \cdot \Bigl( \bigl\| \nabla^{(n)} \Gauss \bigr\|_{L^2_x} + \bigl\| \bfD^{(n)} \phi \bigr\|_{L^2_x} \Bigr).
        \end{align*}
    When the derivatives are instead ``balanced'', namely the terms 
        \begin{align*}
            \bfD^{(a)} \phi \cdot \bfD^{(b)}  \phi \cdot \bfD^{(c)} \phi , \qquad \nabla^{(a)} \Gauss \cdot \bfD^{(b)} \phi, 
        \end{align*}
    where $a, b, c \neq 0$, we use H\"older's inequality, the magnetic interpolation inequality \eqref{GN-interpolation} and the usual Gagliardo-Nirenberg interpolation,
        \begin{align*}
            \bigl\| \text{balanced}  \bigr\|_{L^2_x} 
                &\lesssim \sum_{\substack{a + b + c = n \\ a, b, c \neq 0}}\bigl\| \bfD^{(a)} \phi \bigr\|_{L^{2n/a}_x} \, \bigl\| \bfD^{(b)} \phi \bigr\|_{L^{2n/b}_x} \, \bigl\| \bfD^{(c)} \phi \bigr\|_{L^{2n/c}_x} + \sum_{\substack{a + b = n \\ a, b \neq 0}}\bigl\| \nabla^{(a)} \Gauss \bigr\|_{L^{2n/a}_x} \, \bigl\| \bfD^{(b)} \phi \bigr\|_{L^{2n/b}_x}  \\
                &\lesssim \Bigl( \bigl\|\Gauss \bigr\|_{L^\infty_x}   + \bigl\| \bfD \phi \bigr\|_{L^2_x} \Bigr) \Bigl( \bigl\| \nabla^{(n)} \Gauss \bigr\|_{L^2_x} + \bigl\| \bfD^{(n)} \phi \bigr\|_{L^2_x} \Bigr).
        \end{align*}
   These are acceptable in view of the right-hand side of the energy estimate \eqref{energyestimate}.
\end{proof}

To propagate the charge density, we apply the energy method to the density-flux identity. This loses a derivative, however the top-order bound \eqref{energyestimate} controls $\bfD \phi$, so the estimates in concert close.

\begin{proposition}[Charge density estimate]\label{prop:charge-density}
    Let $(A, \phi)$ be a smooth solution on $[0, T] \times \R^2$ to the Chern-Simons-Schr\"odinger equation \eqref{CSS-amp-schro}, then 
        \begin{equation}\label{eq:apriori-charge}
            \frac{d}{dt} \Bigl( \big\| 1 - |\phi|^2 \big\|_{L^2_x}^2 \Bigr)
                \lesssim \| \phi\|_{L^\infty_x}^2 \| \bfD \phi \|_{L^2_x}^2. 
        \end{equation}
\end{proposition}

\begin{proof}
    Differentiating the charge density in $t$, using the covariant Schr\"odinger equation, and integrating-by-parts appropriately, we obtain
        \begin{align*}
            \partial_t \tfrac12(1 - |\phi|^2)^2
                &= - (1 - |\phi|^2) \, \Re(\overline \phi \bfD_t \phi) \\
                &=  (1 - |\phi|^2) \Im (\overline \phi \bfD^j \bfD_j \phi) \\
                &= (1 - |\phi|^2) \, \partial^j \Im(\overline \phi \bfD_j \phi) \\
                &= \partial^j \Bigl( (1 - |\phi|^2) \Im(\overline \phi \bfD_j \phi) \Bigr) + \Re(\overline \phi \bfD^j \phi) \Im(\overline \phi \bfD_j \phi). 
        \end{align*}
    Integrating on $\R^2$ and applying the divergence theorem yields the differential inequality.
\end{proof}

\begin{remark}
   One can view the top-order energy estimate \eqref{energyestimate} as propagating ``frequency-localisation", which effectively confers ``finite-speed of propagation'' for $\phi$ as a solution to a Schr\"odinger equation and thereby control over lower-order quantities such as the charge density as in \eqref{apriori-charge}. It is instructive to compare with the estimates for the linear Schr\"odinger flow on Zhidkov spaces, c.f. \cite{Gallo2004} and \cite[Section 2.2]{Gerard2006}, which play a similar role in the well-posedness theory for the Gross-Pitaevskii equation. 
\end{remark}

\section{Paradifferential formulation of \eqref{CSS-amp-schro}}\label{sec:paradiff}

Rather than directly working with the Chern-Simons-Schr\"odinger equation, we conduct the bulk of our analysis on a paradifferential-style decomposition of \eqref{CSS-amp-schro} induced by the \textit{dynamic covariant heat equation} on $[0, T] \times \R^2 \times [0, \underline s]$, which consists of the covariant heat equation \eqref{covariant-heat}-\eqref{covariant-heat2} augmented by an equation for the $s$-dynamics of the electric potential $A_t$, 
    \begin{align}
            \bigl( \bfD_s - \bfD^\ell \bfD_\ell \bigr) \phi 
                &= 0, \tag{dP-$\phi$} \label{eq:d-covariant-heat}\\
            F_{s\mu} 
                &= \partial^\ell F_{\ell \mu}. \tag{dP-$A$} \label{eq:d-covariant-heat2}
    \end{align}
The system of equations is covariant under the gauge transformation 
	\begin{align*}
        A_\bfa 
            &\mapsto A_\bfa + \partial_\bfa \chi, \\
        \phi 
            &\mapsto e^{i \chi} \phi,
    \end{align*}
for any real scalar field $\chi : [0, T] \times \R^2 \times [0, \underline s] \to \R$. Given a configuration $(\overline A, \overline{\vphantom{A}\phi})$ on $[0, T] \times \R^2$ solving the Chern-Simons-Schr\"odinger equation, we impose it as heat-temporal initial data for the dynamic covariant heat equation \eqref{d-covariant-heat}-\eqref{d-covariant-heat2} and solve forward in $s$. Altogether, we consider a system of equations which we will refer to as the \textit{paradifferential formulation of the Chern-Simons-Schr\"odinger equation},
	\begin{equation}\tag{dP-CSS}\label{eq:dP-CSS}
        \begin{split}
        \bigl( \bfD_s - \bfD^\ell \bfD_\ell \bigr) \phi
				&= 0 
                \hspace{12em}
                \raisebox{-.7\normalbaselineskip}[0pt][0pt]
                {%
                    on $[0, T] \times \R^2 \times [0, \underline s]$,
                } \\
		F_{s\mu}
				&= \partial^\ell F_{\ell \mu} \\
        \Bigl( i \bfD_t + \bfD^j \bfD_j + \tfrac\lambda2 (1 - |\phi|^2) \Bigr) \phi 
				&= 0
                \hspace{12em}
                \raisebox{-.7\normalbaselineskip}[0pt][0pt]
                {%
                    on $[0, T] \times \R^2 \times \{0\}.$
                }\\
		F_{tj}
				&= - \epsilon_{jk} \Big(\mathfrak \Im(\overline \phi \bfD^k \phi) + \partial^\ell F\indices{_\ell^k}\Big)
        \end{split}
     \end{equation}

\begin{remark}
    The idea of using a geometric heat flow as a tool for smoothing out a geometric dispersive equation was introduced by Tao \cite{Tao2004} in his work on global regularity for wave maps from $\R^2$ to hyperbolic space. Various authors have successfully adopted the approach for other problems, e.g. Bejenaru-Ionescu-Kenig-Tataru \cite{BejenaruEtAl2011} for Schr\"odinger maps, Oh \cite{Oh2014a, Oh2015a}, Gavrus \cite{Gavrus2022}, and Oh-Tataru \cite{OhTataru2019, OhTataru2020, OhTataru2021, OhTataru2022} for the hyperbolic Yang-Mills equation, A reader familiar with the literature will note that our approach slightly differs from the traditional implementation, which would dictate that we use the abelian Higgs gradient flow as opposed to the covariant heat equation. Nevertheless, there do not seem to be any advantages to using the gradient flow for studying the well-posedness of the Chern-Simons-Schr\"odinger equation. 
\end{remark}

Of course, before we can make use of the paradifferential formulation of the Chern-Simons-Schr\"odinger equation, we should address the existence of solutions arising from solutions to the original Chern-Simons-Schr\"odinger equation. Paralleling the notion of smooth solution to \eqref{CSS-amp-schro} from Definition \ref{def:smooth-1}, we define

\begin{definition}\label{def:smooth-dP}
    We say that a configuration $(A, \phi)$ on $[0, T] \times \R^2 \times [0, \underline s]$ is a \textit{smooth solution} to the paradifferential formulation of the Chern-Simons-Schr\"odinger equation if it solves the system \eqref{dP-CSS} and satisfies the following regularity properties: 
    \begin{itemize}        
        \item The configuration $(A_x(t,s ), \phi(t,s)) \in \frE^\infty$ is smooth for each $t \in [0, T]$ and $s \in [0, \underline s]$, and obeys the uniform bounds for each positive integer $N \in \N$,
            \[
                \sup_{s \in [0, \underline s]} \sup_{t \in [0, T]} \cE^N [A_x(t,s), \phi(t, s)]
                    < \infty. 
            \]

        \item The heat-temporal initial data satisfies $(\overline A_x(t), \overline{\vphantom{A}\phi}(t)) \in \frE(\underline s)$ for each $t \in [0, T]$. Moreover, it forms a continuous flow into the space of finite-energy configurations, i.e. $(\overline A_x, \overline{\vphantom{A}\phi}) \in C_t ([0, T] \to \frE(\underline s))$. 
        
        \item The gauge potential has regularity $A_t, A_s, \nabla A_x \in C^\infty_{t, s} H^\infty_x (\R^2)$. 
        
        \item The scalar field has regularity $\phi \in C^\infty_{t, x, s} ([0, T] \times \R^2 \times [0, \underline s])$.
    \end{itemize}
\end{definition}

Just as before, one should not take too much stock into all the details of the definition. The main takeway is that when we say $(A, \phi)$ is a smooth solution on $[0, T] \times \R^2 \times [0, \underline s]$, it satisfies the equations \eqref{dP-CSS} in a classical sense, we can rigorously justify all the ensuing calculations, and the parameter $\underline s$ controls the size of the abelian Higgs energy and the Gauss tension field. Indeed, recalling the definition of $\frE(\underline s)$, we have 
    \[
        1 + \cE_{\text{AH}} [\overline A, \overline{\vphantom{A}\phi}] +\bigl\| \nabla \overline\Gauss \bigr\|_{L^2_x}  
            \lesssim \underline s^{-1}. 
    \]

\begin{proposition}[Existence of smooth solutions for \eqref{dP-CSS}]\label{prop:dP-exist}
    Let $(\overline A, \overline{\vphantom{A} \phi})$ be a smooth solution on $[0, T] \times \R^2$ to the Chern-Simons-Schr\"odinger equation \eqref{CSS-amp-schro} with Schr\"odinger-temporal initial data $(\overline A, \overline{\vphantom{A}\phi})_{|t = 0} \in \frE(\underline s)$. Then there exists a smooth solution $(A, \phi)$ on $[0, T] \times \R^2 \times [0, \underline s]$ to the paradifferential formulation of the Chern-Simons-Schr\"odinger equation \eqref{dP-CSS} with heat-temporal initial data $(A, \phi)_{s = 0} = (\overline A, \overline{\vphantom{A}\phi})$. 
\end{proposition}

\begin{proof}
    Imposing the DeTurck gauge in $s$,  
        \[
            A_s 
                = \partial^\ell A_\ell,
        \]
    it follows from a straight-forward calculation that the dynamic Yang-Mills heat flow \eqref{d-covariant-heat2} reduces to the linear heat equation for the electromagnetic potential, 
        \[
            (\partial_s - \Delta) A_\mu 
                = 0.
        \]
    For each $t$, we first solve in $s$ for the electromagnetic potential, and then the electromagnetic heat equation for the scalar field \eqref{d-covariant-heat} with $(\overline A, \overline{\vphantom{A} \phi})$ as heat-temporal initial data. Since $(\overline A, \overline{\vphantom{A} \phi})$ is a smooth solution to \eqref{CSS-amp-schro}, we by definition know that $(\overline A, \overline{\vphantom{A} \phi}) \in C_t ([0, T] \to \frE)$ and, by conservation of the abelian Higgs energy and transport of the Gauss tension field, we have $(\overline A_x (t), \overline{\vphantom{A}\phi}(t)) \in \frE(\underline s)$ for all $t$. The ensuing configurations $(A_x (t, s), \phi(t, s)) \in \frE^\infty$ arising from the dynamic covariant heat flow are smooth uniformly in $t$ and $s$ by virtue of the covariant energy estimates \eqref{charge-sobolev}, \eqref{Dphi-sobolev-high}, \eqref{gauss-sobolev-high}. The remaining desired regularity properties can easily be read-off from the equations of motion and the calculations from Appendix \ref{sec:smooth-configurations} -- we leave their verification to the reader. 
\end{proof}

\begin{remark}
    Another approach is to impose the caloric gauge and apply Lemma \ref{lem:lwp-heat}, which yields a solution to the covariant-heat equation \eqref{covariant-heat}-\eqref{covariant-heat2} for each fixed $t$. From here, it would remain to construct the electric potential $A_t$, which can be done by formally solving the linear heat equation for the electric field \eqref{Fmunu-heat} and then the transport equation for $A_t$ arising from \eqref{d-covariant-heat2} in the caloric gauge. 
\end{remark}

Our main goal herein is to derive a Schr\"odinger equation for the covariant Hessian $(s \bfD^{(2)}) \phi$, which will play as the leading protagonist in the analysis of the Chern-Simons-Schr\"odinger equation. More precisely, we want to think of the covariant Hessian as high-frequency wave evolving along a linear electromagnetic Schr\"odinger equation with potentials localised to lower-frequencies, and forced by terms we may regard as perturbative on the right,
    \[
        \underbrace{\bigl( i \bfD_t + \bfD^j \bfD_j \bigr)}_{\text{low frequency}} \underbrace{(s \bfD^{(2)}) \phi}_{\text{high frequency}} = \text{perturbative}.
    \]
We will make the above heuristic rigorous in Section \ref{sec:tension}-\ref{sec:energy}. To facilitate the analysis, we also build a supporting cast of covariant parabolic equations for the electromagnetic field, along with the \textit{Schr\"odinger} and \textit{Amp\'ere tension fields}, defined respectively by 
     \begin{align*}
        \CSS
            &:= \left(i \bfD_t + \bfD^j \bfD_j + \tfrac\lambda2 (1 - |\phi|^2)\right) \phi,\\
        \Amp_j
            &:= F_{tj} + \epsilon_{jk} \left( \Im(\overline \phi \bfD^k \phi) + \partial^\ell F\indices{_\ell^k} \right).
    \end{align*}
These measure the extent to which a configuration fails to satisfy the Chern-Simons-Schr\"odinger equation \eqref{CSS-amp-schro} -- it follows from the definitions that a solution to the system \eqref{dP-CSS} satisfies 
    \[
        (\Amp, \CSS)_{|s = 0} = (0, 0).
    \]

\begin{proposition}[Equations of motion for \eqref{dP-CSS}]
    Let $(A, \phi)$ be a configuration on $[0, T] \times \R^2 \times [0, \underline s]$ solving the paradifferential formulation of the Chern-Simons-Schr\"odinger equation \eqref{dP-CSS}. The covariant Hessian satisfies the following covariant Schr\"odinger equation,
        \begin{equation}\label{eq:Ds-schrodinger}
            \begin{split}
                \bigl( i \bfD_t + \bfD^j \bfD_j \bigr) (s\bfD^{(2)}) \phi 
                    &= (s \bfD^{(2)}) \phi +  \phi \cdot \phi \cdot (s\bfD^{(2)}) \phi + \phi \cdot (s^\frac12\bfD) \phi \cdot (s^\frac12 \bfD) \phi \\
                        &\qquad + \phi \cdot (s^\frac12 \nabla) s^\frac12 F_{tx} + (s^\frac12 \bfD) \phi \cdot s^\frac12 F_{tx} \\
                        &\qquad + \phi \cdot (s^\frac12 \nabla)^{(2)} F_{12} + (s^\frac12 \bfD) \phi \cdot (s^\frac12 \nabla) F_{12} + (s \bfD^{(2)}) \phi \cdot F_{12}\\
                        &\qquad + (s^\frac12 \bfD)^{(2)} \CSS. 
            \end{split}
            \end{equation}
    Furthermore, we have the following covariant heat equations for
        \begin{itemize}
            \item the electromagnetic field,
                \begin{equation}\label{eq:Fmunu-heat}
                    (\partial_s - \Delta) F_{\mu\nu} 
                        = 0,
                \end{equation}

            \item the Ampere tension field,
            \begin{equation}\label{eq:Amp-heat}
                (\partial_s - \Delta) \Amp
                    =  \bfD \phi \cdot \bfD^{(2)} \phi + F_{12} \cdot \phi \cdot \bfD \phi,
            \end{equation}
           
            \item the Schr\"odinger tension field,
            \begin{equation}\label{eq:CSS-heat}
                \bigl( \bfD_s - \bfD^\ell \bfD_\ell \bigr) \CSS
                    =  \phi \cdot \bfD \phi \cdot \bfD \phi + \nabla F_{12} \cdot \bfD \phi + F_{tx} \cdot \bfD \phi.
            \end{equation}
        \end{itemize}
\end{proposition}

\begin{proof}[Covariant Hessian equation \eqref{Ds-schrodinger}]
    Applying $\bfD^{(2)}$ to the definition of the Schr\"odinger tension field, commuting covariant derivatives and applying the product rule \eqref{product} wherever appropriate, 
        \begin{align*}
            \bigl( i \bfD_t + \bfD^j \bfD_j \bigr) \bfD^{(2)} \phi 
                &= - \tfrac\lambda2(1 - |\phi|^2) \bfD^{(2)} \phi  + \phi \cdot \phi \cdot \bfD^{(2)} \phi + \phi \cdot \bfD \phi \cdot \bfD \phi\\
                    &\qquad  + \Bigl[ i \bfD_t + \bfD^j \bfD_j , \bfD^{(2)} \Bigr] \phi \\
                    &\qquad + \bfD^{(2)} \CSS.
        \end{align*}
    Using the commutator identities \eqref{commute}, \eqref{laplace-commute}, the commutator may be expressed as
        \begin{align*}
            \Bigl[ i \bfD_t + \bfD^j \bfD_j , \bfD^{(2)} \Bigr] \phi
                &= \Bigl[ i \bfD_t + \bfD^j \bfD_j , \bfD \Bigr] \bfD\phi + \bfD \Bigl[ i \bfD_t + \bfD^j \bfD_j , \bfD \Bigr] \phi  \\
                &= F_{t x} \cdot \bfD \phi + \nabla F_{t x} \cdot \phi + F_{12} \cdot \bfD^{(2)} \phi + \nabla F_{12} \cdot \bfD \phi + \nabla^{(2)} F_{12} \cdot \phi.
        \end{align*}
    Multiplying the equation by $s$ and distributing powers appropriately gives the desired equation. 
\end{proof}

\begin{proof}[Electromagnetic field equation \eqref{Fmunu-heat}]
    By the Bianchi identity \eqref{bianchi} and \eqref{d-covariant-heat2}, 
        \begin{align*}
            \partial_s F_{\mu\nu} 
                &= \partial_\mu F_{s\nu} - \partial_\nu F_{s\mu} \\
                &=\partial^\ell ( \partial_\mu F_{\ell \nu} - \partial_\nu F_{\ell \mu}) \\
                &= \Delta F_{\mu \nu}, 
        \end{align*}
    as desired. 
\end{proof}

\begin{proof}[Ampere tension field equation \eqref{Amp-heat}]
    Applying $\partial_s - \Delta$ to the definition of the Ampere tension field, 
        \begin{align*}
            (\partial_s - \Delta) \Amp_j 
                &= (\partial_s - \Delta) \bigl( F_{tj} + \epsilon_{jk} \Im(\overline \phi \bfD^k \phi) + \partial_j F_{12} \bigr)\\
                &= (\partial_s -\Delta) \epsilon_{jk} \Im(\overline \phi \bfD^k \phi),
        \end{align*}
    since the electromagnetic field satisfies the linear heat equation \eqref{Fmunu-heat}. It remains to see whether the terms generated by the charge current are acceptable. Applying $\bfD_s$ to the charge current, and using the covariant heat equation \eqref{d-covariant-heat} and the equation for $\bfD \phi$ \eqref{Dphi-heat}, we compute via the product rule \eqref{product}
        \begin{align*}
            \partial_s \Im(\overline \phi \bfD_j \phi)
                &= \Im(\overline{\bfD_s \phi} \bfD_j \phi) + \Im (\overline \phi \bfD_s \bfD_j \phi) \\
                &= \Im\bigl( \overline{\bfD^\ell \bfD_\ell \phi} \bfD_j \phi \bigr) + \Im\bigl(\overline \phi \bfD^\ell \bfD_\ell \bfD_j \phi \bigr) - 2 F\indices{_j^\ell} \Re(\overline \phi \bfD_\ell \phi) \\
                &= \Delta \Im(\overline \phi \bfD_j \phi) - 2 \Im\bigl(\overline{\bfD^\ell \phi} \bfD_\ell \bfD_j \phi \bigr) - 2 F\indices{_j^\ell} \Re(\overline \phi \bfD_\ell \phi).
        \end{align*}
    These terms are all acceptable in view of the right-hand side of the equation \eqref{Amp-heat}.
\end{proof}

\begin{proof}[Schr\"odinger tension field equation \eqref{CSS-heat}]
    Commuting $\bfD_s - \bfD^\ell \bfD_\ell$ into the definition of the Schr\"odinger tension field, we see by virtue of \eqref{d-covariant-heat} that $\CSS$ obeys a covariant heat equation forced by the commutator between the covariant parabolic operator and covariant Schr\"odinger operator, 
        \begin{align*}
            \bigl( \bfD_s - \bfD^\ell \bfD_\ell \bigr) \CSS
                &= \Bigl[ \bfD_s - \bfD^\ell \bfD_\ell, i \bfD_t + \bfD^j \bfD_j + \tfrac\lambda2 (1 - |\phi|^2) \Bigr] \phi  .
        \end{align*}
    Separating out the commutator by (bi-)linearity, we compute
        \begin{align*}
             \bigl[ \bfD_s - \bfD^\ell \bfD_\ell, i \bfD_t \bigr] \phi
                &= \bigl( F_{st} - \partial^\ell F_{\ell t} + F\indices{_t^\ell} \bfD_\ell  \bigr) \phi \\
                &= 2 F\indices{_t^\ell} \, \bfD_\ell \phi,\\
            \bigl[ \bfD_s - \bfD^\ell \bfD_\ell , \bfD^j \bfD_j \bigr]\phi
                &= -i \bigl( \partial^j F_{sj} + 2 F\indices{_s^j} \bfD_j \bigr) \phi \\
                &= - 2i \partial^\ell F\indices{_\ell^j} \bfD_j \phi, \\
            \bigl[ \bfD_s - \bfD^\ell \bfD_\ell, \tfrac12 (1 - |\phi|^2) \bigr] \phi 
                &= \phi \, (\partial_s - \Delta) \tfrac12 (1 - |\phi|^2) - \partial^\ell \tfrac12 (1 - |\phi|^2) \, \bfD_\ell \phi \\
                &= \phi \overline{\bfD^\ell \phi} \bfD_\ell \phi + \Re(\overline \phi \bfD^\ell \phi) \bfD_\ell \phi,
        \end{align*}
    using the parabolic commutator identity \eqref{parabolic-commute}, the Laplacian commutator identity \eqref{laplace-commute} and \eqref{d-covariant-heat2} in the first two calculations, and recalling the equation for the charge density \eqref{charge-heat} in the last calculation. These terms are all acceptable in view of the right-hand side of the equation \eqref{CSS-heat}.
\end{proof}

\section{Parabolic smoothing of electromagnetic and tension fields}\label{sec:tension}

In this section, we analyse the parabolic equations of motion arising from the paradifferential formulation of the Chern-Simons-Schr\"odinger equation. Before delving into the estimates, it is instructive to perform a dimensional analysis in order to get a handle on the expected sizes of various dynamical variables. Starting with the tension fields, applying Duhamel's formula to the equations \eqref{Gauss-heat-prelim}, \eqref{Amp-heat}, \eqref{CSS-heat} yield
    \begin{align*}
        \Gauss 
            &\approx (\partial_s - \Delta)^{-1} (\bfD \phi \cdot \bfD \phi) + \text{linear homogeneous flow},\\ 
        \Amp
            &\approx (\partial_s - \Delta)^{-1} (\bfD \phi \cdot \bfD^{(2)} \phi) + \text{lower-order terms},\\
        \CSS 
            &\approx (\partial_s - \Delta)^{-1} (\phi \cdot \bfD \phi \cdot \bfD \phi + \bfD \phi \cdot \Amp + \bfD \phi \cdot \nabla \Gauss).
    \end{align*}
Recall that parabolic scaling dictates $\bfD, \nabla \approx s^{-1/2}$ and $(\partial_s - \Delta)^{-1} \approx s$ and $\| - \|_{L^p_x} \approx s^{1/p}$. Consequently, $\| \bfD \phi \|_{L^2_x} \approx \| \phi \|_{L^\infty_x}$, so, working at energy regularity, we regard the amplitude of the scalar field as size $O(1)$. Following these rules, we arrive at the heuristic bounds, 
    \begin{align*}
         \bigl\| \Gauss \bigr\|_{L^\infty_x} 
            &\lessapprox s^0,\\
        \bigl\| \Amp \bigr\|_{L^\infty_x}  
            &\lessapprox s^{-\frac12}, \\
         \bigl\| \CSS \bigr\|_{L^\infty_x} 
            &\lessapprox s^0.
    \end{align*}
Using the definitions of the tension fields, we can also predict the sizes of the electromagnetic field, 
    \begin{align*}
        \bigl\| F_{12} \bigr\|_{L^\infty_x} 
            &\lessapprox s^0, \\
        \bigl\| F_{tx} \bigr\|_{L^\infty_x} 
            &\lessapprox s^{-\frac12}. 
    \end{align*}
Our goal herein is to make these heuristics rigorous. 

We begin with some ``low-frequency'' bounds on the electromagnetic field. Throughout, it is useful to remember that $\underline s > 0$ was chosen such that 
    \[
        1 + \bigl\| \nabla \overline \Gauss \bigr\|_{L^2_x}^2 + \cE_{\mathrm{AH}} [\overline A, \overline{\vphantom{A}\phi}]
            \lesssim \underline s^{-1}. 
    \]

\begin{proposition}
    Let $(A, \phi)$ be a smooth solution on $[0, T] \times \R^2 \times [0, \underline s]$ to the paradifferential formulation of the Chern-Simons-Schr\"odinger equation \eqref{dP-CSS}. Then, for each non-negative integer $n \in \N_0$ and $2 \leq p \leq \infty$, we have the following:
    \begin{itemize}
        \item \textup{(Magnetic field uniform bound).}
            \begin{equation}\label{eq:F12-Linfty}
                \bigl\| (s^\frac12 \nabla)^{(n)} F_{12}  (s)\bigr\|_{L^\infty_{x}} 
                    \lesssim_n  \bigl\|  \phi (s)\bigr\|_{L^\infty_x}^2 + \bigl\| \Gauss \bigr\|_{L^\infty_x} + \underline s^{-1}.
            \end{equation}   

        \item \textup{(Magnetic field $L^p$-bound).}
            \begin{equation}\label{eq:nabla-F12-Lp}
                \bigl\| s^{\frac12 - \frac1p} (s^\frac12 \nabla)^{(n)} \nabla F_{12} \bigr\|_{L^p_x}
                    \lesssim_n \bigl\| \phi(s) \bigr\|_{L^\infty_{x}} \bigl\| \overline{\bfD \phi} \bigr\|_{L^2_x}  + \underline s^{-1}.
            \end{equation}

        \item \textup{(Electric field $L^p$-bound).}
            \begin{equation}\label{eq:Fmunu-Lp}
                 \bigl\| s^{\frac12 - \frac1p} (s^{\frac12} \nabla)^{(n)} F_{tx} (s) \bigr\|_{L^p_{x}}
                    \lesssim_n \bigl\| \phi(s) \bigr\|_{L^\infty_{x}} \bigl\| \overline{\bfD \phi} \bigr\|_{L^2_x}  + \underline s^{-1}. 
            \end{equation}

    \end{itemize}
\end{proposition}

\begin{proof}[Proof of magnetic field uniform bound \eqref{F12-Linfty}]
    It suffices to prove the bounds for $n = 0$ by parabolic smoothing \eqref{standard-heat} applied to the heat equation \eqref{Fmunu-heat} for the magnetic field. Recalling the definition of the Gauss tension field and its uniform bound \eqref{Gauss-Linfty},
        \begin{align*}
            \bigl\| F_{12} \bigr\|_{L^\infty_x}
                &\lesssim 1 + \bigl\| \phi \bigr\|_{L^\infty_x}^2 + \bigl\| \Gauss \bigr\|_{L^\infty_x} \\
                &\lesssim 1 + \bigl\|  \phi \bigr\|_{L^\infty_x}^2 + \bigl\| \Gauss \bigr\|_{L^\infty_x} + \bigl\| \overline{\bfD \phi}\bigr\|_{L^2_x}^2,
        \end{align*}
    which is of acceptable form. 
\end{proof}    

\begin{proof}[Proof of magnetic field $L^p$-bound \eqref{nabla-F12-Lp}]
    By parabolic smoothing \eqref{standard-heat} for the heat equation \eqref{Fmunu-heat}, it suffices to prove the bounds for $n = 0$ and $p = 2$. Recalling the definition of the Gauss tension field,
        \begin{align*}
            \bigl\| \nabla F_{12} \bigr\|_{L^2_x}
                &\lesssim \bigl\| \phi \bigr\|_{L^\infty_x} \bigl\| \bfD \phi \bigr\|_{L^2_x} + \bigl\| \nabla \Gauss \bigr\|_{L^2_x} \\
                &\lesssim \bigl\| \phi \bigr\|_{L^\infty_x} \bigl\| \overline{\bfD \phi} \bigr\|_{L^2_x} + \bigl\| \nabla \overline \Gauss \bigr\|_{L^2_x} + \bigl\| \overline{\bfD \phi} \bigr\|_{L^2_x}^2,
        \end{align*}
    by parabolic smoothing \eqref{Dphi-sobolev} and the Gauss tension field smoothing \eqref{Gauss-sobolev}. These terms are acceptable. 
\end{proof}

\begin{proof}[Proof of electric field $L^p$-bound \eqref{Fmunu-Lp}]
    By parabolic smoothing \eqref{standard-heat} for the heat equation \eqref{Fmunu-heat}, it suffices to prove the bounds for $n = 0$ and $p = 2$. Recalling the definition of the Amp\'ere tension field,
        \begin{align*}
            \bigl\| s^\frac12 F_{tx} \bigr\|_{L^\infty_x}
                &\lesssim \bigl\| F_{tx} \bigr\|_{L^2_x} \\
                &\lesssim \bigl\| \Amp \bigr\|_{L^2_x} + \bigl\| \phi \bigr\|_{L^\infty_x} \bigl\| \bfD \phi \bigr\|_{L^2_x} + \bigl\| \nabla F_{12} \bigr\|_{L^2_x} \\
                &\lesssim \bigl\| \Amp \bigr\|_{L^2_x} + \bigl\| \phi \bigr\|_{L^\infty_{x}} \bigl\| \overline{\bfD \phi} \bigr\|_{L^2_x}  + \bigl\| \nabla \overline \Gauss \bigr\|_{L^2_x} + \bigl\| \overline{\bfD \phi} \bigr\|_{L^2_x}^2
        \end{align*}
    by parabolic smoothing \eqref{Dphi-sobolev} and the previous estimate \eqref{nabla-F12-Lp}. The last two three terms are acceptable, so it remains to estimate the Amp\'ere tension field. We defer this to the subsequent Proposition \ref{prop:tension-env}, namely the envelope bound \eqref{amp-env}. Choosing the frequency envelope as in, say, Proposition \ref{prop:env-stable}, would complete the proof by the slowly-varying property and finite-energy.  
\end{proof}

\begin{proposition}\label{prop:tension-env}
    Let $(A, \phi)$ be a smooth solution on $[0, T] \times \R^2 \times [0, \underline s]$ to the paradifferential formulation of the Chern-Simons-Schr\"odinger equation \eqref{dP-CSS}. Fix $t \in [0, T]$, and suppose that $\ttc(s) \in L^2_{ds/s} ([0, \underline s])$ is a $\frE$-frequency envelope for $(A (t), \phi(t))$. Then for each non-negative integer $n \in \N_0$, we have the following:
    \begin{itemize}
        \item \textup{(Amp\'ere tension field bound).}
            \begin{equation}\label{eq:amp-env}
                \bigl\| \Amp (t, s) \bigr\|_{L^2_x}
                    \lesssim \bigl\| \overline{\bfD \phi} \bigr\|_{L^2_x} \cdot \ttc(s).
            \end{equation}
        
        \item \textup{(Electric field bound).}
            \begin{equation}\label{eq:Ftx-env}
                \bigl\| (s^\frac12 \nabla)^{(n + 1)} F_{tx} (t, s) \bigr\|_{L^2_x}
                    \lesssim_n \Bigl( 1 + \bigl\| \phi(t, s) \bigr\|_{L^\infty_x} + \bigl\| \overline{\bfD \phi} \bigr\|_{L^2_x} \Bigr) \cdot \ttc(s). 
            \end{equation}

         \item \textup{(Schr\"odinger tension field bound).}
            \begin{equation}\label{eq:css-env}
                \bigl\| s^{-\frac12} (s^\frac12 \bfD)^{(n)} \CSS (t, s) \bigr\|_{L^2_x}
                    \lesssim_n \Bigl( \| \phi(t, s) \|_{L^\infty_x}  \bigl\| \overline{\bfD \phi} \bigr\|_{L^2_x}  + \underline s^{-1} \Bigr) \cdot \ttc(s).
            \end{equation}
    
    \end{itemize}
\end{proposition}

\begin{proof}[Proof of Amp\'ere tension field bound \eqref{amp-env}]
    Using the equation \eqref{Amp-heat}, we decompose the Amp\'ere tension field between a higher-order term and lower-order term, 
            \begin{align*}
                \Amp 
                    := \HOT + \LOT,
            \end{align*}
        where 
            \begin{align}
                \HOT
                    &:=  (\partial_s - \Delta)^{-1} (\bfD \phi \cdot \bfD^{(2)} \phi), \label{eq:amp-heat-1}\\ 
                \LOT
                    &:= (\partial_s - \Delta)^{-1} ( \phi \cdot F_{12} \cdot \bfD \phi) \label{eq:amp-heat-2}.   
            \end{align}
        By the triangle inequality, it suffices to prove \eqref{amp-env} for these two terms. 
        
        We use Duhamel's formula, $(L^1 \to L^2)$-parabolic smoothing \eqref{standard-heat}, and H\"older, 
            \begin{align*}
                \bigl\| \HOT \bigr\|_{L^2_x}
                    &\lesssim \int_0^s \bigl\| e^{(s - s') \Delta}  (\bfD \phi \cdot \bfD^{(2)} \phi)\bigr\|_{L^2_x} ds' \\
                    &\lesssim \int_0^s \Bigl( \frac{s'}{s - s'}\Bigr)^\frac12 \bigl\| |s'|^\frac12 \cdot \bfD \phi \cdot \bfD^{(2)} \phi \bigr\|_{L^1_x} \frac{ds'}{s'} \\
                    &\lesssim \int_0^s \Bigl( \frac{s'}{s - s'}\Bigr)^\frac12 \cdot \bigl\| \overline{\bfD \phi} \bigr\|_{L^2_x} \cdot \ttc(s') \frac{ds'}{s'} \\
                    &\lesssim \int_0^s \Bigl( \frac{s'}{s - s'}\Bigr)^\frac12 \Bigl( \frac{s}{\, s'\,} \Bigr)^{\frac1{1000}} \cdot \bigl\| \overline{\bfD \phi} \bigr\|_{L^2_x} \cdot \ttc(s) \frac{ds'}{s'} \lesssim \bigl\| \overline{\bfD \phi} \bigr\|_{L^2_x} \cdot \ttc(s)
            \end{align*}
        using parabolic smoothing \eqref{Dphi-sobolev} and the slowly-varying property. Similarly, 
            \begin{align*}
                \bigl\|\LOT \bigr\|_{L^2_x} 
                    &\lesssim \int_0^s \bigl\| e^{(s - s') \Delta} (\phi \cdot F_{12} \cdot \bfD \phi) \bigr\|_{L^2_x} \, ds' \\
                    &\lesssim \int_0^s \Bigl( \frac{s'}{s - s'} \Bigr)^\frac12 \bigl\| |s'|^\frac12 \cdot \phi \cdot F_{12} \cdot \bfD \phi \bigr\|_{L^1_x} \frac{ds'}{s'} \\
                    &\lesssim \int_0^s \Bigl( \frac{s'}{s - s'} \Bigr)^\frac12 |s'|^\frac13 \cdot \bigl\| |s'|^\frac16 \phi \bigr\|_{L^\infty_x} \cdot \bigl\| \overline{\bfD \phi} \bigr\|_{L^2_x} \cdot \ttc(\underline s) \frac{ds'}{s'} \\
                    &\lesssim  \int_0^s \Bigl( \frac{s'}{s - s'} \Bigr)^\frac12 \Bigl( \frac{s'}{\,\underline s\,} \Bigr)^\frac13 \cdot \bigl\| \overline{\bfD \phi} \bigr\|_{L^2_x} \cdot \Bigl( \frac{\,\underline s\,}{s} \Bigr)^{\frac1{1000}} \cdot \ttc(s) \frac{ds'}{s'} \lesssim \bigl\| \overline{\bfD \phi} \bigr\|_{L^2_x} \cdot \ttc(s).
            \end{align*}
        using parabolic smoothing \eqref{inefficient-smooth}, \eqref{F12-sobolev}, \eqref{Dphi-sobolev}, and the slowly-varying property. 
\end{proof}

\begin{proof}[Proof of electric field bound \eqref{Ftx-env}]
    By parabolic smoothing \eqref{standard-heat} for the linear heat equation \eqref{Fmunu-heat}, it suffices to show the bound for $n = 0$. We write
        \begin{align*}
            \nabla F_{tx} (s)
                &=  e^{\frac12 s \Delta} \nabla F_{tx} (\tfrac12 s) \\
                &= e^{\frac12 s \Delta} \nabla \Bigl( \Amp + \phi \cdot \bfD \phi + \nabla F_{12} \Bigr)(\tfrac12 s) \\
                &= \nabla e^{\frac12 s \Delta} \Amp (\tfrac12 s) + e^{\frac12 s \Delta} \Bigl( \bfD \phi \cdot \bfD \phi + \phi \cdot \bfD^{(2)} \phi + \nabla^{(2)} F_{12} \Bigr) (\tfrac12 s),
        \end{align*}
    using the fact that the electric field satisfies the linear heat equation \eqref{Fmunu-heat} in the first line, the definition of the Amp\'ere tension field in the second line, and the product rule in the last line. We individually estimate each term in the last line. For the first term, 
        \[
            \bigl\| (s^\frac12 \nabla) e^{\frac12 s \Delta} \Amp (\tfrac12 s) \bigr\|_{L^2_x}
                \lesssim \bigl\| \overline{\bfD \phi} \bigr\|_{L^2_x} \cdot \ttc(\tfrac12 s). 
        \]
    by parabolic smoothing \eqref{standard-heat} and the Amp\'ere tension field bound \eqref{amp-env}. For the second term,
        \begin{align*}
            \bigl\| s^\frac12 e^{\frac12 s \Delta} (\bfD \phi \cdot \bfD \phi)(\tfrac12 s) \bigr\|_{L^2_x}
                &\lesssim \bigl\| (s^\frac12\nabla) (\bfD \phi \cdot \bfD \phi) (\tfrac12 s)\bigr\|_{L^1_x} \\
                &\lesssim \bigl\| \overline{\bfD \phi} \bigr\|_{L^2_x} \cdot \ttc(\tfrac12 s),
        \end{align*}
    by Sobolev embedding, Cauchy-Schwarz, and parabolic smoothing \eqref{Dphi-sobolev}. For the third and last terms,  
        \begin{align*}
            \bigl\| s^\frac12 e^{\frac12 s \Delta} (\phi \cdot \bfD^{(2)} \phi)(\tfrac12 s) \bigr\|_{L^2_x} + \bigl\| s^\frac12 e^{\frac12 s \Delta} (\nabla^{(2)} F_{12})(\tfrac12 s) \bigr\|_{L^2_x}
                &\lesssim \Bigl( 1 + \bigl\| \phi \bigr\|_{L^\infty_x} + \bigl\| \overline{\bfD \phi} \bigr\|_{L^2_x} \Bigr)\cdot \ttc(\tfrac12 s). 
        \end{align*}
    using the magnetic field bound \eqref{F12-env}. Collecting the previous calculations and appealing to the slowly-varying property, $\ttc(s) \sim \ttc(\tfrac12 s)$, we conclude the desired bound. 
\end{proof}

\begin{proof}[Proof of Schr\"odinger tension field bound \eqref{css-env}]
    We claim that the right-hand side of the covariant heat equation for the Schr\"odinger tension field obeys the bound
        \begin{align}\label{eq:css-rhs}
            \bigl\| s^\frac12 (s^\frac12 \bfD)^{(\leq n)} (\text{R.H.S.\eqref{CSS-heat}}) \bigr\|_{L^2_x}
                &\lesssim \bigl( \| \phi \|_{L^\infty_x} \bigl\| \overline{\bfD \phi} \bigr\|_{L^2_x}   + \underline s^{-1} \bigr) \cdot \ttc(s).
        \end{align}
    Assuming the claim, we can prove the desired bound by applying parabolic smoothing. Indeed, in the case $n = 0$, the maximum principle \eqref{maximum-principle} and Duhamel's formula applied to \eqref{CSS-heat} yield
        \begin{align*}
            \| s^{-\frac12} \CSS (s) \|_{L^2_x}
                &\lesssim \int_0^s \Bigl( \frac{s'}{s} \Bigr)^\frac12 \bigl\| {s'}^\frac12 \cdot \bigl(\text{R.H.S.\eqref{CSS-heat}}\bigr) (s') \bigr\|_{L^2_x} \, \frac{ds'}{s'} \\
                &\lesssim \int_0^s \Bigl( \frac{s'}{s} \Bigr)^\frac12 \cdot \Bigl( \| \phi(s') \|_{L^\infty_x} \bigl\| \overline{\bfD \phi} \bigr\|_{L^2_x}  + \underline s^{-1} \Bigr)  \cdot \ttc(s') \, \frac{ds'}{s'} \\
                &\lesssim \int_0^s  \Bigl( \frac{s'}{s} \Bigr)^{\frac12 - \frac1{1000}} \cdot \Bigl( \| \phi(s) \|_{L^\infty_x} \bigl\| \overline{\bfD \phi} \bigr\|_{L^2_x}  +  \log^\frac12 \Bigl( \frac{s}{s'}\Bigr) \cdot \underline s^{-1} + \underline s^{-1} \Bigr) \cdot \ttc(s) \, \frac{ds'}{s'}\\
                &\lesssim \Bigl(  \| \phi (s) \|_{L^\infty_x} \bigl\| \overline{\bfD \phi} \bigr\|_{L^2_x} + \underline s^{-1}  \Bigr) \cdot \ttc(s),
        \end{align*}
    using the Brezis-Gallou\"et-type inequality \eqref{brezis-0} and the slowly-varying property. For higher-derivatives, we apply the abstract dyadic smoothing \eqref{dyadic-smooth} to the equation \eqref{CSS-heat} for the Schr\"odinger tension field, 
        \begin{align*}
            \bigl\| s^{-\frac12} (s^\frac12 \bfD)^{(n)} \CSS (s) \bigr\|_{L^2_x} 
                &\lesssim \bigl\| s^{-\frac12} \CSS(s) \bigr\|_{L^2_x} + \bigl\| {s'}^\frac12 \cdot ({s'}^\frac12 \bfD)^{(\leq n)} \bigl(\text{R.H.S.\eqref{CSS-heat}}\bigr) (s')\bigr\|_{L^\infty_{s'} L^2_x ([\frac12 s, s])} \\
                &\lesssim \Bigl( \| \phi \|_{L^\infty_x} \bigl\| \overline{\bfD \phi} \bigr\|_{L^2_x}  + \underline s^{-1} \Bigr) \cdot \ttc(s)
        \end{align*}
    as desired. 

    It remains to establish the claim. We compute via the product rule
         \begin{align*}
            (s^\frac12 \bfD)^{(n)}\bigl(\text{R.H.S.\eqref{CSS-heat}}\bigr)
                &= \sum_{a + b + c = n} (s^\frac12 \bfD)^{(c)} \phi \cdot (s^\frac12 \bfD)^{(a)} \bfD \phi \cdot (s^\frac12 \bfD)^{(b)} \bfD \phi \\
                    &\qquad +  \sum_{a + b = n}  (s^\frac12 \nabla)^{(a)} \sfQ \cdot (s^\frac12 \bfD)^{(b)} \bfD \phi,
         \end{align*}
    abbreviating the electromagnetic field terms by $\sfQ \in \{ \nabla F_{12}, F_{tx}\}$. The best terms are the ``balanced'' ones, where derivatives are evenly distributed, 
         \[
            (s^\frac12 \bfD)^{(c)} \phi \cdot (s^\frac12 \bfD)^{(a)} \bfD \phi \cdot (s^\frac12 \bfD)^{(b)} \bfD \phi,
         \]
    for $a + b + c = n$ with $c \neq 0$. Placing the first factor in $L^\infty$ and estimating the remaining two factors in $L^2$ via Lemma \ref{lem:bilinear}, we arrive at 
         \begin{align*}
            \bigl\| s^\frac12 (-) \bigr\|_{L^2_x}
                &\lesssim \bigl\| (s^\frac12 \bfD)^{(c)} \phi \bigr\|_{L^\infty_x} \bigl\| (s^\frac12 \bfD)^{(a)} \bfD \phi \bigr\|_{L^2_x} \bigl\| s^{-\frac12} (s^\frac12 \bfD)^{(b)} (s \bfD^{(2)}) \phi \bigr\|_{L^2_x} \\
                &\lesssim \bigl\| \overline{\bfD \phi} \bigr\|_{L^2_x}^2 \cdot \ttc(s),
         \end{align*}
    using parabolic smoothing \eqref{Dphi-sobolev} for the first two factors and the envelope estimate for higher-derivatives of the Hessian \eqref{no-loss} for the last factor. The remaining terms are ``unbalanced'', i.e. we will always see at least one undifferentiated factor of $\phi$, 
         \[
            \phi \cdot  (s^\frac12 \bfD)^{(a)} \bfD \phi \cdot  (s^\frac12 \bfD)^{(b)} \bfD \phi,  \qquad (s^\frac12 \nabla)^{(a)} \sfQ \cdot  (s^\frac12 \bfD)^{(b)} \bfD \phi,
         \]
    for $a + b = n$. Estimating $\phi$ in $L^\infty$ and the remaining factors in $L^2$ via Lemma \ref{lem:bilinear}, 
         \begin{align*}
            \bigl\| s^\frac12 ( - ) \bigr\|_{L^2_x} 
                &\lesssim \Bigl( \|\phi \|_{L^\infty_x} \bigl\| (s^\frac12 \bfD)^{(a)} \bfD \phi \bigr\|_{L^2_x}  + \| (s^\frac12 \nabla)^{(a)} \sfQ \|_{L^2_x}\Bigr) \bigl\| s^{-\frac12} (s^\frac12 \bfD)^{(b)} (s \bfD^{(2)}) \phi \bigr\|_{L^2_x} \\
                    &\qquad + \bigl\| (s^\frac12 \bfD)^{(b)} \bfD \phi \bigr\|_{L^2_x} \bigl\| (s^\frac12 \nabla)^{(a + 1)} \sfQ  \bigr\|_{L^2_x} \\
                &\lesssim \Bigl( \| \phi \|_{L^\infty_x} \bigl\| \overline{\bfD \phi} \bigr\|_{L^2_x} + \underline s^{-1} \Bigr) \cdot \ttc(s),
         \end{align*}
    using the envelope bounds for the derivatives of the covariant Hessian \eqref{no-loss} and the electromagnetic field \eqref{F12-env}, \eqref{Ftx-env} for the ``high-frequency'' factors, and parabolic smoothing \eqref{Dphi-sobolev} and the electromagnetic field bounds \eqref{nabla-F12-Lp}, \eqref{Fmunu-Lp} for the ``low-frequency'' factors. This completes the proof of the claim and thereby the proposition.
\end{proof}

Our last order of business is to prove difference estimates for the electromagnetic field and tension fields with respect to the $\pzcE^0$-metric, which we abbreviate by  
    \begin{align*}
        \dist_{\pzcE^0}
            &\sim  \bigl\| \delta (s^\frac12 \bfD)^{(\leq 2027)} (s^\frac12 \bfD)\phi \bigr\|_{L^2_{\frac{ds}{s}} L^2_x} + \bigl\| \delta \overline \Gauss \bigr\|_{L^2_x} \\
                &\qquad + \bigl\| \delta \underline \bfD^{(\leq 100)} \underline{\bfD \phi} \bigr\|_{L^2_x} +  \bigl\| \delta \underline A \bigr\|_{L^4_x} + \bigl\| \nabla \delta \underline A \bigr\|_{H^{100}_x} + \bigl\| \delta \underline \phi \bigr\|_{L^\infty_x}.
    \end{align*}
In a word, interpolating via Lemma \ref{lem:low-high-inter} between the high-frequency decay alotted by Proposition \ref{prop:tension-env} and summability of low-frequencies, we can show $\log^{1/2}$-Lipschitz-type bounds. 

\begin{proposition}\label{prop:tension-diff}
    Let $(A, \phi)$, $(A', \phi')$ be smooth solutions on $[0, T] \times \R^2 \times [0, \underline s]$ to the paradifferential formulation of the Chern-Simons-Schr\"odinger equation \eqref{dP-CSS} in caloric-temporal gauge \eqref{CTC}. Then for any $0 < s_* < \tfrac1{1000} \underline s$, and each non-negative integer $n \in \N_0$, we have:
    \begin{itemize}
        \item \textup{(Amp\'ere tension field difference bound, higher-order terms).}
            \begin{equation}\label{eq:Amp-diff}
            \begin{split}
                \bigl\| s^\frac12 \delta \HOT \bigr\|_{L^\infty_s L^2_x} + \bigl\| s^\frac12 \delta \HOT \bigr\|_{L^2_{\frac{ds}s} L^2_x}
                    &\lesssim_{\underline s} \log^\frac12 \Bigl( \frac{\underline s}{s_*} \Bigr) \cdot  \Bigl( |s_*|^\frac12 + \dist_{\pzcE^0} \Bigr)  .
            \end{split}
            \end{equation}

        \item \textup{(Amp\'ere tension field difference bound, lower-order terms).}
            \begin{equation}\label{eq:Amp-diff-lot}
            \begin{split}
                \bigl\| s^\frac12 \delta \LOT \bigr\|_{L^\infty_s L^2_x} + \bigl\| s^\frac12 \delta \LOT \bigr\|_{L^2_{\frac{ds}s} L^2_x}
                    &\lesssim_{\underline s} \dist_{\pzcE^0}.
            \end{split}
            \end{equation}

        \item \textup{(Electric field $L^2_{ds/s} L^2_x$-difference bound).}
            \begin{equation}\label{eq:Ftx-diff}
                \bigl\| s^\frac12 (s^\frac12 \nabla)^{(n)} \delta F_{tx} \bigr\|_{L^2_{\frac{ds}{s}} L^2_x} 
                    \lesssim_{\underline s, n} \log^\frac12 \Bigl( \frac{\underline s}{s_*} \Bigr) \cdot\Bigl( |s_*|^\frac12 + \dist_{\pzcE^0} \Bigr) .
            \end{equation}

        \item \textup{(Electric field difference bound at $s = \underline s$)}
            \begin{equation}\label{eq:Ftx-diff-low}
                \bigl\| \nabla^{(n)} \delta \underline F_{tx} \bigr\|_{L^2_x}
                    \lesssim_{\underline s, n}\log^\frac12\Bigl( \frac{\underline s}{s_*} \Bigr) \Bigl( |s_*|^\frac12 + \dist_{\pzcE^0} \Bigr).
            \end{equation}

        \item \textup{(Schr\"odinger tension field difference bound).}
            \begin{equation}\label{eq:CSS-diff}
            \begin{split}
                \bigl\| \delta (s^\frac12 \bfD)^{(\leq 2027)} \CSS \bigr\|_{L^\infty_{s} L^2_x} &+ \bigl\|\delta (s^\frac12 \bfD) (s^\frac12 \bfD)^{(\leq 2027)}\CSS  \bigr\|_{L^2_{\frac{ds}{s}} L^2_x}\\
                    &\lesssim_{\underline s} \Bigl(1 + \bigl\| \overline A_x \bigr\|_{L^4_x}\Bigr) \cdot \log^\frac12 \Bigl(\frac{\underline s}{s_*} \Bigr) \cdot \Bigl( |s_*|^\frac12 + \dist_{\pzcE^0} \Bigr). 
            \end{split}
            \end{equation}
    \end{itemize}
\end{proposition}

\begin{proof}[Proof of Amp\'ere tension field difference bound, higher-order terms \eqref{Amp-diff}]
    Applying the weighted energy estimate \eqref{abs-heat-weight} with $\gamma = \tfrac12$ to the heat equation \eqref{amp-heat-1} for $\delta \HOT$, 
        \begin{align*}
            \bigl\| s^\frac12 \delta \HOT \bigr\|_{L^\infty_s L^2_x} 
                &\lesssim \bigl\| s^\frac12 \delta \HOT \bigr\|_{L^2_{\frac{ds}{s}} L^2_x} + \bigl\| s^{1 + \frac12} \delta (\bfD \phi \cdot \bfD^{(2)} \phi) \bigr\|_{L^1_{\frac{ds}{s}} L^2_x} \\
                &\lesssim \bigl\| s^\frac12 \delta \HOT \bigr\|_{L^2_{\frac{ds}{s}} L^2_x} + \bigl\| \overline{\bfD \phi}\bigr\|_{L^2_x} \bigl\| \delta (s^\frac12 \bfD)^{(\leq 1)} (s^\frac12 \bfD) \phi \bigr\|_{L^2_{\frac{ds}{s}} L^2_x} , 
        \end{align*}
    using Cauchy-Schwarz and parabolic smoothing \eqref{Dphi-sobolev}, \eqref{Dphi-strichartz-heat} to estimate the right-hand side. Thus, it remains to estimate the difference $s^\frac12 \delta \Amp$ in $L^2_{ds/s} L^2_x$. We argue by interpolating between decay of high-frequencies and summability of low-frequencies -- a l\'a Lemma \ref{lem:low-high-inter}, it suffices to show  
        \begin{align*}
            \bigl\| s^\frac12 \HOT (s)\bigr\|_{L^2_x}
                &\lesssim_{\underline s} s^\frac12 ,\\
            \bigl\| s^\frac12 \delta \HOT \bigr\|_{L^2_{\frac{ds}{s}} L^2_x ([s_*, \underline s])}
                &\lesssim_{\underline s}\log^\frac12 \Bigl( \frac{\underline s}{s_*} \Bigr) \cdot \Bigl( |s_*|^\frac12 + \dist_{\pzcE^0}\Bigr).
        \end{align*}
    The former is an immediate consequence of the envelope bound \eqref{amp-env}. For the latter, we write the difference using Duhamel's formula and divide the analysis between two cases, 
        \begin{align*}
            s^\frac12 \delta \HOT (s)
                &= s^\frac12 \Bigl( \int_0^{\frac12 s_*} + \int_{\frac12 s_*}^s \Bigr) e^{(s - s') \Delta} \delta (\bfD \phi \cdot \bfD^{(2)} \phi) (s') \,ds'\\
                &=: \mathrm{I} + \mathrm{II}. 
        \end{align*}
    
        The first term on the right is a $\text{high}\times\text{high} \to \text{low}$-type term -- more precisely, the regime $s' < \tfrac12 s_* < \tfrac12 s$. We leverage parabolic smoothing in $s$, writing $e^{(s - s') \Delta} = e^{\frac12 s \Delta} e^{(\frac12 s - s') \Delta}$, and decay of high-frequencies, 
        \begin{align*}
            \|\mathrm{I} \|_{L^2_x}
                &\lesssim \int_0^{\frac12 s_*} \bigl\| \bfD \phi \cdot \bfD^{(2)} \phi \bigr\|_{L^1_x} \, ds' \\
                &\lesssim \int_0^{\frac12 s_*} |s'|^\frac12 \bigl\| \bfD \phi \bigr\|_{L^2_x} \bigl\| ({s'}^\frac12 \bfD) \bfD \phi \bigr\|_{L^2_x} \frac{ds'}{s'} \\
                &\lesssim \bigl\| \overline{\bfD \phi} \bigr\|_{L^2_x}^2 \cdot |s_*|^\frac12, 
        \end{align*}
    using the smoothing estimate \eqref{standard-heat} in the first line, Cauchy-Schwarz in the second, and parabolic smoothing \eqref{Dphi-sobolev} in the third. Putting this in $L^2_{ds/s} ([s_*, \underline s])$, we see another factor of $\log^{1/2} (\underline s/s_*)$, which is acceptable. 
    
    The second term on the right is a $\text{medium}\times\text{medium} \to \text{low}$-type term, where both $s$ and $s'$ reside in the regime $[s_*, \underline s]$. Here we argue by almost-orthogonality, 
        \begin{align*}
            \| \mathrm{II} \|_{L^2_x} 
                &\lesssim \int_{\frac12 s_*}^s s^\frac12  \bigl\| e^{(s - s') \Delta} (\bfD \phi \cdot \delta \bfD^{(2)} \phi + \delta \bfD \phi \cdot \bfD^{(2)} \phi)\bigr\|_{L^2_x} ds' \, \\
                &\lesssim \int_{\frac12 s_*}^s \Bigl( \frac{s}{s - s'} \Bigr)^\frac12 \Bigl( \bigl\|  \bfD \phi \bigr\|_{L^2_x} \bigl\| \delta ({s'}^\frac12 \bfD)^{(2)} \phi \bigr\|_{L^2_x} + \bigl\| \delta ({s'}^\frac12 \bfD) \phi \bigr\|_{L^2_x} \bigl\| ({s'}^\frac12 \bfD) \bfD \phi \bigr\|_{L^2_x} \Bigr) \frac{ds'}{s'}\\
                &\lesssim \int_{\frac12 s_*}^s \Bigl( \frac{s}{s - s'} \Bigr)^\frac12 \bigl\| \overline{\bfD \phi} \bigr\|_{L^2_x} \bigl\| \delta ({s'}^\frac12 \bfD)^{(\leq 1)} (s^\frac12 \bfD) \phi \bigr\|_{L^2_x} \frac{ds'}{s'} ,
        \end{align*}
    using smoothing in $(s - s')$ and Cauchy-Schwarz in the second line, and parabolic smoothing \eqref{Dphi-sobolev} of $\bfD \phi$ in the third. Viewing this as an integral operator acting on functions on $[\frac12 s_*, \underline s]$, we estimate the kernel, 
        \begin{align*}
            \sup_{s \in [s_*, \underline s]} \int_{\frac12 s_*}^s \Bigl( \frac{s}{s - s'} \Bigr)^\frac12 \frac{ds'}{s'}
                &= \sup_{s \in [s_*, \underline s]} \int_{s_*/2s}^1 \Bigl( \frac{1}{1 - s"} \Bigr)^\frac12 \frac{ds"}{s"} \sim \log \Bigl( \frac{\underline s}{s_*} \Bigr),\\
            \sup_{s' \in [\frac12 s_*, \underline s]} \int_{s'}^{\underline s} \Bigl( \frac{s}{s - s'} \Bigr)^\frac12 \frac{ds}{s} 
                &=\sup_{s' \in [\frac12 s_*, \underline s]} \int_{1}^{\underline s/s'} \Bigl( \frac{1}{1 - s"}\Bigr)^\frac12 \frac{ds"}{s"} \lesssim 1.
         \end{align*}
    Thus, by Schur's test, 
         \begin{align*}
            \bigl\| s^\frac12 \delta \HOT \bigr\|_{L^2_{\frac{ds}{s}} L^2_x ([s_*, \underline s])} 
                &\lesssim \bigl\| \overline{\bfD \phi} \bigr\|_{L^2_x} \cdot \log^\frac12 \Bigl( \frac{\underline s}{s_*} \Bigr) \cdot \bigl\| \delta(s^\frac12 \bfD)^{(\leq 1)} (s^\frac12 \bfD) \phi \bigr\|_{L^2_{\frac{ds}{s}} L^2_x}.
         \end{align*}
    Collecting the previous inequalities, we conclude the desired bound. 
\end{proof}

\begin{proof}[Proof of Amp\'ere tension field difference bound, lower-order terms \eqref{Amp-diff-lot}]
    It suffices to show the stronger bound
        \[
            \bigl\| s^{\frac14} \delta \LOT \bigr\|_{L^\infty_s L^2_x}
                \lesssim_{\underline s} \dist_{\pzcE^0}.
        \]
Writing this term using Duhamel's formula for \eqref{amp-heat-2} and then applying parabolic smoothing \eqref{standard-heat}, 
        \begin{align*}
             \bigl\| \delta \LOT \bigr\|_{L^2_x}
                &\lesssim \int_0^s \bigl\| \delta (\phi \cdot F_{12} \cdot \bfD \phi) (s') \bigr\|_{L^2_x} \frac{ds'}{s'} \\
                &\lesssim \int_0^s \Bigl( \frac{s'}{s - s'} \Bigr)^\frac14 \bigl\| \delta (\phi \cdot {s'}^\frac14 F_{12} \cdot ({s'}^\frac12 \bfD) \phi) (s')\bigr\|_{L^{4/3}_x} \, \frac{ds'}{s'}.
        \end{align*}
We divide into three cases depending on where the difference operator falls, and then estimating the remaining factors by parabolic smoothing for the magnetic field \eqref{F12-sobolev}, $\bfD \phi$ \eqref{Dphi-sobolev}, and $\phi$. When the difference falls on $\phi$, we estimate it by \eqref{phi-diff}, 
        \eqref{phi-diff}
        \begin{align*}
            |s'|^\frac34 \cdot \bigl\| \delta \phi \cdot F_{12} \cdot \bfD \phi \bigr\|_{L^{4/3}}
                &\lesssim \bigl\| \delta \phi \bigr\|_{L^2 + L^\infty} \bigl\| {s'}^\frac14 F_{12} \bigr\|_{L^4}  \bigl\| ({s'}^\frac12 \bfD) \phi \bigr\|_{L^\infty \cap L^2}\\
                &\lesssim_{\underline s} \Bigl( 1 + \log^\frac12 \Bigl( \frac{\underline s}{s'} \Bigr) \Bigr)  \cdot \dist_{\pzcE^0}. 
        \end{align*}
When the difference falls on the magnetic field, we estimate it by \eqref{F12-diff-lot},
        \begin{align*}
            |s'|^\frac34 \cdot \bigl\| \phi \cdot  \delta F_{12} \cdot \bfD\phi \bigr\|_{L^{4/3}}
                &\lesssim |s'|^\frac14 \bigl\| {s'}^\frac18 \phi \bigr\|_{L^\infty} \bigl\| {s'}^\frac18 \delta F_{12} \bigr\|_{L^2}  \bigl\| {s'}^\frac14 \bfD \phi \bigr\|_{L^4}\\
                &\lesssim_{\underline s} |s'|^{\frac14} \dist_{\pzcE^0}. 
        \end{align*}
When the difference falls on $\bfD \phi$, we simply write 
         \begin{align*}
            |s'|^\frac34 \cdot \bigl\| \phi \cdot  F_{12} \cdot \delta \bfD \phi \bigr\|_{L^{4/3}}
                &\lesssim \bigl\| \phi \bigr\|_{L^\infty} \bigl\| {s'}^\frac14 F_{12} \bigr\|_{L^4}  \bigl\| \delta ({s'}^\frac12 \bfD) \phi \bigr\|_{L^2}\\
                &\lesssim_{\underline s} \Bigl( 1 + \log^\frac12 \Bigl( \frac{\underline s}{s'} \Bigr) \Bigr)  \cdot \bigl\| \delta ({s'}^\frac12 \bfD) \phi \bigr\|_{L^2},
        \end{align*}
and then take Cauchy-Schwarz in $ds'/s'$ to see the $\pzcE^0$-metric. Altogether, these calculations yield a bound of the form 
        \begin{align*}
            \bigl\| \delta \LOT \bigr\|_{L^2_x} 
                &\lesssim_{\underline s} \Bigl(1 + \log \Bigl( \frac{\underline s}{s} \Bigr)\Bigr) \cdot \dist_{\pzcE_0}.
        \end{align*}
Multiplying by $s^{1/4}$ handles the logarithmic divergence near $s = 0$. This completes the proof. 
\end{proof}

\begin{proof}[Proof of electric field $L^2_{ds/s} L^2_x$-difference bound \eqref{Ftx-diff}]
    By parabolic smoothing \eqref{standard-heat} for the linear heat equation \eqref{Fmunu-heat}, it suffices to consider $n = 0$. Interpolating via Lemma \ref{lem:low-high-inter}, the desired bound follows provided we show decay of high-frequencies $s \leq s_*$ and summability of low-frequencies $s \geq s_*$, 
        \begin{align*}
            \bigl\| s^\frac12 F_{tx} (s)\bigr\|_{L^2_x}
                &\lesssim_{\underline s} s^\frac12 \log^\frac12 \Bigl( \frac{\underline s}{\, s\,} \Bigr) ,\\
            \bigl\| s^\frac12 F_{tx} \bigr\|_{L^2_{\frac{ds}{s}} L^2_x([s_*, \underline s])}
                &\lesssim \log^\frac12 \Bigl( \frac{\underline s}{s_*} \Bigr) \cdot \Bigl( \dist_{\pzcE^0} +  |s_*|^\frac12  \Bigr). 
        \end{align*}
    The high frequency decay follows as an immediate consequence of the electric field bound \eqref{Fmunu-Lp} and the Brezis-Gallou\"et-type inequality \eqref{brezis}. For the latter, we recall the definition of the Amp\'ere tension field 
        \begin{align*}
            s^\frac12 \delta F_{tx} 
                = s^\frac12 \delta \Bigl( \Amp + \phi \cdot \bfD \phi + \nabla F_{12}\Bigr) .
        \end{align*}
    The Amp\'ere tension field difference was controlled in \eqref{Amp-diff}-\eqref{Amp-diff-lot}, while the magnetic field difference was controlled in \eqref{F12-diff}. For the contribution of $\delta \cdot (\phi \cdot \bfD \phi)$, we write 
        \begin{align*}
            \bigl\| s^\frac12 \delta (\phi \cdot \bfD \phi)\bigr\|_{L^2_{\frac{ds}{s}} L^2_x([s_*, \underline s])}
                &\lesssim \bigl\| \phi \bigr\|_{L^\infty_{s, x} ([s_*, \underline s])} \bigl\| \delta (s^\frac12 \bfD) \phi \bigr\|_{L^2_{\frac{ds}{s}} L^2_x} + \bigl\| \delta \phi \bigr\|_{L^\infty_s (L^2 + L^\infty)_x ([s_*, \underline s])} \bigl\| (s^\frac12 \bfD) \phi \bigr\|_{L^2_{\frac{ds}{s}} (L^\infty \cap L^2)_x} \\
                &\lesssim \log^{\frac12} \Bigl( \frac{\underline s}{s_*} \Bigr) \bigl\| \delta (s^\frac12 \bfD)^{(\leq 1)} (s^\frac12 \bfD) \phi \bigr\|_{L^2_{\frac{ds}{s}} L^2_x} + \bigl\| \delta \underline \phi \bigr\|_{L^\infty_x},
        \end{align*}
    by the Brezis-Gallou\"et-type bound \eqref{brezis}, its difference analogue \eqref{phi-diff}, and parabolic smoothing \eqref{Dphi-sobolev}, \eqref{Dphi-strichartz-heat}. Collecting the previous bounds completes the proof. 
\end{proof}

\begin{proof}[Proof of electric field difference bound at $s = \underline s$ \eqref{Ftx-diff-low}]
    By parabolic smoothing \eqref{standard-heat} for the linear heat equation \eqref{Fmunu-heat}, 
        \begin{align*}
            \bigl\| \nabla^{(n)} \delta \underline F_{tx} \bigr\|_{L^2_x} 
                &\lesssim_n \bigl\| \delta\underline F_{tx} (\tfrac12 \underline s)\bigr\|_{L^2_x}. 
        \end{align*}
    Expanding via the definition of the Amp\'ere tension field,
        \begin{align*}
             \delta F_{tx} (\tfrac12 \underline s)
                =  \delta \Bigl( \Amp + \phi \cdot \bfD \phi + \nabla F_{12}\Bigr) (\tfrac12 \underline s). 
        \end{align*}
    The difference of the Amp\'ere tension field was estimated in \eqref{Amp-diff}-\eqref{Amp-diff-lot}, while the difference of the magnetic field as estimated in \eqref{F12-diff}. When the difference falls on $\phi$, we use \eqref{phi-diff} and parabolic smoothing \eqref{Dphi-sobolev}. When the difference falls on $\bfD \phi$, we use \eqref{Dphi-s-diff-2} and estimate $\phi$ by \eqref{brezis-0}. Since we are working on the dyadic scale $[\underline s/2, \underline s]$, the ensuing extra negative powers of $s$ and logarithmic factors are harmless. 
\end{proof}

\begin{proof}[Proof of Schr\"odinger tension field difference bound \eqref{CSS-diff}]
    By the envelope bound for the Schr\"odinger tension field \eqref{css-env} and the Brezis-Gallou\"et-type bound \eqref{brezis},  
        \begin{align*}
            \bigl\| (s^\frac12 \bfD)^{(n)} \CSS (s) \bigr\|_{L^2_x}
                &\lesssim_{\underline s, \|} s^{\frac12} \cdot \Bigl( \| \phi(t, s) \|_{L^\infty_x} \bigl\| \overline{\bfD \phi} \bigr\|_{L^2_x}  + \underline s^{-1} \Bigr) \cdot \ttc(s) \\
                &\lesssim_{\underline s} 1 .
        \end{align*}
    Thus, applying the abstract difference bound from Lemma \ref{lem:diff-smooth} to the equation for the Schr\"odinger tension field \eqref{CSS-heat}, 
        \begin{align*}
            \bigl\| \delta (s^\frac12 \bfD)^{(\leq 2027)} \CSS \bigr\|_{L^\infty_s L^2_x} &+ \bigl\| \delta (s^\frac12 \bfD)^{(\leq 2027)} (s^\frac12 \bfD) \CSS \bigr\|_{L^2_{\frac{ds}{s}} L^2_x} \\
                &\lesssim_{\underline s} \bigl\| \delta \overline A_x \bigr\|_{L^4_x} + \bigl\| \nabla \delta  \overline A_x  \bigr\|_{L^2_x}  + \bigl\| \overline A_x \bigr\|_{L^4_x} \bigl\| \delta \overline A_x \bigr\|_{L^4_x} \\
                    &\qquad+ \bigl\|s \cdot \delta (s^\frac12 \bfD)^{(\leq 2027)} (\text{R.H.S.\eqref{CSS-heat}}) \bigr\|_{L^1_{\frac{ds}{s}} L^2_x}\\
                &\lesssim_{\underline s} \Bigl(1 + \bigl\| \overline A_x \bigr\|_{L^4_x} \Bigr) \cdot \log^\frac12\Bigl( \frac{\underline s}{s_*} \Bigr) \cdot \Bigl( |s_*|^\frac12 + \dist_{\pzcE^0} \Bigr) \\
                    &\qquad+ \bigl\|s \cdot \delta (s^\frac12 \bfD)^{(\leq 2027)} (\text{R.H.S.\eqref{CSS-heat}}) \bigr\|_{L^1_{\frac{ds}{s}} L^2_x},
        \end{align*}
    estimating the difference of the magnetic potential via \eqref{F12-diff}\footnote{Actually, if one inspects the proof of Lemma \ref{lem:diff-smooth}, we can put more powers of $s$ on the magnetic field and estimate it using the lower-order estimate \eqref{F12-diff-lot}. }and \eqref{Ax-diff}. Thus, it remains to show
        \begin{align*}
            \bigl\|s \cdot \delta (s^\frac12 \bfD)^{(\leq 2027)} (\text{R.H.S.\eqref{CSS-heat}}) \bigr\|_{L^1_{\frac{ds}{s}} L^2_x}
                &\lesssim_{\underline s} \log^\frac12 \Bigl( \frac{\underline s}{s_*} \Bigr) \cdot \Bigl( |s_*|^\frac12 +  \dist_{\pzcE^0} \Bigr).
        \end{align*}
    To prove this estimate, we interpolate between high and low frequencies using Lemma \ref{lem:low-high-inter}, reducing the task to showing 
       \begin{align*}
            \bigl\|s \, (s^\frac12 \bfD)^{(\leq 2027)} (\text{R.H.S.\eqref{CSS-heat}})\bigr\|_{L^2_x}
                &\lesssim_{\underline s} s^\frac12 \log^\frac12 \Bigl( \frac{\underline s}{s} \Bigr),\\
            \bigl\|s \, \delta (s^\frac12 \bfD)^{(\leq 2027)} (\text{R.H.S.\eqref{CSS-heat}}) \bigr\|_{L^1_{\frac{ds}{s}} L^2_x ([s_*, \underline s])}
                &\lesssim_{\underline s} \log^\frac12 \Bigl( \frac{\underline s}{s_*} \Bigr) \cdot \Bigl( |s_*|^\frac12 +  \dist_{\pzcE^0} \Bigr).
       \end{align*}

    For the high-frequency bound, recall that, during the proof of the envelope bound \eqref{css-env}, we showed that 
       \begin{align*}
            \bigl\|s \, (s^\frac12 \bfD)^{(\leq 2027)} (\text{R.H.S.\eqref{CSS-heat}})\bigr\|_{L^2_x}
                &\lesssim_{\underline s} s^\frac12 \Bigl( 1 +  \| \phi(s) \|_{L^\infty_x} \Bigr).
       \end{align*}
    Thus the high-frequency bound for $s \leq s_*$ follows from the Brezis-Gallou\"et-type inequality \eqref{brezis}. 
    
    For the low-frequency bound, we compute 
       \begin{align*}
            s \cdot \delta (s^\frac12 \bfD)^{(n)} (\text{R.H.S.\eqref{CSS-heat}})
                &= \sum_{a + b = n} \delta \phi \cdot (s^\frac12 \bfD)^{(a + 1)} \phi \cdot (s^\frac12 \bfD)^{(b + 1)} \phi \\
                &\qquad +  \sum_{a + b + c = n} (s^\frac12 \bfD)^{(c)} \phi \cdot (s^\frac12 \bfD)^{(a + 1)} \phi \cdot \delta (s^\frac12 \bfD)^{(b + 1)} \phi \\
                    &\qquad + \sum_{a + b = n} (s^\frac12 \nabla)^{(a)} \delta s^\frac12 \sfQ \cdot (s^\frac12 \bfD)^{(b + 1)} \phi  \\
                    &\qquad + \sum_{a + b = n} (s^\frac12 \nabla)^{(a)}  s^\frac12 \sfQ \cdot \delta (s^\frac12 \bfD)^{(b + 1)} \phi \\
                &=: \mathrm{I}_n + \mathrm{II}_n + \mathrm{III}_n + \mathrm{IV}_n,
       \end{align*}
    where $\sfQ \in \{ \nabla F_{12}, F_{tx}\}$. For the first term, we estimate the difference of the scalar field in $L^2 + L^\infty$ via \eqref{phi-diff}, and the remaining terms in $L^2_{ds/s} (L^\infty \cap L^2)_x$ and $L^2_{ds/s} L^\infty_x$ via parabolic smoothing \eqref{Dphi-sobolev}, \eqref{Dphi-strichartz-heat}, 
       \begin{align*}
            \bigl\| \mathrm{I}_n \bigr\|_{L^1_{\frac{ds}{s}} L^2_x ([s_*, \underline s])}
                &\lesssim \sum_{a + b = n}\bigl\| \delta \phi \bigr\|_{L^\infty_s (L^2 + L^\infty)_x ([s_*, \underline s])} \bigl\| (s^\frac12 \bfD)^{(a + 1)} \phi \bigr\|_{L^2_{\frac{ds}{s}} (L^\infty \cap L^2)_x} \bigl\| (s^\frac12 \bfD)^{(b + 1)} \phi \bigr\|_{L^2_{\frac{ds}{s}} L^\infty_x} \\
                &\lesssim_{\underline s} \log^\frac12 \Bigl( \frac{\underline s}{s_*} \Bigr) \cdot \Bigl( \bigl\| \delta \underline \phi \bigr\|_{L^\infty_x} + \bigl\| \delta (s^\frac12 \bfD)^{(2)} \phi \bigr\|_{L^2_{\frac{ds}{s}} L^2_x} \Bigr) ,
       \end{align*}
    which is acceptable. For the second term, we make analogous moves,
       \begin{align*}
            \bigl\| \mathrm{II}_n \bigr\|_{L^1_{\frac{ds}{s}} L^2_x ([s_*, \underline s])}
                &\lesssim \sum_{a + b + c = n} \bigl\| (s^\frac12 \bfD)^{(a)}\phi \bigr\|_{L^\infty_s (L^2 + L^\infty)_x ([s_*, \underline s])} \bigl\| (s^\frac12 \bfD)^{(a + 1)} \phi \bigr\|_{L^2_{\frac{ds}{s}} (L^\infty \cap L^2)_x} \bigl\| \delta (s^\frac12 \bfD)^{(b + 1)} \phi \bigr\|_{L^2_{\frac{ds}{s}} L^2_x}  \\
                &\lesssim_{\underline s} \log^\frac12 \Bigl(\frac{\underline s}{s_*} \Bigr) \cdot \bigl\| \delta (s^\frac12 \bfD)^{(\leq n)} (s^\frac12 \bfD) \phi \bigr\|_{L^2_{\frac{ds}{s}} L^2_x}
       \end{align*}
    where we bound $a = 0$ via the Brezis-Gallou\"et-type inequality \eqref{brezis}. For the third term, we estimate the electromagnetic field differences via \eqref{F12-diff} and \eqref{Ftx-diff} and $\bfD \phi$ via parabolic smoothing \eqref{Dphi-sobolev}, \eqref{Dphi-strichartz-heat} 
       \begin{align*}
            \bigl\| \mathrm{III}_n \bigr\|_{L^1_{\frac{ds}{s}} L^2_x} 
                &\lesssim \sum_{a + b = n} \bigl\| (s^\frac12 \nabla)^{(a)} \delta s^\frac12 \sfQ \bigr\|_{L^2_{\frac{ds}{s}} L^2_x} \bigl\| (s^\frac12 \bfD)^{(b + 1)} \phi \bigr\|_{L^2_{\frac{ds}{s}} L^\infty_x}  \\
                &\lesssim_{\underline s} \log^\frac12 \Bigl( \frac{\underline s}{s_*} \Bigr) \cdot \Bigl( |s_*|^\frac12 + \dist_{\pzcE^0} \Bigr) .
       \end{align*}
    For the fourth term, we estimate the electromagnetic field via \eqref{nabla-F12-Lp} and \eqref{Fmunu-Lp}, and the ensuing contribution of $\phi(s)$ via the Brezis-Gallou\"et-type inequality \eqref{brezis}, 
       \begin{align*}
            \bigl\| \mathrm{IV}_n \bigr\|_{L^1_{\frac{ds}{s}} L^2_x ([s_*, \underline s])} 
                &\lesssim \sum_{a + b = n} \bigl\| (s^\frac12 \nabla)^{(a)}  s^\frac12 \sfQ \bigr\|_{L^2_{\frac{ds}{s}} L^\infty_x ([s_*, \underline s])} \bigl\| \delta (s^\frac12 \bfD)^{(b + 1)} \phi \bigr\|_{L^2_{\frac{ds}{s}} L^2_x} \\
                &\lesssim_{\underline s} \log^\frac12 \Bigl( \frac{\underline s}{s_*} \Bigr)\cdot  \bigl\| \delta (s^\frac12 \bfD)^{(\leq n)} (s^\frac12 \bfD)\phi \bigr\|_{L^2_{\frac{ds}{s}} L^2_x} .
       \end{align*}
    Taking $n \leq 2027$, these terms are all acceptable.
\end{proof}

\begin{numremark}\label{rem:hihilo}
    The Amp\'ere tension field difference bound \eqref{Amp-diff} is analogous to the electric potential difference bound \cite[Lemma 5.3]{Lim2018} in the setting of the Jackiw-Pi Chern-Simons-Schr\"odinger equation \eqref{JP-CSS} in Coulomb gauge. Unfortunately, in our setting, there seems to be an extra logarithmic divergence arising from the term 
        \[
            (\partial_s - \Delta)^{-1} (\bfD \phi \cdot \delta \bfD^{(2)} \phi).
        \]
    One can view this as a $\text{high}\times \text{high} \to \text{low}$ paraproduct. If we try to estimate this in a very weak norm, e.g. $L^1$, we obtain by Cauchy-Schwarz
        \[
            \bigl\| (\partial_s - \Delta)^{-1} (\bfD \phi \cdot \delta \bfD^{(2)} \phi) \bigr\|_{L^1_x}
                \leq \int_0^s \bigl\| \bfD \phi \bigr\|_{L^2_x} \bigl\| \delta ({s'}^\frac12 \bfD)^{(2)} \phi \bigr\|_{L^2_x} \frac{ds'}{s'}
        \]
    The difference is $L^2_{ds/s}$-summable, however the factor of $\bfD \phi$ is merely $L^\infty_s$. Furthermore, there does not seem to be any room to leverage Strichartz estimates, as all the powers of $s'$ have been ``used up".
\end{numremark}

\section{The caloric-temporal-Coulomb gauge}\label{sec:gauge}

As alluded to in the introduction, we conduct the bulk of our analysis of the Chern-Simons-Schr\"odinger equation in an auxiliary choice of gauge. Continuing on the paradifferential analogy, it is prudent to choose a gauge which renders the electromagnetic potential as ``low-frequency'' so as to give the left-hand side of the equation for the covariant Hessian \eqref{Ds-schrodinger} a ``low-high paraproduct'' structure. We begin with imposing the caloric gauge $A_s = 0$, which we used in Section \ref{subsec:LP-topology} to give a Littlewood-Paley-type characterisation the topology on the space of finite-energy configurations. This fixes the gauge-freedom up to transformations independent of the heat-time $s$, and renders the dynamic Yang-Mills heat flow \eqref{d-covariant-heat2} into ODEs in $s$. Since parabolic smoothing, as established in Sections \ref{subsec:LP-sobolev} and \ref{sec:tension}, implies that the geometric quantities are smooth at $s = \underline s$, we leverage the regularity by imposing the temporal gauge $A_t = 0$ at $s = \underline s$. Altogether, we arrive at the \textit{caloric-temporal gauge} 
    \begin{equation}\label{eq:CTC}\tag{C-T}
    \begin{split}
        A_s 
            &= 0
            \qquad \text{on $[0, T] \times \R^2 \times [0, \underline s]$,}
            \notag 
            \\
        {A_t}
            &= 0
            \qquad \text{on $[0, T] \times \R^2 \times \{\underline s\}$.}
    \end{split}
    \end{equation}
There is a residual gauge-freedom in $(t,s)$-independent transformations; this was eliminated by imposing Coulomb gauge \eqref{coulomb-data} at $t = 0$ and $s = 0$. However, to better leverage the parabolic smoothing, it is more convenient to impose the Coulomb gauge at $s = \underline s$, 
    \begin{equation}\tag{$\underline{\mathrm{Coul}}$}\label{eq:coulomb-sbar}
		\partial^j {A_j} = 0 \qquad \text{on $\{0\} \times \R^2 \times \{\underline s\}$}.
	\end{equation}
Of course, since the curl-free part of the magnetic field is transported \eqref{A-transport} along the Yang-Mills heat flow in caloric gauge, the conditions \eqref{coulomb-data} and \eqref{coulomb-sbar} are actually equivalent under \eqref{CTC}. 

\begin{remark}
    Our choice of caloric-temporal gauge takes inspiration from the work of Oh \cite{Oh2014a, Oh2015a} on the hyperbolic Yang-Mills equation on $\R^{1 + 3}$, along with its later use by Gavrus \cite{Gavrus2022}. Indeed, similar to these works, our setting is ``sub-critical'' in some sense, so the gauge-choice implicitly distinguishes a particular choice of ``low-frequency'', $s = \underline s$ in our notation.  
\end{remark}

Before turning to the analysis of the caloric-temporal-Coulomb gauge, we should first briefly address the existence of a smooth solution to paradifferential formulation of the Chern-Simons-Schr\"odinger equation satisfying these gauge conditions. By Proposition \ref{prop:dP-exist}, there exists a smooth solution, so it remains to find a gauge-transformation which sends it to a smooth solution satisfying \eqref{CTC} and \eqref{coulomb-sbar}. 

\begin{definition}
    We say that $\chi : [0, T] \times \R^2 \times [0, \underline s] \to \R$ is \textit{regular gauge-transformation} if it has the following regularity properties: 
    \begin{itemize}
        \item $\nabla_{t, x, s} \chi \in C^\infty_{t, s} H^\infty_x ([0, T] \times \R^2 \times [0, \underline s])$.
        
        \item $\chi \in C^\infty_{t, x,s, \loc} ([0, T] \times \R^2 \times [0, \underline s])$.
    \end{itemize}
\end{definition}

Comparing with Definition \ref{def:smooth-dP}, it is clear that a regular gauge-transformation sends a smooth solution $(\widetilde A, \widetilde \phi)$ to another smooth solution $(\widetilde A + d\chi, e^{i \chi}\widetilde \phi)$. We can construct the gauge-transformation to caloric-temporal-Coulomb gauge by solving a system of linear transport and elliptic equations, where the regularity can easily be read-off. 

\begin{proposition}\label{prop:ctc-exists}
    Let $(\widetilde A, \widetilde \phi)$ be a configuration on $[0, T] \times \R^2 \times [0, \underline s]$ solving the paradifferential formulation of the Chern-Simons-Schr\"odinger equation \eqref{dP-CSS}. Then 
        \begin{align*}
            A_\bfa 
                &:= \widetilde A_\bfa + \partial_\bfa \chi, \\
            \phi 
                &:= e^{i \chi} \widetilde\phi ,
        \end{align*}
    where 
         \[
            \chi 
                := -(\partial_s)^{-1} A_s - (\partial_t)^{-1} \underline A_t - \Delta^{-1} \partial^\ell {\underline {A}_\ell}_{|t = 0} ,
        \]
    solves \eqref{dP-CSS} in caloric-temporal gauge \eqref{CTC} with initial data in Coulomb gauge \eqref{coulomb-sbar}. Moreover, if $(\widetilde A, \widetilde \phi)$ is a smooth solution, then $\chi$ is a regular gauge transformation, and thereby $(A, \phi)$ is also a smooth solution. 
\end{proposition}

Once we are in caloric-temporal-Coulomb gauge, we can recover the electromagnetic potential at any given $(t, s)$ by solving a system of transport equations in $s$, with initial data at $s = \underline s$, then solving a transport equation in $t$, with initial data at $t = 0$, and finally solving an elliptic equation at $t = 0$ and $s = \underline s$. To be more precise, we will decouple the divergence-free and curl-free parts of the magnetic potential.

\begin{proposition}[Electromagnetic potentials under caloric-temporal-Coulomb gauge]\label{prop:A-ctc}
    Let $(A, \phi)$ be a smooth solution on $[0, T] \times \R^2 \times [0, \underline s]$ to the paradifferential formulation of the Chern-Simons-Schr\"odinger equation \eqref{dP-CSS} in caloric-temporal gauge \eqref{CTC} with initial data in Coulomb gauge \eqref{coulomb-sbar}. Then the electric potential obeys:  
    \begin{itemize}
        \item transport equation in $s$,
            \begin{equation}\label{eq:At-transport}
            \begin{split}
                \partial_s A_t 
                    &= \nabla F_{tx},\\
                {A_t}_{|s = \underline s}
                    &= 0,
            \end{split}
            \end{equation}
    \end{itemize}
    curl-free part of the magnetic potential obeys:
    \begin{itemize}
        \item transport equation in $s$,
            \begin{equation}\label{eq:Ax-transport-cf}
            \begin{split}
                \partial_s \PP^\cf A_x 
                    &= 0,\\
                {\PP^\cf A_x}_{|s = \underline s}
                    &= \PP^\cf \underline A_x, 
            \end{split}
            \end{equation}

        \item transport equation in $t$ at $s = \underline s$, 
            \begin{equation}\label{eq:Ax-transport-t}
            \begin{split}
                \partial_t \PP^\cf \underline A_x 
                    &= \underline F_{tx},\\
                \PP^\cf {\underline A_x}_{|t = 0}
                    &= 0.
            \end{split}
            \end{equation} 
    \end{itemize}
    and the divergence-free part of the magnetic potential obeys: 
    \begin{itemize}
        \item transport equation in $s$,
            \begin{equation}\label{eq:Ax-transport-df}
            \begin{split}
                \partial_s \PP^\df A_x 
                    &= \nabla F_{12},\\
                {\PP^\df A_x}_{|s = \underline s}
                    &= \PP^\df \underline A_x,
            \end{split}
            \end{equation}
        
        \item modified Biot-Savart law at $s = \underline s$, 
            \begin{equation}\label{eq:Ax-elliptic}
                \PP^\df \underline A_x
                    = \nabla^{-1} \Bigl(  (1 - |\underline\phi|^2) \cdot \underline F_{12} +\underline{\bfD \phi} \cdot \underline{\bfD \phi} \Bigr) + \PP^\df \bigl( \underline \phi \cdot \underline{\bfD \phi} \bigr).
            \end{equation}
    \end{itemize}

\end{proposition}

The transport equations \eqref{At-transport}, \eqref{Ax-transport-cf}, \eqref{Ax-transport-df} in $s$ for the electromagnetic potential are obvious from writing the (dynamic) Yang-Mills heat flow \eqref{d-covariant-heat2} in caloric-temporal gauge \eqref{CTC}, and then projecting to the curl- and divergence-free parts in the equation for the magnetic potential. The transport equation \eqref{Ax-transport-t} in $t$ for the magnetic potential is simply the projection to the divergence-free part of the curvature component $\underline F_{tx}$ written in temporal gauge at $s = \underline s$ and Coulomb gauge \eqref{coulomb-sbar} at $t = 0$. We defer the proof of the modified Biot-Savart law \eqref{Ax-elliptic} for the divergence-free part of the magnetic potential to Lemma \ref{lem:vorticity}.

Integrating the transport equations in $s$, we write the electromagnetic potentials using the fundamental theorem of calculus as a sum of the contribution at $s = \underline s$, which should be regarded as ``low-frequency'' and thereby smooth, and the contribution from the heat-times $[s, \underline s]$, which represents the ``high-frequency'' part -- schematically, we have 
    \begin{align*}
        \text{high-frequency part of } A_t 
            &\approx (\partial_s)^{-1} \nabla F_{tx}, \\
        \text{high-frequency part of } \PP^\df A_x 
            &\approx (\partial_s)^{-1} \nabla F_{12}.
    \end{align*}
Recalling that parabolic scaling dictates $\nabla \approx s^{-1/2}$ and $(\partial_s)^{-1} \approx s$, along with the heuristic bounds on the electromagnetic field from Section \ref{sec:tension}, we predict that the electromagnetic potentials have sizes
    \begin{align*}
        \| \text{high-frequency part of } A_t  \|_{L^\infty_x} 
            &\lessapprox s^0, \\
        \| \text{high-frequency part of } \PP^\df A_x  \|_{L^\infty_x} 
            &\lessapprox s^{\frac12},
    \end{align*}
In practice, we will not need sharp bounds, settling instead for the less-efficient collection of bounds, 

\begin{proposition}[Bounds for electromagnetic potentials in caloric-temporal-Coulomb gauge]\label{thm:ctc-bounds}
    Let $(A, \phi)$ be a smooth solution on $[0, T] \times \R^2 \times [0, \underline s]$ to the paradifferential formulation of the Chern-Simons-Schr\"odinger equation \eqref{dP-CSS} in caloric-temporal gauge \eqref{CTC} with initial data in Coulomb gauge \eqref{coulomb-sbar}. Then the following hold:
    \begin{itemize}
        \item \textup{(Bounds for curl-free part of magnetic potential).}
            \begin{equation}\label{eq:Ax-ctc-cf}
                \bigl\| \PP^\cf A_x \bigr\|_{L^\infty_{t, s} H^{100}_x} 
                    \lesssim_{\underline s} T.
            \end{equation}

        \item \textup{(Bounds for divergence-free part of magnetic potential).} 
            \begin{equation}\label{eq:Ax-ctc-df}
                \bigl\| \PP^\df A_x \bigr\|_{L^\infty_{t, s} L^4_x} + \bigl\| \nabla \PP^\df A_x \bigr\|_{L^\infty_{t, s} H^{1/2}_x}
                    \lesssim_{\underline s} 1.
            \end{equation}

        \item \textup{(Uniform bounds for magnetic potential).}
            \begin{equation}\label{eq:Ax-ctc}
                \bigl\| A_x \bigr\|_{L^\infty_{t, x, s}} + \bigl\| \partial^j A_j \bigr\|_{L^\infty_{t, x, s}} 
                    \lesssim_{\underline s} \langle T \rangle.
            \end{equation}

        \item \textup{(Uniform bounds for electric potential).} 
            \begin{equation}\label{eq:At-ctc}
                \bigl\|  A_t(s) \bigr\|_{L^\infty_{t, x}}
                    \lesssim_{\underline s} 1 + \log^{2} \Bigl( \frac{\,\underline s\,}{s} \bigr) .
            \end{equation}
    \end{itemize}
\end{proposition}

\begin{remark}
    As suggested from dimensional analysis, the top-order bounds in Theorem \ref{thm:ctc-bounds} are not the most efficient. Indeed, for each $\gamma > 0$, one can replace the Sobolev bound for the magnetic potential in \eqref{Ax-ctc-df} with
        \[
            \bigl\| \nabla \PP^\df A_x \bigr\|_{L^\infty_{t, s} H^{1- \gamma}_x}
                \lesssim_{\gamma, \underline s} 1
        \]
    by working harder. 
\end{remark}

\begin{remark}
    We can replace the $L^4$-norm in \eqref{Ax-ctc-df} with any $L^p$-norm with $2 < p \leq \infty$. At the endpoint $p = 2$, our use of Hardy-Littlewood-Sobolev breaks down, though we can use a logarithmic refinement, such as the weak $L^2$-norm, which is consistent with the expected $r^{-1}$-asymptotics for the low-frequency part of $A_x$. 
\end{remark}

It is easy to see that the uniform bound on the magnetic potential \eqref{Ax-ctc} follow from the Sobolev bounds on its curl-free and divergence-free parts \eqref{Ax-ctc-cf}-\eqref{Ax-ctc-df} along with Sobolev embedding. The estimates for the curl-free part of the magnetic potential \eqref{Ax-ctc-cf} and the electric potential \eqref{At-ctc} follow from integrating their respective transport equations and using parabolic smoothing,

\begin{proof}[Proof of bound on curl-free part of magnetic potential \eqref{Ax-ctc-cf}]
    Integrating the transport equations \eqref{Ax-transport-cf} and \eqref{Ax-transport-t}, 
        \[
            \PP^\cf A_x (t, s)
                = \PP^\cf (\partial_t)^{-1} \underline F_{tx}.
        \]
    Estimating the right-hand side, 
        \begin{align*}
            \bigl\| \PP^\cf  (\partial_t)^{-1} \underline F_{tx} \bigr\|_{L^\infty_t H^{100}_x}
                &\lesssim \int_0^T \bigl\| F_{tx} (\underline s) \bigr\|_{H^{100}_x} \, dt \\
                &\lesssim_{\underline s} \int_0^T \bigl\| F_{tx} (\tfrac12 \underline s) \bigr\|_{L^2_x} \, dt \\
                &\lesssim_{\underline s} \int_0^T \Bigl(\| \Amp \|_{L^2_x} +  \| \phi \|_{L^\infty_x} \| \bfD \phi \|_{L^2_x} + \| \nabla F_{12} \|_{L^2_x} \Bigr) (\tfrac12 \underline s) \, dt \\
                &\lesssim_{\underline s} T,
        \end{align*}
    using parabolic smoothing \eqref{standard-heat} for the linear heat equation \eqref{Fmunu-heat} in the second line, the definition of the Amp\'ere tension field in the third line, and the smoothing bounds \eqref{inefficient-smooth}, \eqref{Dphi-sobolev}, \eqref{Gauss-sobolev} in the last line. 
\end{proof}

\begin{proof}[Proof of electric potential bound \eqref{At-ctc}]
    Integrating the transport equation \eqref{At-transport},
            \[
                A_t (s) 
                    = -(s \partial_s)^{-1} (s \nabla F_{tx}). 
            \]
    Rewriting the right-hand side using the dynamic Yang-Mills heat flow \eqref{d-covariant-heat2} and estimating, 
        \begin{align*}
            \bigl\| (s\partial_s)^{-1} (s \nabla F_{tx}) \bigr\|_{L^\infty_x}
                &\lesssim \int_s^{\underline s} \bigl\| {s'}^\frac12 ({s'}^\frac12 \nabla) F_{tx}  \bigr\|_{L^\infty_x} \frac{ds'}{s'} \\
                &\lesssim \int_s^{\underline s} \Bigl( \| \phi(s') \|_{L^\infty_x}^2  + \underline s^{-1}\Bigr) \frac{ds'}{s'} \\
                &\lesssim \int_s^{\underline s} \underline s^{-1} \Bigl( 1 + \log\Bigl( \frac{\underline s}{\, s\,} \Bigr) \Bigr) \frac{ds'}{s'}\\
                &\lesssim_{\underline s} 1 + \log^{2} \Bigl( \frac{\,\underline s\,}{s} \bigr) 
        \end{align*}
    using the electric field bound \eqref{Fmunu-Lp} in the second line, the Brezis-Gallou\"et-type inequality \eqref{brezis} in the third line, and integrating appropriately. 
\end{proof}

It remains to bound the divergence-free part of the magnetic potential \eqref{Ax-ctc-df}. Thanks to conservation of the abelian Higgs energy and transport of the Gauss tension field, it actually suffices to prove the bound for a general caloric extension $(A, \phi) \in \pzcE^1 ([0, \underline s])$. We begin by decomposing, 
    \[
        \PP^\df A_x 
            = \underbrace{\PP^\df \underline A_x}_{\text{``low-frequency''}} + \quad \underbrace{\bigl(\PP^\df A_x - \PP^\df \underline A_x\bigr)}_{\text{``high-frequency''}}. 
    \]

\begin{proposition}\label{prop:Ax-coulomb}
    Let $(A, \phi) \in \mathpzc E^1([0, \underline s])$ be a caloric extension on $\R^2 \times [0, \underline s]$. Then we have:
    \begin{itemize}
        \item \textup{(Bound for ``high-frequency'' part).}
            \begin{equation}\label{eq:Ax-coulomb-high}
                \bigl\| \PP^\df A - \PP^\df \underline A \bigr\|_{L^\infty_{s} H^{3/2}_x}
                    \lesssim_{\underline s} 1.
        \end{equation}
        
        \item \textup{(Bound for ``low-frequency'' part).}
            \begin{equation}\label{eq:Ax-coulomb-low}
            \bigl\| \PP^\df \underline A \bigr\|_{L^4} + \bigl\| \nabla \PP^\df \underline A \bigr\|_{H^{100}_x} 
                \lesssim_{\underline s} 1. 
            \end{equation}
        
        \item \textup{(Bound for magnetic potential).}
            \begin{equation}\label{eq:Ax-coulomb}
            \bigl\| \PP^\df A \bigr\|_{L^\infty_s L^4_x} + \bigl\| \nabla \PP^\df A \bigr\|_{L^\infty_s H^{1/2}_x} 
                \lesssim_{\underline s} 1. 
            \end{equation}
    \end{itemize}
\end{proposition}

The bound \eqref{Ax-coulomb} follows from the previous two bounds \eqref{Ax-coulomb-high}-\eqref{Ax-coulomb-low}, the triangle inequality, and Sobolev embedding. Observe that \eqref{Ax-transport-df} continues to hold, so integrating and using parabolic smoothing yields

\begin{proof}[Proof of high-frequency bound \eqref{Ax-coulomb-high}]
Integrating the transport equation \eqref{Ax-transport-df}, 
    \[
        \PP^\df A_x (s) - \PP^\df \underline A_x 
            = - (s \partial_s)^{-1} (s F_{sx}).
    \]
    Estimating the right-hand side, 
    \begin{align*}
        \bigl\| (s \partial_s)^{-1} (s F_{sx}) \bigr\|_{L^\infty_{t, s} H^{3/2}_x}
            &\lesssim \int_0^{\underline s} s^{\frac14} \bigl\| (s^{\frac34} \langle \nabla \rangle^{\frac32}) \nabla F_{12} \bigr\|_{L^\infty_t L^2_x} \frac{ds}{s} \\
            &\lesssim \int_0^{\underline s} s^\frac14 \bigl\| \nabla F_{12} \bigr\|_{L^\infty_t L^2_x}\frac{ds}{s} \\
            &\lesssim \int_0^{\underline s} s^\frac14 \Bigl( \| \phi \|_{L^\infty_x} \| \bfD \phi \|_{L^2_x} + \| \nabla \Gauss \|_{L^2_x} \Bigr) \frac{ds}{s}  \\
            &\lesssim_{\underline s} 1. 
        \end{align*}
    using parabolic smoothing \eqref{standard-heat} for the linear heat equation \eqref{Fmunu-heat} in the second line, the definition of the Gauss tension field in the third line, and the smoothing bounds \eqref{inefficient-smooth}, \eqref{Dphi-sobolev}, \eqref{Gauss-sobolev} in the last line. 
\end{proof}

To prove \eqref{Ax-coulomb-low}, a first-take would be to use the Biot-Savart law to re-express the divergence-free part of the magnetic potential in terms of the magnetic field,
    \[
        \PP^{\df} A_x
            = \nabla^{-1} F_{12}.
    \]
The Biot-Savart law is favourable for estimating $A$ at high-frequencies, particularly its derivatives $\nabla^{(n)} \nabla A$, however, the formula does not seem suitable for estimating low-frequencies due to the slow spatial decay of the magnetic field, as it is only square-integrable for generic finite-energy configurations. The most then which can be read-off at bottom-order from the Biot-Savart law, by Hardy-Littlewood-Sobolev inequality, is a bounded mean oscillation bound -- far from the desired $L^4$-bound. To overcome this issue, we instead use the modified Biot-Savart formula \eqref{Ax-elliptic}, which should be thought of as a decomposition into a ``topological'' part and a ``non-topological'' part, 
    \[
        \PP^\df A_x
            = \underbrace{\nabla^{-1}\Bigl(  (1 - |\phi|^2) \cdot F_{12} +{\bfD \phi} \cdot {\bfD \phi} \Bigr)}_{\substack{\text{topologically non-trivial}\\\text{$r^{-1}$-asymptotics}}} + \underbrace{\PP^\df (\phi \cdot \bfD \phi)}_{\substack{\text{topologically trivial}\\\text{$L^2$-integrable}}}.
    \]
For further discussion on this perspective, see the remarks surrounding Lemma \ref{lem:vorticity} in Appendix \ref{app:conservation}.

\begin{proof}[Proof of low-frequency bound \eqref{Ax-coulomb-low}]
    Using the Biot-Savart law and parabolic smoothing \eqref{F12-sobolev} for the linear heat equation \eqref{F12-heat}, we obtain the top-order Sobolev bound
        \begin{align*}
            \bigl\| \nabla \PP^\df \underline A_x \bigr\|_{H^{100}_x} 
                &\lesssim \bigl\| \underline F_{12} \bigr\|_{H^{100}_x} \lesssim_{\underline s} 1. 
        \end{align*}
    For the $L^4$-bound, we appeal to the modified Biot-Savart law \eqref{Ax-elliptic}, which is of the form 
        \[
            \PP^\df \underline A_x 
                = \nabla^{-1} (\underline\sfU \cdot \underline\sfU) + \PP^\df (\underline \phi \cdot \underline{\bfD \phi}), \qquad \sfU = \tfrac12 (1 - |\underline \phi|^2) \text{ or } F_{12} \text{ or } \bfD \phi. 
        \]
    For the first term on the right,
        \begin{align*}
            \bigl\| \nabla^{-1} (\underline \sfU \cdot \underline \sfU) \bigr\|_{L^4_x}
                &\lesssim \bigl\| \underline \sfU \cdot \underline \sfU \bigr\|_{L^{4/3}_x} \\
                &\lesssim \bigl\| \underline \sfU \bigr\|_{L^2_x} \bigl\| \underline \sfU \bigr\|_{L^4_x} \\
                &\lesssim_{\underline s} 1,
        \end{align*}
    using Hardy-Littlewood-Sobolev in the first line, H\"older in the second line, and parabolic smoothing \eqref{F12-sobolev}, \eqref{Dphi-sobolev}, \eqref{charge-sobolev} in the last line. For the second term on the right, 
        \begin{align*}
            \bigl\| \PP^\df (\underline \phi \cdot \underline{\bfD \phi})\bigr\|_{L^4_x}
                &\lesssim \bigl\| \underline \phi \bigr\|_{L^\infty_x} \bigl\| \underline{\bfD \phi} \bigr\|_{L^4_x} \\
                &\lesssim_{\underline s} 1,
        \end{align*}
    by the smoothing bounds \eqref{Dphi-sobolev}, \eqref{inefficient-smooth}. Collecting the previous estimates completes the proof. 
\end{proof}

\begin{remark}
    As a trivial consequence of \eqref{Ax-coulomb}, any finite-energy configuration $(A, \phi) \in \frE(\underline s)$ on $\R^2$ in Coulomb gauge obeys
        \[
           \| A \|_{L^4_x} 
                    \lesssim_{\underline s} 1.
        \]
    Interestingly, it is not clear how to directly deduce this bound from the Biot-Savart law. This seems, at least superficially, similar to a scenario encountered in the context of wave maps on $\R^2$ \cite[Remark 4.5]{Tao2009}.
\end{remark}

We proved some difference bounds for the magnetic potential in Proposition \ref{prop:Ax-diff}. We supplement those by showing how differences of the magnetic potential are propagated under the caloric-temporal gauge, and also proving a fixed-$t$ difference bound for the electric potential. These will come to play in Section \ref{sec:difference}. 

\begin{proposition}[Difference bounds for electromagnetic potentials in caloric-temporal-Coulomb gauge]
    Let $(A, \phi)$, $(A', \phi')$ be smooth solutions on $[0, T] \times \R^2 \times [0, \underline s]$ to the paradifferential formulation of the Chern-Simons-Schr\"odinger equation \eqref{dP-CSS} in caloric-temporal gauge \eqref{CTC}. Then, for $s_* \leq \tfrac1{1000} \, \underline s$, we have:
    \begin{itemize}
        \item \textup{(Propagation of magnetic potential differences).}
            \begin{equation}\label{eq:Ax-diff-trans}
                \frac{d}{dt} \bigl\| \delta \nabla \underline A_x \bigr\|_{H^{100}_x}^2
                    \lesssim \log\Bigl( \frac{\underline s}{s_*} \Bigr)^\frac12 \cdot\Bigl( |s_*|+ |\dist_{\pzcE^0}|^2 \Bigr),
            \end{equation}
        and 
            \begin{equation}\label{eq:Ax-diff-trans-2}
                \frac{d}{dt} \bigl\|\delta  \underline A_x \bigr\|_{L^4_x}^2 
                     \lesssim \log\Bigl( \frac{\underline s}{s_*} \Bigr)^\frac12 \cdot\Bigl( |s_*| + |
                     \dist_{\pzcE^0}|^2 \Bigr).
            \end{equation}

        \item \textup{(Difference bound for electric potential).}
            \begin{equation}\label{eq:At-diff}
                \bigl\| \delta A_t \bigr\|_{L^\infty_s L^2_x}
                    \lesssim_{\underline s} \log^\frac12 \Bigl( \frac{\underline s}{s_*} \Bigr) \cdot\Bigl( |s_*|^\frac12 + \dist_{\pzcE^0} \Bigr) .
            \end{equation}
    \end{itemize}

\end{proposition}

\begin{proof}[Proof of \eqref{Ax-diff-trans}, \eqref{Ax-diff-trans-2}]
    Writing in temporal gauge, we arrive at the transport equation
        \begin{align*}
            \partial_t \delta \underline A_x 
                = \delta\underline F_{tx}. 
        \end{align*}
    We can immediately estimate the right-hand side via the electric field difference bound \eqref{Ftx-diff-low} and Sobolev embedding, 
        \[
            \| \underline \delta F_{tx} \|_{L^4} + \| \nabla \delta \underline F_{tx} \|_{H^{100}}
                \lesssim \log\Bigl( \frac{\underline s}{s_*} \Bigr)^\frac12 \cdot\Bigl( |s_*|^\frac12 + \dist_{\pzcE^0} \Bigr),
        \]
    as desired. 
\end{proof}

\begin{proof}[Proof of difference bound for electric potential \eqref{At-diff}]
    Integrating the transport equation \eqref{At-transport},
        \[
            \delta A_t (s) 
                = - (s \partial_s)^{-1} (s \nabla \delta F_{tx} ) (s).
        \]
    Estimating the right-hand side, 
        \begin{align*}
            \bigl\| (s \partial_s)^{-1} (s\nabla \delta F_{tx})(s) \bigr\|_{L^2_x}
                &\lesssim  \Bigl\| \int_{s}^{\underline s}  {s'}^\frac12 ({s'}^\frac12 \nabla) \delta  F_{tx} \frac{ds'}{s'} \Bigr\|_{L^2_x} \\
                &\lesssim \bigl\| s^\frac12 (s^\frac12 \nabla)^{(\leq 2)} \delta F_{tx} \bigr\|_{L^2_{\frac{ds}{s}} L^2_x} \\
                &\lesssim_{\underline s} \log^\frac12 \Bigl( \frac{\underline s}{s_*} \Bigr) \cdot\Bigl( |s_*|^\frac12 + \dist_{\pzcE^0} \Bigr) ,
        \end{align*}
    by almost orthogonality (c.f. Lemma \ref{lem:orthogonal}) and the electric field difference bound \eqref{Ftx-diff}.  
\end{proof}

\section{The linear paradifferential equation}\label{sec:strichartz}

Continuing on the paradifferential analogy, we want to abstractly view the covariant Hessian $(s \bfD^{(2)}) \phi(s)$ as a high-frequency wave evolving with respect to a Schr\"odinger equation with low-frequency coefficients. To this end, given $(A, \phi)$ a smooth solution on $[0, T] \times \R^2 \times [0, \underline s]$ to the paradifferential formulation of the Chern-Simons-Schr\"odinger equation \eqref{dP-CSS}, we introduce the \textit{linear paradifferential equation} as the family of linear electromagnetic Schr\"odinger equations indexed by $s \in [0, \underline s]$,
    \begin{equation}\label{eq:paralinear-abstract}\tag{Paralin}
    \begin{split}
        \bigl( i \bfD_t + \bfD^j \bfD_j \bigr) u (s)
            &= \sfN (s) ,\\
        u(s)_{|t = 0}
            &= u^{\mathrm{in}} (s),
    \end{split}
    \end{equation}
where $\bfD_\mu = \partial_\mu - i A_\mu (s)$ and $u(s), \sfN(s) : [0, T] \times \R^2 \to \C$. Then, in view of \eqref{Ds-schrodinger}, we can recast the covariant Hessian as a solution to \eqref{paralinear-abstract}.

One should view the solution $u(s)$ as localised to frequency $s^{-1/2}$, and the potential $A_\mu (s)$ as localised to lower-frequencies. We will make the former heuristic rigorous by showing that the linear paradifferential equation \eqref{paralinear-abstract} propagates $s$-weighted Sobolev-type norms, while we achieve the latter by imposing the caloric-temporal gauge \eqref{CTC} on $[0, T] \times \R^2 \times [0, \underline s]$ and Coulomb gauge \eqref{coulomb-sbar} at $t = 0$ and $s = \underline s$. To bound the growth of energy along the flow, we introduce the \textit{control parameter}, 
    \[
        \cC^\sharp(t, s)
            :=  \bigl\| \overline\Gauss (t) \bigr\|_{L^\infty_x} + \bigl\| \phi(t, s) \bigr\|_{L^\infty_x}^2 + \underline s^{-1}.
    \]

\begin{proposition}\label{prop:strichartz}
    Let $(A, \phi)$ be a smooth solution on $[0, T] \times \R^2 \times [0, \underline s]$ to the paradifferential formulation of the Chern-Simons-Schr\"odinger equation \eqref{dP-CSS}, and suppose $u (s): [0, T] \times \R^2 \to \C$ solves the linear paradifferential equation \eqref{paralinear-abstract}. Then the following estimates hold:   
    \begin{itemize}
        \item \textup{(Energy estimate).} For each non-negative integer $n \in \N_0$, we have
            \begin{equation}\label{eq:paralinear-energy}
            \frac{d}{dt} \bigl\| (s^\frac12 \bfD)^{(n)} u \bigr\|_{L^2_x}^2
                \lesssim \bigl\| (s^\frac12 \bfD)^{(n)} u \bigr\|_{L^2_x} \bigl\| (s^\frac12 \bfD)^{(n)} \sfN \bigr\|_{L^2_x} + \cC^\sharp \cdot \bigl\| (s^\frac12 \bfD)^{(\leq n)} u \bigr\|_{L^2_x}^2 .
            \end{equation}

        \item \textup{(Strichartz estimate).} Suppose furthermore that $(A, \phi)$ satisfies the caloric-temporal gauge \eqref{CTC} with initial data in Coulomb gauge \eqref{coulomb-sbar}. Then there exists a constant $C > 0$ such that, for Strichartz admissible $(p, q)$, i.e. $2 \leq p, q \leq \infty$ such that $\tfrac2p + \tfrac2q = 1$ and $(p, q) \neq (2, \infty)$, we have 
            \begin{equation}\label{eq:paralinear-str}
            s^{\frac1{2p}} \bigl\|  u \bigr\|_{L^p_t L^q_x}
                \lesssim_{\underline s, p, T} \Bigl( \bigl\| (s^\frac12 \bfD)^{(\leq 1)} u^{\mathrm{in}} \bigr\|_{L^2_x} + \bigl\| (s^\frac12 \bfD)^{(\leq 1)}\sfN \bigr\|_{L^1_t L^2_x} \Bigr) \exp \Bigl( \int_0^T C \cdot\cC^\sharp (t) \, dt \Bigr).
            \end{equation}
    \end{itemize}
\end{proposition}

\begin{remark}
    The weight $s^{\frac1{2p}}$ on the left-hand side of \eqref{paralinear-str} represents a $\tfrac1p$-derivative loss compared to the classical Strichartz estimate for the constant-coefficient Schr\"odinger equation. 
\end{remark}

We begin by proving the energy estimate for the linear paradifferential equation. Amusingly, by working in a ``covariant'' fashion, i.e. commuting using magnetic derivatives $\bfD$ rather than usual derivatives $\nabla$, we can prove the estimate without having to fix a gauge for the underlying background $(A, \phi)$ solving the flow \eqref{dP-CSS}. 

\begin{lemma}[Low-high paraproduct commutator]\label{lem:low-high}
    Let $(A, \phi)$ be a smooth solution on $[0, T] \times \R^2 \times [0, \underline s]$ to the paradifferential formulation of the Chern-Simons-Schr\"odinger equation \eqref{dP-CSS}, and suppose that $u : \R^2 \to \C$ is a sufficiently regular scalar field. Then for each positive integer $n \in \N$, we have 
        \begin{equation}\label{eq:low-high-bound}
            \Bigl\| \bigl[ i \bfD_t + \bfD^j \bfD_j , (s^{\frac12}\bfD)^{(n)} \bigr] u  \Bigr\|_{L^2_x}
                \lesssim \cC^\sharp \cdot \bigl\| (s^{\frac12} \bfD)^{(\leq n)} u  \bigr\|_{L^2_x}.
        \end{equation}
\end{lemma}

\begin{proof}
    We claim the following commutator identity, 
        \begin{equation}\label{eq:low-high-commutator}
            \begin{split}
            \Bigl[i \bfD_t + \bfD^j \bfD_j , \bfD^{(n)}\Bigr]  
                &= \sum_{a + b = n - 1} \nabla^{(a)} F_{t x} \cdot \bfD^{(b)} + \sum_{a + b = n} \nabla^{(a)} F_{12} \cdot \bfD^{(b)}.
            \end{split}
        \end{equation}
    This follows by induction on $n$. Indeed, for $n = 1$, we use \eqref{commute}, \eqref{laplace-commute} to obtain
        \begin{align*}
            \Bigl[i \bfD_t + \bfD^j \bfD_j , \bfD \Bigr] 
                &= F_{tx} + F_{12} \cdot \bfD  + \nabla F_{12} .
        \end{align*}
    Assuming the result for $n - 1$, we write the commutator for $n$ as 
        \begin{align*}
            \Bigl[i \bfD_t + \bfD^j \bfD_j , \bfD^{(n)}\Bigr]  
                &= \Bigl[i \bfD_t + \bfD^j \bfD_j , \bfD\Bigr] \bfD^{(n - 1)} + \bfD \Bigl[i \bfD_t + \bfD^j \bfD_j, \bfD^{(n - 1)} \Bigr] .
        \end{align*}
    Then the desired identity clearly follows from the base case $n = 1$ and the inductive hypothesis. 

    We then multiply the commutator identity \eqref{low-high-commutator} by $s^{\frac{n}{2}}$ and distribute powers of $s$ appropriately into the right-hand side. For the first two terms on the right-hand side of the commutator identity \eqref{low-high-commutator}, we apply the estimate \eqref{Fmunu-Lp} for the electric field,  
        \begin{align*}
            \bigl\| s^\frac12 (s^{\frac12} \nabla)^{(a)} F_{t x}\cdot (s^{\frac12} \bfD)^{(b)} u \bigr\|_{L^2_x}
                &\lesssim \bigl\| s^\frac12 (s^{\frac12} \nabla)^{(a)} F_{t x} \bigr\|_{L^\infty_x} \bigl\| (s^{\frac12} \bfD)^{(b)} u \bigr\|_{L^2_x} \\
                &\lesssim \cC^\sharp \cdot \bigl\| (s^\frac12 \bfD)^{(b)} u \bigr\|_{L^2_x}  
        \end{align*}
    for $a + b = n - 1$, and similarly \eqref{F12-Linfty} for the magnetic field, 
        \begin{align*}
            \bigl\| (s^{\frac12} \nabla)^{(a)} F_{12} \cdot (s^\frac12\bfD)^{(b)} u \bigr\|_{L^2_x} 
                &\lesssim \bigl\| (s^{\frac12} \nabla)^{(a)} F_{12} \bigr\|_{L^\infty_x} \, \bigl\|  (s^\frac12\bfD)^{(b)} u \bigr\|_{L^2_x} \\
                &\lesssim \cC^\sharp  \cdot \bigl\|  (s^\frac12\bfD)^{(b)} u \bigr\|_{L^2_x},
        \end{align*}
    for $a + b = n$. Collecting these calculations completes the proof. 
\end{proof}

\begin{proof}[Proof of energy estimate \eqref{paralinear-energy}]
    Commuting $(s^{\frac12} \bfD)^{(n)}$ through the linear paradifferential equation \eqref{paralinear-abstract},
        \[
            \bigl( i \bfD_t + \bfD^j \bfD_j \bigr) (s^\frac12\bfD)^{(n)} u 
                = \Bigl[ i \bfD_t + \bfD^j \bfD_j , (s^{\frac12} \bfD)^{(n)} \Bigr] u + (s^{\frac12} \bfD)^{(n)} \sfN.
        \]
    Applying the abstract energy estimate \eqref{abstract-energy} to the equation above, we obtain the differential inequality
        \begin{align*}
            \frac{d}{dt} \bigl\|  (s^\frac12\bfD)^{(n)} u \bigr\|_{L^2_x}^2  
                &\lesssim \Bigl\langle  (s^{\frac12} \bfD)^{(n)} u  , \bigl[ i \bfD_t + \bfD^j \bfD_j , (s^{\frac12} \bfD)^{(n)} \bigr] u \Bigr\rangle_{L^2_x} + \Bigl\langle (s^\frac12\bfD)^{(n)} u  , (s^\frac12\bfD)^{(n)} \sfN \Bigr\rangle_{L^2_x}\\
                &\lesssim \cC^\sharp \cdot \bigl\| (s^\frac12 \bfD)^{(\leq n)} u \bigr\|_{L^2_x}^2 + \bigl\| (s^\frac12 \bfD)^{(n)} u \bigr\|_{L^2_x} \bigl\| (s^\frac12 \bfD)^{(n)} \sfN \bigr\|_{L^2_x} , 
        \end{align*}
    using Cauchy-Schwarz and the low-high paraproduct commutator bound \eqref{low-high-bound} to estimate the first term on the right-hand side. This completes the proof.
\end{proof}

It remains to establish the Strichartz estimate for the linear paradifferential equation. One might naively attempt to regard the equation as a perturbation of the constant-coefficient Schr\"odinger equation, for which we can use the classical Strichartz estimate \cite[Chapter 2.3]{Tao2006},

\begin{lemma}[Strichartz estimate for constant-coefficient Schr\"odinger equation]\label{lem:strichartz-classical}
    Let $u: [0, T] \times \R^2 \to \C$ be a solution to the linear Schr\"odinger equation,
        \begin{equation}\label{eq:constant-coeff}\tag{$e^{i t\Delta}$}
        \begin{split}
             \big(i \partial_t + \Delta\big) u
                &= \sfN ,\\
            u_{|t = 0}
                &= u^{\mathrm{in}},
        \end{split}
        \end{equation}
    then, for Strichartz admissible $(p, q)$, the following estimate holds, 
        \[
            \| u \|_{L^p_t L^q_x}
                \lesssim \| u^{\mathrm{in}} \|_{L^2_x} + \| \sfN \|_{L^1_t L^2_x}.
        \]
\end{lemma}

Unfortunately, at least at face-value, this strategy cannot work due to loss of derivatives. Expanding the magnetic Laplacian, \eqref{paralinear-abstract} may be written as
    \[
        (i\partial_t + \Delta) u - 2i A^j \bfD_j u - \bigl(i \partial^j A_j + A_t - A^j A_j\bigr) u
            = \sfN.
    \]
Placing the most dangerous term, magnetic drift $A \cdot \bfD \phi$, on the right, we formally estimate
    \[
        \bigl\| A^j (s) \, \bfD_j u (s) \bigr\|_{L^1_t L^2_x([0, T])} 
            \approx |s|^{-\frac12} \cdot |T| \cdot \bigl\| A \bigr\|_{L^\infty_{t, x, s} ([0, T])} \cdot \bigl\| u (s) \bigr\|_{L^\infty_t L^2_x([0, T])}.
    \]
On unit time-scales $T \sim 1$, this move loses a derivative in the estimate. Nevertheless, we can ammeliorate the derivative-loss by restricting to short time-scales $T \lesssim s^{1/2}$, which, in view of the paradifferential analogy, we refer as the semi-classical time-scale. Our strategy then for proving the Strichartz estimate for the linear paradifferential equation consists of two steps: first, we divide the unit time-scale $[0, T]$ into sub-intervals $J \subseteq [0, T]$ on the semi-classical time-scale; second, we prove Strichartz estimates without derivative-loss on these sub-intervals $J$, where we can regard \eqref{paralinear-abstract} as a perturbation of \eqref{constant-coeff}.

\begin{lemma}[Strichartz estimate on short time-scales]
    Let $(A, \phi)$ be a smooth solution on $[0, T] \times \R^2 \times [0, \underline s]$ to the paradifferential formulation of the Chern-Simons-Schr\"odinger equation \eqref{dP-CSS} in caloric-temporal gauge \eqref{CTC} with initial data in Coulomb gauge \eqref{coulomb-sbar}, and suppose that $u (s): [0, T] \times \R^2 \to \C$ solves the linear paradifferential equation \eqref{paralinear-abstract}. If $(p, q)$ are Strichartz admissible and $J \subseteq [0, T]$ is a sub-interval of length $|J| \sim s^{1/2}$, then 
        \begin{equation}\label{eq:semi-classical-strichartz}
            \| u\|_{L^p_t L^q_x (J)} 
                \lesssim_{\underline s, p, T} \bigl\| (s^\frac12 \bfD)^{(\leq 1)} u \bigr\|_{L^\infty_t L^2_x} + \| \sfN \|_{L^1_t L^2_x}.
        \end{equation}
\end{lemma}

\begin{proof}
    We recast the linear paradifferential equation as a perturbation of the constant-coefficient Schr\"odinger equation, placing the variable-coefficient terms on the right-hand side, 
        \[
            (i \partial_t + \Delta) u 
                = 2i A^j \, \bfD_j u + \big( i \partial^j A_j + A_t - A^j A_j \big) u + \sfN \qquad \text{on $J \times \R^2 \times \{s\}$},
        \]
    and apply the classical Strichartz estimate from Lemma \ref{lem:strichartz-classical}, the triangle inequality, and recalling $|J| \sim s^{1/2}$, 
        \begin{align*}
            \| u  \|_{L^p_t L^q_x (J)}
                &\lesssim \| u \|_{L^\infty_t L^2_x} +  \bigl\|2i A^j \, \bfD_j u + \big( i \partial^j A_j + A_t - A^j A_j \big) u + \sfN \bigr\|_{L^1_t L^2_x(J)} \\
                &\lesssim \| u \|_{L^\infty_t L^2_x} +  \| \sfN \|_{L^1_t L^2_x} + s^\frac12 \Bigl( \bigl\| A_x \cdot \bfD u \bigr\|_{L^\infty_t L^2_x} + \bigl\| A_t \, u \bigr\|_{L^\infty_t L^2_x} + \bigl\| (\partial^j A_j - A^j A_j) \, u \bigr\|_{L^\infty_t L^2_x} \Bigr). 
        \end{align*}
    It remains to show that the last three terms in the second line are controlled by $(s^\frac12 \bfD)^{(\leq 1)} u$. Distributing the powers of $s$ appropriately and estimating $L^2_x \to L^\infty_x \times L^2_x$,  
        \begin{align*}
            s^\frac12 \bigl\| A_x \cdot \bfD u \bigr\|_{L^\infty_t L^2_x}
                &\lesssim \bigl\| A_x \bigr\|_{L^\infty_{t, x, s}} \bigl\| (s^\frac12 \bfD) u \bigr\|_{L^\infty_t L^2_x}  \\
                &\lesssim_{\underline s, T} \bigl\| (s^\frac12 \bfD) u \bigr\|_{L^\infty_t L^2_x},
        \end{align*}
    using the uniform bound on the magnetic potential \eqref{Ax-ctc},
        \begin{align*}
            s^\frac12 \bigl\| A_t\, u \bigr\|_{L^\infty_t L^2_x} 
                &\lesssim \bigl\| s^\frac12 A_t \bigr\|_{L^\infty_{t, x, s}} \| u \|_{L^\infty_t L^2_x} \\
                &\lesssim_{\underline s} \| u \|_{L^\infty_t L^2_x},
        \end{align*}
    using the uniform bound on the electric potential \eqref{At-ctc},
        \begin{align*}
            s^\frac12 \bigl\| (i \partial^j A_j + A^j A_j) u \bigr\|_{L^\infty_t L^2_x} 
                &\lesssim \Bigl( \bigl\| \partial^j A_j \bigr\|_{L^\infty_{t, x, s}} + \bigl\| A_x \bigr\|_{L^\infty_{t, x, s}}^2 \Bigr) \| u\|_{L^\infty_t L^2_x} \\
                &\lesssim_{\underline s, T} \| u \|_{L^\infty_t L^2_x}, 
        \end{align*}
    using the uniform bounds on the magnetic potential \eqref{Ax-ctc}. These complete the proof.
\end{proof}

\begin{proof}[Proof of Strichartz estimate \eqref{paralinear-str}]
    Partition the time interval $[0, T]$ into $O (s^{-1/2})$-many sub-intervals $J \subseteq [0, T]$ with size $|J| \sim s^{1/2}$. We estimate the Strichartz norm on the full interval by those on the sub-intervals,
        \begin{align*}
            \| u  \|_{L^p_t L^q_x ([0, T])} 
                &\leq \Big( \sum_{\substack{J \subseteq [0, T] \text{ disjoint and } \, |J| \sim s^{\frac12}}} \| u \|_{L^p_t L^q_x (J)}^p \Big)^{1/p}\\
                &\lesssim_T s^{-\frac{1}{2p}} \sup_{\substack{J \subseteq [0, T] \text{ and } |J| \sim s^{\frac12}}} \| u \|_{L^p_t L^q_x (J)}.
        \end{align*}
    Combining the short-time Strichartz estimate \eqref{semi-classical-strichartz} with the energy estimate \eqref{paralinear-energy} and Gronwall's inequality, 
        \begin{align*}
            \bigl\|  u  \bigr\|_{L^p_t L^q_x (J)} 
                &\lesssim_{\underline s, p, T} \bigl\| (s^\frac12 \bfD)^{(\leq 1)} u  \bigr\|_{L^\infty_t L^2_x} + \bigl\| \sfN  \bigr\|_{L^1_t L^2_x} \\ 
                &\lesssim_{\underline s, p, T} \Bigl( \bigl\| (s^\frac12 \bfD)^{(\leq 1)} u  \bigr\|_{L^\infty_t L^2_x} + \bigl\| (s^\frac12 \bfD)^{(\leq 1)}\sfN \bigr\|_{L^1_t L^2_x} \Bigr) \exp \Bigl( \int_0^T C \cdot \overline C (t) \, dt \Bigr).
        \end{align*}
    Collecting the previous two calculations completes the proof. 
\end{proof}

\begin{remark}
    The Strichartz estimate \eqref{paralinear-str} and the proof herein is analogous to those of Lim \cite[Proposition 3.7]{Lim2018} in the context of the Jackiw-Pi Chern-Simons-Schr\"odinger equation -- the reader should find it instructive to compare the two arguments. In particular, Lim also divides the time-interval into frequency-dependent sub-intervals, which is a well-known idea in the dispersive PDE literature; a comprehensive overview may be found, among other things, in the article of Burq-G\'erard-Tzvetkov \cite{BurqEtAl2004a}.
\end{remark}

\section{Covariant estimates for \eqref{dP-CSS}}\label{sec:energy}

Armed with the machinery developed thus far, we can finally show global regularity for smooth solutions to the Chern-Simons-Schr\"odinger equation as promised in Section \ref{sec:high-LWP}. In fact, we will prove a more refined statement, namely that the energy profile of the initial data propagates under the flow, and that the solution obeys covariant Strichartz-type bounds. The former roughly implies that energy cannot move from low-to-high frequencies too quickly, thereby conferring a kind of compactness for the Chern-Simons-Schr\"odinger flow. We will use this fact in the concluding Section \ref{sec:cts} when establishing continuity of the data-to-solution map. The latter implies that low-high interactions along the flow are more regular than expected from, say, Sobolev embedding, and accordingly plays a key role in closing the control parameter bound.  

We begin by stating the estimates on unit-time scales, 

\begin{theorem}[Covariant estimates for smooth solutions]\label{thm:dP-est}
    Let $(A, \phi)$ be a smooth solution on $[0, T] \times \R^2 \times [0, \underline s]$ to the paradifferential formulation of the Chern-Simons-Schr\"odinger equation \eqref{dP-CSS}, and suppose that $T \leq T_*$, where $T_* \equiv T_*(\| \overline\Gauss \|_{L^\infty}, \underline s) > 0$, and $\epsilon_* \ll 1$ is a fixed constant. If $\ttc(s) \in L^2_{ds/s} ([0, \underline s])$ is a $\frE$-frequency envelope for $(A, \phi)_{|t = 0} = (A^{\mathrm{in}}, \phi^{\mathrm{in}})$, then we have the following bounds:
    \begin{itemize}
        \item \textup{(Control parameter bound).} 
            \begin{equation}\label{eq:covar-control}
                \int_0^T \cC^\sharp(t, s) \, dt 
                    \leq \epsilon_*. 
            \end{equation}
    
        \item \textup{(Covariant energy estimate).} 
            \begin{equation}\label{eq:covar-env}
                \sup_{t \in [0, T]} \sqrt{\tte^1 [A(t, s), \phi(t, s)]}
                    \leq 2 \ttc(s). 
            \end{equation}

        \item \textup{(Covariant Strichartz estimate).} For Strichartz admissible $(p, q)$, 
            \begin{equation}\label{eq:covar-str}
            \bigl\| s^{\frac{1}{2p} - \frac12} (s^\frac12 \bfD)^{(\leq 99)} (s \bfD^{(2)}) \phi (s) \bigr\|_{L^p_t L^q_x}
                \lesssim_{\underline s, p} \ttc(s).
            \end{equation}

    \end{itemize}
\end{theorem}

\begin{remark}
    Compared to Sobolev embedding, the Strichartz estimate allows one to exchange time-averaging for (almost) half a derivative, or, in the language of parabolic scaling, one-quarter of an $s$-weight. Indeed, the parabolic smoothing estimate \eqref{Dphi-sobolev} yields 
        \[
            \bigl\| s^{- \frac1q} (s \bfD^{(2)}) \phi (s) \bigr\|_{L^q_x}
                \lesssim \ttc(s),
        \]
    while, massaging the exponents slightly, the Strichartz estimate \eqref{covar-str} gives
        \[
            \bigl\| s^{\frac12 (\frac12 - \frac1q)} s^{-\frac12} (s \bfD^{(2)}) \phi (s) \bigr\|_{L^p_t L^q_x}
                \lesssim \ttc(s). 
        \]
    Taking $q \to \infty$, the $s$-weight of the Strichartz estimate saves (up to) $s^{1/4}$ compared to parabolic smoothing.  
\end{remark}

Leveraging conservation of the abelian Higgs energy and transport of the Gauss tension field along the Chern-Simons-Schr\"odinger flow, we can iterate the unit time-scale bounds \eqref{covar-control}, \eqref{covar-env}, \eqref{covar-str}, and appeal to parabolic smoothing \eqref{no-loss} to obtain the following global-in-time \textit{a priori} bounds, 

\begin{corollary}[Global covariant bounds]\label{cor:global-2}
    Let $(A, \phi)$ be a smooth solution on $[0, \infty) \times \R^2 \times [0, \underline s]$ to the paradifferential formulation of the Chern-Simons-Schr\"odinger equation \eqref{dP-CSS}, and suppose that $\kappa \equiv \kappa (\| \overline\Gauss \|_{L^\infty_x}, \underline s) > 0$. If $\ttc(s) \in L^2_{ds/s} ([0, \underline s])$ is a $\frE$-frequency envelope for $(A, \phi)_{|t = 0} = (A^{\mathrm{in}}, \phi^{\mathrm{in}})$, then the following bounds hold:
        \begin{itemize}
            \item \textup{(Control parameter bound).} 
                \begin{equation}\label{eq:global-control}
                    \int_0^T \cC^\sharp(t, s) \, dt 
                        \lesssim \kappa  T.
                \end{equation}
            
            \item \textup{(Covariant energy estimate).}
                \begin{equation}
                    \bigl\| s^{-\frac12} (s^\frac12 \bfD)^{(\leq 2026)} (s \bfD^{(2)}) \phi (t, s) \bigr\|_{L^2_x}
                        \lesssim \ttc(s) \exp(\kappa t).
                \end{equation}

             \item \textup{(Covariant Strichartz estimate).} For Strichartz admissible $(p, q)$, 
            \begin{equation}\label{eq:covar-str-global}
            \bigl\| s^{\frac{1}{2p} - \frac12} (s^\frac12 \bfD)^{(\leq 99)} (s \bfD^{(2)}) \phi (s) \bigr\|_{L^p_t L^q_x}
                \lesssim_{\underline s, p} \ttc(s) \exp(\kappa t).
            \end{equation}
        \end{itemize}
\end{corollary}

Turning to the proof of Theorem \ref{thm:dP-est}, we proceed by a bootstrap argument, assuming there exists a time $T$ such that 
    \begin{itemize}
        \item (Control parameter bound bootstrap).
            \begin{equation}\label{eq:BS-contr}
                \int_0^T \cC^\sharp(t, s) \, dt 
                    \leq 100 \, \epsilon_*.
            \end{equation}

        \item (Covariant energy estimate bootstrap).
            \begin{equation}\label{eq:covar-env-BS}
                \sup_{t \in [0, T]}\sqrt{\tte^1 [A, \phi] (t, s)} 
                    \leq 100 \,  \ttc(s). 
            \end{equation}
    \end{itemize}
To get the continuity argument started, observe that since $(A, \phi)$ is a smooth solution, there exists $T \ll 1$ such that \eqref{BS-contr} holds and, by virtue of Proposition \ref{prop:env-stable}, \eqref{covar-env-BS} also holds. With these assumptions at hand, our goal is to prove the covariant Strichartz estimate \eqref{covar-str} and close the bootstrap by improving \eqref{BS-contr}, \eqref{covar-env-BS} to \eqref{covar-control}, \eqref{covar-env}. The strategy is to view the covariant Hessian as a solution to the linear paradifferential equation and treat the right-hand side as perturbative, in the sense that 

\begin{proposition}[Estimating the perturbative terms]\label{prop:non-linear}
    Let $(A, \phi)$ be a smooth solution on $[0, T] \times \R^2 \times [0, \underline s]$ to the paradifferential formulation of the Chern-Simons-Schr\"odinger equation \eqref{dP-CSS}. Fix $t \in [0, T]$, and suppose $\ttc(s) \in L^2_{ds/s} ([0, \underline s])$ is a $\frE$-frequency envelope for $(A(t), \phi(t))$. Then the right-hand side of the covariant Schr\"odinger equation for the covariant Hessian satisfies, for each non-negative integer $n \in \N_0$, 
        \begin{equation}
            \bigl\| s^{-\frac12} (s^\frac12 \bfD)^{(n)} \bigl(\textup{R.H.S.\eqref{Ds-schrodinger}}\bigr)(t) \bigr\|_{L^2_x} 
                \lesssim_n \cC^\sharp(t) \cdot \ttc(s).
        \end{equation}
\end{proposition}

Taking Proposition \ref{prop:non-linear} for granted (that is, deferring the proof to the end of the section), we can use the analysis of the linear paradifferential equation from Section \ref{sec:strichartz} to complete the proof of Theorem \ref{thm:dP-est},

\begin{proof}[Proof of Theorem \ref{thm:dP-est}]
    We begin with improving the bootstrapped energy estimate \eqref{covar-env-BS} to \eqref{covar-env}. Recall that the Gauss tension field is transported along the Chern-Simons-Schr\"odinger flow \eqref{gauss-transport-intro}, so $\overline\Gauss(t) \equiv \overline\Gauss(0)$. Combined with the energy estimate for the covariant Hessian, where we multiply its equation \eqref{Ds-schrodinger} by $s^{-1/2}$ and apply the abstract energy estimate \eqref{abstract-energy},
    \begin{align*}
        \sup_{t \in [0, T]} \tte^1 [A(t), \phi(t)]
            &\leq \tte^1 [A^{\mathrm{in}}, \phi^{\mathrm{in}}]  + \Bigl|\int_0^T \bigl\langle  s^{-\frac12} (s \bfD^{(2)}) \phi , s^{-\frac12} \text{R.H.S.\eqref{Ds-schrodinger}} \bigr\rangle_{L^2_x} \, dt \Bigr| \\
            &\leq  \tte^1 [A^{\mathrm{in}}, \phi^{\mathrm{in}}]   + \int_0^T  \bigl\| s^{-\frac12} (s \bfD^{(2)}) \phi \bigr\|_{L^2_x} \bigl\| s^{-\frac12} \text{R.H.S.\eqref{Ds-schrodinger}} \bigr\|_{L^2_x} \, dt \\
            &\leq |\ttc(s)|^2 + |\ttc(s)|^2 \cdot O\Bigl( \int_0^T  \cC^\sharp(t) \, dt \Bigr) \\
            &\leq |2\ttc(s)|^2
    \end{align*}
using the bootstrapped energy estimate \eqref{covar-env-BS} and the non-linear estimate from Proposition \ref{prop:non-linear} in the third line and the bootstrapped control parameter bound \eqref{BS-contr} in the last line.

The next order of business is the covariant Strichartz bound \eqref{covar-str}. Commuting $s^{-\frac12} (s^\frac12 \bfD)^{(\leq 99)}$ through the equation \eqref{Ds-schrodinger} for the covariant Hessian, we obtain the covariant Schr\"odinger equation
    \begin{align*}
        \bigl( i \bfD_t + \bfD^j \bfD_j \bigr) s^{-\frac12} (s^\frac12 \bfD)^{(\leq 99)} (s \bfD^{(2)}) \phi 
            &= s^{-\frac12} (s^\frac12 \bfD)^{(\leq 99)} \bigl(\text{R.H.S.\eqref{Ds-schrodinger}} \bigr) \\
                &\qquad + \Bigl[ i \bfD_t + \bfD^j \bfD_j, (s^\frac12 \bfD)^{(\leq 99)} \Bigr]s^{-\frac12} (s \bfD^{(2)}) \phi .
    \end{align*}
We interpret the equation as the linear paradifferential equation \eqref{paralinear-abstract} for $u(s) := s^{-\frac12} (s^\frac12 \bfD)^{(\leq 99)} (s \bfD^{(2)}) \phi$ with perturbative right-hand side. Breaking from gauge-covariance and transforming to caloric-temporal gauge \eqref{CTC} with initial data in Coulomb gauge \eqref{coulomb-sbar}, we can now apply the Strichartz estimate \eqref{paralinear-str} from Proposition \ref{prop:strichartz},
    \begin{align*}
        \bigl\| s^{\frac1{2p}} u(s) \bigr\|_{L^p_t L^q_x}
            &\lesssim_{\underline s, p, T} \bigl\| (s^\frac12 \bfD)^{(\leq 1)} u^{\mathrm{in}} (s) \bigr\|_{L^2_x} + \bigl\| s^{-\frac12} (s^\frac12 \bfD)^{(\leq 100)} \bigl( \text{R.H.S.\eqref{Ds-schrodinger}}\bigr) \bigr\|_{L^1_t L^2_x} \\
                &\qquad \qquad + \Bigl\| (s^\frac12 \bfD)^{(\leq 1)} \Bigl[ i \bfD_t + \bfD^j \bfD_j, (s^\frac12 \bfD)^{(\leq 99)} \Bigr]s^{-\frac12} (s \bfD^{(2)}) \phi \Bigr\|_{L^1_t L^2_x}\\ 
            &\lesssim \ttc(s) \Bigl( 1 + \int_0^T \cC^\sharp(t) \, dt \Bigr) \\
            &\lesssim \ttc(s),
    \end{align*}
using Proposition \ref{prop:non-linear} to estimate the right-hand side of \eqref{Ds-schrodinger} and Lemma \ref{lem:low-high} to estimate the commutator terms, the frequency envelope bound alotted by the energy estimate \eqref{covar-env} and parabolic smoothing \eqref{no-loss}, and the control parameter bootstrap \eqref{BS-contr} in the last line. 

To conclude the proof, it remains to improve the assumed bound on the control parameter \eqref{BS-contr} to \eqref{covar-control}. It suffices to show an estimate of the form 
    \[
         \int_0^T \cC^\sharp (t) \, dt 
            \leq T^{0+} \cdot C\bigl(\| \overline \Gauss^{\mathrm{in}} \|_{L^\infty_x}, \underline s\bigr), 
    \]
as then taking $T_* \ll 1$ completes the proof. Unraveling the definition of the control parameter, we have
    \[
        \int_0^T \cC^\sharp (t) \, dt 
            \lesssim T \cdot \Bigl(  \bigl\| \overline \Gauss^{\mathrm{in}} \bigr\|_{L^\infty_x} + \underline s^{-1} \Bigr) + \bigl\| \phi\bigr\|_{L^2_t L^\infty_x ([0, T])}^2,
    \]
so it remains to estimate the $L^2_t L^\infty_x$-norm of the scalar field. We decompose the analysis between low- and high-frequencies, estimating $|\phi|$ using the fundamental theorem of calculus in $s$, the diamagnetic inequality \eqref{diamagnetic}, and the covariant heat equation \eqref{d-covariant-heat} that
    \begin{align*}
        |\phi| 
            &\leq |\underline \phi| + \int_0^{\underline s} |(s \bfD_s) \phi(s)| \frac{ds}{s} \\
            &\leq |\underline \phi| + \int_0^{\underline s} |(s \bfD^{(2)}) \phi(s)| \, \frac{ds}{s}.
    \end{align*}
For the ``low-frequency" part, we use the Sobolev-type inequality \eqref{inefficient-smooth},
    \begin{align*}
        \bigl\| \underline \phi \bigr\|_{L^2_t L^\infty_x} 
            &\lesssim_{\underline s} T^{\frac12}.
    \end{align*}
For the ``high-frequency part",  
    \begin{align*}
        \int_0^{\underline s} \bigl\| (s \bfD^{(2)}) \phi \bigr\|_{L^2_t L^\infty_x} \frac{ds}{s}
            &\lesssim T^{\frac14} \cdot \int_0^{\underline s} s^{\frac18} \bigl\| s^{\frac18 - \frac12} (s \bfD^{(2)}) \phi \bigr\|_{L^4_{t, x}}^\frac12 \bigl\| s^{\frac18 - \frac12} (s^\frac12 \bfD) (s \bfD^{(2)}) \phi \bigr\|_{L^4_{t, x}}^\frac12 \frac{ds}{s} \\
            &\lesssim_{\underline s} T^{\frac14}
    \end{align*}
using a combination of H\"older's inequality in $t$, the $W^{1, 4} (\R^2) \hookrightarrow L^\infty (\R^2)$ interpolation inequality, and the diamagnetic inequality \eqref{diamagnetic} in the first line, then the $L^4_{t, x}$-Strichartz bound \eqref{covar-str} in the second line. Taking into account the previous estimates, we can improve \eqref{BS-contr} to \eqref{covar-control}, as desired. 
\end{proof}

It remains then to show that the right-hand side of the equation for the covariant Hessian is perturbative,

\begin{proof}[Proof of Proposition \ref{prop:non-linear}]
    Recall that 
    \begin{align*}
        \text{R.H.S.\eqref{Ds-schrodinger}} 
            &= (s \bfD^{(2)}) \phi +  \phi \cdot \phi \cdot (s\bfD^{(2)}) \phi + \phi \cdot (s^\frac12\bfD) \phi \cdot (s^\frac12 \bfD) \phi \\
                &\qquad + \phi \cdot (s^\frac12 \nabla) s^\frac12 F_{tx} + (s^\frac12 \bfD) \phi \cdot s^\frac12 F_{tx} \\
                &\qquad + \phi \cdot (s^\frac12 \nabla)^{(2)} F_{12} + (s^\frac12 \bfD) \phi \cdot (s^\frac12 \nabla) F_{12} + (s \bfD^{(2)}) \phi \cdot F_{12}\\
                &\qquad + (s^\frac12 \bfD)^{(2)} \CSS. 
    \end{align*}
To estimate, the rules of the game are as follows: (\texttt{i}) undifferentiated factors of $\phi$ are estimated in $L^\infty$, (\texttt{ii}) terms with at least two derivatives of the scalar field, one derivative of the electric field, or two derivatives of the magnetic field are estimated in $L^2$ by the frequency envelope. We win by showing each term obeys an estimate of the form 
    \begin{align*}
        \bigl\| s^{-\frac12} (s^\frac12 \bfD)^{(n)} (-) \bigr\|_{L^2_x}
            &\lesssim \cC^\sharp \cdot \ttc(s).
    \end{align*}

Following these rules, we can estimate the terms of the form 
    \[
        (s \bfD^{(2)}) \phi, \qquad \phi \cdot \phi \cdot (s\bfD^{(2)}) \phi,\qquad \phi \cdot (s^\frac12 \nabla) s^\frac12 F_{tx} , \qquad \phi \cdot (s^\frac12 \nabla)^{(2)} F_{12}, \qquad (s \bfD^{(2)}) \phi \cdot F_{12},
    \]
by using the electromagnetic field bounds \eqref{F12-env}, \eqref{Ftx-env} and Proposition \ref{prop:no-loss} to handle the derivatives falling on the ``differentiated'' factors, and parabolic smoothing \eqref{Dphi-strichartz-heat} for when derivatives fall on $\phi$, and the bound \eqref{F12-Linfty} for when derivatives fall on $F_{12}$. The terms of the form  
    \[
        \phi \cdot (s^\frac12\bfD) \phi \cdot (s^\frac12 \bfD) \phi, \qquad (s^\frac12 \bfD) \phi \cdot s^\frac12 F_{tx} , \qquad (s^\frac12 \bfD) \phi \cdot (s^\frac12 \nabla) F_{12},
    \]
do not satisfy rule (\texttt{ii}), though recall that we can use Lemma \ref{lem:bilinear} to ``see'' another derivative, reducing to the previous case. The last term to estimate is the Schr\"odinger tension field -- in this case, it is easy to see that our estimate \eqref{css-env} from the previous section is acceptable. 
\end{proof}

\section{Difference estimates and uniqueness}\label{sec:difference}

To use PDE jargon, we regard the Chern-Simons-Schr\"odinger equation as a dispersive perturbation of a quasi-linear hyperbolic equation, and aim to compare solutions to the Chern-Simons-Schr\"odinger equation in a weaker topology than the energy topology. Somewhat heterodoxically, we only show that the data-to-solution map is weakly H\"older continuous, rather than the weak Lipschitz continuity in a typical well-posedness argument. This is somewhat necessary given the limitations of our approach (see Remark \ref{rem:hihilo-2}). On the other hand, our difference bound has the advantage that it does not rely on any auxiliary bounds -- it propagates at the level of the energy class and thus furnishes unconditional uniqueness of solutions in $C_t ([0, T] \to \frE)$. Furthermore, weak Lipschitz bounds are not, strictly speaking, necessary for closing a well-posedness argument, and indeed our much milder statement will suffice when proving strong continuity of the data-to-solution map in Section \ref{sec:cts}. 

The bulk of the analysis herein is dedicated to studying weak continuity of the data-to-solution map; see Section \ref{subsec:diff}. We conclude in Section \ref{subsec:unique} with some brief remarks concerning unconditional uniqueness. 

\subsection{Difference estimates in caloric-temporal gauge}\label{subsec:diff}

For convenience, let us recall the construction of the $\pzcE^0$-distance functional from Definition \ref{def:dist-2},
    \begin{align*}
        |\dist_{\pzcE^0} (t)|^2
            &\equiv \bigl\| \delta (s^\frac12 \bfD)^{(\leq 2027)} (s^\frac12 \bfD) \phi (t) \bigr\|_{L^2_{\frac{ds}{s}} L^2_x}^2 + \bigl\| \delta \overline\Gauss(t) \bigr\|_{L^2_x}^2 \\
                &\qquad + \bigl\| \delta \underline \bfD^{(\leq 100)} \underline{\bfD \phi} (t) \bigr\|_{L^2_x}^2 + \bigl\| \delta \underline A (t)\bigr\|_{L^4_x}^2 + \bigl\| \delta \underline A (t) \bigr\|_{H^{100}_x}^2 + \bigl\| \underline \phi (t) \bigr\|_{L^\infty_x}^2 .
    \end{align*}
We aim to show that, on the ``bounded" set of data $\frE_\df (\underline s)$, the data-to-solution map for the Chern-Simons-Schr\"odinger equation in caloric-temporal gauge is H\"older continuous with respect to the $\pzcE^0$-distance functional. 

\begin{theorem}[Very weak H\"older continuity of data-to-solution]\label{thm:weak-cty}
    Let $(A, \phi)$, $(A', \phi')$ be smooth solutions on $[0, T] \times \R^2 \times [0, \underline s]$ to the paradifferential formulation of the Chern-Simons-Schr\"odinger equation \eqref{dP-CSS} in  caloric-temporal gauge \eqref{CTC} with initial data in Coulomb gauge \eqref{coulomb-sbar}. If the $\pzcE^0$-distance of the initial data satisfies 
        \[
            |\dist_{\pzcE^0} (0)|^2 
                \leq \tfrac{1}{9000} ,
        \]
    then there exists a time-scale $T_* \equiv T_* (\underline s, \| \overline \Gauss \|_{L^\infty_x}, \dist_{\pzcE^0} (0)) > 0$ satisfying $|T_*|^{-1} = o(1)$ as $\dist_{\pzcE^0} (0) \to 0$, and a constant $\kappa \equiv \kappa(\underline s, \| \Gauss \|_{L^\infty_x}) > 0$ such that, if $T \leq T_*$, we have
        \begin{equation}\label{eq:E0-weak}
            \dist_{\pzcE^0} (t)
                \leq \dist_{\pzcE^0} (0)^{\exp (- \kappa \langle t \rangle^2)}.
        \end{equation}
\end{theorem}

The strategy amounts to observing that Gronwall's inequality can ``tolerate'' a logarithm. More generally, uniqueness continues to hold for ODE which  satisfy Osgood's criterion \cite[Chapter 3.1]{BahouriEtAl2011}.

\begin{lemma}[Non-linear Gronwall's inequality]\label{lem:gronwall}
    Let $f(t) : [0, T] \to [0, \infty)$ be a non-negative continuous function obeying the inequality
        \begin{equation}\label{eq:gronwall-assume}
           f (t)
            \leq f(0) + \int_0^t \cK(t') \cdot \log \bigl( \tfrac{1}{f(t')} \bigr) \cdot f(t') \, dt',
        \end{equation}
    for some non-negative integrable function $\cK : [0, T] \to [0, \infty]$. If $f(0) \leq \tfrac1{9000}$ and $\int_0^T \cK(t) \, dt \leq \log\log (1/f(0))$, then the following double-exponential bound holds for all $t \in [0, T]$, 
        \begin{equation}\label{eq:gronwall}
            f(t) 
                \leq f(0)^{\exp (- \int_0^t \cK(t') \, dt')}.
        \end{equation}
\end{lemma}

\begin{proof}
    Set 
        \[
            I(t) 
                := f(0) + \int_0^t \cK (t')\, \log\bigl( \tfrac1{f(t')} \bigr) f(t') \, dt'.
        \]
    By a regularisation and continuity argument, we may assume that $0 < I(t) \leq e^{-1}$ for all $t$. This function obeys the differential inequality
        \begin{align*}
            \frac{d}{dt} I(t)
                &\leq \cK \, \langle t \rangle \, \log\bigl( \tfrac1{f(t)} \bigr) f(t)  \\
                &\leq \cK(t)\, \log \bigl( \tfrac{1}{I(t)} \bigr) \, I(t),
        \end{align*}
    using the assumption \eqref{gronwall-assume} and monotonicity of $x \log(1/x)$ when $x < e^{-1}$. Repeated application of the chain rule implies the differential inequality, 
        \begin{align*}
            \frac{d}{dt} \log \log \bigl( \tfrac{1}{I(t)} \bigr)
                \geq -\cK(t) .
        \end{align*}
    Integrating, applying the fundamental theorem of calculus, and undoing the logarithms yields the double-exponential-type bound \eqref{gronwall} with $I(t)$ in place of $f(t)$. Of course, $f(t) \leq I(t)$ by assumption \eqref{gronwall-assume}, and $f(0) = I(0)$ by construction, so we conclude the desired bound for $f(t)$ follows as well. The restriction the time-scale $T$ arises from our assumption on the upper bound $I(t) \leq e^{-1}$, so setting 
        \[
            f(0)^{\exp(- \int_0^t \cK(t') \, dt')}
                = e^{-1}
        \]
    and rearranging gives the time-scale. 
\end{proof}

Lemma \ref{lem:gronwall} reduces the task of proving Theorem \ref{thm:weak-cty} to showing the differential inequality\footnote{Strictly speaking, this will not quite be true, since we will have a differential inequality for $\tfrac{d}{dt} |\delta\underline \phi(x)|^2$ uniformly in $x \in \R^2$ rather than $\tfrac{d}{dt} \|\delta  \underline \phi \|_{L^\infty_x}^2$. Nevertheless, the argument below carries through with minor modifications -- the integral inequality \eqref{gronwall-assume} is robust.}
    \begin{equation}\label{eq:diff-ineq}
        \frac{d}{dt} |\dist_{\pzcE^0} (t)|^2 
            \lesssim_{\underline s} \Bigl( \langle t \rangle + \bigl\| \overline \Gauss \bigr\|_{L^\infty_x} \Bigr) \cdot \log \Bigl( \frac{\underline s}{s_*} \Bigr) \cdot \Bigl( |s_*| + |\dist_{\pzcE^0} (t)|^2 \Bigr),
    \end{equation}
for all $0 < s_* < \tfrac1{1000} \underline s$. Indeed, formally taking $s_* := |\dist_{\pzcE^0}|^2$ and integrating the ensuing differential inequality gives an integral inequality of the form \eqref{gronwall-assume} for the $\pzcE^0$-distance functional. 

It remains to establish (at least, morally) the differential inequality \eqref{diff-ineq}; we give a full proof of Theorem \ref{thm:weak-cty} at the end of the section. To propagate the differences of $\bfD \phi$ and its derivatives, we argue by the energy method, 

\begin{lemma}[Difference equations of motion, I]
    Let $(A, \phi), (A', \phi')$ be configurations on $[0, T] \times \R^2 \times [0, \underline s]$ to the paradifferential formulation of the Chern-Simons-Schr\"odinger equation \eqref{dP-CSS} in caloric-temporal gauge \eqref{CTC}. Then, for each positive integer $n \in \N$, we have the magnetic Schr\"odinger equations,
        \begin{equation}\label{eq:Dphi-diff-t}
        \begin{split}
            \bigl( i \bfD_t + \bfD^j \bfD_j \bigr) \delta (s^\frac12 \bfD)^{(n)} \phi
                &= \bigl[ \delta, i \bfD_t + \bfD^j \bfD_j \bigr] (s^\frac12 \bfD)^{(n)} \phi \\
                    &\qquad + \delta (s^\frac12 \bfD)^{(n)} \phi + \sum_{a + b + c = n} (s^\frac12 \bfD)^{(a)} \phi \cdot (s^\frac12 \bfD)^{(b)} \phi \cdot \delta (s^\frac12 \bfD)^{(c)} \phi \\
                    &\qquad + \delta \sum_{a + b = n - 1} (s^\frac12 \bfD)^{(a)} \phi \cdot (s^\frac12 \nabla)^{(b)} s^\frac12 F_{tx} \\
                    &\qquad + \delta \sum_{a + b = n} (s^\frac12 \bfD)^{(a)} \phi \cdot (s^\frac12 \nabla)^{(b)} F_{12} \\
                    &\qquad + \delta (s^\frac12 \bfD)^{(n)} \CSS .
        \end{split}
        \end{equation}
\end{lemma}

\begin{proof}
    Commuting $(s^\frac12 \bfD)^{(n)}$ through the definition of the Schr\"odinger tension field,
        \begin{align*}
            (i \bfD_t + \bfD^j \bfD_j) (s^\frac12 \bfD)^{(n)} \phi 
                &= \bigl[i \bfD_t + \bfD^j \bfD_j , (s^\frac12 \bfD)^{(n)}  \bigr] \phi \\
                &\qquad + \sum_{a + b + c = n} (s^\frac12 \bfD)^{(a)} \phi \cdot (s^\frac12 \bfD)^{(b)} \phi \cdot (s^\frac12 \bfD)^{(c)}\phi +  (s^\frac12 \bfD)^{(n)} \CSS\\
                &= \sum_{a + b = n - 1} (s^\frac12 \nabla)^{(a)} s^\frac12 F_{t x} \cdot (s^\frac12\bfD)^{(b)}\phi + \sum_{a + b = n} (s^\frac12\nabla)^{(a)} F_{12} \cdot (s^\frac12 \bfD)^{(b)} \phi \\
                &\qquad  + \sum_{a + b + c = n} (s^\frac12 \bfD)^{(a)} \phi \cdot (s^\frac12 \bfD)^{(b)} \phi \cdot (s^\frac12 \bfD)^{(c)} \phi +  (s^\frac12 \bfD)^{(n)} \CSS,
        \end{align*}
    applying the product rule \eqref{product} in the first line, and recalling the commutator identity \eqref{low-high-commutator} in the second. We then commute $\delta$ through the equation. 
\end{proof}

\begin{proposition}
    Let $(A, \phi)$, $(A', \phi')$ be smooth solutions on $[0, T] \times \R^2 \times [0, \underline s]$ to the paradifferential formulation of the Chern-Simons-Schr\"odinger equation \eqref{dP-CSS} in caloric-temporal gauge \eqref{CTC} with initial data in Coulomb gauge \eqref{coulomb-sbar}. Then, for any $0 < s_* < \tfrac1{1000} \underline s$, we have: 
    \begin{itemize}
        \item \textup{(Propagation of high-frequency differences).}
             \begin{equation}\label{eq:high-diff}
                \frac{d}{dt}\bigl\| \delta (s^\frac12 \bfD)^{(\leq 2027)} (s^\frac12 \bfD) \phi \bigr\|_{L^2_{\frac{ds}{s}}L^2_x}^2  
                    \lesssim_{\underline s} \Bigl( \bigl\| \overline \Gauss \bigr\|_{L^\infty_x} + \langle t \rangle \Bigr) \cdot \log\Bigl( \frac{\underline s}{s_*} \Bigr) \cdot\Bigl( |s_*| + |\dist_{\pzcE^0}|^2 \Bigr).
            \end{equation} 
        \item \textup{(Propagation of low-frequency differences).}
            \begin{equation}\label{eq:low-diff}
                \frac{d}{dt} \bigl\| \delta \underline\bfD^{(\leq 100)} \underline{\bfD \phi} \bigr\|_{L^2_x}^2
                    \lesssim_{\underline s} \langle t \rangle \cdot \log\Bigl( \frac{\underline s}{s_*} \Bigr) \cdot\Bigl( |s_*| + |\dist_{\pzcE^0}|^2 \Bigr).
            \end{equation}
    \end{itemize}

\end{proposition}

\begin{proof}[Proof of high-frequency difference bound \eqref{high-diff}]
    Applying the abstract energy estimate \eqref{abstract-energy} to the Schr\"odinger equation for the differences \eqref{Dphi-diff-t} and square-summing in $s$, our task reduces to showing 
        \begin{align*}
            \bigl\| \text{R.H.S.\eqref{Dphi-diff-t}} \bigr\|_{L^2_{\frac{ds}{s}} L^2_x}
                \lesssim \Bigl(\langle t \rangle + \bigl\| \overline \Gauss \bigr\|_{L^\infty_x} \Bigr) \cdot \log \Bigl( \frac{\underline s}{s_*} \Bigr)  \cdot \Bigl( |s_*|^\frac12 + \dist_{\pzcE^0} \Bigr),
        \end{align*}
    for each $n = 1, \dots, 2028$. 

    Expanding out the commutator using the identities \eqref{commute-diff} and \eqref{commute-diff-laplace}, we see the terms
        \[
            \delta A_x \cdot (s^\frac12 \bfD)^{(n)} \bfD \phi, \qquad \delta A_t \cdot (s^\frac12 \bfD)^{(n)} \phi , \qquad \delta \partial^\ell A_\ell \cdot (s^\frac12 \bfD)^{(n)} \phi, \qquad \delta A_x \cdot A_x \cdot (s^\frac12 \bfD)^{(n)} \phi. 
        \]
    The first term is the most dangerous -- the extra derivative falling on $(s^\frac12 \bfD)^{(n)} \phi$ forces us to place it in $L^2_{ds/s} L^2_x$ and thereby the difference of the magnetic potentials in $L^\infty_x$ in order to control it at energy regularity. For the remaining terms, we can use the estimate $L^\infty_s L^2_x \times L^2_{ds/s} L^\infty_x \to L^2_{ds/s} L^2_{x}$. Altogether,
        \begin{align*}
            \bigl\| [\delta, i \bfD_t + \bfD^j \bfD_j] (s^\frac12 \bfD)^{(n)} \phi \bigr\|_{L^2_{\frac{ds}{s}} L^2_x}  
                &\lesssim \bigl\| \delta A_x \bigr\|_{L^\infty_{s, x}} \bigl\| (s^\frac12 \bfD)^{(n)} \bfD \phi \bigr\|_{L^2_{\frac{ds}{s}} L^2_x} \\
                    &+ \Bigl( \bigl\| \delta A_t \bigr\|_{L^\infty_s L^2_x} + \bigl\| \delta \partial^\ell \underline A_\ell \bigr\|_{L^2_x} + \bigl\| \delta A_x \bigr\|_{L^\infty_{s} L^4_x} \bigl\| A_x \bigr\|_{L^\infty_{s} L^4_x}\Bigr) \bigl\| (s^\frac12 \bfD)^{(n)} \phi \bigr\|_{L^2_{\frac{ds}{s}} L^\infty_x}\\
                &\lesssim_{\underline s} \langle t \rangle \cdot \log\Bigl( \frac{\underline s}{s_*}\Bigr)  \cdot \Bigl( |s_*|^\frac12 + \dist_{\pzcE^0} \Bigr),
        \end{align*}
    estimating the differences of the magnetic potential by \eqref{Ax-diff}, \eqref{Ax-diff-linfty}, the difference of the electric potential by \eqref{At-diff}, the growth of the magnetic potential by \eqref{Ax-ctc-cf}-\eqref{Ax-ctc-df}, and the derivatives of $\phi$ by parabolic smoothing \eqref{Dphi-sobolev} and \eqref{Dphi-strichartz-heat}. Thus, these terms are acceptable. 
    
    For the remaining terms on the right-hand side of \eqref{Dphi-diff-t}, we interpolate between high-frequency decay and low-frequency summability a l\'a Lemma \ref{lem:low-high-inter}, aiming to show that they satisfy estimates of the form
        \begin{align}
            \bigl\| (-)(s)  \bigr\|_{L^2_x}
                &\lesssim_{\underline s} \Bigl( 1 + \bigl\| \overline \Gauss \bigr\|_{L^\infty_x} \Bigr) \cdot s^\frac12 \cdot \log \Bigl( \frac{\underline s}{s} \Bigr) \label{eq:high-decay},\\
            \bigl\| - \bigr\|_{L^2_{\frac{ds}{s}} L^2_x ([s_*, \underline s])}
                &\lesssim_{\underline s} \Bigl(\langle t \rangle + \bigl\| \overline \Gauss \bigr\|_{L^\infty_x} \Bigr) \cdot \log \Bigl( \frac{\underline s}{s_*} \Bigr)  \cdot \Bigl( |s_*|^\frac12 + \dist_{\pzcE^0} \Bigr) \label{eq:low-sum} .
        \end{align}
    These terms more or less are exactly ones encountered in the proof of Proposition \ref{prop:non-linear}, so they obey a bound of the form 
        \[  
             \bigl\| s^{-\frac12} (-)(s)  \bigr\|_{L^2_x}
                \lesssim_{\underline s} \Bigl( 1 + \| \phi(s)\|_{L^\infty_x}^2 + \bigl\| \overline \Gauss \bigr\|_{L^\infty_x} \Bigr)  .
        \]
    Estimating $\phi(s)$ by the Brezis-Gallou\"et-type inequality \eqref{brezis} yields \eqref{high-decay}. It remains to show \eqref{low-sum}

    The contribution of the Schr\"odinger tension field difference may be estimated by \eqref{CSS-diff}. For the remaining terms, we divide into three scenarios depending on where the difference operator $\delta$ falls: (\texttt{i}) differences of $\bfD \phi$ and its derivatives, (\texttt{ii}) differences of the undifferentiated scalar field $\phi$, and (\texttt{iii}) differences of the electromagnetic field $F_{\mu\nu}$. The key thing to keep track of is the half-powers of $\log(\underline s/s_*)$ appearing. 

    When (\texttt{i}) the difference falls on $\bfD\phi$ and its derivatives, i.e. the terms of the form
        \[
            \delta (s^\frac12 \bfD)^{(a)} \phi \cdot (s^\frac12 \bfD)^{(b)} \phi \cdot (s^\frac12 \bfD)^{(c)} \phi, \qquad \delta (s^\frac12 \bfD)^{(a)} \phi \cdot (s^\frac12 \nabla)^{(b)} F_{12}, \qquad \delta (s^\frac12 \bfD)^{(a)} \phi \cdot (s^\frac12 \nabla)^{(b)} F_{tx},
        \]
    where $a \neq 0$ and $a \leq n$, we place the difference in $L^2_{ds/s} L^2_x$ and the remaining factors in $L^\infty_{s, x}$,
        \begin{align*}
            \bigl\|  \delta (s^\frac12 \bfD)^{(a)} \phi \cdot (-) \bigr\|_{L^2_{\frac{ds}{s}} L^2_x ([s_*, \underline s])}
                &\lesssim  \bigl\|  \delta (s^\frac12 \bfD)^{(a)} \phi \bigr\|_{L^2_{\frac{ds}{s}} L^2_x}  \bigl\| (-) \bigr\|_{L^\infty_{s, x} ([s_*, \underline s])} \\
                &\lesssim \dist_{\pzcE^0} \cdot   \bigl\| - \bigr\|_{L^\infty_{s, x} ([s_*, \underline s])} .
        \end{align*}    
    Differentiated factors of $\phi$ are estimated by parabolic smoothing \eqref{Dphi-sobolev}, and the electromagnetic field by \eqref{F12-Linfty}, \eqref{Fmunu-Lp}. We will see at most 2 factors of $\phi(s)$ in $L^\infty_x$, so applying the Brezis-Gallou\"et-type inequality \eqref{brezis},
        \begin{align*}
            \bigl\| (s^\frac12 \bfD)^{(b)} \phi \cdot (s^\frac12 \bfD)^{(c)} \phi \bigr\|_{L^\infty_{s, x} ([s_*, \underline s])} &+  \bigl\| (s^\frac12 \nabla)^{(b)} F_{12} \bigr\|_{L^\infty_{s, x} ([s_*, \underline s])}  +  \bigl\| (s^\frac12 \nabla)^{(b)} F_{tx}\bigr\|_{L^\infty_{s, x} ([s_*, \underline s])} \\
                &\lesssim_{\underline s} 1 + \bigl\| \overline \Gauss \bigr\|_{L^\infty_x} + \bigl\| \phi(s_*)\bigr\|_{L^\infty_{x}}^2\\
                &\lesssim_{\underline s} \bigl\| \overline \Gauss \bigr\|_{L^\infty_x}  + \log \Bigl( \frac{\underline s}{s_*} \Bigr).
        \end{align*}
    Altogether, these terms obey a bound of the form \eqref{low-sum}. 

    When (\texttt{i}) the difference falls on the undifferentiated scalar field, i.e. the terms of the form
        \[
            \delta \phi \cdot (s^\frac12 \bfD)^{(a)} \phi \cdot (s^\frac12 \bfD)^{(b)} \phi, \qquad \delta \phi \cdot (s^\frac12 \nabla)^{(n)} F_{12}, \qquad \delta \phi \cdot (s^\frac12 \nabla)^{(n - 1)} s^\frac12 F_{tx},
        \]
    we place the difference in $L^\infty_s (L^2 + L^\infty)_x$, and then use \eqref{phi-diff}, incurring one factor of $\log^{1/2} (\underline s/s_*)$, and the remaining factors in $L^2_{ds/s} (L^\infty \cap L^2)_x$,
        \begin{align*}
            \bigl\| \delta \phi \cdot (-) \bigr\|_{L^2_{\frac{ds}{s}} L^2_x}   
                &\lesssim \bigl\| \delta \phi \bigr\|_{L^\infty_s (L^2 + L^\infty)_x ([s_*, \underline s])}\bigl\| - \bigr\|_{L^2_{\frac{ds}{s}} (L^\infty \cap L^2)_x ([s_*, \underline s])} \\
                &\lesssim_{\underline s}  \log^\frac12 \Bigl( \frac{\underline s}{s_*} \Bigr)  \cdot \Bigl( |s_*|^\frac12 + \dist_{\pzcE^0} \Bigr) \cdot \bigl\| - \bigr\|_{L^2_{\frac{ds}{s}} (L^\infty \cap L^2)_x ([s_*, \underline s])}.
        \end{align*}
    There is more room when placing the remaining factors in $L^2_{ds/s} L^2_x$, so let us consider the case $L^2_{ds/s} L^\infty_x$. In the first term, we will see at most one undifferentiated factor of $\phi$, which we handle by the Brezis-Gallou\"et-type bound \eqref{brezis}, thereby incurring another factor of $\log^{1/2}(\underline s/s_*)$. The differentiated factors can be estimated by parabolic smoothing \eqref{Dphi-sobolev} and parabolic Strichartz \eqref{Dphi-strichartz-heat}. Thus
        \begin{align*}
            \bigl\| (s^\frac12 \bfD)^{(a)} \phi \cdot (s^\frac12 \bfD)^{(b)} \phi\bigr\|_{L^2_{\frac{ds}{s}} L^\infty_x ([s_*, \underline s])} 
                &\lesssim \bigl\| \phi (s_*) \bigr\|_{L^\infty} \bigl\| (s^\frac12 \bfD)^{(n)} \phi \bigr\|_{L^2_{\frac{ds}{s}} L^\infty_x} \\
                    &\qquad + \sum_{a + b  \neq 0}\bigl\| (s^\frac12 \bfD)^{(a)} \phi \bigr\|_{L^\infty_{s, x}} \bigl\| (s^\frac12 \bfD)^{(b)} \phi\bigr\|_{L^2_{\frac{ds}{s}} L^\infty_x} \\
            &\lesssim_{\underline s} \log^\frac12 \Bigl( \frac{\underline s}{s_*} \Bigr).
        \end{align*}
    Altogether, we see at most 2 factors of $\log^{1/2} (\underline s/s_*)$, which is acceptable for a bound of the form \eqref{low-sum}. For the second and third terms, we can use parabolic Strichartz \eqref{parabolic-strichartz} for the heat equation \eqref{Fmunu-heat} to estimate the electromagnetic field by $\nabla F_{12}(s_*)$ and $F_{tx}(s_*)$ in $L^2_x$, which we can bound by \eqref{nabla-F12-Lp} and \eqref{Fmunu-Lp} respectively. These incur one factor of $\phi(s_*)$ in $L^\infty$, and thereby one factor of $\log^{1/2}(\underline s/s_*)$ when using \eqref{brezis}. Thus,
        \begin{align*}
            \bigl\| (s^\frac12 \nabla)^{(n)} F_{12} \bigr\|_{L^2_{\frac{ds}{s}} L^\infty_x ([s_*, \underline s])} + \bigl\| (s^\frac12 \nabla)^{(n - 1)} s^\frac12 F_{tx} \bigr\|_{L^2_{\frac{ds}{s}} L^\infty_x ([s_*, \underline s])}  
                &\lesssim \bigl\| \nabla F_{12} (s_*) \bigr\|_{L^2_x} + \bigl\| F_{tx} (s_*) \bigr\|_{L^2_x}\\
                &\lesssim_{\underline s} \| \phi(s_*)\|_{L^\infty_x} \\
                &\lesssim \log^\frac12 \Bigl( \frac{\underline s}{s_*} \Bigr).
        \end{align*}
    Altogether, these two terms will see only 2 factors of $\log^{1/2}(\underline s/s_*)$, which is acceptable for a bound of the form \eqref{low-sum}. 
  
    Finally, when (\texttt{iii}) the difference falls on the electromagnetic field, i.e. the terms of the form
        \[
            (s^\frac12 \bfD)^{(a)} \phi \cdot (s^\frac12 \nabla)^{(b)} \delta F_{12}, \qquad (s^\frac12 \bfD)^{(a)} \phi \cdot (s^\frac12 \nabla)^{(b)} s^\frac12 \delta F_{tx},
        \]
    we place the difference in $L^2_{ds/s} L^2_x$, and the remaining factors in $L^\infty_{s,x}$. Estimating the magnetic field by \eqref{F12-diff} and electric field by \eqref{Ftx-diff} incurs a $\log^{1/2} (\underline s/s_*)$ factor, so
        \begin{align*}
            \bigl\| (s^\frac12 \bfD)^{(a)} \phi \cdot (-)\bigr\|_{L^2_{\frac{ds}{s}} L^2_x ([s_*, \underline s])} 
                &\lesssim \bigl\| (s^\frac12 \bfD)^{(a)} \phi \bigr\|_{L^\infty_{s, x}([s_*, \underline s])} \bigl\| (-)\bigr\|_{L^2_{\frac{ds}{s}} L^2_x}\\
                &\lesssim \bigl\| (s^\frac12 \bfD)^{(a)} \phi \bigr\|_{L^\infty_{s, x}([s_*, \underline s])} \cdot \log^\frac12 \Bigl( \frac{\underline s}{s_*} \Bigr) \cdot \Bigl( |s_*|^\frac12 + \dist_{\pzcE^0} \Bigr). 
        \end{align*}
    The most dangerous case is $a = 0$, where we need to estimate the undifferentiated scalar field by  Brezis-Gallou\"et \eqref{brezis}, incurring another factor of $\log^{1/2} (\underline s/s_*)$. If $a \neq 0$, we can use parabolic smoothing \eqref{Dphi-sobolev},
        \begin{align*}
            \bigl\| (s^\frac12 \bfD)^{(a)} \phi \bigr\|_{L^\infty_{s, x}([s_*, \underline s])}
                &\lesssim_{\underline s}  \log^\frac12 \Bigl( \frac{\underline s}{s_*} \Bigr). 
        \end{align*}
    Altogether, we see at most two factors of $\log^{1/2} (\underline s/s_*)$, which is acceptable for \eqref{low-sum}. 
\end{proof}

\begin{proof}[Proof of low-frequency difference bound \eqref{low-diff}]
    The proof proceeds more or less the same as that of \eqref{high-diff}, however the argument is strictly easier since we are conducting the analysis at a fixed low-frequency $s = \underline s$. We leave the details to the reader. 
\end{proof}

We view the configuration $(\underline A, \underline \phi)$ at $s = \underline s$ as evolving by transport in $t$. 

\begin{lemma}[Difference equations of motion, II]
    Let $(A, \phi), (A', \phi')$ be configurations on $[0, T] \times \R^2 \times [0, \underline s]$ to the paradifferential formulation of the Chern-Simons-Schr\"odinger equation \eqref{dP-CSS} in caloric-temporal gauge \eqref{CTC}. Then $(\underline A, \underline \phi)$ satisfy 
        \begin{itemize}
            \item \textup{(Transport equation for $\underline A$).}
                \begin{equation}\label{eq:transport-sbar2}
                    \partial_t \, \delta \underline A_x 
                        = \delta \underline F_{tx}. 
                \end{equation}
                
            \item \textup{(Transport equation for $\underline \phi$).}
                \begin{equation}\label{eq:transport-sbar}
                    \partial_t \, \delta \underline \phi 
                        = \delta \underline \phi + \underline \phi \cdot \underline \phi \cdot \delta \underline \phi + \delta \underline \bfD^{(2)} \underline \phi + \delta \underline \CSS.
                    \end{equation}
        \end{itemize}
\end{lemma}

\begin{proof}
    These are simply the definitions of the electric field and Schr\"odinger tension field written at $s = \underline s$ in temporal gauge $\underline A_t = 0$. 
\end{proof}

We already propagated difference bounds for $\underline A_x$ using \eqref{transport-sbar2} in \eqref{Ax-diff-trans}-\eqref{Ax-diff-trans-2}, so it remains to propagate difference bounds for $\underline \phi$ using \eqref{transport-sbar}. 

\begin{proposition}[Propagation of $\delta \underline \phi$ in $L^\infty$]
     Let $(A, \phi)$, $(A', \phi')$ be smooth solutions on $[0, T] \times \R^2 \times [0, \underline s]$ to the paradifferential formulation of the Chern-Simons-Schr\"odinger equation \eqref{dP-CSS} in caloric-temporal gauge \eqref{CTC} with initial data in Coulomb gauge \eqref{coulomb-sbar}. Then, for any $0 < s_* < \tfrac1{1000} \underline s$ and $x \in \R^2$, we have
            \begin{equation}\label{eq:Linfty-diff}
            \frac{d}{dt}|\delta \underline \phi(x)|^2
                \lesssim_{\underline s} \langle t \rangle \cdot \log\Bigl( \frac{\underline s}{s_*} \Bigr) \cdot\Bigl( |s_*| + |\dist_{\pzcE^0}|^2 \Bigr).
             \end{equation}
\end{proposition}

\begin{proof}[Proof of $\delta \underline \phi$ bound \eqref{Linfty-diff}]
    Estimating the right-hand side of the transport equation \eqref{transport-sbar}, recall the scalar field $\underline \phi$ is uniformly bounded thanks to \eqref{brezis}, so it remains to estimate the last two terms, namely 
        \[
            \delta \underline\bfD^{(2)} \underline\phi, \qquad \delta \underline \CSS.
        \]  
    Writing $\sfU \in \{ \bfD^{(2)} \phi, \CSS \}$, 
    \begin{align*}
        \bigl\| \delta \underline\sfU \bigr\|_{L^\infty_x}
            &\lesssim \bigl\| \delta \underline \sfU \bigr\|_{L^2_x} + \bigl\| \underline \bfD^{(2)} \delta \underline \sfU \bigr\|_{L^2_x} \\
            &\lesssim \bigl\| \delta \underline \bfD^{(\leq 2)}\underline \sfU \bigr\|_{L^2_x} + \bigl\| [\delta, \underline \bfD^{(2)}]\underline \sfU \bigr\|_{L^2_x}  \\
            &\lesssim \bigl\| \delta \underline \bfD^{(\leq 2)}\underline \sfU \bigr\|_{L^2_x} + \bigl\| \delta \underline A_x \bigr\|_{L^4_x} \bigl\| \underline{\bfD \sfU} \bigr\|_{L^4_x} + \Bigl( \bigl\| \delta \nabla \underline A_x \bigr\|_{L^2_x} +  \bigl\| \delta \underline A_x \bigr\|_{L^4_x} \bigl\|  \underline A_x \bigr\|_{L^4_x}  \Bigr) \bigl\| \underline \sfU \bigr\|_{L^\infty_x} \\
            &\lesssim_{\underline s} \Bigl(1 + \bigl\| \underline A_x \bigr\|_{L^4_x}\Bigr) \cdot \log^\frac12 \Bigl( \frac{\underline s}{s_*} \Bigr) \cdot \Bigl( |s_*|^\frac12 + \dist_{\pzcE^0} \Bigr)
    \end{align*}
    by Sobolev embedding \eqref{GNinfty} in the first line, commuting $\delta$ through the Hessian in the second and third lines and then applying H\"older appropriately, estimating $\underline{\sfU}$ using parabolic smoothing \eqref{Dphi-sobolev} and the Schr\"odinger tension field bound \eqref{css-env}, and finally estimating differences of the Schr\"odinger tension field by \eqref{CSS-diff}. The magnetic potential grows at the rate $\langle t \rangle$ by \eqref{Ax-ctc}, completing the proof. 
\end{proof}

\begin{proof}[Proof of Theorem \ref{thm:weak-cty}]
    By the analyses from Sections \ref{subsec:LP-topology} and \ref{subsec:weak}, the $\pzcE^0$-distance functional is continuous along the flows of two smooth solutions to the Chern-Simons-Schr\"odinger equation. Thus, we can invoke continuity argument and assume that the $\pzcE^0$-distance remains small, e.g.
    \[
        \sup_{t \in [0, T]} |\dist_{\pzcE^0} (t)|^2
            \leq \tfrac{1}{5000}.
    \]
    Fix $\epsilon \ll \tfrac{1}{5000}$, then for each $t \in [0, T]$, we can choose 
        \[
            s_*
                := \underline s \cdot \Bigl( |\dist_{\pzcE^0}|^2 + \epsilon \Bigr).
        \]
    Inserting this into the right-hand sides of the differential inequalities for the differences of magnetic derivatives \eqref{Dphi-s-diff-1}, \eqref{Dphi-s-diff-2}, the differences of $(\underline A, \underline \phi)$ \eqref{Ax-diff-trans}, \eqref{Ax-diff-trans-2}, \eqref{Linfty-diff}, recalling the trivial differential inequality for differences of the Gauss tension field, and integrating in $t$, we arrive at the integral inequality \eqref{gronwall-assume} with
        \begin{align*}
            f(t) 
                &:= |\dist_{\pzcE^0}|^2 + \epsilon ,\\
            \cK(t) 
                &:= \kappa (\underline s, \| \overline \Gauss \|_{L^\infty_x}) \cdot \langle t \rangle,
        \end{align*}
    where $\kappa \equiv \kappa (\underline s, \| \overline \Gauss \|_{L^\infty_x})  > 0$ is some large constant. Applying the non-linear Gronwall inequality, Lemma \ref{lem:gronwall}, and taking $\epsilon \to 0$ gives the double-exponential bound. Clearly the time-scale for which it holds may be extended to $t \to \infty$ as one takes the initial $\pzcE^0$-metric to zero. This completes the proof. 
\end{proof}

\begin{numremark}\label{rem:hihilo-2}
    The logarithmic losses in the preceding argument arise from two types of terms, namely ``low-high'' and ``high-high-to-low''. For the low-high terms, e.g. 
        \[  
            \phi \cdot \phi \cdot \delta (s^\frac12 \bfD) \phi, \qquad \delta \phi \cdot \phi \cdot (s^\frac12 \bfD) \phi,
        \]
    one can avoid the logarithmic divergences by replacing the Brezis-Gallou\"et bounds with Strichartz bounds. When the difference falls on the electric field, e.g.  
        \[
            \phi \cdot \delta s^\frac12 F_{tx},
        \]
    there is an unavoidable logarithmic divergence when estimating a $\text{high} \times \text{high} \to \text{low}$-type interaction arising from the Amp\'ere tension field, as there is no room to interpolate using Strichartz, see Remark \ref{rem:hihilo}. Thus, it does not seem possible to prove a weak Lipschitz bound relying solely on the tools herein. Nevertheless, we expect the issue is merely technical; \eqref{CSS-amp-schro} should resemble, at least in some high-frequency limit sense, the Jackiw-Pi model \eqref{JP-CSS}, where well-posedness holds using a contraction mapping argument in $H^{0+}$ by exploiting bilinear estimates and the null structure in the equation \cite{LiuEtAl2014}.
\end{numremark}

\subsection{Unconditional uniqueness of finite-energy solutions}\label{subsec:unique}

After some careful bookkeeping, the astute reader will realise that smoothness is far more regularity than needed to prove Theorem \ref{thm:weak-cty}. Indeed, the weak difference bound \eqref{E0-weak} propagates at the level of the energy regularity, and essentially solely relies on energy or transport-type estimates from Sections \ref{sec:littlewoodpaley}, \ref{sec:paradiff}, \ref{sec:tension} and \ref{sec:gauge}. The only place where we needed smoothness was Section \ref{sec:energy}, where we used it to bootstrap the control parameter bound -- here, we substituted the control parameter bound for the Brezis-Gallou\"et-type bound, which holds for fairly rough configurations. 

These observations lead us to claim a Kato-style unconditional uniqueness for finite-energy solutions, 

\begin{theorem}[Unconditional uniqueness in $\frE$]\label{thm:unique}
    Let $(A, \phi)$, $(A', \phi') \in C_t ([0, T] \to \frE)$ solve the Chern-Simons-Schr\"odinger equation \eqref{CSS-amp-schro} in DeTurck gauge \eqref{DeTurck} with initial data in Coulomb gauge \eqref{coulomb-data}. If the solutions admit the same initial data, i.e. $(A, \phi)(0) = (A', \phi')(0)$, then they coincide, i.e. $(A, \phi)(t) = (A', \phi')(t)$ for all $t \in [0, T]$. 
\end{theorem}

\begin{remark}
    Theorem \ref{thm:unique} is analogous to the $H^1$-unconditional uniqueness result of Huh \cite{Huh2013} for Jackiw-Pi's Chern-Simons-Schr\"odinger equation \eqref{JP-CSS} in Coulomb gauge. 
\end{remark}

This furnishes the unconditional uniqueness part of Theorem \ref{thm:gwp}. To prove Theorem \ref{thm:unique}, there are two technical points we have to verify: first, the prior analysis was conducted on the paradifferential formulation of the Chern-Simons-Schr\"odinger equation in caloric-temporal gauge, so we will need to extend the solution to $[0, T] \times \R^2 \times [0, \underline s]$ and make a gauge-transform to make use of it, and second, we need to verify that the $\pzcE^0$-distance functional is continuous to justify the use of Gronwall. 

We give a brief sketch, and leave the details to the reader. Set $\underline s > 0$ such that\footnote[1]{We take supremum in $t$ here because \textit{a priori} these solutions need not conserve energy. } 
    \[
        \sup_{t \in [0, T]} \Bigl( 1 + \cE_{\mathrm{AH}} [A,\phi] + \| \nabla \Gauss \|_{L^2}^2 \Bigr) + \sup_{t \in [0, T]} \Bigl( 1 + \cE_{\mathrm{AH}} [A',\phi'] + \| \nabla \Gauss' \|_{L^2}^2 \Bigr)
            \ll \underline s^{-1},
    \]
and, following Proposition \ref{prop:dP-exist}, we extend $(A, \phi)$ and $(A', \phi')$ to solutions to \eqref{dP-CSS} on $[0, T] \times \R^2 \times [0, \underline s]$ by solving the dynamic covariant heat equation \eqref{d-covariant-heat}-\eqref{d-covariant-heat2} in the DeTurck gauge,
    \[
        A_s
            = \partial^\ell A_\ell. 
    \]
This yields solutions to \eqref{dP-CSS} which are in ``DeTurck-DeTurck gauge''. By Proposition \ref{thm:topology}, the ensuing configurations form continuous flows into the space of caloric extensions $C_t ([0, T] \to \pzcE^1 ([0, \underline s]))$. Observe that the electromagnetic potential $A_\mu$ solves the linear heat equation in $s$ in this gauge, so one can follow the argument of Lemma \ref{lem:weak} to show that the $\pzcE^0$-distance functional is continuous in $t$ for these solutions. Following Proposition \ref{prop:ctc-exists}, we convert to solutions $(\widetilde A, \widetilde \phi)$ and $(\widetilde A', \widetilde \phi')$ in caloric-temporal gauge by making the gauge-transformation 
    \[
        \chi 
            := - (\partial_s)^{-1} A_s - (\partial_t)^{-1} \underline A_t. 
    \]
It remains to verify continuity of the $\pzcE^0$-distance between $(\widetilde A, \widetilde \phi)$ and $(\widetilde A', \widetilde \phi')$. This amounts to showing that $\chi$ is sufficiently regular. For the contribution of the electric potential at $s = \underline s$, recall that $A_t$ solves the linear heat equation in $s$ under DeTurck gauge in $s$, and the heat-temporal initial data is given by $\partial^j \overline A_j$ under DeTurck gauge in $t$, so this part of the transformation is in fact smooth, $(\partial_t)^{-1} \underline A_t \in C^\infty_t H^\infty_x ([0, T] \times\R^2)$ by parabolic smoothing. For the contribution of $A_s$, we rewrite under DeTurck gauge in $s$, leveraging almost orthogonality and parabolic smoothing, 
    \begin{align*}
        \| \nabla^{(2)} (\partial_s)^{-1} A_s \|_{L^2_x}
            &\lesssim \Bigl\| \int_s^{\underline s} ( {s'}^\frac12 \nabla)^{(2)} \partial^\ell A_\ell (s')\frac{ds'}{s'} \Bigr\|_{L^2_x} \\
            &\lesssim \bigl\| ( s^\frac12 \nabla)^{(\leq 2)} ( {s'}^\frac12 \nabla) \partial^\ell A_\ell  \bigr\|_{L^2_{\frac{ds}{s}} L^2_x ([s, \underline s])} \\
            &\lesssim \| \nabla \overline A_x \|_{L^2_x}.
    \end{align*}
An easier argument holds at bottom order. Altogether, $\chi \in C^0_t H^2_x ([0, T] \times \R^2)$, so by a similar argument as in Lemma \ref{lem:gauge-diff}, we can conclude continuity of the $\pzcE^0$-distance functional. This is enough regularity to justify Theorem \ref{thm:weak-cty}, so we can conclude uniquenes of finite-energy solutions of the Chern-Simons-Schr\"odinger equation in DeTurck gauge. 

\section{Continuity of the flow map}\label{sec:cts}

To complete the proof of Theorem \ref{thm:gwp}, it remains to construct solutions $(A, \phi) \in C_t ([0, \infty) \to \frE)$ to the Chern-Simons-Schr\"odinger equation in DeTurck gauge with initial data $(A^{\mathrm{in}}, \phi^{\mathrm{in}}) \in \frE_\df$, and showing that the ensuing data-to-solution map is continuous with respect to the $\frE$-topology. In the process, we will also show that the analogous statement holds in caloric-temporal gauge; altogether, 

\begin{theorem}\label{thm:cts}
    The following statements hold for the Chern-Simons-Schr\"odinger equation \eqref{CSS-amp-schro} in either caloric-temporal gauge \eqref{CTC} or DeTurck gauge \eqref{DeTurck}:
    \begin{itemize}
        \item \textup{(Rough solutions as limits of smooth solutions).} Let $(A^{\mathrm{in}}, \phi^{\mathrm{in}}) \in \frE_\df$ be a finite-energy configuration on $\R^2$ in Coulomb gauge \eqref{coulomb-data}. Then there exists a unique solution $(A, \phi) \in C_t ([0, \infty) \to \frE)$ with initial data $(A, \phi)_{|t = 0} = (A^{\mathrm{in}}, \phi^{\mathrm{in}})$. The solution may be realised as the limit of smooth solutions $\{ (A^{[n]}, \phi^{[n]})  \}_n$ with respect to the $\frE$-metric, i.e., for each $T > 0$, 
            \[
                \lim_{n \to \infty} \sup_{t \in [0, T]} \dist_\frE \bigl( (A^{[n]}, \phi^{[n]}) , (A, \phi) \bigr) = 0 .
            \]
        
        \item \textup{(Continuity of data-to-solution map).} For each $T > 0$, the data-to-solution map 
            \begin{align*}
                \frE_\df 
                    &\longrightarrow C_t ([0, T] \to \frE), \\
                (A^{\mathrm{in}}, \phi^{\mathrm{in}}) 
                    &\longmapsto (A, \phi), 
            \end{align*}
        is continuous. 
    \end{itemize}
    
\end{theorem}

As an immediate corollary of constructing rough solutions as limits of smooth solutions, the rough solutions inherit the desired qualitative properties from the class of smooth solutions, namely continuity of the flow $(A, \phi) \in C_t ([0, T] \to \frE)$, conservation of the abelian Higgs energy and the topological degree. 

We first prove Theorem \ref{thm:cts} for the caloric-temporal gauge, then transfer the result to the DeTurck gauge by estimating the gauge-transformation. That is, we regard the data-to-solution map for DeTurck gauge as a composition of that for caloric-temporal gauge, and the gauge-transformation map, 
    \begin{align*}
		\begin{matrix}
			\frE_\df \\
			\text{initial data in \eqref{coulomb-data}}
		\end{matrix}
		\qquad\overset{\text{data-to-sol.}}{\longrightarrow}\qquad
		\begin{matrix}
			C_t ([0, T] \to \frE)\\
			\textup{solution in \eqref{CTC}} 
		\end{matrix}
		\qquad\overset{\text{gauge trans.}}{\longrightarrow}\qquad
		\begin{matrix}
			C_t ([0, T] \to \frE)\\
			\textup{solution in \eqref{DeTurck}}.
		\end{matrix}
	\end{align*}
It is a standard fact of life that if each individual map is continuous, then their composition is also continuous. Thus, we divide the analysis of the equation in DeTurck gauge between these two tasks. 

\subsection{Data-to-solution map in caloric-temporal gauge}\label{subsec:data-to-sol}

We construct solutions in $C_t ([0, \infty) \to \frE)$ to the Chern-Simons-Schr\"odinger equation in caloric-temporal gauge with initial data in $\frE_\df$ and prove continuity of the ensuing data-to-solution map following the Ifrim-Tataru roadmap \cite{IfrimTataru2022}. The continuous dependence argument follows more-or-less the same lines as the existence argument \textit{mutatis mutandis}, so we focus on the construction, which proceeds as follows:
    \begin{enumerate}[label=(\tt\roman*)]
        \item Regularising the initial data (Corollary \ref{cor:smooth-data}). 
        
        \item Constructing smooth solutions on $[0, \infty) \times \R^2$ arising from the regularised data (Theorems \ref{thm:smooth-LWP}, \ref{thm:global-1}, Propositions \ref{prop:dP-exist}, \ref{prop:ctc-exists}): 
        
        \item Proving uniform frequency-envelope bounds for smooth solutions (Theorem \ref{thm:dP-est}, Corollary \ref{cor:global-2}).
        
        \item Proving continuity of the data-to-solution map with respect to a weaker topology (Theorem \ref{thm:weak-cty}).
        
        \item Interpolate between the frequency-envelope bounds and the weak-difference bounds to conclude convergence in the energy topology. 
    \end{enumerate}
It remains to establish the last step in the roadmap. Taking stock of the preceding steps, we have smooth solutions $(A^{[n]}, \phi^{[n]})$ on $[0, \infty) \times \R^2 \times [0, \underline s]$ to the Chern-Simons-Schr\"odinger equation \eqref{CSS-amp-schro} in the caloric-temporal gauge \eqref{CTC} with initial data in Coulomb gauge \eqref{coulomb-data} and satisfying 
    \[
        \lim_{n \to \infty} \dist_{\pzcE^1 ([0, \underline s])} \bigl( (A^{\mathrm{in}}, \phi^{\mathrm{in}}), (A^{[n], \mathrm{in}}, \phi^{[n], \mathrm{in}}) \bigr) 
            = 0, \qquad \lim_{n \to \infty} \| \Gauss^{[n]} \|_{L^\infty_x} \leq \| \Gauss \|_{L^\infty_x}.
    \]
By Proposition \ref{prop:env-stable}, there exists a sequence of $\frE$-frequency envelopes $\ttc^{[n]} (s)$ for $(A^{[n],\mathrm{in}}, \phi^{[n], \mathrm{in}})$ converging to a $\frE$-frequency envelope $\ttc(s)$ for $(A^{\mathrm{in}}, \phi^{\mathrm{in}})$,
    \begin{equation}\label{eq:env-converge}
        \lim_{n \to \infty}\| \ttc^{[n]} (s)- \ttc(s) \|_{L^2_{\frac{ds}{s}}} = 0.
    \end{equation}
By Corollary \ref{cor:converge-low}, the initial data also converge with respect to the weaker $\pzcE^0$-distance functional, so, combined with Theorem \ref{thm:weak-cty}, the sequence of smooth solutions is Cauchy with respect to the $\pzcE^0$-distance. Thus, we can find a limiting weak solution $(A, \phi)$ with the initial data $(A, \phi)_{|t = 0} = (A^{\mathrm{in}}, \phi^{\mathrm{in}})$, and, for each $T > 0$, 
    \begin{equation}\label{eq:E0-converge}
        \lim_{n \to \infty}\sup_{t \in [0, T]} \dist_{\pzcE^0 ([0, \underline s])} \bigl( (A, \phi), (A^{[n]}, \phi^{[n]}) \bigr)  = 0.
    \end{equation}

We claim that this convergence also holds with respect to the $C_t \frE$-metric -- equivalently, by Proposition \ref{thm:topology}, we want to show that 
    \begin{equation}\label{eq:E1-converge}
        \lim_{n \to \infty} \sup_{t \in [0, T]}\dist_{\pzcE^1 ([0, \underline s])} \bigl( (A, \phi), (A^{[n]}, \phi^{[n]}) \bigr)
            = 0.
    \end{equation}
For brevity, denote $\delta^{[n]}$ for the difference operator between fields involving $(A, \phi)$ and $(A^{[n]}, \phi^{[n]})$. Since the Gauss tension field evolves by transport \eqref{gauss-transport-intro} under the Chern-Simons-Schr\"odinger flow, we trivially have convergence of the Gauss tension field for all time. Furthermore, the $\pzcE^0$-convergence along with the magnetic potential difference bounds \eqref{F12-diff}, \eqref{Ax-diff} imply convergence of $\overline A$. Altogether, 
        \begin{align*}
            \lim_{n \to \infty} \sup_{t \in [0, T]} \bigl\| \delta^{[n]} \overline\Gauss \bigr\|_{H^1} + \bigl\| \delta^{[n]} \overline A \bigr\|_{L^4 \cap \dot H^1} 
                = 0.
        \end{align*}
It remains to show convergence of the differences of the covariant Hessians along with their derivatives. To this end, we introduce a frequency cut-off $s_* \in [0, \underline s]$ and split the analysis between high-frequencies $[0, s_*]$ and low-frequencies $[s_*, \underline s]$. We control the contribution of the high-frequencies using the respective frequency envelope bounds from Corollary \ref{cor:global-2}, while for the low-frequencies, we use the trivial inequality $|s|^{-1/2} \leq |s_*|^{-1/2}$ so as to bound this contribution by the ``weak'' distance functional, 
        \begin{align*}
            \bigl\| s^{-\frac12} \delta^{[n]}(s^\frac12 \bfD)^{(\leq 2026)} (s \bfD^{(2)}) \phi \bigr\|_{L^2_{\frac{ds}{s}} L^2_x}
                &\leq \bigl\| s^{-\frac12} \delta^{[n]} (s^\frac12 \bfD)^{(\leq 2026)} (s \bfD^{(2)}) \phi \bigr\|_{L^2_{\frac{ds}{s}} L^2_x ([0, s_*])} \\
                    &\qquad + \bigl\| s^{-\frac12} \delta^{[n]} (s^\frac12 \bfD)^{(\leq 2026)} (s \bfD^{(2)}) \phi \bigr\|_{L^2_{\frac{ds}{s}} L^2_x ([s_*, \underline s])}\\
                &\lesssim \| \ttc(s) \|_{L^2_{\frac{ds}{s}} ([0, s_*])} + \| \ttc^{[n]} (s) \|_{L^2_{\frac{ds}{s}} ([0, s_*])} \\
                    &\qquad + |s_*|^{-\frac12} \cdot \dist_{\pzcE^0} \bigl( (A, \phi), (A^{[n]}, \phi^{[n]} )\bigr)
        \end{align*}
    Passing to the limit $n \to \infty$, the contribution of the low-frequency part of the difference vanishes due \eqref{E0-converge}, while the envelopes converge \eqref{env-converge},
        \[
            \lim_{n \to \infty} \sup_{t \in [0, T]}\bigl\| s^{-\frac12} \delta (s^\frac12 \bfD)^{(\leq 2026)} (s \bfD^{(2)}) \phi \bigr\|_{L^2_{\frac{ds}{s}} L^2_x} 
                \lesssim \| \ttc(s) \|_{L^2_{\frac{ds}{s}} ([0, s_*])}.
        \]
    Taking $s_* \to 0$, the right-hand side vanishes by monotone convergence theorem. This completes the proof of \eqref{E1-converge}. 

\subsection{Transformation from caloric-temporal to DeTurck gauge}\label{subsec:CTC-DT}

We transfer the well-posedness in caloric-temporal gauge to DeTurck gauge by estimating the ensuing gauge-transformation. 

\begin{proposition}[Transformation from \eqref{CTC} to \eqref{DeTurck}]\label{prop:ctc-to-dt}
    Let $(A, \phi)$ be a configuration on $[0, T] \times \R^2 \times [0, \underline s]$ solving the paradifferential formulation of the Chern-Simons-Schr\"odinger equation \eqref{dP-CSS} in caloric-temporal gauge \eqref{CTC} with initial data in Coulomb gauge \eqref{coulomb-data}. Set
        \begin{align*}
            \widetilde A_\bfa
                &:= A_\bfa + \partial_\bfa \chi, \\
            \widetilde \phi 
                &:= e^{i \chi} \phi,
        \end{align*}
    where 
        \[
            \chi 
                := - (\partial_t - \Delta)^{-1} (\overline A_t - \partial^j \overline A_j).
        \]
    Then $(\widetilde A, \widetilde \phi)$ solves the Chern-Simons-Schr\"odinger equation \eqref{CSS-amp-schro} in DeTurck gauge \eqref{DeTurck} with initial data in Coulomb gauge \eqref{coulomb-data}. Furthermore, the gauge transformation $(A, \phi) \mapsto (\widetilde A, \widetilde \phi)$ is continuous on $C_t ([0, T] \to \frE)$.
\end{proposition}

It will be convenient to record the following general lemma, 

\begin{lemma}[Stability of $\frE$-metric under $H^2$-gauge-transformations]\label{lem:gauge-diff}
    Let $(A, \phi), (A', \phi') \in \frE$ be finite-energy configurations on $\R^2$, and suppose that $\chi, \chi' \in H^2 (\R^2)$ are gauge-transformations. Define the gauge-transformed finite-energy configurations $(\widetilde A ,\widetilde \phi), (\widetilde A', \widetilde \phi') \in \frE$ by 
        \begin{equation*}
            \begin{split}
                \widetilde A_j 
                &:= A_j + \partial_j \chi, \\
            \widetilde \phi 
                &:= e^{i \chi} \phi,
            \end{split}
            \qquad\qquad
            \begin{split}
                \widetilde A_j' 
                &:= A_j' + \partial_j \chi', \\
            \widetilde \phi' 
                &:= e^{i \chi'} \phi'.
            \end{split}
        \end{equation*}
    Then 
        \[
            \dist_\frE \bigl( (\widetilde A, \widetilde \phi), (\widetilde A', \widetilde \phi') \bigr)
                \lesssim \dist_\frE \bigl( (A, \phi), (A', \phi') \bigr) + \Bigl( 1 + \cE_{\mathrm{AH}} [A, \phi] + \cE_{\mathrm{AH}} [A', \phi'] \Bigr)^\frac12 \bigl\| \delta \chi \bigr\|_{H^2}.
        \]
\end{lemma}

\begin{proof}
    We recall the definition of the metric below for convenience, 
        \begin{align*}
			\dist_\frE \bigl( (A, \phi), (A', \phi') \bigr) 
				&:= \bigl\| \delta \bfD \phi \bigr\|_{L^2} + \bigl\| \delta \Gauss \bigr\|_{H^1} + \bigl\| \delta \tfrac12(1 - |\phi|^2) \bigr\|_{L^2} + \bigl\| \delta A \bigr\|_{L^4 \, \cap\, \dot H^1} + \bigl\| \delta \phi \bigr\|_{L^2 + L^\infty}.
		\end{align*}
    The second and third terms on the right are manifestly gauge-invariant, so it remains to control the gauge-transformed differences of the magnetic derivatives, and of the configuration itself. A calculation reveals
        \begin{align*}
            \delta \widetilde \phi 
                &= e^{i \chi} \Bigl( \delta \phi + e^{i \delta \chi} \phi \Bigr), \\
            \delta \widetilde{\bfD \phi} 
                &= e^{i \chi} \Bigl( \delta \bfD \phi + e^{i \delta \chi} \bfD \phi \Bigr).
        \end{align*}
    The difference of the gauge-transformed scalar field therefore satisfies 
        \begin{align*}
            \bigl\| \delta \widetilde \phi  \bigr\|_{L^2 + L^\infty} 
                &\lesssim \bigl\| \delta \phi \bigr\|_{L^2 + L^\infty} + \bigl\| e^{i \delta \chi} \bigr\|_{L^\infty} \bigl\| \phi \bigr\|_{L^2 + L^\infty} \\
                &\lesssim  \bigl\| \delta \phi \bigr\|_{L^2 + L^\infty} + \bigl\| \delta \chi \bigr\|_{H^2} \Bigl( 1 + \bigl\| 1 - |\phi|^2 \bigr\|_{L^2_x} \Bigr),
        \end{align*}
    using the triangle inequality in the first line, Lipschitz continuity of the exponential and Sobolev embedding to estimate $e^{i \delta \chi}$, and the $(L^2 + L^\infty)$-bound from Lemma \ref{lem:gerard} to handle the norm of $\phi$. By a similar circle of ideas for the difference of the magnetic derivatives, 
        \begin{align*}
            \bigl\| \delta \widetilde{\bfD \phi} \bigr\|_{L^2}
                &\lesssim \bigl\| \delta \bfD \phi \bigr\|_{L^2} + \bigl\| e^{i \delta \chi} \bigr\|_{L^\infty} \bigl\| \bfD \phi \bigr\|_{L^2} \\
                &\lesssim \bigl\| \delta \bfD \phi \bigr\|_{L^2} + \bigl\| \delta \chi \bigr\|_{H^2} \bigl\| \bfD \phi \bigr\|_{L^2}.
        \end{align*}
    The gauge-transformation acts on the magnetic potential linearly, so 
        \[
            \delta \widetilde A 
                = \delta A + \nabla \delta \chi. 
        \]
    A straightforward application of the triangle inequality and Sobolev embedding yields 
        \begin{align*}
            \bigl\|\delta \widetilde A \bigr\|_{L^4 \, \cap \, \dot H^1}
                &\lesssim \bigl\| \delta A \bigr\|_{L^4 \, \cap \, \dot H^1} + \bigl\| \nabla \delta \chi \bigr\|_{L^4 \, \cap \, \dot H^1} \\
                &\lesssim \bigl\| \delta A \bigr\|_{L^4 \, \cap \, \dot H^1} + \bigl\| \delta \chi \bigr\|_{H^2}.
        \end{align*}
    Collecting the previous inequalities, we conclude the proof. 
\end{proof}

It will also be useful to record the following energy estimate for the inhomogeneous linear heat equation with zero initial data, 

\begin{lemma}[Energy estimate for heat equation]
    Let $\sfN : [0, T] \times \R^2 \to \C$ be sufficiently regular, then 
    \begin{equation}\label{eq:heat-energy-estimate}
        \bigl\| (\partial_t - \Delta)^{-1} \nabla \sfN \bigr\|_{L^\infty_t L^2_x} + \bigl\| \nabla(\partial_t - \Delta)^{-1} \nabla \sfN \bigr\|_{L^2_{t, x}}
            \lesssim \bigl\| \sfN \bigr\|_{L^2_{t, x}}.
        \end{equation}
\end{lemma}

\begin{proof}
    Writing the equation,
        \[
            (\partial_t - \Delta)  (\partial_t - \Delta)^{-1} \nabla \sfN 
                = \nabla \sfN, 
        \]
    the result follows by integration-by-parts on both the left and right. 
\end{proof}

\begin{proof}[Proof of Proposition \ref{prop:ctc-to-dt}]
    There are two tasks at hand: first, showing that if $(A,\phi) \in C_t ([0, T] \to \frE)$, then $(\widetilde A, \widetilde \phi) \in C_t ([0, T] \to \frE)$, second, showing that $(A, \phi) \mapsto (\widetilde A, \widetilde \phi)$ is continuous on $C_t ([0, T] \to \frE)$. By Lemma \ref{lem:gauge-diff}, it suffices to show $\chi \in C_t H^2 ([0, T] \times \R^2)$ and
        \begin{align*}
            \bigl\| \delta \chi \bigr\|_{L^\infty_t H^2_x} 
                &= o(1) , \qquad \text{as $\sup_t \dist_\frE \to 0$.}
        \end{align*}
    We focus on showing the difference bound above; showing $\chi \in C_t H^2 ([0, T] \times \R^2)$ is similar. 

    By Corollary \ref{cor:converge-low} and Proposition \ref{thm:topology}, it is useful to recall that the $\pzcE^0$- and $\pzcE^1$-metrics also vanish, 
        \[
            \sup_t \dist_{\pzcE^0}, \quad \sup_t \dist_{\pzcE^1} = o(1), \qquad \text{as $\sup_t \dist_\frE \to 0$.}
        \]
    By a standard parabolic energy estimate for the (inhomogeneous) linear heat equation \eqref{heat-energy-estimate}, 
        \begin{align*}
            \bigl\| \delta \chi \bigr\|_{L^\infty_t H^2_x}
                &\lesssim \bigl\| \delta \partial^j \overline A_j \bigr\|_{L^1_t H^2_x} + \bigl\| \langle \nabla \rangle \delta \overline A_t \bigr\|_{L^2_{t, x}}. 
        \end{align*}
    Since the curl-free part of $A_x$ is transported \eqref{Ax-transport-cf} along the Yang-Mills heat flow in caloric gauge, the first term on the right is immediately controlled by the $\pzcE^0$-distance. To handle the contribution of the electric potential, we integrate its transport equation \eqref{At-transport} in $s$ under the caloric-temporal gauge. We leverage almost-orthogonality, expanding out
        \begin{align*}
            \bigl\| \nabla \delta \overline A_t \bigr\|_{L^2_{t, x}}^2
                &= \Bigl\| \int_0^{\underline s} (s^\frac12 \nabla)^{(2)} \delta F_{tx} (s) \, \frac{ds}{s} \Bigr\|_{L^2_{t, x}}^2 \\
                &= 2 \int_0^{\underline s} \int_{s"}^{\underline s} \langle ({s'}^\frac12 \nabla)^{(2)} \delta F_{tx} (s'), ({s"}^\frac12 \nabla)^{(2)} \delta F_{tx} (s") \rangle_{L^2_{t, x}} \frac{ds'}{s'} \frac{ds"}{s"} \\
                &= - 2 \int_0^{\underline s} \int_{s"}^{\underline s} \Bigl( \frac{s"}{s'} \Bigr)^{\frac12} \langle ({s'}^\frac12 \nabla)^{(3)} \delta F_{tx} (s'), ({s"}^\frac12 \nabla) \delta F_{tx} (s") \rangle_{L^2_{t, x}} \frac{ds'}{s'} \frac{ds"}{s"},
        \end{align*}
    symmetrising the integrals in $s'$ and $s"$ in the second line and integrating-by-parts in the third. Applying Cauchy-Schwarz and Schur's test (compare with Lemma \ref{lem:orthogonal}), followed by parabolic smoothing \eqref{parabolic-strichartz} for the linear heat equation \eqref{Fmunu-heat}, 
        \begin{align*}
            \bigl\| \nabla \delta\overline A_t \bigr\|_{L^2_{t, x}}
                &\lesssim \bigl\| (s^\frac12 \nabla)^{(\leq 2)} (s^\frac12 \nabla) \delta F_{tx} \bigr\|_{L^2_{\frac{ds}{s}} L^2_{t, x}} \\
                &\lesssim \bigl\| \delta \overline F_{tx} \bigr\|_{L^2_{t, x}}. 
        \end{align*}
    By Amp\'ere's law and the definition of the Gauss tension field, 
        \[
            \delta \overline F_{tx} 
                = \delta \overline \phi \cdot \overline{\bfD \phi} + \phi \cdot \delta \overline{\bfD \phi} + \nabla \delta \overline \Gauss. 
        \]
    We need to estimate this in $L^2_{t, x}$. The difference of the Gauss tension field is estimated by the $\frE$-metric. When the difference falls on the differentiated scalar field, we can put the difference in $L^\infty_t L^2_x$, which is controlled by the $\frE$-metric, and the scalar field in $L^2_t L^\infty_x$, which is bounded by the control parameter \eqref{control-bound}. The most challenging term term to estimate is when the difference falls on the scalar field; we place it in $L^2_t L^\infty_x$ and the derivative in $L^\infty_t L^2_x$. By the fundamental theorem of calculus in $s$, we may decompose the difference of the scalar field into 
        \begin{align*}
            \bigl\| \delta \overline \phi \bigr\|_{L^2_t L^\infty_x}
                &\lesssim_T \bigl\| \delta \underline \phi \bigr\|_{L^\infty_{t, x}} + \int_0^{\underline s} \bigl\| \delta (s \bfD^{(2)}) \phi \bigr\|_{L^2_t L^\infty_x} \frac{ds}{s} .
        \end{align*}
    The difference at $s = \underline s$ is controlled by the $\pzcE^1$-metric. For the high-frequency term, we interpolate between Strichartz estimate and the $\pzcE^1$-metric. That is, by Gagliardo-Nirenberg interpolation (and the diamagnetic inequality \eqref{diamagnetic}), and H\"older in $t$, 
        \begin{align*}
            \int_0^{\underline s}\bigl\| \delta (s \bfD^{(2)}) \phi \bigr\|_{L^2_t L^\infty_x} \frac{ds}{s}
                &\lesssim_T \int_0^{\underline s} s^{\frac1{12}} \bigl\| s^{-\frac12} \delta (s \bfD^{(2)}) \phi\bigr\|_{L^\infty_t L^2_x}^{\frac13} \bigl\| s^{\frac18} s^{-\frac12} (s^\frac12 \bfD) \delta (s \bfD^{(2)}) \phi \bigr\|_{L^4_{t, x}}^{\frac23} \frac{ds}{s} \\
                &\lesssim_{T, \underline s} o(1), 
        \end{align*}
    estimating the difference in $L^4_{t, x}$ using Strichartz \eqref{covar-str-global}, and the difference in $L^\infty_t L^2_x$ using \eqref{heat-lip}. 
\end{proof}

\begin{remark}
    One can probably show that the gauge-transformation is in fact Lipschitz continuous by refining our crude interpolation argument with Strichartz-type bounds for the differences. 
\end{remark}

\begin{remark}
    It is not difficult to show that $\chi \in L^\infty_t H^{5/2-}_x ([0, T] \times \R^2)$, so in fact the magnetic potential in DeTurck gauge has regularity $\nabla A_x \in L^\infty_t H^{1/2 -} ([0, T] \times \R^2)$.
\end{remark}

\appendix

\section{Conservation laws}\label{app:conservation}

The Chern-Simons-Schr\"odinger equation obeys the phase-rotation symmetry $\phi \mapsto e^{i \theta_0 }\phi$ for $\theta_0 \in \R$. The infinitesimal generator of this symmetry is $\phi \mapsto i \phi$, so, in light of Noether's theorem, multiplying the covariant Schr\"odinger equation by $i \phi$ and integrating-by-parts appropriately yields a density-flux identity relating the charge density to the charge current. Amusingly, Ampere's law implies a parallel density-flux identity relating the magnetic field to the charge current -- subtracting the two identities yields the transport equation for the Gauss tension field \eqref{gauss-transport-intro},

\begin{proposition}[Transport of Gauss tension field]\label{prop:CSS-Gauss}
    Let $(A, \phi)$ be a configuration on $[0, T] \times \R^2$ obeying the Chern-Simons-Schr\"odinger equation \eqref{CSS-amp-schro}. Then the magnetic field and charge density obey the density-flux identities,
        \begin{align}\label{eq:density-flux1}
             \partial_t F_{12} 
                &= \partial^j \Im(\overline \phi \bfD_j \phi),\\
            \partial_t \tfrac12(1 - |\phi|^2)
                &= \partial^j \Im(\overline \phi \bfD_j \phi).\label{eq:density-flux2}
        \end{align}
    In particular, the Gauss tension field obeys the transport equation
        \begin{equation}\label{eq:t-gauss}
            \partial_t \Gauss = 0. 
        \end{equation}
\end{proposition}

\begin{proof}
    From the Bianchi identity \eqref{bianchi} and Ampere's law from \eqref{CSS-amp-schro}, we compute the density-flux identity \eqref{density-flux1} for the curvature,
        \begin{align*}
            \partial_t F_{12} 
                &= \partial_1 F_{t2} - \partial_2 F_{t1} \\
                &= \partial_1 \left( \Im(\overline \phi \bfD_1 \phi) - \partial_2 F_{12} \right) - \partial_2 \left( - \Im(\overline \phi \bfD_2 \phi) - \partial_1 F_{12} \right)\\
                &= \Im(\overline \phi \bfD^j \bfD_j \phi) = \partial^j \Im(\overline \phi \bfD_j \phi).
        \end{align*}
    Using the covariant Schr\"odinger equation from \eqref{CSS-amp-schro}, we compute the density-flux identity \eqref{density-flux2} for the charge, 
        \begin{align*}
            \partial_t \tfrac12(1 - |\phi|^2)
                &= - \Re(\overline \phi \bfD_t \phi) \\
                &= \Im(\overline \phi \bfD^j \bfD_j \phi) = \partial^j \Im(\overline \phi \bfD_j \phi).
        \end{align*}
    Subtracting the two equations, we conclude that the Gauss tension field satisfies \eqref{t-gauss}.
\end{proof}

The conservation of the abelian Higgs energy arises from the time-translation invariance of the equation $(A, \phi)(t) \mapsto (A, \phi)(t + t_0)$ for $t_0 \in \R$. The infinitesimal generator of this symmetry is $(A, \phi) \mapsto (\partial_t A, \partial_t \phi)$. Thus, using these as multipliers should yield an appropriate density-flux identity, namely conservation of the abelian Higgs energy, 

\begin{proposition}[Conservation of abelian Higgs energy]\label{prop:conserve}
    Let $(A, \phi)$ be a configuration on $[0, T] \times \R^2$ satisfying the Chern-Simons-Schr\"odinger equation \eqref{CSS-amp-schro}. Then 
        \begin{equation}\label{eq:energy-density-flux}
            \tfrac12 \partial_t \Bigl( \overline{\bfD^j \phi} \bfD_j \phi + |F_{12}|^2 + \tfrac\lambda4 (1 -|\phi|^2)^2 \Bigr)
                = \partial^j \Re(\overline{\bfD_t \phi} \bfD_j \phi) - \epsilon\indices{_j^k} \partial_k (F\indices{_t^j} F_{12}).
        \end{equation}
    Furthermore, if $(A, \phi)$ is a smooth solution, then the abelian Higgs energy is conserved under the flow, 
        \begin{equation}\label{eq:energy-conservation}
            \cE_{\mathrm{AH}} [A(t), \phi(t)] = \cE_{\mathrm{AH}} [A(0), \phi(0)] , \qquad \text{for all $t \in [0, T]$}.
        \end{equation}
\end{proposition}

\begin{proof}
    Working in a gauge-covariant fashion, we multiply the covariant Schr\"odinger equation by $\bfD_t \phi$,
	\begin{align*}
		0 
			&= \Re \Bigl(\overline{\bfD_t \phi} (i \bfD_t \phi + \bfD^j \bfD_j + \tfrac\lambda2 (1 - |\phi|^2) \phi)\Bigr)\\
			&= \partial^j \Re(\overline{\bfD_t \phi} \bfD_j \phi) - \Re(\overline{\bfD^j \bfD_t \phi} \bfD_j \phi) - \partial_t \tfrac12 \tfrac\lambda4 (1 - |\phi|^2)^2 \\
            &= \partial^j \Re(\overline{\bfD_t \phi} \bfD_j \phi) - \Re(\overline{\bfD_t \bfD^j \phi} \bfD_j \phi) - F\indices{_t^j} \Im(\overline \phi \bfD_j \phi) - \partial_t \tfrac12 \tfrac\lambda4 (1 - |\phi|^2)^2 \\
            &= \partial^j \Re(\overline{\bfD_t \phi} \bfD_j \phi) - F\indices{_t^j} \Im(\overline \phi \bfD_j \phi) - \tfrac12 \partial_t \Bigl( \overline{\bfD^j \phi} \bfD_j \phi - \tfrac\lambda4 (1 - |\phi|^2)^2 \Bigr),
	\end{align*}
differentiating-by-parts and commuting covariant derivatives appropriately, and similarly multiplying Ampere's law by $F_{tj}$ appropriately, 
    \begin{align*}
        0 
            &= F\indices{_t^j} \Bigl( \epsilon\indices{_j^k} F_{tk} - \Im(\overline \phi \bfD_j \phi) + \epsilon\indices{_j^k} \partial_k F_{12} \Bigr) \\
            &= - F\indices{_t^j} \Im(\overline \phi \bfD_j \phi) + \epsilon\indices{_j^k} \partial_k (F\indices{_t^j} F_{12}) - \epsilon\indices{_j^k} \partial_k F\indices{_t^j} F_{12} \\
            &= - F\indices{_t^j} \Im(\overline \phi \bfD_j \phi) + \epsilon\indices{_j^k} \partial_k (F\indices{_t^j} F_{12}) + \partial_t F_{12} \, F_{12}\\
            &= - F\indices{_t^j} \Im(\overline \phi \bfD_j \phi) + \epsilon\indices{_j^k} \partial_k (F\indices{_t^j} F_{12}) + \tfrac12 \partial_t \bigl(|F_{12}|^2 \bigr),
    \end{align*}
differentiating-by-parts in the second line and using the Bianchi identity \eqref{bianchi} in the third line. Subtracting the previous two calculations, we obtain the density-flux relation \eqref{energy-density-flux}. 

To prove the conservation of the abelian Higgs energy \eqref{energy-conservation}, fix a smooth cut-off $\chi \in C^\infty_c (\R^2)$ such that $\chi \equiv 1$ in a neighborhood of the origin, and test the identity \eqref{energy-density-flux} against its rescaling $\chi(\tfrac{x}{R})$ for $R > 0$. Integrating-by-parts on $[0, t] \times \R^2$ gives the schematic identity
    \begin{align*}
        \int_{\R^2} \mathfrak e[A(t'), \phi(t')] \, \chi\bigl( \tfrac{x}{R} \bigr) \ dx \Big|_{t = 0}^{t' = t}
            = \frac1R \int_0^t \int_{\R^2} (\nabla \chi)\bigl( \tfrac{x}{R} \bigr) \Bigl( \bfD_t \phi \cdot \bfD \phi + F_{tx} \cdot F_{12} \Bigr) \, dx dt',
    \end{align*}
abbreviating $\mathfrak e[A(t), \phi (t)]$ for the energy density. Taking $R \to \infty$, the expression on the left converges to the differences of the energies by the dominated convergence theorem, 
    \[
        \lim_{R \to \infty} \text{L.H.S.} = \cE_{\text{AH}}[A(t), \phi(t)] - \cE_{\text{AH}}[A^{\mathrm{in}}, \phi^{\mathrm{in}}].
    \]
On the right-hand side, an application of Cauchy-Schwarz and passing to the limit yields 
    \[
        \Big| \text{R.H.S.} \Big|
            \lesssim_t \frac{1}R \sup_{t' \in [0, t]}\Big( \| \bfD_t \phi \|_{L^2_x} \| \bfD \phi\|_{L^2_x} + \| F_{tx}\|_{L^2} \| F_{12} \|_{L^2_x}  \Big) \overset{R \to \infty}{\longrightarrow} 0.
    \]
Indeed, it is easy to check that the supremum is $O_t (1)$ for smooth solutions using the equation and Sobolev embedding \eqref{inefficientLinfty}. 
\end{proof}

Define the \textit{vorticity} of a configuration $(A, \phi)$ on $\R^2$ as the curl of the magnetic potential plus the charge current, 
   \[
        \omega (A, \phi)
            := \epsilon_{jk} \partial^j \Bigl(A^k + \Im(\overline \phi \bfD^k \phi)\Bigr).
    \]
The vorticity captures the topologically non-trivial part of the configuration. Indeed, Czubak-Jerrard \cite[Lemma 2.3]{CzubakJerrard2014} showed that the \textit{topological degree} of a finite-energy configuration, defined by 
    \[
        \deg[A , \phi]
            := \tfrac{1}{2\pi}\int_{\R^2} \omega(A, \phi)\, dx,
    \]
is integer-valued. Furthermore, it is not difficult to see from their proof that the topological degree agrees with the usual notion of degree for the \textit{boundary map at infinity} $\phi^\infty : \mathbb S^1 \to \mathbb S^1$, 
    \[
        \phi^\infty (\theta)
            := \lim_{r \to \infty} \phi(r e^{i \theta}),
        \]
when it is well-defined and continuous. 

At first glance, this definition of the topological degree seems different from ours from the introduction, however, a short calculation (namely the identity \eqref{vorticity-identity}) shows that they are nonetheless the same. We write the vorticity this way because it makes the derivation of the modified Biot-Savart law \eqref{Ax-elliptic} rather immediate. 

\begin{lemma}\label{lem:vorticity}
    Let $(A, \phi) \in \frE$ be a finite-energy configuration on $\R^2$. Then 
        \begin{enumerate}
            \item \textup{(Vorticity identity).} 
                \begin{equation}\label{eq:vorticity-identity}
                    \begin{split}
                        \omega(A, \phi)
                            &= (1 - |\phi|^2) F_{12} + 2 \Im(\overline{\bfD_1 \phi} \bfD_2 \phi).
                    \end{split}
                    \end{equation}

            \item \textup{(Modified Biot-Savart law).} 
                \begin{align}
                    \PP^\df A 
                        &= \nabla^{-1} \omega(A, \phi) + \PP^\df (\phi \cdot \bfD \phi ) .\label{eq:BS-low}
                \end{align}
        \end{enumerate}
\end{lemma}

\begin{proof}
\leavevmode
    \begin{enumerate}
        \item Differentiating and commuting appropriately yields
        \begin{align*}
            \omega(A, \phi)
                &= F_{12} + \epsilon_{jk} \partial^j \Im(\overline \phi \bfD^k \phi) \\
                &= F_{12} + \epsilon_{jk} \left( \Im(\overline{\bfD^j \phi} \bfD^k \phi) + \Im(\overline \phi \bfD^j \bfD^k  \phi) \right) \\
                &= (1 - |\phi|^2) F_{12} + 2 \Im(\overline{\bfD_1 \phi} \bfD_2 \phi),
        \end{align*}
    as desired.

        \item We write 
                \[
                    \PP^\df A 
                        = \PP^\df (A + \Im(\overline \phi \bfD \phi)) - \PP^\df \Im(\overline \phi \bfD \phi),
                \]
            and use the Biot-Savart law to re-express the first term on the right-hand side. 
    \end{enumerate}
\end{proof}

Amusingly, the vorticity of an arbitrary configuration $(A, \phi)$ on $[0, T] \times \R^2$ obeys a suitable density-flux relation. Integrating thereby yields conservation of the topological degree along the $1$-parameter family of configurations $(A(t), \phi(t))$ on $\R^2$, provided suitable regularity assumptions. In particular, smooth solutions to the Chern-Simons-Schr\"odinger equation conserve their topological degree along the flow. 

\begin{proposition}[Conservation of topological degree]\label{lem:vorticity-II}
    Let $(A, \phi)$ be a configuration on $[0, T] \times \R^2$. Then the vorticity obeys the density-flux identity
        \begin{equation}\label{eq:difference-flux}
            \partial_t \omega(A, \phi)
                    = \epsilon_{jk} \partial^j \Bigl( (1 - |\phi|^2) F\indices{_t^k}  +  2 \Im(\overline{\bfD_t \phi} \bfD^k \phi)\Bigr).
        \end{equation} 
     Furthermore, if $(A, \phi)$ is a smooth solution to the Chern-Simons-Schr\"odinger equation, then the topological degree is conserved under the flow, 
        \begin{equation}\label{eq:top-conserve}
            \deg[A(t_1), \phi(t_1)] = \deg[A(t_2), \phi(t_2)] \qquad \text{for all $t_1, t_2 \in [0, T]$}.
        \end{equation}
\end{proposition}

\begin{proof}
    By the Bianchi identity \eqref{bianchi},
            \begin{align*}
                \partial_t \epsilon_{jk} \partial^j A 
                    = \partial_t F_{12} 
                    = \epsilon_{jk} \partial^j F\indices{_t^k} .
            \end{align*}
        Differentiating the charge current with respect to $h$ and commuting appropriately gives
            \begin{align*}
                \partial_t \epsilon_{jk} \partial^j \Im(\overline \phi \bfD^k \phi) 
                    &= \epsilon_{jk} \partial^j \Big( \Im(\overline{\bfD_t \phi} \bfD^k \phi) + \Im(\overline \phi \bfD_t \bfD^k \phi) \Big)\\
                    &= \epsilon_{jk} \partial^j \Big( \Im(\overline{\bfD_t \phi} \bfD^k \phi) + \Im(\overline \phi \bfD^k \bfD_t \phi) - |\phi|^2 F\indices{_t^k} \Big)\\
                    &= \epsilon_{jk} \partial^j \Big( 2\Im(\overline{\bfD_t \phi} \bfD^k \phi)  - |\phi|^2 F\indices{_t^k} \Big).
            \end{align*}
        Collecting the two calculations yields \eqref{difference-flux}. 

        The proof of conservation of the topological degree \eqref{top-conserve} is similar to that of conservation of the abelian Higgs energy -- simply integrate the density-flux relation \eqref{difference-flux} appropriately. 
\end{proof}

\begin{remark}
    The density-flux relation for the vorticity \eqref{difference-flux} is completely algebraic and holds for any configuration on $[0, T] \times \R^2$. Moreover, the conservation of the topological degree \eqref{top-conserve} holds as long as $(A, \phi)$ has uniformly bounded abelian Higgs energy and finite $L^2$-length, i.e. 
    \[
        \sup_{t \in [0, T]} \cE_{\textup{AH}} [A(t), \phi(t)] + \int_0^T \| \bfD_t \phi  \|_{L^2_x} + \| F_{tx} \|_{L^2_x} \, dt < \infty.
    \]
    One can think of this as the statement that two configurations $(A(0), \phi(0))$ and $(A(T), \phi(T))$ have the same topological degree if they admit a ``reasonable'' homotopy. 
\end{remark}

\section{Spaces of configurations}\label{app:example}
\subsection{Examples of smooth configurations}\label{subapp:example}

To illustrate the complexity of the space of finite-energy configurations, we adapt an example due to G\'erard \cite{Gerard2006} in the context of the non-gauged Ginzburg-Landau energy to our setting. Let $\rho \in C^\infty_c ([0, \infty) \to \R)$ be a radial test function such that $\rho(r) = 1$ for $|r| \leq 1$. Then, for $n \in \Z$, and parameters $0 \leq \alpha \leq \tfrac{1}{1000}$ and $0 \leq \gamma < \tfrac12$, we define the configurations $(A^{(n, \alpha,\gamma)}, \phi^{(n, \alpha,\gamma)})$ by 
		\begin{align*}
			A_r (r, \theta)
				&= 0,\\
			A_\theta (r, \theta) 
				&= (1 - \rho (r)) \bigl( n + \alpha \cos (\theta)\bigr) ,\\
			\phi (r, \theta)
				&=  (1 - \rho(r))  \exp\bigl( i \log^\gamma (2 + r) \bigr) \exp\bigl( i \bigl( n\theta + \alpha \sin (\theta)\bigr) \bigr),
		\end{align*}
where $(r, \theta)$ are the usual polar coordinates on $\R^2$. 

\begin{proposition}\label{prop:example}
    $(A^{(n, \alpha,\gamma)}, \phi^{(n, \alpha,\gamma)}) \in \frE^\infty \cap C^\infty (\R^2)$ are smooth configurations. 
\end{proposition}

\begin{proof}
    For brevity, we suppress the dependence on the parameters and write $(A, \phi) \equiv (A^{(n, \alpha,\gamma)}, \phi^{(n, \alpha,\gamma)})$. It is easy to see by inspection that $(A, \phi) \in C^\infty (\R^2)$, so it remains to verify $(A, \phi) \in \frE^\infty$. 
    
    Computing the magnetic gradient in polar coordinates, 
        \begin{align*}
            \bfD_r \phi 
                &= \partial_r \phi \\
                &= \Bigl( - \partial_r \rho + i \gamma (1 - \rho) \frac{\log^{\gamma - 1} (2 + r)}{2 + r} \Bigr) \exp\bigl( i \log^{\gamma}(2 + r) \bigr) \exp\bigl( i \bigl( n\theta + \alpha \sin (\theta)\bigr) \bigr),\\ 
            \bfD_\theta \phi 
                &= (\partial_\theta - i A_\theta) \phi \\
                &= i (1 - \rho) \Bigl( \partial_\theta (n \theta + \alpha \sin(\theta)) - (n + \alpha \cos (\theta))\Bigr)\exp\bigl( i \bigl( n\theta + \alpha \sin (\theta)\bigr) \bigr) = 0.
        \end{align*}
    Then, since $\rho$ is smooth, compactly supported, and $(2 + r)^{-1} \log^{\gamma - 1} \in L^2 (r dr)$ whenever $\gamma < \tfrac12$, we see that $\bfD \phi \in L^2 (\R^2)$. More generally, observe that for any combination of magnetic derivatives, the contributions of the angular derivatives $\bfD_\theta$ vanish, while consecutive applications of the radial derivatives $\bfD_r$ remain smooth and can only improve decay. Thus, $\bfD^{(n)} \phi \in L^2 (\R^2)$ for any positive integer $n \in \N$. 
    
    For the curvature and charge density, 
        \begin{align*}
            F_{r \theta} 
                &= \partial_r A_\theta \\
                &= - \partial_r \rho \cdot (n + \alpha \cos (\theta)), \\
            \tfrac12(1 - |\phi|^2) 
                &= \tfrac12(1 - |1 - \rho|^2) \\
                &= \rho - \tfrac12 \rho^2,
        \end{align*}
    so it is easy to see that $F_{12}, \tfrac12 (1 - |\phi|^2) \in C^\infty_c (\R^2)$ by choice of $\rho$. 
\end{proof}

\begin{proposition}
    $(A^{(n, 0,\gamma)}, \phi^{(n, 0,\gamma)})$ is in Coulomb gauge. 
\end{proposition}

\begin{proof}
    We calculate the divergence of the magnetic potential in polar coordinates, 
        \[
            \partial^j A_j = \Bigl( \tfrac1r + \partial_r \Bigr) A_r + \frac{1}{r^2} \partial_\theta A_\theta = -\alpha (1 - \rho) \sin(\theta). 
        \]  
    Thus, $A$ is in Coulomb gauge if and only if $\alpha = 0$. 
\end{proof}

We use the family of examples $\{ (A^{(n, \alpha, \gamma)}, \phi^{(n,\alpha, \gamma)})\}_{n, \alpha, \gamma}$ to illustrate that the space of finite-energy configurations cannot be effectively captured by any familiar function spaces. More precisely, we show that 
    \[
        \frE_\phi 
            := \bigl\{ \phi \in H^1_\loc (\R^2) : (A, \phi) \in \frE \bigr\},
    \]
can neither be embedded in an affine $L^p$-space, nor a magnetic Sobolev space. 

\begin{proposition}
\leavevmode
    \begin{enumerate}
        \item If $\gamma \neq \gamma'$, then for each $1 \leq p < \infty$ we have 
            \[
				\phi^{(n, \alpha, \gamma)} - \phi^{(n, \alpha, \gamma')} \not\in L^p (\R^2).
			\]
        In particular, $\frE_\phi \not\subseteq \mathring \phi + L^p (\R^2)$ for any function $\mathring \phi : \R^2 \to \C$. 
        
        \item If $n \neq 0$, then
            \[
                \nabla \phi^{(n, \alpha, \gamma)} \not\in L^2 (\R^2).
            \] 
        In particular, $\frE_\phi \not\subseteq \dot H^1 (\R^2)$. 
        
        \item If $\alpha \neq \alpha'$, then 
            \[
				\bigl(\nabla - i A^{(n, \alpha, \gamma)}\bigr) \phi^{(n, \alpha', \gamma)} \not\in L^2 (\R^2).
			\]
        In particular, $\frE_\phi \not\subseteq \dot H^1_{\mathring A} (\R^2)$ for any connection $1$-form $\mathring A : \R^2 \to \R^2$.
    \end{enumerate}
\end{proposition}

\begin{proof}
\leavevmode
    \begin{enumerate}
        \item This is due to non-resonant phase oscillation. Indeed, we calculate, 
            \begin{align*}
                \bigl|\phi^{(n, \alpha, \gamma)} - \phi^{(n, \alpha, \gamma')}\bigr|
                    = |1 - \rho| \cdot \bigl|\exp(i\log^\gamma (2 + r)) - \exp(i\log^{\gamma'} (2 + r))\bigr|,
            \end{align*}
        so, unless $\gamma = \gamma'$, we can easily find infinitely many unit-length intervals where the amplitude of the difference is of size, say, at least $c > 0$. 

        \item We calculate 
            \[  
                |\partial_\theta \phi^{(n, \alpha, \gamma)} |
                    = |1 - \rho| \cdot |n + \alpha \cos (\theta)|,
            \]
        If $n \neq 0$, the calculation above shows that the gradient obeys the pointwise bound 
            \[
                \bigl|\nabla \phi^{(n, \alpha, \gamma)}\bigr|
                    \gtrsim \tfrac1r 
            \]
        in the region where $|\cos(\theta)| \ll 1$ and $r \gg 1$, so it is not $L^2$-integrable. 

        \item We calculate
            \begin{align*}
                \bigl|(\partial_\theta - i A^{(n, \alpha, \gamma)}_\theta) \phi^{(n, \alpha', \gamma)} \bigr| 
                    &= |1 - \rho| \cdot \Bigl| \partial_\theta (n \theta   \alpha \sin (\theta)) - (n + \alpha' \cos (\theta)) \Bigr|\\
                    &= |1 - \rho| \cdot |\alpha - \alpha'| \cdot |\cos(\theta)|.
            \end{align*}
        If $\alpha \neq \alpha'$, the calculation above shows that the gradient obeys the pointwise bound 
            \[
                \bigl| (\nabla - i A^{(n, \alpha, \gamma)}) \phi^{(n, \alpha', \gamma)} \bigr| 
                    \gtrsim \tfrac1r, 
            \]
        in the region where $|\cos(\theta)| \ll 1$ and $r \gg 1$, so it is not $L^2$-integrable. 

        Suppose now, towards a contradiction, that we can find a fixed connection $1$-form $\mathring A$ such that $\mathring \bfD \phi^{(n, \alpha, \gamma)}, \mathring \bfD \phi^{(n, \alpha', \gamma)} \in L^2(\R^2)$, where $\mathring \bfD = \nabla - i \mathring A$. However, we can write  
            \begin{align*}
                (\nabla - i A^{(n, \alpha, \gamma)}) \phi^{(n, \alpha', \gamma)}
                    &= \mathring \bfD \phi^{(n, \alpha', \gamma)} + i \bigl(\mathring A - A^{(n, \alpha, \gamma)}\bigr) \phi^{(n, \alpha', \gamma)},\\
                    &= \mathring \bfD \phi^{(n, \alpha', \gamma)} + i \Bigl(   (\nabla - i A^{(n, \alpha, \gamma)}) \phi^{(n, \alpha, \gamma)} - \mathring \bfD \phi^{(n, \alpha, \gamma)}\Bigr)\exp \Bigl( i (\alpha' - \alpha) \sin(\theta) \Bigr).
            \end{align*}
        By the assumptions and Proposition \ref{prop:example}, the right-hand side is in $L^2 (\R^2)$, a contradiction. 
    \end{enumerate}
\end{proof}

\subsection{Space of smooth configurations}
\label{sec:smooth-configurations}

In this section, we show that a smooth configuration $(A, \phi) \in \frE^\infty$ in Coulomb gauge \eqref{coulomb-data} is smooth in the usual sense, that is, $(A, \phi) \in C^\infty (\R^2)$. Recall that configuration is smooth if
    \[
        \bigl\| \tfrac12 (1 - |\phi|^2) \bigr\|_{L^2_x} + \bigl\| \bfD^{(n)} \phi \bigr\|_{L^2_x} + \bigl\| \nabla^{(n)} F_{12} \bigr\|_{L^2_x} 
            < \infty
    \]
for all non-negative integers $n \in \N$.

\begin{proposition}[$A$ is smooth]
    Let $(A, \phi) \in \frE^\infty$ be a smooth configuration on $\R^2$ in Coulomb gauge \eqref{coulomb-data}. Then $A \in C^\infty (\R^2)$, and, for each non-negative integer $n \in \N_0$, 
        \begin{equation}\label{eq:A-Cinfty}
            \bigl\| \nabla^{(n)} A \bigr\|_{L^\infty_x}
                \lesssim_{\| \frac12 (1 - |\phi|^2) \|_{L^2}, \| F_{12} \|_{H^{n + 1}}, \| \bfD^{(\leq 1)}\bfD \phi \|_{L^2}} 1.
        \end{equation}
\end{proposition}

\begin{proof}
    For each positive integer $n \in \N$, it follows from Sobolev embedding and the Coulomb gauge \eqref{coulomb-data} that 
        \[
            \| \nabla^{(n - 1)} \nabla A \|_{L^\infty}
                \lesssim \| F_{12} \|_{H^{n + 1}}.
        \]
    This proves \eqref{A-Cinfty} at top-order. At bottom-order, we interpolate,
        \[
            \| A \|_{L^\infty_x}
                \lesssim \| A \|_{L^4} + \| F_{12} \|_{H^1}.
        \]
    To estimate $A$ in $L^4$, we use the modified Biot-Savart law \eqref{BS-low} to write $A = \nabla^{-1} \omega (A, \phi) + \phi \cdot \bfD \phi$. For the contribution of the charge current, 
        \begin{align*}
            \| \phi \cdot \bfD \phi \|_{L^4} 
                &\lesssim \| \phi \|_{L^\infty} \| \bfD \phi \|_{L^4} \\
                &\lesssim 1 + \| \tfrac12 (1 - |\phi|^2) \|_{L^2}^2 + \| \bfD^{(\leq 1)} \bfD \phi \|_{L^2}^2
        \end{align*}
    estimating $\phi$ by the $L^\infty$-bound \eqref{inefficientLinfty}, and its derivative by interpolation \eqref{GN-interpolation}. For the contribution of the vorticity, 
        \begin{align*}
            \| \nabla^{-1} \omega (A, \phi) \|_{L^4}
                &\lesssim \| \omega(A, \phi) \|_{L^{4/3}} \\
                &\lesssim \| \bfD \phi \|_{L^2} \| \bfD\phi \|_{L^4} + \| \tfrac12 (1 - |\phi|^2) \|_{L^2} \| F_{12} \|_{L^4}\\
                &\lesssim  \| \bfD^{(\leq 1)} \bfD \phi \|_{L^2}^2 + \| \tfrac12 (1 - |\phi|^2) \|_{L^2} \| F_{12} \|_{H^1}. 
        \end{align*}
    using Hardy-Littlewood-Sobolev in the first line, H\"older's inequality applied to the identity \eqref{vorticity-identity} in the second, and Sobolev embedding in the last. Collecting the inequalities above, we conclude \eqref{A-Cinfty} at bottom-order. This completes the proof.  
\end{proof}

\begin{proposition}[$\phi$ is smooth]
     Let $(A, \phi) \in \frE^\infty$ be a smooth configuration on $\R^2$ in Coulomb gauge \eqref{coulomb-data}. Then $\phi \in C^\infty (\R^2)$, and, for each non-negative integer $n \in \N_0$, 
        \begin{equation}\label{eq:phi-Cinfty}
            \bigl\| \nabla^{(n)} \phi \bigr\|_{L^\infty_x}
                \lesssim_{\| \bfD^{(\leq n + 1)} \bfD \phi \|_{L^2}, \| \frac12(1 - |\phi|^2) \|_{L^2}, \| F_{12} \|_{H^{n}}} 1.
        \end{equation}
\end{proposition}

\begin{proof}
    The case $n = 0$ is a consequence of the Sobolev-type estimate \eqref{inefficientLinfty}. From here, we proceed inductively, assuming the result holds up to $n$, and aiming to prove the estimate for $n + 1$. To facilitate the analysis, observe the following schematic formula, 
        \[
            \nabla^{(n + 1)} \phi 
                = \bfD^{(n + 1)} \phi + (\nabla^{(\leq n)} A)^{(\leq n)} \nabla^{(\leq n)} \phi,
        \]
    where $\nabla^{(\leq n)}$ denotes differentiation up to order $n$, and $u^{(\leq n)}$ denotes a product consisting of up to $n$-many terms. It is not difficult to see the formula holds by an inductive application of the definition of the covariant derivative $\bfD = \nabla - i A$ and the product rule. It follows from the formula that 
        \begin{align*}
            \bigl\| \nabla^{(n + 1)} \phi \bigr\|_{L^\infty_x} 
                &\lesssim \bigl\| \bfD^{(n + 1)} \phi \bigr\|_{L^\infty_x} + \bigl\| \nabla^{(\leq n)} A \bigr\|_{L^\infty_x}^{\leq n} \, \bigl\| \nabla^{(\leq n)} \phi \bigr\|_{L^\infty_x} \\
                &\lesssim_{\| \bfD^{(\leq n + 1)} \bfD\phi \|_{L^2}, \| \frac12(1 - |\phi|^2) \|_{L^2}, \| F_{12} \|_{H^{n + 1}}} \bigl\| \bfD^{(n + 1)} \phi \bigr\|_{L^2_x}^{\frac12} \bigl\| \bfD^{(n + 3)} \phi \bigr\|_{L^2_x}^{\frac12} + 1 \\
                &\lesssim_{\| \bfD^{(\leq n + 2)} \bfD\phi \|_{L^2}, \| \frac12(1 - |\phi|^2) \|_{L^2}, \| F_{12} \|_{H^{n + 1}}} 1,
        \end{align*}
    using Sobolev embedding \eqref{GNinfty}, the $C^k$-estimates for the magnetic potential \eqref{A-Cinfty}, and the inductive hypothesis. This completes the proof. 
\end{proof}

\subsection{Comparing magnetic derivatives}
\label{subsec:lip-N}

To facilitate the high-regularity local well-posedness in Appendix \ref{app:picard}, we prove a Lipschitz-type property for the non-linear map
    \[
        (A, \Phi) \mapsto \Bigl( \tfrac12(1 - |\Phi|^2), F_{12}, \bfD_A \Phi \Bigr),
    \]
and also for higher-order magnetic derivatives,
    \[
        (A, \Phi) \mapsto \bfD_A^{(n)} \Phi.
    \]

\begin{proposition}[Lipschitz continuity of abelian Higgs energy]\label{prop:magnetic-lip-1}
\leavevmode
    \begin{itemize}
        \item Let $\Phi, \Psi \in (L^4 + L^\infty) (\R^2)$ such that $\Phi - \Psi \in H^1 (\R^2)$, then 
        \begin{equation}\label{eq:charge-perturb}
            \bigl\| |\Phi|^2 - |\Psi|^2 \bigr\|_{L^2_x} 
                \lesssim \Bigl( \| \Phi \|_{L^4_x + L^\infty_x} + \| \Psi \|_{L^4_x + L^\infty_x} \Bigr) \bigl\| \Phi - \Psi \bigr\|_{H^1_x}.
            \end{equation}
            Furthermore, fixing $\mathring{\vphantom{A} \Phi} \in (L^4 + L^\infty) (\R^2)$, perturbations by $u, v \in H^1 (\R^2)$ satisfy 
            \begin{equation}\label{eq:charge-perturb-2}
                \bigl\| |\mathring{\vphantom{A} \Phi} + u|^2 - |\mathring{\vphantom{A} \Phi} + v|^2 \bigr\|_{L^2_x} 
                    \lesssim \Bigl( \| \mathring{\vphantom{A} \Phi} \|_{L^4_x + L^\infty_x} + \| u \|_{H^1_x} + \| v \|_{H^1_x} \Bigr) \| u - v \|_{H^1_x}.
            \end{equation} 

        \item Let $A, B \in H^1_\loc (\R^2)$ such that $\nabla (A - B) \in L^2 (\R^2)$, then 
            \begin{equation}\label{eq:field-perturb}
                \bigl\| (F_A)_{12} - (F_B)_{12} \bigr\|_{L^2_x} 
                    \lesssim \| \nabla (A - B) \|_{L^2_x}. 
            \end{equation}
            Furthermore, fixing $\mathring A \in H^1_\loc (\R^2)$, perturbations by $a, b \in H^1 (\R^2)$ satisfy 
            \begin{equation}\label{eq:field-perturb-2}
                \bigl\| (F_{\mathring A + a})_{12} - (F_{\mathring A + b})_{12} \bigr\|_{L^2_x} 
                    \lesssim \| a - b \|_{H^1_x}. 
            \end{equation}

        \item Let $A - B, \Phi - \Psi \in H^1 (\R^2)$ 
            \begin{equation}\label{eq:DA-perturb}
                \bigl\| \bfD_A \Phi - \bfD_B \Psi \bigr\|_{L^2_x} 
                    \lesssim \Bigl( 1 + \| \Phi \|_{L^4_x + L^\infty_x} + \| B \|_{L^4_x} \Bigr) \| (A, \Phi) - (B, \Psi) \|_{H^1_x}.
            \end{equation}
            Furthermore, fixing $(\mathring A, \mathring{\vphantom{A} \Phi}) \in H^1_\loc (\R^2)$, perturbations by $(a, u), (b, v) \in H^1 (\R^2)$ satisfy 
            \begin{equation}\label{eq:DA-perturb-2}
                \bigl\| \bfD_{\mathring A + a} (\mathring{\vphantom{A} \Phi} + u) - \bfD_{\mathring A + b} (\mathring{\vphantom{A} \Phi} + v) \bigr\|_{L^2_x} 
                    \lesssim \Bigl( \| \mathring{\vphantom{A} \Phi} \|_{L^4_x + L^\infty_x} + \| \mathring A \|_{L^4} + \| u \|_{L^4_x} + \| b \|_{L^4_x} \Bigr) \| (a, u) - (b, v) \|_{H^1_x}.
            \end{equation}
    \end{itemize}
\end{proposition}

\begin{proof}
    The magnetic field bounds \eqref{field-perturb}-\eqref{field-perturb-2} are obvious by linearity, and the bounds \eqref{charge-perturb-2} and \eqref{DA-perturb-2} are immediate consequences of the previous \eqref{charge-perturb} and \eqref{DA-perturb} by the triangle inequality and Sobolev embedding. Writing 
        \[
            |\Phi|^2 - |\Psi|^2 = (|\Phi| + |\Psi|)(|\Phi| - |\Psi|),
        \] 
    and 
        \begin{align*}
            \bfD_A \Phi - \bfD_B \Psi 
                = \nabla(\Phi - \Psi) -i (A - B) \Phi - i B (\Phi - \Psi),
        \end{align*}
    the bounds \eqref{charge-perturb} and \eqref{charge-perturb} follow from H\"older's inequality, placing the coefficients in $L^4 + L^\infty$ and the differences in $L^4 \cap L^2$, and the Sobolev embedding. 
\end{proof}

For higher-orders of magnetic derivatives, we need to handle multi-linear expressions of increasing complexity. Nonetheless, the main point one should keep in mind is that, given magnetic potentials $A \equiv A_j \, dx^j$ and $B \equiv B_j \, dx^j$, the principal parts of the corresponding magnetic derivatives are the same, 
    \[
        \bfD_A^{(n)} = \bfD_B^{(n)} + \text{lower-order terms}.
    \]
With a bit of combinatorics and interpolation, we can handle the lower-order terms inductively. To facilitate the analysis, we expand out the magnetic derivatives by hand,

\begin{lemma}[Schematic magnetic derivative identity]
    Let $A = A_j \, dx^j$ and $B = B_j \, dx^j$ be connection $1$-forms, and let $\Phi$ be a complex scalar field. Then, for each non-negative integer $n \in \N_0$,  
        \begin{equation}\label{eq:schematic-D}
            \bfD_A^{(n + 1)} \Phi 
                = \bfD_B^{(n + 1)} \Phi + \sum_{j + k = n} \sum_{a + b = k} \bigl( \nabla^{(\leq a)} (A - B) \bigr)^{j + 1} \cdot \bfD_B^{(b)} \Phi
        \end{equation}
    where $\nabla^{(\leq n)}$ denotes differentiation up to order $n$.
\end{lemma}

\begin{proof}
    Clearly the expression holds for $n = 0$ by definition $\bfD_A = \bfD_B - i (A - B)$. Assume for induction then \eqref{schematic-D}, then by the product rule, 
        \begin{align*}
            \bfD^{(n + 2)}_A \Phi 
                &= (\bfD_B - i (A - B)) \Bigl( \bfD_B^{(n + 1)} \Phi + \sum_{j + k = n} \sum_{a + b = k} \bigl( \nabla^{(\leq a)} (A - B) \bigr)^{j + 1} \cdot \bfD_B^{(b)} \Phi \Bigr) \\
                &= \bfD_B^{(n + 2)} \Phi + \sum_{j + k = n + 1} \sum_{a + b = k} \bigl( \nabla^{(\leq a)} (A - B) \bigr)^{j + 1} \cdot \bfD_B^{(b)} \Phi.
        \end{align*}
    as desired. 
\end{proof}

\begin{remark}
    In the schematic formula \eqref{schematic-D}, one should think of the index $j$ as the number of factors of $A$, the index $k$ as the number of derivatives, and the indices $a$ and $b$ denoting the distribution of derivatives between the magnetic potentials and the scalar field respectively. Note the expression $(\nabla^{(\leq a)} (A - B))^{j + 1}$ is rather inefficient with regards to the distribution of derivatives, but it will suffice for our purposes, as it succinctly captures the worst possible scenario where all the derivatives fall on one factor. 
\end{remark}

As a first application, we can show that the magnetic Sobolev norms are equivalent, provided that the differences of the magnetic potentials are sufficiently regular,

\begin{lemma}[Equivalence of magnetic Sobolev norms]
    Let $n \geq 2$ be an integer, $A = A_j \, dx^j$ and $B = B_j \, dx^j$ be connection $1$-forms, and $\Phi \in H^n_\loc (\R^2)$. 
    \begin{itemize}
        \item If $A - B \in H^{n} (\R^2)$, then 
            \begin{equation}\label{eq:equiv-HN}
                \Bigl( 1 + \| A - B \|_{H^{n - 1}}^{n} \Bigr)^{-1}  \| \Phi \|_{H^n_B} \lesssim \| \Phi \|_{H^n_A} \lesssim \Bigl( 1 + \| A - B \|_{H^n}^{n} \Bigr) \| \Phi \|_{H^n_B}.
            \end{equation}

        \item If $A \in L^4 (\R^2)$ and $\nabla A \in H^{n - 2} (\R^2)$, then 
            \begin{equation}\label{eq:equiv-HN-2}
                \Bigl( 1 + \| A \|_{L^4}^{n} + \| \nabla A \|_{H^{n - 2}}^n \Bigr)^{-1}  \| \Phi \|_{H^n} \lesssim \| \Phi \|_{H^n_A} \lesssim \Bigl( 1 + \| A \|_{L^4}^{n} + \| \nabla A \|_{H^{n - 2}}^n \Bigr) \| \Phi \|_{H^n}.
            \end{equation}
    \end{itemize}

\end{lemma}

\begin{proof}[Proof of \eqref{equiv-HN}]
    Appropriately applying H\"older and Sobolev embedding, e.g. \eqref{GNinfty} and \eqref{GN-interpolation}, 
        \begin{align*}
            \bigl\| \bfD_A \Phi \bigr\|_{L^2} 
                &\lesssim \bigl\| \bfD_B \Phi \bigr\|_{L^2} + \bigl\| A - B \bigr\|_{L^4} \bigl\| \Phi \bigr\|_{L^4} \\
                &\lesssim \Bigl( 1 + \bigl\| A - B \bigr\|_{H^1} \Bigr) \bigl\| \Phi \bigr\|_{H^1_B}, \\ 
            \bigl\| \bfD_A^{(2)} \Phi \bigr\|_{L^2} 
                &\lesssim \bigl\| \bfD_B^{(2)} \Phi \bigr\|_{L^2} + \bigl\| A - B \bigr\|_{L^4} \bigl\| \bfD_B \Phi \bigr\|_{L^4} + \bigl\| \nabla(A - B) \bigr\|_{L^2} \bigl\| \Phi \bigr\|_{L^\infty} + \bigl\| A - B \bigr\|_{L^4}^2 \|\Phi \|_{L^\infty} \\
                &\lesssim \Bigl( 1 + \bigl\| A - B \bigr\|_{H^1}^2 \Bigr) \bigl\| \Phi \bigr\|_{H^2_B}.
        \end{align*}
    The reverse inequality is a symmetric calculation. We can conclude \eqref{equiv-HN} for $n = 2$. 

    For higher-orders, assume for induction the result holds up to $n$. Again, by symmetry, it suffices to prove one direction of the inequality \eqref{equiv-HN}. We estimate the schematic formula \eqref{schematic-D}, 
        \begin{align*}
            \bigl\| \bfD^{(n + 1)}_A \Phi \bigr\|_{L^2} 
                &\lesssim \bigl\| \bfD_B^{(n + 1)} \Phi  \bigr\|_{L^2} + \sum_{j + k = n} \sum_{a + b = k} \Bigl\| \bigl( \nabla^{(\leq a)} (A - B) \bigr)^{j + 1} \cdot \bfD_B^{(b)} \Phi \Bigr\|_{L^2}.
        \end{align*}
    It remains to estimate the sum; we work case-by-case depending on how the derivatives are distributed. When all the derivatives fall on the magnetic potentials, i.e. $b = 0$, we may estimate using H\"older and Sobolev embedding \eqref{GNinfty}
        \begin{align*}
            \sum_{j + k = n}\Bigl\| \bigl( \nabla^{(\leq k)} (A - B) \bigr)^{j + 1} \cdot \Phi \Bigr\|_{L^2}
                &\lesssim \| A - B \|_{H^n}^{n + 1} \|\Phi \|_{L^\infty} \\
                &\lesssim \| A - B \|_{H^n}^{n + 1} \| \bfD^{(2)}_B \Phi \|_{L^2}. 
        \end{align*}
    When the derivatives are distributed between the two, i.e. $b \neq 0$, then the number of derivatives hitting the difference of the magnetic potentials is at most $n - 1$. Thus, putting both factors in, say, $L^4$, 
        \begin{align*}
            \sum_{j + k = n} \sum_{\substack{a + b = k \\ b \neq 0}} \Bigl\| \bigl( \nabla^{(\leq a)} (A - B) \bigr)^{j + 1} \cdot \bfD_B^{(b)} \Phi \Bigr\|_{L^2}.
                &\lesssim \sum_{j + k = n} \sum_{\substack{a + b = k \\ b \neq 0}}  \bigl\| (\nabla^{(\leq a)} (A - B))^{j + 1} \bigr\|_{L^4} \bigl\| \bfD^{(b)}_B \Phi \bigr\|_{L^4} \\
                &\lesssim \| A - B \|_{H^n}^{n} \| \Phi \|_{H^{n + 1}_B},  
        \end{align*}
    using Sobolev embedding wherever appropriate. 
\end{proof}

\begin{proof}[Proof of \eqref{equiv-HN-2}]
    Following the proof of \eqref{equiv-HN} and taking $B \equiv 0$, it is not difficult to see the result holds for $n = 2$. For higher-orders, assume for induction the result holds up to $n$. By symmetry, it suffices to prove one direction of the inequality \eqref{equiv-HN-2}. We estimate the schematic formula \eqref{schematic-D}, 
        \begin{align*}
            \bigl\| \bfD^{(n + 1)}_A \Phi \bigr\|_{L^2} 
                &\lesssim \bigl\| \nabla^{(n + 1)} \Phi  \bigr\|_{L^2} + \sum_{j + k = n} \sum_{a + b = k} \Bigl\| \bigl( \nabla^{(\leq a)} A \bigr)^{j + 1} \cdot \nabla^{(b)} \Phi \Bigr\|_{L^2}.
        \end{align*}
    When all the derivatives fall on the scalar field, i.e. $a = 0$, using H\"older and Sobolev embedding yields
        \begin{align*}
            \sum_{j + k = n} \Bigl\| A^{j + 1} \cdot \nabla^{(k)} \Phi \Bigr\|_{L^2}
                &\lesssim \Bigl( 1 + \| A \|_{L^4}^{n + 1} \Bigr) \| \Phi\|_{H^{n + 1}}.
        \end{align*}
    When the derivatives are distributed between the two, i.e. $a \neq 0$, we can put both factors in, say, $L^4$ and estimate using Sobolev embedding
        \begin{align*}
            \sum_{j + k = n}\sum_{\substack{a + b = k \\ a \neq 0}} \Bigl\| \bigl( \nabla^{(\leq a)} A \bigr)^{j + 1} \cdot \nabla^{(b)} \Phi \Bigr\|_{L^2}
                &\lesssim  \sum_{j + k = n}\sum_{\substack{a + b = k \\ a \neq 0}} \bigl\| \bigl( \nabla^{(\leq a)} A \bigr)^{j + 1} \bigr\|_{L^4} \bigl\| \nabla^{(b)} \Phi \bigr\|_{L^4} \\
                &\lesssim \Bigl( 1 + \| \nabla A \|_{H^n}^{n + 1} \Bigr) \| \Phi \|_{H^{n + 1}}.
        \end{align*}
    This completes the proof. 
\end{proof}

\begin{proposition}[Lipschitz continuity of higher-order derivatives]\label{prop:magnetic-lip-2}
    Let $n \in \N_0$ be a non-negative integer, $A \equiv A_j \, dx^j$ and $B \equiv B_j \, dx^j$ be connection $1$-forms such that $A, B \in L^4 (\R^2)$ and $\nabla A, \nabla B \in H^{n - 1} (\R^2)$, and let $\Phi \in X^{n + 1}_A (\R^2)$ and $\Psi \in X^{n + 1}_B (\R^2)$. Suppose also that $A - B \in H^n (\R^2)$ and $\Phi - \Psi \in H^{n + 1} (\R^2)$, then  
        \begin{equation}\label{eq:magn-perturb}
            \bigl\| \bfD_A^{(n + 1)} \Phi - \bfD_B^{(n + 1)}\Psi \bigr\|_{L^2}
                \leq C_n \cdot \Bigl( \|A - B \|_{H^n} + \| \Phi - \Psi \|_{H^{n + 1}} \Bigr),
        \end{equation}
    where the constant depends on 
        \[
            C_n 
                \equiv C_n \bigl( \| A \|_{L^4}, \| \nabla A \|_{H^{n - 1}}, \| A - B \|_{H^{n}}, \| \Phi \|_{X^n_A}\bigr).
        \]
    Furthermore, fixing $(\mathring A, \mathring{\vphantom{A} \Phi})$ satisfying $\mathring A \in L^4 (\R^2)$ and $\nabla \mathring A \in H^{n - 1} (\R^2)$ and $\mathring{\vphantom{A} \Phi} \in H^n_{\mathring A} (\R^2)$, perturbations by perturbations by $(a, u), (b, v) \in (H^n \times H^{n + 1}) (\R^2)$ satisfy 
        \begin{equation}\label{eq:magn-perturb-2}
            \bigl\| \bfD_{\mathring A + a}^{(n + 1)} (\mathring{\vphantom{A} \Phi} + u) - \bfD_{\mathring A + b}^{(n + 1)} (\mathring{\vphantom{A} \Phi} + v) \bigr\|_{L^2} 
                \lesssim C_n\cdot \|(a, u) - (b, v) \|_{H^n \times H^{n + 1}},
        \end{equation}
    where the constant depends on 
        \[
            C_n \equiv C_n \bigl( \| \mathring A \|_{L^4}, \| \nabla \mathring A \|_{H^{n - 1}}, \| \mathring{\vphantom{A} \Phi} \|_{X^{n}_{\mathring A}}, \| (a, u) \|_{H^n \times X^{n}}, \| (b, v) \|_{H^n \times X^{n}}  \bigr).
        \]
\end{proposition}

\begin{proof}
    The bound \eqref{magn-perturb-2} is an consequence of \eqref{magn-perturb}, taking 
        \[
            (A, \Phi) = (\mathring A, \mathring{\vphantom{A} \Phi}) + (a, u), \qquad (B, \Psi) = (\mathring A, \mathring{\vphantom{A} \Phi}) + (b, v),
        \]
    and then carefully estimating the constants appropriately. We leave this easy verification to the reader. 

    By slightly modifying the proof of \eqref{DA-perturb}, placing the scalar field in $L^\infty$ rather than $L^4 + L^\infty$, we may obtain the $n = 0$ case. From here, we proceed by induction, assuming the result up to $n - 1$. Using the schematic identity \eqref{schematic-D}, we may rewrite everything in terms of magnetic derivatives with respect to $B$, which yields 
        \begin{align*}
            \bfD_A^{(n + 1)} \Phi - \bfD_B^{{(n + 1)}} \Psi 
                &= \bfD_B^{(n + 1)} (\Phi - \Psi) + \sum_{j + k = n} \sum_{a + b = k} \bigl( \nabla^{(\leq a)} (A - B) \bigr)^{j + 1} \cdot \bfD_B^{(b)} \Phi.
        \end{align*}
    Putting the above in $L^2$, the first term on the right is acceptable in view of \eqref{equiv-HN}. For the second term on the right, when all the derivatives fall on the magnetic potentials, i.e. $b = 0$, we may estimate using H\"older and Sobolev embedding
        \begin{align*}
            \sum_{j + k = n} \bigl\| \bigl( \nabla^{(\leq k)} (A - B) \bigr)^{j + 1} \cdot \Phi \bigr\|_{L^2}
                &\lesssim \bigl\| A - B \bigr\|_{H^{n}}^{n + 1} \| \Phi \|_{L^\infty}.
        \end{align*}
    When all the derivatives fall on the scalar field, i.e. $a = 0$, we similarly have
        \begin{align*}
            \sum_{j +k = n} \bigl\| (A - B)^{j + 1} \cdot \bfD^{(k)}_B\Phi\bigr\|_{L^2}
                &\lesssim \sum_{j + k = n}\| A - B \|_{L^\infty_x}^{j + 1} \| \bfD^{(k)}_B \Phi \|_{L^2_x} \\
                &\lesssim \Bigl( 1 + \| A - B \|_{H^n_x}^{n + 1} \Bigr) \Bigl( \sum_{k = 1}^n  \| \bfD^{(k)}_B \Phi \|_{L^2_x}\Bigr) \\
                &\lesssim C_n \Bigl( 1 + \| A - B \|_{H^n_x}^{n + 1} \Bigr) \Bigl( \sum_{k = 1}^n  \| \bfD^{(k)}_A \Phi \|_{L^2_x}\Bigr),
        \end{align*}
    using the induction hypothesis in the last line. When at least one derivative falls on both the magnetic potentials and the scalar fields, i.e. $a, b \neq 0$, then at most $n - 1$ derivatives falls on any given factor. Placing both in $L^4$ and using Sobolev embedding, 
        \begin{align*}
            \sum_{j + k = n} \sum_{\substack{a + b = k \\ a, b \neq 0}} \Bigl\|  \bigl( \nabla^{(\leq a)} (A - B) \bigr)^{j + 1} \cdot \bfD_B^{(b)} \Phi \Bigr\|_{L^2} 
                &\lesssim \sum_{j + k = n} \sum_{\substack{a + b = k \\ b \neq 0}}  \bigl\| (\nabla^{(\leq a)} (A - B))^{j + 1}\bigr\|_{L^4}\bigl\| \bfD_B^{(b)} \Phi \bigr\|_{L^{4}} \\
                &\lesssim \Bigl( 1 +  \| A - B \|_{H^{n}}^{n + 1} \Bigr) \sum_{k = 1}^n \bigl\| \bfD^{(k)}_B \Phi \bigr\|_{L^2}\\
                &\lesssim C_n \Bigl( 1 + \| A - B \|_{H^n_x}^{n + 1} \Bigr) \Bigl( \sum_{k = 1}^n  \| \bfD^{(k)}_A \Phi \|_{L^2_x}\Bigr),
        \end{align*}
    again using the induction hypothesis in the last line. Collecting the previous inequalities, we conclude the result. 
\end{proof}

\section{Existence of smooth solutions}\label{app:picard}

To initiate the analysis of the Chern-Simons-Schr\"odinger equation, we took on a technical debt in Section \ref{sec:high-LWP} by assuming the existence of smooth solutions. Being fiscally responsible authors, we conclude the article by finally paying our dues and construct such solutions arising from smooth initial data. We also want the time-interval of existence to depend only on a higher-order covariant energy, which we combined with the covariant energy estimates to reduce the global-in-time regularity of smooth solutions to the continuation criterion \eqref{high-blow-up}. 
 
Herein, we will fix a smooth background configuration and view the solutions as residing within a scale of affine Sobolev spaces. That is, given a smooth configuration $(\mathring A, \mathring{\vphantom{A}\phi}) \in \frE^\infty$ in Coulomb gauge \eqref{coulomb-data}, we consider initial data in the set
    \[
        \mathfrak I^N 
            := \Bigl\{ (A, \phi) \in H^N_{\loc} (\R^2) : (A, \phi) \in (\mathring A, \mathring{\vphantom{A}\phi}) + H^N (\R^2) \text{ and } \Gauss \in H^N (\R^2) \Bigr\},
    \]
and aim to construct solutions in the affine Sobolev space\footnote{Amusingly, while the initial data set $\mathfrak I^N$ is not a closed with respect to the $H^N$-metric, we can check that the Chern-Simons-Schr\"odinger flow nonetheless propagates the regularity $\Gauss \in H^N (\R^2)$ \textit{a posteriori} by virtue of the transport equation \eqref{gauss-transport-intro}.}
    \[
        \mathfrak S^N ([0, T])
            := (\mathring A, \mathring{\vphantom{A}\phi}) + L^\infty_t H^N_x ([0, T] \times \R^2).
    \]
The reader should not take too much stock in these function spaces, as they will serve a mostly qualitative role in our analysis herein. Namely, by working within an affine Sobolev space, we have access to standard functional analysis tools upon passing to the equation for the perturbation of the background configuration. Furthermore, we can use the results of Appendix \ref{app:example} to justify the covariant and magnetic energy estimates. For example, Proposition \ref{prop:magnetic-lip-2} implies that if $(A, \phi) \in \frI^N$, then 
    \[
        \cE^N [A, \phi] 
            < \infty.
    \]

The following existence theory within this scale of function spaces completes the proof of Theorem \ref{thm:smooth-LWP},

\begin{theorem}[$\mathfrak I^{9000}$-existence theory for \eqref{CSS-amp-schro}]\label{thm:picard}
    Let $(\mathring A, \mathring{\vphantom{A}\phi}) \in \frE^\infty$ be a smooth configuration in Coulomb gauge \eqref{coulomb-data}, and suppose that $(A^{\mathrm{in}}, \phi^{\mathrm{in}}) \in \mathfrak I^{9000}$. Set
        \begin{equation}\label{eq:high-data}
            \sfE 
                := \cE^{9000} [A^{\mathrm{in}}, \phi^\mathrm{in}].
        \end{equation}
    Then there exists a positive time $T \equiv T(\sfE)$ and a solution $(A, \phi) \in \mathfrak S^{9000}([0, T])$ to the Chern-Simons-Schr\"odinger equation in DeTurck gauge \eqref{CSS-DT-intro} with initial data $(A, \phi)_{|t = 0} = (A^{\mathrm{in}}, \phi^{\mathrm{in}})$. 

    Furthermore, if $(A^{\mathrm{in}}, \phi^{\mathrm{in}}) \in \frI^\infty$, then $(A, \phi)$ is a smooth solution to the Chern-Simons-Schr\"odinger equation. 
\end{theorem}

Upgrading the existence theory in $\mathfrak S^{9000} ([0, T])$ to the existence of smooth solutions follows by standard persistence of regularity arguments, then using the equation and Sobolev embedding to read-off regularity in $t$ and $x$. We leave the technical verifications to the reader -- one may find the lemmas transfering bounds for magnetic derivatives to non-magnetic derivatives in Appendix \ref{app:example} useful. Our focus herein then will be the construction of solutions in $\mathfrak S^{9000} ([0, T])$ by a Picard iteration scheme. 

Rather than directly working with the Schr\"odinger-heat system \eqref{CSS-DT-intro}, it is convenient to decouple the curl- and divergence-free parts of the magnetic potential,

\begin{lemma}[Equations of motion for $A$ in DeTurck gauge]\label{lem:CSS-DT}
    Let $(A, \phi)$ be a configuration on $[0, T] \times \R^2$ solving the Chern-Simons-Schr\"odinger equation \eqref{CSS-amp-schro} in DeTurck gauge \eqref{DeTurck}. Then the magnetic potential satisfies the heat equation
        \begin{equation}\label{eq:A-DT}
            (\partial_t - \Delta) A_j 
                = - \epsilon_{jk} \Im(\overline \phi \bfD^k \phi) + \Bigl( \Re(\overline \phi \bfD_j \phi) + \partial_j \Gauss \Bigr) + \epsilon_{\ell j} \Bigl( \Re(\overline \phi \bfD^\ell \phi) + \partial^\ell \Gauss \Bigr).
        \end{equation}
    Furthermore, the curl-free and divergence-free parts satisfy the schematic equations 
        \begin{align}
            (\partial_t - \Delta) \PP^\cf A 
                &= \PP^\cf (\phi \cdot \bfD \phi) + \nabla \Gauss,\label{eq:A-DT-cf}\\ 
            \Delta \PP^\df A 
                &= \nabla^\perp \bigl( \tfrac12 (1 - |\phi|^2) + \Gauss \bigr). \label{eq:A-DT-df}
        \end{align}
\end{lemma}

\begin{proof}
    Using the DeTurck gauge, we can write $\partial^\ell F_{\ell j} = \Delta A_j - \partial_j A_t$. Inserting this into Amp\'ere's law, 
        \begin{align*}
            (\partial_t - \Delta) A_j 
                &= F_{tj} - \partial^\ell F_{\ell j} \\
                &= - \epsilon_{jk} \Im(\overline \phi \bfD^k \phi) - \partial_j F_{12} - \partial^\ell F_{\ell j} \\
                &= - \epsilon_{jk} \Im(\overline \phi \bfD^k \phi) + \Bigl( \Re(\overline \phi \bfD_j \phi) + \partial_j \Gauss \Bigr) + \epsilon_{\ell j} \Bigl( \Re(\overline \phi \bfD^\ell \phi) + \partial^\ell \Gauss \Bigr),
        \end{align*}
    rewriting instances of the magnetic field in terms of the scalar field and the Gauss tension field. Projecting onto curl-free vector fields, we obtain the desired parabolic equation \eqref{A-DT-cf} for the curl-free part. The elliptic equation \eqref{A-DT-df} for the divergence-free part is precisely the Biot-Savart law. 
\end{proof}

When designing the Picard iteration scheme, we should keep in mind that the divergence-free part of the magnetic potential is more regular than the curl-free part, and that the Gauss tension field has better regularity than expected from the definition. By virtue of Lemma \ref{lem:CSS-DT}, we can recast a solution $(A, \phi)$ to the Chern-Simons-Schr\"odinger equation \eqref{CSS-amp-schro} in DeTurck gauge \eqref{DeTurck} with initial data $(A, \phi)_{|t = 0} = (A^{\mathrm{in}}, \phi^{\mathrm{in}})$ as a fixed-point of the non-linear mapping 
    \begin{align*}
        \mathfrak S^{9000} ([0, T])
            &\longrightarrow \mathfrak S^{9000} ([0, T]) ,\\
        (A, \phi) 
            &\longmapsto (\widehat A, \widehat \phi),
    \end{align*}
where 
    \[
        \widehat A 
            := \widehat{A^\df} + \widehat{A^{\cf}},
    \]
and the variables $(\widehat{A^\df}, \widehat{A^\cf}, \widehat \phi)$ are defined as the solution to the linear Schr\"odinger-heat-elliptic system 
    \begin{equation}\label{eq:iterationeq}
        \begin{split}
            \Bigl( i \partial_t + \partial^\ell \widehat A_\ell + \bigl( \bfD_{\widehat A} \bigr)^j \bigl( \bfD_{\widehat A}\bigr)_j + \tfrac\lambda2 (1 - |\phi|^2) \Bigr) \widehat\phi
                &= 0 ,\\
            (\partial_t - \Delta) \widehat{A^\cf} 
                &= \PP (\phi \cdot \bfD \phi) + \nabla \Gauss,\\
            \Delta \widehat{A^\df} 
                &= \nabla^\perp \bigl( \tfrac12(1 - |\phi|^2) + \Gauss \bigr),
        \end{split}
    \end{equation} 
subject to initial data 
    \[
        (\widehat{A^\cf}, \widehat \phi)_{|t = 0}
            = (\PP^\cf A^{\mathrm{in}}, \phi^{\mathrm{in}}).
    \]    
We abuse notation, fixing the Gauss tension field here as that for the initial data $\Gauss := \Gauss^{\mathrm{in}}$. 

Define the Picard iterates by the recursive algorithm
    \begin{align*}
        (\widehat A^{[n]}, \widehat\phi^{[n]}) 
            &:= (\widehat{\widehat A^{[n - 1]}}, \widehat{\widehat\phi^{[n - 1]}}),\\
        (\widehat A^{[0]}, \widehat \phi^{[0]})
            &:= (A^{\mathrm{in}}, \phi^{\mathrm{in}}).
    \end{align*}
We will establish the following properties for the Picard iterates: 
    \begin{enumerate}[label=\tt\arabic*.]
        \item \textit{A priori} uniform magnetic energy bounds; defining the magnetic energy functional
            \[
                \slashed \cE^{9000} [A, \phi]
                    := \bigl\| \tfrac12 (1 - |\phi|^2) \bigr\|_{L^2_x}^2 + \sum_{n = 1}^{9000} \bigl\| \bfD_A^{(n)} \phi \bigr\|_{L^2_x}^2,
            \]
        there exists a uniform time $T \ll_\sfE 1$ such that
            \[
                \sup_{t \in [0, T]} \slashed\cE^{9000}[\widehat A^{[n]}, \widehat \phi^{[n]}]
                    \leq 2 \sfE, \qquad \text{for all $n \in \N_0$.}
            \]
        This is the subject of Section \ref{sec:high-energy}.

        \item Existence in the affine Sobolev space $\mathfrak S^{9000} ([0, T])$ and uniform affine bounds 
            \[
                \bigl\| (\widehat A^{[n]}, \widehat \phi^{[n]}) - (\mathring A, \mathring{\vphantom{A}\phi}) \bigr\|_{L^\infty_t H^{9000}_x}
                    \lesssim_{\sfE, (A^{\mathrm{in}}, \phi^{\mathrm{in}}), (\mathring A, \mathring{\vphantom{A}\phi})} 1, \qquad \text{ for all $n \in \N_0$}.
            \]
        This is the subject of Section \ref{sec:high-qualitative}. 
        
        \item Cauchy convergence with respect to the $L^\infty_t (L^2 \times H^1)_x$-topology,
            \[
                \lim_{n, m \to \infty} \bigl\| (\widehat A^{[n]}, \widehat \phi^{[n]}) - (\widehat A^{[m]}, \widehat \phi^{[m]}) \bigr\|_{L^\infty_t (L^2 \times H^1)_x}
                    = 0.
            \]
        This is the subject of Section \ref{sec:high-contract}. 
    \end{enumerate}
Passing to the limit as $n \to \infty$ furnishes a fixed point of \eqref{iterationeq} and thereby a solution to the Chern-Simons-Schr\"odinger equation \eqref{CSS-amp-schro} in DeTurck gauge \eqref{DeTurck}.

\subsection{Uniform \textit{a priori} bounds}\label{sec:high-energy}

Mimicking the covariant energy estimates from Section \ref{sec:high-LWP}, we aim to show that there exists a uniform time $T \ll_\sfE 1$ such that the magnetic energies of the Picard iterates are uniformly bounded. It will also be convenient to prove some bounds for the magnetic potential. 

\begin{proposition}[Energy estimates]\label{prop:iterate-energy}
    Let $(A, \phi)$ be a configuration on $[0, T] \times \R^2$ satisfying
        \begin{equation}\label{eq:prev-iterate}
            \Bigl( i \partial_t + \bigl( \bfD_A\bigr)^j \bigl( \bfD_A \bigr)_j + \mathsf V \Bigr) \phi 
                = 0,
        \end{equation}
    for some real-valued potential $\mathsf V :[0, T] \times \R^2 \to \R$, and the energy bound 
        \begin{equation}\label{eq:prev-iterate-energy}
            \sup_{t \in [0, T]} \slashed\cE^{9000} [A(t), \phi(t)] 
                \leq 2 \sfE. 
        \end{equation}
    Then there exists a sufficiently small time $T \ll_\sfE 1$ such that the solution $(\widehat A, \widehat \phi)$ on $[0, T] \times \R^2$ to the system \eqref{iterationeq} satisfies 
        \begin{itemize}
            \item the magnetic energy bound 
                \begin{equation}\label{eq:iterate-energy}
                      \sup_{t \in [0, T]} \slashed\cE^{9000} [\widehat A(t), \widehat \phi(t)] 
                \leq 2 \sfE. 
                \end{equation}

            \item the top-order Sobolev bound for the magnetic potential, 
                \begin{equation}\label{eq:iterate-A-HN-top}
                    \bigl\| \nabla  \widehat A \bigr\|_{L^\infty_t H^{8999}_x} + \bigl\| \nabla^{(2)} \widehat A \bigr\|_{L^2_t H^{8999}_x}
                \lesssim 1 + \bigl\| \nabla A^{\mathrm{in}} \bigr\|_{H^{8999}_x},
                \end{equation}

            \item the bottom-order affine Sobolev estimate for the magnetic potential, 
                \begin{equation}\label{eq:iterate-A-HN-bot}
                    \bigl\| \widehat A - \mathring A\bigr\|_{L^\infty_t L^2_x} 
                        \lesssim 1 + \bigl\| A^{\mathrm{in}} - \mathring A \bigr\|_{L^2_x}.
                \end{equation}
                \end{itemize}
\end{proposition}

To facilitate the magnetic energy estimates, it is useful to remember that each Picard iterate $\widehat \phi^{[n]}$ for $n \geq 1$ solves an electromagnetic Schr\"odinger equation, so the subsequent iterate $\widehat A^{[n + 1]}$ obeys an analogue of the Gauss constraint equation \eqref{t-gauss}. More precisely, we have the transport equations 

\begin{lemma}[Transport equations]
    Let $(A, \phi)$ be a configuration on $[0, T] \times \R^2$ and let $(\widehat A, \widehat \phi)$ be a solution to the system \eqref{iterationeq}. Then the charge density of $\widehat \phi$ satisfies the density-flux relation
        \begin{align}\label{eq:iterate-charge-transport}
            \partial_t \tfrac12 (1 - |\widehat \phi|^2) 
                = \nabla ( \widehat \phi \cdot \bfD_{\widehat A} \widehat \phi).
        \end{align}
    Furthermore, if $(A, \phi)$ satisfies \eqref{prev-iterate}, then the divergence-free part of $\widehat A$ obeys the transport equation,
        \begin{equation}\label{eq:iterate-transport}
            \partial_t \widehat{A^\df}
                = \PP^\df (\phi \cdot \bfD_A \phi).
        \end{equation}
\end{lemma}

\begin{proof}
    Differentiating the charge density and using the Schr\"odinger equation for $\widehat \phi$ and the product rule, 
        \begin{align*}
            \partial_t \tfrac12(1 - |\widehat \phi|^2) 
                &= - \Re\bigl( \overline{\widehat \phi} \partial_t \widehat \phi \bigr) \\
                &= \Im\Bigl(\overline{\widehat \phi} \bigl( \bfD_{\widehat A} \bigr)^j \bigl( \bfD_{\widehat A} \bigr)_j \widehat \phi \Bigr) \\
                &= \partial^j \Im\Bigl(\overline{\widehat \phi}  \bigl( \bfD_{\widehat A} \bigr)_j \widehat \phi \Bigr) = \nabla \bigl( \widehat \phi \cdot \bfD_{\widehat A} \widehat \phi \bigr),
        \end{align*}
    yielding \eqref{iterate-charge-transport}. 
    
    Differentiating the elliptic equation for the divergence-free part in \eqref{iterationeq}, using the Schr\"odinger equation for $\phi$ and the product rule, 
        \begin{align*}
            \partial_t  \widehat{A^\df} 
                &= \nabla^{-1} \Re(\overline \phi \partial_t \phi) \\
                &=\nabla^{-1} \Im\Bigl( \overline \phi \bigl( \bfD_{A} \bigr)^j \bigl( \bfD_{A}\bigr)_j \phi \Bigr) \\
                &= \nabla^{-1} \partial^j \Im\Bigl( \overline \phi (\bfD_A)_j \phi \Bigr) = \PP^\df (\overline \phi \cdot \bfD_A \phi), 
        \end{align*}
    yielding \eqref{iterate-transport}
\end{proof}

\begin{remark}
    For the first Picard iterate $(\widehat A^{[1]}, \widehat \phi^{[1]})$, instead of \eqref{iterate-transport} we can use the simpler transport equation 
        \[
            \partial_t \widehat A^{[1]} = 0. 
        \]
\end{remark}

As with the covariant energy estimates from Section \ref{sec:high-LWP}, the Schr\"odinger equation in \eqref{iterationeq} has favourable energy structure when one commutes with magnetic derivatives rather than usual derivatives, since the divergence-free part of the magnetic potential has better regularity compared to the curl-free part.  

\begin{lemma}[Higher-order equations of motion]
    Let $(A, \phi)$ be a configuration on $[0, T] \times \R^2$ satisfying \eqref{prev-iterate}, and suppose $(\widehat A, \widehat \phi)$ is a solution to the system \eqref{iterationeq}. Then, for each positive integer $n \in \N$, 
    \begin{itemize}
        \item the scalar field obeys the Schr\"odinger equation,
            \begin{equation}\label{eq:iterate-schro}
            \begin{split}
                \Bigl( i \partial_t + \partial^\ell \widehat A_\ell + \bigl( \bfD_{\widehat A} \bigr)^j \bigl( \bfD_{\widehat A}\bigr)_j + \tfrac\lambda2 (1 - |\phi|^2) \Bigr) \bfD_{\widehat A}^{(n)} \widehat\phi
                    &= \bfD_{\widehat A}^{(n)} \widehat \phi + \sum_{a + b + c = n} \PP \bigl(\bfD_A^{(a)} \phi \cdot \bfD_A^{(b)} \phi\bigr) \cdot \bfD_{\widehat A}^{(c)} \widehat \phi \\
                        &\qquad + \sum_{a + b = n} \PP \bigl(\nabla^{(a)} \Gauss\bigr) \cdot \bfD_A^{(b)} \widehat \phi,
            \end{split}
            \end{equation}
        where $\PP$ is a Fourier multiplier of order zero,
        
        \item the curl-free part of the magnetic potential obeys the heat equation,
            \begin{equation}\label{eq:iterate-heat}
                (\partial_t - \Delta) \nabla^{(n)}  \widehat{A^\cf} 
                    = \nabla \Bigl( \sum_{a + b = n} \bfD_A^{(a)} \phi \cdot \bfD_A^{(b)} \phi + \nabla^{(n)} \Gauss \Bigr),
            \end{equation}
        
        \item the divergence-free part of the magnetic potential obeys
            \begin{equation}\label{eq:iterate-elliptic}
            \begin{split}
                \nabla  \widehat{A^\df} 
                    &= \tfrac12(1 - |\phi|^2) + \Gauss, \\ 
                \nabla^{(n)}  \widehat{A^\df} 
                    &= 
                        \sum_{a + b = n - 1} \bfD^{(a)}_A \cdot \bfD^{(b)}_A \phi + \nabla^{(n - 1)} \Gauss, \qquad \text{for $n \geq 2$.}
            \end{split}
            \end{equation}

    \end{itemize}
\end{lemma}

\begin{proof}[Derivation of Schr\"odinger equation \eqref{iterate-schro}]
    Commuting,  
        \begin{align*}
            \Bigl[ \bfD_{\widehat A} , \bigl( \bfD_{\widehat A} \bigr) (\bfD_{\widehat A}) + \tfrac\lambda2 (1 - |\phi|^2) \Bigr]
                &= \widehat F_{12} \cdot \bfD_{\widehat A} + \nabla \widehat F_{12} + \phi \cdot \bfD_A \phi \\
                &= \bfD_{\widehat A} + \phi \cdot \phi \cdot \bfD_{\widehat A} + \Gauss \cdot \bfD_{\widehat A} + \nabla \Gauss + \phi \cdot \bfD_A \phi,
            \end{align*}
    using the elliptic equation for the divergence-free part to rewrite the magnetic field. The commutator with the magnetic time derivative is given by  
        \begin{align*}
            \Bigl[ (\bfD_{\widehat A})_k, i \partial_t + \partial^\ell \widehat A_\ell \Bigr]
                &= \partial_t \widehat A_k -  \partial_k \partial^\ell \widehat A_\ell  \\
                &= (\partial_t - \Delta) (\PP^\cf_k \widehat A + \PP^\df_k \widehat A) + \epsilon_{jk} \partial^j \widehat F_{12} \\
                &= \PP (\phi \cdot \bfD_A \phi) + \PP \nabla \Gauss,
            \end{align*}
    using the heat equation for the curl-free part, the elliptic equation for the divergence-free part from \eqref{iterationeq}, and the transport equation \eqref{iterate-transport} for the divergence-free part to re-express its time-derivative.

    With these two commutator identities at hand, we can easily show \eqref{iterate-schro} by induction. Indeed, the case $n = 1$ is immediate from commuting $\bfD_{\widehat A}$ into \eqref{iterationeq}; proceeding inductively, we apply $\bfD_{\widehat A}$ to the equation for $n$. The right-hand side is acceptable due to the product rule, while the commutators from the left are also easily seen to be of the desired form by the calculations above. 
\end{proof}

We are now in a position to prove Proposition \ref{prop:iterate-energy}. 

\begin{proof}[Proof of magnetic energy bound \eqref{iterate-energy}]
    The proof will follow a similar argument to that of the covariant energy estimate \eqref{EN-bound-1}, plus estimating the $L^\infty$-norms of the scalar fields by the Sobolev-type inequality \eqref{inefficientLinfty}. For these reason, we will be brief with writing out all the details. Making the bootstrap assumption
        \begin{equation}\label{eq:BS-iterate}
            \sup_{t \in [0, T]}\slashed\cE^N [\widehat A(t), \widehat \phi(t)] 
                \leq 10\sfE,
        \end{equation}
    our goal is to improve the bound to $2\sfE$. Throughout, we denote $C(\sfE)$ for some polynomial expression in $\sfE$, the definition of which may change line-by-line.
    
    Multiplying the transport equation \eqref{iterate-charge-transport} by the charge density, and integrating-by-parts, we obtain
        \begin{equation}\label{eq:iterate-charge-est}
        \begin{split}
            \bigl\| \tfrac12 (1 - |\widehat \phi|^2) \bigr\|_{L^\infty_t L^2_x}^2
                &\leq \bigl\| \tfrac12(1 - |\phi^{\mathrm{in}}|^2) \bigr\|_{L^2_x}^2 + T \cdot \bigl\| \widehat \phi \bigr\|_{L^\infty_{t, x}}^2 \bigl\| \bfD_{\widehat A} \widehat \phi \bigr\|_{L^\infty_t L^2_x}^2 \\
                &\leq \bigl\| \tfrac12(1 - |\phi^{\mathrm{in}}|^2) \bigr\|_{L^2_x}^2 + T \cdot C(\sfE) ,
        \end{split}
        \end{equation}
    using the bootstrap assumption \eqref{BS-iterate} and the Sobolev-type inequality \eqref{inefficientLinfty} in the second line. 

    To control magnetic energy, we apply the abstract energy estimate \eqref{abstract-energy} to the Schr\"odinger equation \eqref{iterate-schro}, and estimate the right-hand side. The linear term is easily treated, while for the non-linear terms, we can use H\"older's inequality, the magnetic interpolation inequalities \eqref{GN-interpolation}, and the Sobolev-type inequality \eqref{inefficientLinfty}, and the bootstrap assumption \eqref{BS-iterate},
        \begin{align*}
            \bigl\| \text{R.H.S. \eqref{iterate-schro}} \bigr\|_{L^2_x} 
                &\lesssim \bigl\| \bfD^{(n)}_{\widehat A} \widehat \phi \bigr\|_{L^2_x} + \sum_{a + b + c = n}\bigl\| \bfD^{(a)}_A \phi \bigr\|_{L^{2n/a}_x} \, \bigl\| \bfD^{(b)}_A \phi \bigr\|_{L^{2n/b}_x} \, \bigl\| \bfD^{(c)}_{\widehat A} \widehat \phi \bigr\|_{L^{2n/c}_x}\\
                &\qquad +  \sum_{a + b = n} \bigl\| \nabla^{(a)} \Gauss  \bigr\|_{L^{2n/a}_x} \bigl\| \bfD_{\widehat A} ^{(b)} \widehat \phi  \bigr\|_{L^{2n/b}_x} \\
                &\lesssim C(\sfE).
        \end{align*}
    Inserting this into the energy estimate \eqref{abstract-energy} for \eqref{iterate-schro},
        \begin{equation} \label{eq:iterate-Dphi-est}
        \begin{split}
            \bigl\| \bfD_{\widehat A}^{(n)} \widehat \phi \bigr\|_{L^\infty_t L^2}^2 
                &\leq \bigl\| \bfD^{(n)}_{A^{\mathrm{in}}} \phi^{\mathrm{in}} \bigr\|_{L^2_x}^2  + T \cdot  \bigl\| \bfD_{\widehat A}^{(n)} \widehat \phi \bigr\|_{L^\infty_t L^2_x} \bigl\| \text{R.H.S. \eqref{iterate-schro}} \bigr\|_{L^\infty_t L^2_x} \\ 
                &\leq \bigl\| \bfD^{(n)}_{A^{\mathrm{in}}} \phi^{\mathrm{in}} \bigr\|_{L^2_x}^2  + T \cdot C(\sfE). 
        \end{split}
        \end{equation}

    Putting together the estimates \eqref{iterate-charge-est} and \eqref{iterate-Dphi-est} for $n = 1, \dots, 9000$, we obtain 
        \begin{align*}
            \sup_{t \in [0, T]}\slashed\cE^{9000} [\widehat A(t), \widehat \phi(t)]
                &\leq \slashed\cE^{9000}  [A^{\mathrm{in}}, \phi^{\mathrm{in}}] + T \cdot C(\sfE) \leq 2 \sfE,
        \end{align*}
    taking $T \ll_\sfE 1$. This closes the bootstrap \eqref{BS-iterate}. 
\end{proof}

\begin{proof}[Proof of magnetic potential bounds \eqref{iterate-A-HN-top}-\eqref{iterate-A-HN-bot}]
    Applying the energy estimate \eqref{heat-energy-estimate} to the heat equation for the curl-free part \eqref{iterate-heat}, using interpolation \eqref{GN-interpolation} and the Sobolev-type estimate \eqref{inefficientLinfty} wherever appropriate, 
        \begin{align*}
            \bigl\| \nabla^{(n)}  \widehat{A^\cf} \bigr\|_{L^\infty_t L^2_x} + \bigl\| \nabla^{(n + 1)}  \widehat{A^\cf} \bigr\|_{L^2_{t, x}}
                &\lesssim \bigl\| \nabla^{(n)} \PP^\cf A^{\mathrm{in}} \bigr\|_{L^2_x} \\
                    &\qquad + \Bigl\| \sum_{a + b = n} \bigl\| \bfD^{(a)}_A \phi  \bigr\|_{L^{2n/a}_x} \bigl\| \bfD_A^{(b)} \phi \bigr\|_{L^{2n/b}_x} + \bigl\| \nabla^{(n)} \Gauss \bigr\|_{L^2_x}\Bigr\|_{L^2_t} \\
                    &\lesssim \bigl\| \nabla^{(n)} \PP^\cf A^{\mathrm{in}} \bigr\|_{L^2_x} \\
                        &\qquad+ T^\frac12 \Bigl( \bigl( \| \phi \|_{L^\infty_{t, x}} + \| \bfD_A \phi \|_{L^\infty_t L^2_x} \bigr) \| \bfD^{(n)}_A \phi \|_{L^\infty_t L^2_x} + \bigl\| \nabla^{(n)} \Gauss \bigr\|_{L^2_x}\Bigr) \\
                    &\lesssim \bigl\| \nabla^{(n)} \PP^\cf A^{\mathrm{in}} \bigr\|_{L^2_x} + T^\frac12 \cdot C(\sfE)
        \end{align*}
    for each $n = 1, \dots, N$. Similarly estimating the expressions \eqref{iterate-elliptic} for the divergence-free part gives 
        \begin{align*}
            \bigl\| \nabla \widehat{A^\df} \bigr\|_{L^\infty_t L^2_x}
                &\lesssim\bigl\| \tfrac12 (1 - |\phi|^2) \bigr\|_{L^\infty_t L^2_x} + \| \Gauss \|_{L^2_x} \\
                &\lesssim C(\sfE), \\
            \bigl\| \nabla^{(n)} \nabla  \widehat{A^\df} \bigr\|_{L^\infty_t L^2_x}
                &\lesssim \sum_{a + b = n} \bigl\| \bfD^{(a)}_A \phi \bigr\|_{L^\infty_t L^{2n/a}_x} \bigl\| \bfD^{(b)}_A \phi \bigr\|_{L^\infty_t L^{2n/b}_x} + \bigl\| \nabla^{(n)} \Gauss \bigr\|_{L^2_x} \\
                &\lesssim \bigl( \| \phi \|_{L^\infty_t L^2_x} + \| \bfD_A \phi \|_{L^\infty_t L^2_x} \bigr) \| \bfD^{(n)}_A \phi \|_{L^\infty_t L^2_x} + \| \nabla^{(n)} \Gauss \|_{L^2_x} \\
                &\lesssim C(\sfE),
        \end{align*}
    for each $n = 1, \dots, 9000$. Collecting these inequalities and taking $T \ll_\sfE 1$ gives the top-order estimate \eqref{iterate-A-HN-top}. 

    Integrating the transport equation \eqref{iterate-transport} and the heat equation \eqref{iterate-heat} in time and estimating the right-hand side using the magnetic energy estimates \eqref{iterate-energy} and Sobolev embedding,
        \begin{align*}
            \bigl\| \widehat A - \mathring A \bigr\|_{L^\infty_t L^2_x} 
                &\lesssim \bigl\| A^{\mathrm{in}} - \mathring A \bigr\|_{L^2_x} +\bigl\| \Delta \PP^{\cf} \widehat A \bigr\|_{L^1_t L^2_x} + \bigl\| \phi \cdot \bfD_A \phi \bigr\|_{L^1_t L^2_x} \\
                &\lesssim  \bigl\| A^{\mathrm{in}} - \mathring A \bigr\|_{L^2_x} + T \cdot \Bigl( \bigl\| \nabla^{(2)} \PP^{\cf} \widehat A \bigr\|_{L^\infty_t L^2_x} +  \| \phi \|_{L^\infty_{t, x}} \| \bfD_A \phi \|_{L^\infty_t L^2_x} \Bigr) \\
                &\lesssim  \bigl\| A^{\mathrm{in}} - \mathring A \bigr\|_{L^2_x} + T \cdot C(\sfE).
        \end{align*}
    using the top-order bound \eqref{iterate-A-HN-top} to control the derivatives of the curl-free part of the magnetic potential. Taking $T \ll_\sfE 1$, we conclude the bottom-order affine Sobolev estimate \eqref{iterate-A-HN-bot}.
\end{proof}

\subsection{Existence of iterates and affine Sobolev bounds} \label{sec:high-qualitative}

The estimates in Proposition \ref{prop:iterate-energy} are \textit{a priori} estimates, that is, assuming the existence of a sufficiently regular solution to the system \eqref{iterationeq}. Nevertheless, it is a matter of principle that quantitative \textit{a priori} estimates for a linear system of equations can be turned into qualitative existence theory. It will be convenient to introduce the additional control parameter
    \[
        \sfA
            := \cE^{9100} [\mathring A, \mathring{\vphantom{A}\phi}] + \bigl\| A^{\mathrm{in}} - \mathring A \bigr\|_{H^{9000}_x}.
    \]
Then we claim that 

\begin{proposition}[Existence of Picard iterates]\label{prop:linear-LWP}
    Let $(A, \phi)$ be a configuration on $[0, T] \times \R^2$ satisfying the energy bound \eqref{prev-iterate-energy}. Then there exists $T \ll_\sfE 1$ and a solution $(\widehat A, \widehat \phi) \in \mathfrak S^{9000} ([0, T])$ to the system \eqref{iterationeq} satisfying 
        \begin{equation}\label{eq:iterate-phi-affine}
            \bigl\| \widehat \phi- \mathring{\vphantom{A}\phi} \bigr\|_{L^\infty_t H^{9000}_x}
                \lesssim_{\sfE, \sfA} \bigl\| \phi^{\mathrm{in}} - \mathring{\vphantom{A}\phi} \bigr\|_{H^{9000}_x}.
        \end{equation}
    and 
        \begin{equation}\label{eq:iterate-A-affine}
            \bigl\| \widehat A - \mathring A\bigr\|_{L^\infty_t H^{9000}_x} + \bigl\| \nabla(\widehat A - \mathring A)\bigr\|_{L^2_t H^{9000}_x}
                \lesssim_{\sfE, \sfA} 1.
        \end{equation}

\end{proposition}

One can easily construct the magnetic potential $\widehat A$ by using the fundamental solution for the linear heat equation and the Laplacian. The bound \eqref{iterate-A-affine} follows from the triangle inequality, the top-order bound \eqref{iterate-A-HN-top}, the bottom-order bound \eqref{iterate-A-HN-bot}, and the bound on the smooth background (c.f. the proof of \eqref{A-Cinfty}), 
        \begin{equation} \label{eq:background-bound}
            \| \mathring A \|_{L^4} + \| \nabla \mathring A \|_{H^{9050}}
                \lesssim_\sfA 1.
        \end{equation}
It remains to show existence of the scalar field $\widehat \phi$ and the bound \eqref{iterate-phi-affine}. We first pass to the equation for the perturbation $u := \widehat \phi - \mathring{\vphantom{A}\phi}$, which is a straightforward calculation, 

\begin{lemma}[Equation for the perturbation]
    Let $(\widehat A, \widehat \phi)$ be the configuration on $[0, T] \times \R^2$ solving the system \eqref{iterationeq}, and set $u := \widehat \phi - \mathring{\vphantom{A}\phi}$. Then the perturbation $u : [0, T] \times \R^2 \to \C$ satisfies the linear inhomogeneous Schr\"odinger equation 
        \begin{equation}\label{eq:perturb-base}
            \Bigl( i \partial_t + \partial^\ell \widehat A_\ell + \bigl( \bfD_{\widehat A} \bigr)^j \bigl( \bfD_{\widehat A}\bigr)_j + \tfrac\lambda2 (1 - |\phi|^2) \Bigr) u 
                = \Bigl( \partial^\ell \widehat A_\ell + \bigl( \bfD_{\widehat A} \bigr)^j \bigl( \bfD_{\widehat A}\bigr)_j + \tfrac\lambda2 (1 - |\phi|^2) \Bigr) \mathring{\vphantom{A}\phi}, 
        \end{equation}
    Furthermore, for each positive integer $n \in \N$, we have the higher-order equation 
        \begin{align}
            \Bigl( i \partial_t + \partial^\ell \widehat A_\ell + \bigl( \bfD_{\widehat A} \bigr)^j \bigl( \bfD_{\widehat A}\bigr)_j + \tfrac\lambda2 (1 - |\phi|^2) \Bigr) \bfD_{\widehat A}^{(n)} u
                &= \bfD_{\widehat A}^{(n)} u + \sum_{a + b + c = n} \PP \bigl(\bfD_A^{(a)} \phi \cdot \bfD_A^{(b)} \phi\bigr) \cdot \bfD_{\widehat A}^{(c)} u \notag\\
                    &\qquad + \sum_{a + b = n} \PP \nabla^{(a)} \Gauss \cdot \bfD_{\widehat A}^{(b)} u \label{eq:perturb-high} \\
                    &\qquad + \bfD^{(n)}_{\widehat A} \Bigl( \partial^\ell \widehat A_\ell + \bigl( \bfD_{\widehat A} \bigr)^j \bigl( \bfD_{\widehat A}\bigr)_j + \tfrac\lambda2 (1 - |\phi|^2) \Bigr) \mathring{\vphantom{A}\phi} ,\notag
        \end{align}
    where $\PP$ is a Fourier multiplier of order zero. 
\end{lemma}

It will be convenient to record the equivalence of non-magnetic Sobolev norms and magnetic Sobolev norms, with implicit constants depending on the control parameter $\sfA$, 

\begin{lemma}[Equivalence of Sobolev norms]
    Suppose $\widehat A \in \mathring A + L^\infty_t H^N_x (\R^2)$, then 
        \begin{align}\label{eq:iterate-equiv}
            \| u \|_{L^\infty_t H^N_{\widehat A}} 
                \sim_{\sfE, \sfA} \|u \|_{L^\infty_t H^N_x} \sim_\sfA \| u \|_{L^\infty_t H^N_{\mathring A}}.
        \end{align} 
\end{lemma}

\begin{proof}
    This follows from \eqref{equiv-HN}-\eqref{equiv-HN-2}, using \eqref{iterate-A-affine} and \eqref{background-bound} to control the implicit constant. 
\end{proof}

We are now in a position to complete the proof of Proposition \ref{prop:linear-LWP}. It suffices to prove the \textit{a priori} estimate \eqref{iterate-phi-affine} for the perturbation $u$, as local existence follows by textbook arguments, e.g. duality \cite[Section 3]{IfrimTataru2022} or the Galerkin method \cite[Chapter 8]{Metivier2008}. Applying the abstract energy estimate \eqref{abstract-energy} to the equation \eqref{perturb-high} yields for each $n = 0, \dots, 9000$, 
        \begin{align*}
            \bigl\| \bfD_{\widehat A}^{(n)} (\widehat \phi - \mathring{\vphantom{A}\phi}) \bigr\|_{L^\infty_t L^2_x}
                &\lesssim \bigl\| \bfD_{\widehat A}^{(n)} (\phi^{\mathrm{in}} - \mathring{\vphantom{A}\phi})\bigr\|_{L^\infty_t L^2_x} + \bigl\| \text{R.H.S \eqref{perturb-high}} \bigr\|_{L^1_t L^2_x}.
        \end{align*}
    Following a similar argument to the proof of the magnetic energy estimate \eqref{iterate-energy}, we can absorb the first three terms on the right-hand side of the equation \eqref{perturb-high} into the left. Combined with \eqref{iterate-equiv},
        \[
            \bigl\| \widehat \phi - \mathring{\vphantom{A}\phi} \bigr\|_{L^\infty_t H^{9000}_x} 
                \lesssim_{\sfA, \sfE} \bigl\| \phi^{\mathrm{in}} - \mathring{\vphantom{A}\phi} \bigr\|_{L^\infty_t H^N_x} + \Bigl\| \Bigl( \partial^\ell \widehat A_\ell + \bigl( \bfD_{\widehat A} \bigr)^j \bigl( \bfD_{\widehat A}\bigr)_j + \tfrac\lambda2 (1 - |\phi|^2) \Bigr) \mathring{\vphantom{A}\phi} \Bigr\|_{L^1_t H^{9000}_{\mathring A}}. 
        \]
    It remains to control the contribution of this forcing term. Rewriting the magnetic Laplacian with respect to $\widehat A$ in terms of that with respect to $\mathring A$ and the differences, 
        \begin{align*}
            \bigl( \bfD_{\widehat A} \bigr)^j \bigl( \bfD_{\widehat A}\bigr)_j \mathring{\vphantom{A}\phi}
                &= \bigl( \bfD_{\mathring A} \bigr)^j \bigl( \bfD_{\mathring A}\bigr)_j \mathring{\vphantom{A}\phi} - 2i (\widehat A - \mathring A)^j \bigl( \bfD_{\mathring A} \bigr)_j \mathring{\vphantom{A}\phi} - 2 (\widehat A - \mathring A)^j \mathring A_j \mathring{\vphantom{A}\phi}- i \partial^j (\widehat A - \mathring A)_j \mathring{\vphantom{A}\phi} - (\widehat A^j \widehat A_j - \mathring A^j \mathring A_j) \mathring{\vphantom{A}\phi}.
        \end{align*}  
    Let us consider the roughest terms, and leave the lower-order terms to the reader. One possibility is that all the derivatives fall on the magnetic Laplacian; this is clearly controlled by $\sfA$, 
        \[
            \bigl\| \bfD_{\mathring A}^{(9000)} \bigl( \bfD_{\mathring A} \bigr)^j \bigl( \bfD_{\mathring A}\bigr)_j \mathring{\vphantom{A}\phi}  \bigr\|_{L^1_t L^2_x}
                \lesssim_{\sfA} T.
        \]
    Another possibility is that all derivatives fall on $\nabla (\widehat A - \mathring A)$. Estimating $\mathring{\vphantom{A}\phi}$ in $L^\infty_x$-norm using \eqref{inefficientLinfty} and the difference using the $L^2_{t, x}$-bound from \eqref{iterate-A-affine},
        \[
            \bigl\| \nabla^{(9000)} \nabla(\widehat A - \mathring A) \cdot \mathring{\vphantom{A}\phi} \bigr\|_{L^1_t L^2_x}
                \lesssim T^{\frac12} \bigl\| \nabla^{(9000)} \nabla(\widehat A - \mathring A) \cdot \mathring{\vphantom{A}\phi} \bigr\|_{L^2_{t, x}} \bigl\| \mathring{\vphantom{A}\phi} \bigr\|_{L^\infty_x} \lesssim_{\sfE, \sfA} T^{\frac12}.
        \]
    The lower-order terms follow a similar circle of ideas. We conclude
        \[
            \Bigl\| \Bigl( \partial^\ell \widehat A_\ell + \bigl( \bfD_{\widehat A} \bigr)^j \bigl( \bfD_{\widehat A}\bigr)_j + \tfrac\lambda2 (1 - |\phi|^2) \Bigr) \mathring{\vphantom{A}\phi} \Bigr\|_{L^1_t H^N_{\mathring A}}
                \lesssim_{\sfA, \sfE} 1. 
        \]
    This completes the proof.

\subsection{Convergence of iteration scheme}
\label{sec:high-contract}

To prove convergence of the Picard iteration, it suffices to prove that the contractive property of the iteration map, 

\begin{proposition}[Contractive property of iteration]\label{prop:high-contract}
    Let $(A, \phi)$ and $(B, \psi)$ be configurations on $[0, T] \times \R^2$ satisfying the equation \eqref{prev-iterate} and the bounds \eqref{prev-iterate-energy}, and suppose $(\widehat A, \widehat \phi)$ and $(\widehat B, \widehat \psi)$ are the corresponding solutions to the system \eqref{iterationeq}. Then there exists a sufficiently small $T \ll_\sfE 1$ such that 
        \begin{equation}
            \bigl\| (\widehat A ,\widehat \phi) - (\widehat B, \widehat \psi) \bigr\|_{L^\infty_t (L^2 \times H^1)_x}
                \leq \tfrac12 \bigl\| (A, \phi) - (B, \psi) \bigr\|_{L^\infty_t (L^2 \times H^1)_x}.
        \end{equation}
\end{proposition}

To facilitate the proof, we derive the equations of motion for the differences of configurations. In line with the energy estimates, it will be convenient to consider the differences of the magnetic derivatives, rather than the non-magnetic derivatives,

\begin{lemma}[Difference equations]
    Let $(A, \phi)$ and $(B, \psi)$ be configurations on $[0, T] \times \R^2$ satisfying \eqref{prev-iterate}, and suppose $(\widehat A, \widehat \phi)$ and $(\widehat B, \widehat \psi)$ are solutions to the system \eqref{iterationeq}. Then 
    \begin{itemize}
        \item the scalar field differences obey the Schr\"odinger equations
            \begin{align}
            \begin{split}
                \Bigl( i \partial_t + \partial^\ell \widehat A_\ell &+ \bigl( \bfD_{\widehat A} \bigr)^j \bigl( \bfD_{\widehat A} \bigr)_j + \tfrac\lambda2 (1 - |\phi|^2) \Bigr) \bigl(\widehat\phi - \widehat\psi\bigr) \\
                    &= \Bigl( \nabla \PP^\cf \bigl(\widehat A - \widehat B\bigr) + \bigl( \widehat A \cdot \widehat A - \widehat B \cdot \widehat B \bigr)+ \bigl( \phi \cdot \phi - \psi \cdot \psi \bigr)\Bigr) \cdot \widehat \psi \\
                    &\qquad + \bigl( \widehat A - \widehat B \bigr) \cdot \nabla \widehat \psi,
            \end{split}\label{eq:iterate-schro-diff}
            \end{align}
            and 
            \begin{equation}
                \begin{split}
                \Bigl( i \partial_t + \partial^\ell \widehat A_\ell &+ \bigl( \bfD_{\widehat A} \bigr)^j \bigl( \bfD_{\widehat A} \bigr)_j + \tfrac\lambda2 (1 - |\phi|^2) \Bigr) \bigl(\bfD_{\widehat A} \widehat \phi - \bfD_{\widehat B} \widehat\psi\bigr)\\
                    &= \Bigl( \nabla \PP^\cf \bigl(\widehat A - \widehat B\bigr) + \bigl( \widehat A \cdot \widehat A - \widehat B \cdot \widehat B \bigr)+ \bigl( \phi \cdot \phi - \psi \cdot \psi \bigr)\Bigr) \cdot \bfD_{\widehat B} \widehat \psi \\
                    &\qquad + \bigl( \widehat A - \widehat B \bigr) \cdot \nabla \bfD_{\widehat B} \widehat\psi + \Bigl( 1 + \phi \cdot \phi + \Gauss \Bigr) \cdot \Bigl( \bfD_{\widehat A}  \widehat \phi - \bfD_{\widehat B} \widehat \psi \Bigr) \\
                    &\qquad + \nabla \Gauss \cdot \bigl( \widehat \phi - \widehat \psi \bigr) + \Bigl( \phi \cdot \bfD_A \phi \cdot \widehat \phi - \psi \cdot \bfD_B \psi \cdot \widehat \psi \Bigr),
            \end{split}\label{eq:iterate-schro-diff-2}
            \end{equation} 

        \item the curl-free differences of the magnetic potentials obey the heat equations
            \begin{align}
                 (\partial_t - \Delta) \PP^\cf \bigl(\widehat A - \widehat B\bigr)
                    &= \PP \Bigl( (\phi - \psi) \cdot \bfD_A \phi + \psi \cdot (\bfD_A \phi - \bfD_B \psi) \Bigr),\label{eq:iterate-heat-diff}
            \end{align}

        \item the divergence-free differences of the magnetic potentials obey
            \begin{align}
            \partial_t \PP^\df \bigl(\widehat A - \widehat B\bigr) 
                &= \PP \Bigl( (\phi - \psi) \cdot \bfD_A \phi + \psi \cdot(\bfD_A \phi - \bfD_B \psi) \Bigr).\label{eq:iterate-transport-diff}
            \end{align}

    \end{itemize}
\end{lemma}

\begin{proof}[Derivation of equations \eqref{iterate-schro-diff}-\eqref{iterate-schro-diff-2} for scalar field differences]
    Taking the differences of the operators on the left-hand sides, 
        \begin{align*}
            \Bigl( \partial^\ell \widehat A_\ell &+ \bigl( \bfD_{\widehat A} \bigr)^j \bigl( \bfD_{\widehat A} \bigr)_j + \tfrac\lambda2 (1 - |\phi|^2) \Bigr) - \Bigl( \partial^\ell \widehat B_\ell + \bigl( \bfD_{\widehat B} \bigr)^j \bigl( \bfD_{\widehat B} \bigr)_j + \tfrac\lambda2 (1 - |\psi|^2) \Bigr) \\
                &= \Bigl( \partial^\ell (\widehat A_\ell - \widehat B_\ell) + i \partial^j (\widehat B_j - \widehat A_j) + (\widehat B^j \widehat B_j - \widehat A^j \widehat A_j) \Bigr)  + 2i \bigl(\widehat B^j - \widehat A^j\bigr) \partial_j + \tfrac\lambda2(|\psi|^2 - |\phi|^2) .
        \end{align*}
    Subtracting the equations \eqref{iterate-schro} for $(\widehat A, \widehat \phi)$ and $(\widehat B, \widehat \psi)$ for $n = 0, 1$ respectively and using the identity above gives the desired equations. 
\end{proof}

\begin{proof}[Derivation of equations \eqref{iterate-heat-diff}-\eqref{iterate-transport-diff} for divergence-free differences]
    These follow from subtracting the respective heat equations \eqref{iterate-heat} and transport equations \eqref{iterate-transport}.
\end{proof}

\begin{proof}[Proof of Proposition \ref{prop:high-contract}]
    It will suffice to show a bound of the form
        \begin{align}\label{eq:iterate-magn-diff}
            \bigl\| (\widehat A, \widehat \phi) - (\widehat B, \widehat \psi) \bigr\|_{L^\infty_t L^2_x} + \bigl\| \bfD_{\widehat A} \widehat \phi - \bfD_{\widehat B} \widehat \psi \bigr\|_{L^\infty_t L^2_x} 
                \lesssim_\cE T \cdot \Bigl( \bigl\|  \phi - \psi \bigr\|_{L^\infty_t L^2_x} + \bigl\| \bfD_A  \phi - \bfD_{B} \psi \bigr\|_{L^\infty_t L^2_x} \Bigr).
        \end{align}
    Indeed, the left-hand side controls the $(L^2 \times H^1)$-difference, writing
        \begin{align*}
            \bigl\| \nabla(\widehat \phi - \widehat \psi) \bigr\|_{L^2}
                &\leq \bigl\| \widehat A \cdot \widehat \phi - \widehat B \cdot \widehat \psi \bigr\|_{L^2} +  \bigl\| \bfD_{\widehat A} \widehat \phi - \bfD_{\widehat B} \widehat \psi \bigr\|_{L^2} \\
                &\leq \bigl\| \widehat A \bigr\|_{L^\infty} \bigl\| \widehat \phi - \widehat \psi \bigr\|_{L^2} + \bigl\| \widehat A - \widehat B \bigr\|_{L^2}\bigl\| \widehat \psi \bigr\|_{L^\infty_x} +  \bigl\| \bfD_{\widehat A} \widehat \phi - \bfD_{\widehat B} \widehat \psi \bigr\|_{L^2}\\
                &\lesssim_\sfE \bigl\| \widehat \phi - \widehat \psi \bigr\|_{L^2}  + \bigl\| \bfD_{\widehat A} \widehat \phi - \bfD_{\widehat B} \widehat \psi \bigr\|_{L^2} + \bigl\| \widehat A - \widehat B \bigr\|_{L^2_x},
        \end{align*}
    and on the right-hand side, we may bound the difference of the magnetic derivatives by
        \begin{align*}
            \| \phi - \psi \|_{L^2} + \bigl\| \bfD_A \phi - \bfD_B \psi \bigr\|_{L^2}
                &\leq \| \phi - \psi\|_{H^1} + \| A \cdot \phi - B \cdot \psi \|_{L^2} \\
                &\leq \| \phi - \psi \|_{H^1} + \| A \|_{L^\infty} \| \phi - \psi \|_{L^2} + \| A - B \|_{L^2} \| \psi \|_{L^\infty} \\ 
                &\lesssim_\sfE \| \phi - \psi \|_{H^1} + \| A - B \|_{L^2},
        \end{align*}
    using energy estimates to bound the coefficients. Collecting these three inequalities, we conclude
        \begin{align*}
            \bigl\| (\widehat A, \widehat \phi) - (\widehat B, \widehat \psi) \bigr\|_{L^\infty_t (L^2 \times H^1)_x}
                \lesssim_\sfE T \cdot \| (A, \phi) - (B, \psi) \|_{L^\infty_t (L^2 \times H^1)_x}.
        \end{align*}
    Taking $T \ll_\sfE 1$ gives the result.             

    To estimate the $L^2$-difference of the magnetic potentials, integrate the transport equation for the difference of the divergence-free parts \eqref{iterate-transport-diff}, and apply the energy estimate \eqref{heat-energy-estimate} to the heat equation for the difference of the curl-free parts \eqref{iterate-heat-diff}. All together, these yield
        \begin{align}
            \bigl\| \widehat A - \widehat B \bigr\|_{L^\infty_t L^2_x} +  \bigl\| \nabla \PP^\cf \bigl( \widehat A - \widehat B \bigr) \bigr\|_{L^2_{t, x}} 
                &\lesssim \bigl\| \text{R.H.S. \eqref{iterate-heat-diff}} \bigr\|_{L^1_t L^2_x} \notag\\
                &\lesssim T \cdot \Bigl( \| \phi - \psi \|_{L^\infty_t L^2_x} \| \bfD_A \phi \|_{L^\infty_{t, x}} + \| \psi \|_{L^\infty_{t, x}} \| \bfD_A \phi - \bfD_B \psi \|_{L^\infty_t L^2_x} \Bigr) \label{eq:iterate-diff-A}\\
                &\lesssim_\sfE T \cdot\Bigl( \| \phi - \psi \|_{L^\infty_t L^2_x} + \| \bfD_A \phi - \bfD_B \psi \|_{L^\infty_t L^2_x} \Bigr) . \notag
        \end{align}
    
    To control the differences of the scalar fields, we apply the energy estimate \eqref{abstract-energy} to the Schr\"odinger equations for the differences of the scalar fields \eqref{iterate-schro-diff} and their magnetic derivatives \eqref{iterate-schro-diff-2},
        \begin{align*}
            \bigl\| \widehat \phi - \widehat \psi \bigr\|_{L^\infty_t L^2_x} + \bigl\| \bfD_{\widehat A} \widehat \phi - \bfD_{\widehat B} \widehat \psi \bigr\|_{L^\infty_t L^2_x}
                &\lesssim \| \text{R.H.S. \eqref{iterate-schro-diff}} \|_{L^1_t L^2_x} + \| \text{R.H.S. \eqref{iterate-schro-diff-2}} \|_{L^1_t L^2_x}.
        \end{align*}
    We estimate the right-hand side by 
        \begin{align*}
             \| \text{R.H.S. \eqref{iterate-schro-diff}} \|_{L^1_t L^2_x} 
                &\lesssim T\cdot \Bigl( \bigl( \| \widehat A \|_{L^\infty_{t, x}} + \| \widehat B \|_{L^\infty_{t, x}} \bigr) \bigl\| \widehat A - \widehat B \bigr\|_{L^\infty_t L^2_x} + \bigl( \| \phi \|_{L^\infty_{t, x}} + \|\psi \|_{L^\infty_{t, x}} \bigr) \| \phi - \psi \|_{L^\infty_t L^2_x} \Bigr) \bigl\| \widehat \psi \bigr\|_{L^\infty_{t, x}} \\
                    &\qquad + T^\frac12 \bigl\| \nabla \PP^\cf ( \widehat A - \widehat B) \bigr\|_{L^2_{t, x}} \bigl\| \widehat \psi \bigr\|_{L^\infty_{t, x}} \\
                    &\qquad + T \bigl\| \widehat A - \widehat B \bigr\|_{L^\infty_t L^2_x} \bigl\| \nabla \widehat \psi \bigr\|_{L^\infty_{t, x}} \\
                &\lesssim_\sfE T \cdot \Bigl(  \| \phi - \psi \|_{L^\infty_t L^2_x} + \| \bfD_A \phi - \bfD_B \psi \|_{L^\infty_t L^2_x} \Bigr),
        \end{align*}
    and similarly,
        \begin{align*}
             \| \text{R.H.S. \eqref{iterate-schro-diff-2}} \|_{L^1_t L^2_x}
                &\lesssim T\cdot \Bigl( \bigl( \| \widehat A \|_{L^\infty_{t, x}} + \| \widehat B \|_{L^\infty_{t, x}} \bigr) \bigl\| \widehat A - \widehat B \bigr\|_{L^\infty_t L^2_x} + \bigl( \| \phi \|_{L^\infty_{t, x}} + \|\psi \|_{L^\infty_{t, x}} \bigr) \| \phi - \psi \|_{L^\infty_t L^2_x} \Bigr) \bigl\| \bfD_{\widehat B} \widehat \psi \bigr\|_{L^\infty_{t, x}} \\
                    &\qquad + T^\frac12 \bigl\| \nabla \PP^\cf (\widehat A - \widehat B)\bigr\|_{L^2_{t, x}}\bigl\| \bfD_{\widehat B} \widehat \psi \bigr\|_{L^\infty_{t, x}} \\
                    &\qquad + T \bigl\| \widehat A - \widehat B \bigr\|_{L^\infty_t L^2_x} \bigl\| \nabla \bfD_{\widehat B} \widehat \psi \bigr\|_{L^\infty_{t, x}} \\
                    &\qquad + T \cdot \Bigl( 1 + \| \phi \|_{L^\infty_{t, x}}^2 + \| \Gauss \|_{L^\infty_x} \Bigr) \bigl\| \bfD_{\widehat A} \widehat \phi - \bfD_{\widehat B} \widehat \psi \bigr\|_{L^\infty_t L^2_x}\\
                    &\qquad + T \cdot \Bigl( \| \nabla \Gauss \|_{L^\infty_x} + \| \psi \|_{L^\infty_{t, x}} \| \bfD_B \psi \|_{L^\infty_{t, x}} \Bigr) \bigl\| \widehat \phi - \widehat \psi \bigr\|_{L^\infty_t L^2_x} \\
                    &\qquad + T \cdot \Bigl( \| \widehat \phi \|_{L^\infty_{t, x}} \| \bfD_B \psi \|_{L^\infty_{t, x}} \| \phi - \psi \|_{L^\infty_t L^2_x} + \| \widehat \phi \|_{L^\infty_{t, x}} \| \phi \|_{L^\infty_{t, x}} \bigl\| \bfD_A \phi - \bfD_B \psi \bigr\|_{L^\infty_t L^2_x}\Bigr) \\
                    &\lesssim_\sfE T \cdot \Bigl( \| \phi - \psi \|_{L^\infty_{t, x}} + \| \bfD_A \phi - \bfD_B \psi \|_{L^\infty_t L^2_x}\Bigr) \\
                        &\qquad+ T \cdot \Bigl( \| \widehat\phi - \widehat\psi \|_{L^\infty_{t, x}} +  \| \bfD_{\widehat A} \widehat \phi - \bfD_{\widehat B} \widehat\psi \|_{L^\infty_t L^2_x}\Bigr) ,
        \end{align*}
    using the $L^2_{t, x}$-bound from \eqref{iterate-diff-A} to control $\nabla \PP^\cf (\widehat A - \widehat B)$, and estimating the terms in $L^\infty_{t, x}$ using an appropriate combination of Sobolev embedding \eqref{GNinfty}, \eqref{inefficientLinfty} and the uniform energy bounds \eqref{prev-iterate-energy}, \eqref{iterate-energy}. Taking $T \ll_\cE 1$, the second term in the last line can be absorbed by the left-hand sides of the previous calculations. In total, the previous three calculations yield 
        \begin{align}\label{eq:iterate-diff-phi}
            \bigl\| \widehat \phi - \widehat \psi \bigr\|_{L^\infty_t L^2_x} + \bigl\| \bfD_{\widehat A} \widehat \phi - \bfD_{\widehat B} \widehat \psi \bigr\|_{L^\infty_t L^2_x}
                \lesssim_\sfE T \cdot \Bigl( \bigl\| \phi - \psi \bigr\|_{L^\infty_t L^2_x} + \bigl\| \bfD_A  \phi - \bfD_{B}\psi \bigr\|_{L^\infty_t L^2_x} \Bigr). 
        \end{align}
    Combining \eqref{iterate-diff-A}-\eqref{iterate-diff-phi} gives the desired inequality \eqref{iterate-magn-diff}.
\end{proof}

\bibliographystyle{alpha}
\bibliography{external/biblio}

\end{document}